\documentclass[a4paper]{amsart}  

\usepackage{amssymb} \usepackage{amscd} \usepackage{times}

\usepackage{amsmath}

\usepackage[french]{babel}

\def\zbb{\mathbb{Z}}  
  
  \def\phi{\varphi}
 \def\p1{{\mathbb{P}^1_\zbb}}

\newcommand{\be} {\begin{equation}}
\newcommand{\ee} {\end{equation}}

\title{Quelques remarques sur les vari\'et\'es, fonctions de Green et formule de Stokes.}

\author{Samy Skander Bahoura}

\address{Equipe d'Analyse Complexe et Geometrie, Universite Pierre et Marie Curie, 2, place Jussieu, 75005, Paris, France.}

\email{ samybahoura@gmail.com}

\begin{document}

\maketitle

\begin{abstract}

On donne quelques remarques sur les surfaces K3, projectifs complexe, r\'eels, le Tore, la topologie des surfaces r\'eelles de dimension 2, la fonction distance au bord, l'orientabilit\'e des boules geodesiques, traces de Sobolev et fonctions de Green et un Th\'eor\`eme de Dualit\'e, le degr\'e topologique de Leray-Schauder et des points fixes topologiques (en particulier, application en dimension 2), les solutions topologiques dans le cas negatif (critique et supercrtique, Eq de courbure scalaire prescrite negatif et Eq. supercritique negatif), la metrique induite sur un ouvert Lipschitzien, les coordonn\'ees geodesiques polaires et la formule de Gauss-Bonnet et sur le Th\'eor\`eme de la masse positive en dimension $\geq 3 $, dans le cas localement conform\'ement plat et non localement conform\'ement plat. Et le flot de Ricci.  Et sur les in\'egalit\'es du type $ \sup u = f(\inf u) $, le volume conforme au sens de la convexit\'e, l'enroulement et la distortion, et sur la th\'eorie de Yang-Mills et la relativit\'e g\'en\'erale et la cosmologie quantique et leur relation avec l'equation de Yamabe et de la courbure prescrite et du type courbure prescrite. Et des obstructions d'astronomie. Et sur les cordes, supercordes et D-branes et la symetrie. Et g\'eom\'etrisation. Et, un algorithme reconnaissant le oui et le non, l'algorithme SAT et un algorithme du probl\`eme $ p=np $.

\end{abstract}

\section{Quelques remarques:}

\smallskip

1) Sur les mesures de Haussdorf et la formule dintegration par parties ( Pourquoi la version Fusco est equivalente a celle de Necas ?):
$ H^{n-1}(\partial \Omega) <+\infty $ alors la mesure $ H^{n-1}_{|\partial \Omega} $ est une mesure de Radon, donc reguliere.

\smallskip

Soit $ A \subset \partial \Omega $, $ H^{n-1}_{|\partial \Omega} $ mesurable, comme la mesure est une mesure de Radon, elle est reguliere:

$$\forall \, \epsilon >0 \, \exists \, K_{\epsilon}, O_{\epsilon}, K_{\epsilon} \subset A \subset O_{\epsilon}, \,\, H^{n-1}_{|\partial \Omega} (O_{\epsilon}-K_{\epsilon}) \leq \epsilon. $$

On a, $ O_{\epsilon}=V_{\epsilon}\cap \partial \Omega $, avec $ V_{\epsilon} $ ouvert de $ {\mathbb R}^n $.

\smallskip

D'apres le Theoreme d'Urysohn;

$$ \exists \,\, f_{\epsilon} \in C^0({\mathbb R}^n, {\mathbb R}^+),\,\, 0 \leq f_{\epsilon} \leq 1,\,\, f_{\epsilon} \equiv 1, \,\, {\rm sur } \,\, K_{\epsilon}, \,\, supp f_{\epsilon} \subset V_{\epsilon}, $$

Soient $ \nu_1, \nu_2 $ les normales (interieures) dans la formulation de Fusco et Necas: 

\smallskip

$ \Omega $ est un ouvert Lipschitzien, on ecrit la formule d'integration par parties, dans les deux formulation, la mesure sur le bord est la mesure de Haussdorf $ H^{n-1}_{|\partial \Omega} $, ceci est du a la formule de l'aire appliquee localement dans des cartes. De plus par Fusco  le bord reduit est egal au bord usuel et $ |D \chi_{\Omega}|= H^{n-1}_{|\partial \Omega} $:

$$\int_{\Omega} div \phi =-\int_{\partial \Omega} \phi \nu_1 dH^{n-1}=- \int_{\partial \Omega} \phi \nu_2 dH^{n-1}, \,\, \forall \,\,\phi \in C^1_c({\mathbb R}^n, {\mathbb R}^n),$$

Avec  $\nu_1, \nu_2 \in L^{\infty}(\partial \Omega) $.

Donc,
$$ \int_{\partial \Omega} \phi (\nu_1-\nu_2) dH^{n-1} = 0, \,\, \forall \,\, \phi \in C^1_c, $$

Comme $ C^1_c $ est dense dans $ C^0_c $, la formule precedente est vraie pour $ \phi \in C^0_c $

\smallskip

On applique cette formule  a $ f_{\epsilon} $, on obtient pour tout mesurable $ A $:

$$ \int_{A} (\nu_1-\nu_2) dH^{n-1} = O(\epsilon) \to 0, $$

Donc, 

$$ \nu_1 =\nu_2, \,\, H^{n-1}_{|\partial \Omega}, p.p $$

Ce qu'on a montr\'e cest que la normale explicite de la formulation Necas  est egale a la normale theorique de la formulation de Fusco.

\smallskip

Dans la formulation de Fusco, la normale est obtenue comme "blow-up" autour d'un point du bord, cela revient \`a prendre la normale \`a la tangente au bord, par des cartes Lipschitziennes, on voit par ce proc\'ed\'e "blow-up", que la normale exterieure est la normale exterieure usuelle. (Car, l'espace tangent "\`a la mesure"("blow-up") est l'espace tangent usuel).

\smallskip

La formule d'integration par parties est bien connue. Pour sa preuve, voir les livres de Fusco (bords avec singularit\'es par recouvrement et suite de compacts $ K_n \to \partial \Omega-\{x_1, x_2,\ldots x_k\} $  dans le cas o\`u il n'y a qu'un nombre fini de points singuliers), Necas et Az\'e.

\smallskip

La formule d'integration par parties sur un domaine singulier en un nombre fini de points en dimension 2 (par exemple):

\smallskip

La formule est locale, on considere $ B(x_j,1/n) $ avec $ x_j $ un point singulier du bord. La formule est locale, par un recouverement et une partition de l'unit\'e on ecrit que (voir le Az\'e, dans chaque partie du recouvrement et on rajoute la boule o\`u il y a la singularit\'e $ B(x_j,1/n)$):

$$ \int_{\Omega} div ((1-\alpha_n) u) dx= \int_{\partial \Omega} (1-\alpha_n) u d\sigma, $$ 

avec $ \alpha_n $ une fonction reguliere et \`a support compact dans $ B(x_j,1/n) $,  $ \alpha_n \equiv 1 $ dans $ B(x_j,1/2n) $. Avec le fait, $ |\nabla \alpha_n|\leq Cn $ et $\mu_L(B(x_j,1/n))\leq C/n^2 $ et $ H_1(\{x_j\})=0 $, la mesure de Hausdorff 1-dimensionnelle (ici, on l'a fait pour le point $ x_j $). Puis on fait tendre $ n\to +\infty $. On obtient la formule de Stokes ou Green-Riemann pour un domaine contenant un nombre fini de singularit\'es au bord.

\smallskip

2) Pourquoi quand on a  $ u \in W^{1,p} \cap C^0(\bar \Omega) $, $ tr(u)=u $ ? on utilise Necas et Brezis:

\smallskip

On a:

\smallskip

$ \Omega $ est un ouvert Lipschitzien $ \Rightarrow $ cartes Lipschitziennes $\Rightarrow $ un operateur de Prolongement $ Pu \in C^0_c({\mathbb R}^n) \cap W^{1,p}({\mathbb R}^n) $, d'apres la consturction de loperateur de Prolongement dans des cartes, $ Pu_{| \bar \Omega} = u $ et,
$ \chi_n (\rho_n*Pu) \to Pu $ dans $ C^0_K({\mathbb R}^n) $ et dans $ W^{1,p} $ o\`u $ K  $ est un compact de $ {\mathbb R}^n $, donc:

$$ ||Tr(\chi_n(\rho_n*Pu)) - Tr(u)||_{L^p(\partial \Omega)} \to 0, $$

Et,

$$ ||Tr(\chi_n(\rho_n*Pu))-u||_{C^0(\partial \Omega)} \to 0, $$

Do\`u,

$$ Tr(u)=u,\,\, H^{n-1}_{|\partial \Omega} p.p, $$

Remarque: Par la formulation de Fusco, 

$$ Tr(u)=u_{\Omega}= \lim_{\rho \to 0} \dfrac{\int_{B_{\rho}(x)} u(y) dy}{|B_{\rho}(x)|} $$

 on voit que si $ u\in C^0(\bar \Omega) $, $ Tr(u)=u $.

\smallskip

3) Pourquoi $ W^{1,1}_0= \{ u\in W^{1,1}, tr(u)=0\} $ avec $ \Omega $ ouvert Lipschitzien. Voir le livre de Necas. On raisonne par rapport aux suites ($ \exists \,\, u_i \in C^{\infty} ({\mathbb R}^n), \, u_i \to u \in W^{1,p}(\Omega) $ pour les traces et la restrictions a des sous ensembles d'ouverts Lipschitzien)).

\smallskip

a) On considere $ \tilde u = u\times 1_{\Omega} \in W^{1,1}({\mathbb R}^n) $ par integration par parties et formule de Stokes.

\smallskip

b) Par des cartes et on multipliant par une fonction test a support compact, localement $\tilde u \in W^{1,1}(Q) $ avec $ Q $ un n-cube "au dela" du bord et nulle avant d'atteindre les cotes du bord qui sont interieurs.

\smallskip

On considere alors $ \tilde u_{\lambda} $, "la translat\'ee" de $ \tilde u $ suivant la direction normale au bord.  le support de $ \tilde u_{\lambda} $ est compact et strictement  \'a l'interieur du domaine. 

\smallskip
 
On a: (on applique ces id\'ees \`a un cube infini dans la direction normale au bord image  (locale) de $ \partial \Omega $, puis on se restreint au n-cube de depart. Sur les bords, la restriction de l'homeorphisme Lipschitzien envoie les traces sur les traces par la formule de l'aire ou coaire. (par convergence de suites sur des ouverts Lipschitziens).).

\smallskip

$ \tilde u_{\lambda} \in W^{1,1}(Q_{-\lambda}) $ avec,

\smallskip

$ \nabla (\tilde u_{\lambda})=(\nabla \tilde u)_{\lambda} $.

\smallskip

puis en utilisant des fonctions tests dans $ Q $  on obtient:

\smallskip

$ \tilde u_{\lambda} \in W^{1,1}(Q) $ pour $ \lambda $ assez petit.

\smallskip

Le support de $ \tilde u_{\lambda} $ est strictement dans $ Q $, on obtient alors, $ \tilde u_{\lambda} \in W^{1,1}_0(Q) $.

\smallskip

Par un theoreme de la moyenne (voir Necas, $ \tilde u $ s'annule en dehors de $ Q $), on obtient $ \tilde u_{\lambda} \to \tilde u $ dans $ W^{1,1} $. Donc, $ \tilde u \in W^{1,1}_0 $ et $ u\in W^{1,1}_0 $.
 
\smallskip

4) Soit $ (M,g) $ une variete riemannienne et soit $ B_r(x) $ une boule geodesique tres petite, telle que l'exponentielle realise un diffeomorphisme sur la boule de $ {\mathbb R}^n $: pourquoi, la normale exterieure $ \nu= \partial_r $ ?

\smallskip

On a: l'image de la base canonique par l'exponentielle est la base canonique:

\smallskip

$ d(exp_x) (\partial_i)=\tilde \partial_i $, on identifie les 1-formes $ dx_i $ et $ d\tilde x_i $.

\smallskip

Soit $ B^0 $ la base canonique de  $ {\mathbb R}^n $ et $ \tilde B^0 $ la base $ \{ \tilde \partial_1,. . ., \tilde \partial_n \} $ la base canonique de la boule geodesique $ B_r(x) $, alors $ d(exp_x) (B^0)= \tilde B^0 $.

\smallskip

La base $ B_{0,r}=\{\partial_{0,r}, \partial_{\theta_1},. . .,\partial_{\theta_{n-1}} \} $ est directe, donc $ det_{B^0} (B_{0,r}) >0 $.

\smallskip

Soit $ \tilde B_r=\{ d(exp_x)(\partial_{0,r}), \tilde \partial_{\theta_1}= d(exp_x)(\partial_{\theta_1}),. . ., \tilde  \partial_{\theta_{n-1}}= d(exp_x)(\partial_{\theta_{n-1}}) \} $, alors:

\smallskip

$ det_{\tilde B^0} (\tilde B_r) = [ det: B_r \to B_{0, r}] o [ det: B_{0, r} \to B_r, \,\,{\rm \,qui \, est \, directe \,}]o [det: B_r \to \tilde B_{0,r}]  >0$

\smallskip

La premiere est $ d(exp_x)^{-1} $ et la derniere $ d(exp_x) $ le determinant selimine, il reste celui de $ B_{0,r} $ vers $ B^0 $ qui est direct.

\smallskip

Donc la base $ \tilde B_r $ est directe par rapport a la base $ \tilde B^0 $.

\smallskip

$ n_{ext} $ est definit precisement comme cela tel que $ (n_{ext}, \tilde \partial_{\theta_1}, . . ., \tilde \partial_{\theta_{n-1}}) $ est directe par rapport a $ \tilde B^0 $.

\smallskip

On peut utiliser la formulation de S. Lang:

\smallskip

Soient $ \omega_{\theta} $ et  $\Omega $ les formes volumes sur la sphere $ S_r(x) $ et la boule $ B_r(x) $ $ n_{ext} $ est tel que son dual $ n_{ext}^* $ verifie:

$$ n_{ext}^* \wedge \omega_{\theta}= \Omega, $$

Cest-a-dire que \`a des constantes multiplicatives pres (determinant des metriques)):

$$ n_{ext}^* \wedge d\tilde \theta_1 \wedge . . .\wedge d \tilde \theta_{n-1} = d\tilde x_1 \wedge...\wedge d\tilde x_n $$

cela veut dire que:

$$ \Omega (n_{ext}, \tilde \partial_{\theta_1}, . . ., \tilde \partial_{\theta_{n-1}}) >0 $$

c'est a dire que la base $ \{n_{ext}, \tilde \partial_{\theta_1},. . ., \tilde \partial_{\theta_{n-1}} \} $ est directe par rapport a la base $ \{ \tilde \partial_1,. . ., \tilde \partial_n \} $ et donc $ n_{ext}=\partial_r $.

\smallskip

Voir aussi le livre de Gallot-Hulin-Lafontaine o\`u la normale exterieure est definie par le produit interieur. (Par rapport \`a la base du bord et l'element de volume du bord) On a le meme resultat.

\smallskip

Par sa construction, la boule geodesique est orientable et orient\'ee. En chaque point, l'exponentielle conserve l'orientation et induit une orientation par raport \`a la carte $ (\Omega, \phi) $ \`a partir de laquelle elle construite. Par un argument de recouvrement, si on raisonne sur un compact $ K $ en utilisant l'exopentielle et cette propri\'et\'e d'orientabilit\'e des boules geodesiques via l'exponentielle, il faudera suppos\'e la vari\'et\'e orientable.

\smallskip

5) Dans la construction de la fonction de Green pour un operateur coercif $ \Delta+a $ voir le livre d'Aubin et le monograph de Frederic Robert, on utilise le produit de convolution et le fait que les traces des fonctions sont nulles jusqu'\`a l'etape:

$$ -\Delta V_x+a V_x=0, \,\, {\rm dans \,} M, V_x=-G_x,\,\, {\rm sur \,} \partial M, $$

Alors la fonction $ G_x \in C^{1,\theta} \cap W^{2, p} (M-\{x\}) $, car la parmetrix $ H_x $ verifie au sens $ C^2_0 $ (Agmon),

$$ -\Delta_{dist} H_x = \delta_x- \Delta_y H_x \in L^{\infty} (M-\{x\})$$

Les $ \Gamma_i $ verifient cette proprit\'e, elles sont nulles au bord car leur support est dans $ B(m, \delta_m/2) $ et le rayon dinjectivit\'e $ =\delta_m \leq d(m, \partial \Omega) $. Puis on utlise des cartes pour se ramener au demi-espace et utiliser les estimations de Agmon-Douglis-Nirenberg.

\smallskip

On peut aussi, utiliser $ G_x $ dans lequation et utiliser les Theoremes de Gilbarg-Trudinger :

$$ \Delta (V_x-\eta G_x)+a(V_x-\eta G_x)= -\Delta (\eta G_x)+a(\eta G_x) \in L^p,\,\, {\rm et} \,\, (V_x-\eta G_x)= 0,\,\, {\rm sur} \,\, \partial M, $$

Avec $ \eta $ une fonction cutoff egale a $ 1 $ dans un voisinage de $ \partial M $ et $ 0 $ au voisinage de $ x $ ( nulle au voisinage de la singularit\'e).

\smallskip

On obtient $ V_x-\eta G_x  \in W^{2,p}(M)$ et comme $ \eta G_x \in W^{2,p}(M) $, alors : $ V_x \in W^{2,p}(M) $.

\smallskip

Remarques:

\smallskip

1) Dans la formulation variationelle on pour une solution $ u \in W^{1,2}_0 $ de :

$$ -\Delta u+ a u = f,\,\, {\rm et} \,\, u=0\,\, {\rm sur }\,\, \partial M,$$

Dans le monograph de Brezis-Marcus-Ponce (ils supposent les solutions $ W^{1,1}_0 $), ils prouvent que ces solutions sont au sens $ C^2_0 $ (Agmon), d'o\`u on peut prouver qu'elles sont $ W^{2,p} $ et qu'on a les estimations a priori dans $ W^{2,p} $ d'Agmon ( m\'ethode des quotients differentiels) ou en sachant qu'elles sont $ W^{2,p} $ qu'on a les estimations a priori par Calderon-Zygmund ecrits dans le Gilbarg-Trudinger.

\smallskip

Dans le cas du Laplacien, la constante dans l'inegalit\'e de Calderon-Zygmund, ne depend pas de $ \Omega $ ce qui fait qu'on peut par coninuit\'e obtenir l'estimation $ W^{2,p}$ d'un operateur general \`a partir du Laplacien. On approche l'operateur par un operateur constant.

\smallskip

2) Les estimations d'Agmon-Douglis-Nirenberg, dans ce cas sont bas\'ees (sur de l'intergartion par parties, par rapport a la variable $ t>0 $) et sur les integrales singulieres de Calderon-Zygmund.

\smallskip

C'est une autre preuve par le potentiel du demi-espace.

\smallskip

Par exemple:

\smallskip

En effet, on considere le probleme suivant:

$$ -L u=0,\,\, u= \phi \, {\rm sur } \,\, \partial \Omega $$ 

On peut supposer $ L $ le Laplacien euclidien, le cas general se ramene au cas constant (voir Gilbarg-Trudinger, par continuit\'e on approche l'operateur general par un operateur constant).

\smallskip

De meme pour Agmon-Douglis-Nirenebrg, ils considerent des operateurs a coefficients constants, et dans le cas du Laplacien, la constante dans linegalit\'e obtenue ne depend pas de $ \Omega $ ce qui fait qu'on peut par coninuit\'e obtenir l'estimation $ W^{2,p}$ d'un operateur general \`a partir du Laplacien. On approche l'operateur par un operateur constant.

\smallskip

Par exemple on ecrit \`a l'ordre 1:

\smallskip

Il suffit, par des cartes, de consider le demi-espace et on utilise le noyau de Poisson $ P(x, t)= \dfrac{t}{(t^2+|x|^2)^{(n+1)/2}} $.

\smallskip

On suppose que $ \phi \in C^2 \cap W^{2, p} $. On ecrit:

$$ \partial_s (P(x-y,t+s)\phi (y, s))= \partial_s P \phi+ P \partial_s \phi, $$

On intergre en $ y $ et $ s $ et on obtient (on utilise la representation integrale de $ u $ en fonction du noyau de Poisson):

$$ u(x,t)= \int_{{\mathbb R}^n} P(x-y, t) \phi(y,0) dy=\int P(x-y, t+T) \phi (y, T)+ \int \int  \partial_t P \phi + \int \int P \partial_t \phi, $$

On a $ \partial_t P \equiv \dfrac{1}{|Q|^{n+1}} $ , $ Q $ un point quelconque, et $ dim ({\mathbb R}^n \times {\mathbb R}) = n+1 $, on utilise l'integrale singuliere de Calderon-Zygmund.

\smallskip

Pour $ \partial_t u, \partial_x u $, on derive sous le signe $ \int $ et puis on integre par parties sur $ {\mathbb R}^n $. On obtient des integrales du type,

 $$ \int \int \partial_t P \partial_t \phi,\,\,  \int \int \partial_t P \partial_x \phi,\,\, {\rm et} \,\, \int \int \partial_x P \partial_t \phi $$.

On elimine $ \int P(x-y, t+T) \phi (y, T) $ par son estimation en $ c(T) \to 0  $ quand $ T \to + \infty $.

\smallskip

On obtient des estimations du type:

$$ ||\partial u ||_{L^p} \leq  C ||\phi||_{W^{1,p}}, $$

On refait la meme chose en derivant une deuxieme fois. 

\smallskip

Comme dans le Gilbarg-Trudinger, on ecrit les solutions en fonction du potentiel du demi-espace (voir chapitre 4 de Gilbarg-Trudinger). Dans le cas d'Agmon-Douglis-Nirenberg, ils considerent le demi-espace.

\smallskip

6) Sur le degr\'e topologique et l'article de De Figueiredo-Lions-Nussbaum: 

\smallskip

On explique les etapes qui permettent de calculer le degr\'e topologique de Leray-Schauder. Ce degr\'e vaut $ -1 \not = 0 $, ce qui veut dire qu'il y a au moins une solution topologique:

\smallskip

Ils utilisent la definition du degr\'e topologique de Leray-Schauder.

$ \exists \, T_0 \in {\mathbb R}^+, deg [x-F(x,T_0) B_{R_2}, 0]=0 $ cela veut dire qu'il n'y a pas de solution a cette equation.

$ x-F(x,t) \not = 0, \, t\in {\mathbb R}^+, x\in \partial B_{R_2} $, cela veut dire quon peut utiliser l'homotopie:

$ deg[x-F(x,t), B_{R_2}, 0]= deg [x-F(x, 0), B_{R_2}, 0]= deg [x-F(x,T_0), B_{R_2}, 0]=0 $,

$ t=0 $, $ F(x,0)=\Phi(x) $ et par Nussbaum et la propri\'et\'e d'excision et d'additivit\'e, on peut d\'ecomposer le degre en 2, la condition est que $ \Phi(x)\not = x $ sur $ \partial B_{R_1} $ et $ \Phi(x)=F(x,0) \not = x $ sur $ \partial B_{R_2} $, alors:

$$ 0=deg [x-\Phi(x), B_{R_2}, 0]= deg [x-\Phi(x), B_{R_1}, 0]+ deg [x-\Phi(x), \{ R_1 < ||x|| < R_2 \}, 0], $$

La derniere condition est celle de l'homotopie avec l'identit\'e et le degre de l'identit\'e est 1.

$ deg[x-\beta \Phi(x), B_{R_1}, 0]= deg [x-\Phi(x), B_{R_1}, 0]= deg [x, B_{R_1}, 0] = 1, $

Donc,

$$ deg [x-\Phi(x), \{ R_1 < ||x||<R_2 \}, 0] =-1, $$

Dans leur exemple, De Figueiredo-Lions-Nussbaum, utilisent l'estimation integrale $ \int_{\Omega} f(u(x)+t) dx \leq C $, $ C $ independante de $ t $ et les estimations a priori pour verifier que les conditions d'applications du degr\'e topologique sont verifi\'ees, dans ce cas les $ t $ et les solutions sont uniform\'ement born\'es.

\smallskip

Concernant l'article de Crandall-Rabinowitz, il traite du cas o\`u en particulier on a un terme nonlineaire du type exponentiel. Grace \`a la condition de stabilit\'e (ils prennent des fonctions tests particulieres), ils prouvent que la masse et le volume ou energie, sont born\'ees (jusqu'a la dimension 10).

\smallskip

La formulation du degr\'e pr\'ec\'edente s'applique aux equations en $ \sinh $ sans poids ou avec poids, pour le laplacien ou l'operateur plus general que le laplacien, en dimension 2 par exemple. Voir l'article de Brezis-Turner et vers la fin de ce print, sur le degr\'e en dimnension 2 et pour des systemes avec singualrit\'es.

\smallskip

7) Sur les surfaces $ K3 $. Soit $ E $ un fiber bundle (fibr\'e vectoriel) et $ M $ une vari\'et\'e Complexe de metrique $ h $. 

\smallskip

Remarquons dabord que cette varit\'et\'e est orientable car les changements de cartes sont holomorphes donc le Jacobien r\'eel est le carr\'e du Jacobien complexe, cest du a l'holomorphie, comme en dimension 2 (equation de Cauchy-Riemann).

\smallskip

Pour ce fibre bundle $ E $, il existe une unique connexion compatible avec l'holomorphie et la metrique, c'est l'analogue complexe de la connexion de Levi-Cevita:

$$  \pi^{0,1} \nabla= \bar \partial $$

et,

$$\nabla h = 0 $$

C'est la connexion hermitienne ou connexion de Chern, semblable a la connexion de Levi-Civita dans le cas Riemannien.

\smallskip

Des qu'on a une connexion, une derivee covariante, on a la courbure, on derive 2 fois, comme dans le cas Riemannien, on a la courbure et donc la courbure de Ricci et la classe de Chern , cest la classe de Ricci.

\smallskip

Connexion de Chern $ \Rightarrow $ symboles de Christofels $ \Rightarrow $ Connexion form, $ \omega $, $ \Rightarrow $ Curvature form,

\smallskip

Dans le cas $ E=\wedge^2 T^*(M) $ les 2-formes alternees:

\smallskip

$ \omega= \omega_i dx^i \otimes (dx^j\wedge dx^k) $, car ici,ce qui joue le role de champs de vecteurs (pour la connexion de Levi-Cevita), c'est les sections sur $ E $ (les vecteurs de base sont les 2-formes alternees), les champs de 2-formes aleternees,

\smallskip

$ d\omega= \Omega_{ij} dx^i\wedge dx^j \otimes dx^l \wedge dx^m $,des qu'on d\'erive encore deuxieme fois on obtient la courbure.

\smallskip

Dans le cas des surfaces $ K3 $, $ E=\wedge^2 T^*(M)=0$, 2 formes alternees:

\smallskip

Si $ E =0 $, la connexion nulle convient pour la derivation sur $ E $ or la connexion de Chern est unique dou $\nabla \equiv 0 $. Donc, il ny a pas de courbure do la classe de Chern est nulle, $ c_1=0 $.

\smallskip

Les groupes de cohomologies sont connues: 

\smallskip

faire attention:

\smallskip

1) ici, c'est les coherent sheaves sur les varietes complexe, il faut prendre la dimension complexe 2 et non reelle pour les groupes de Cohomologie du Sheaf et pour la dualite de Serre, cela va de $ 0 $ \`a $ 2 $, par contre pour les nombres de Betti et la theorie de Hodge, on prend la dimension reelle 4, avec la convention $ H^r= \sum_{\{p+q=r, 0 \leq p, q \leq 2 \}} H^{p,q} $, $ p, q $ varient selon la dimension complexe, quand on inverse $ p $ et $ q $ on a par la conjuguaison complexe , la meme dimension. 

\smallskip

2) Par la cohomologie et dualit\'e de Dolbeault, $ H^{0,1} $ est isomorphe  \`a la cohomologie de Cech $ H^1(S=K3, \Omega^0)=H^1(S=K3, O_S) $ car les $0$-formes sont les fonctions et on prend la cohomologie dans le faisceau de fonctions $ O_S $ et par definition $ H^1(S, O_S)=0 $ d'ou $ h^{0,1} = dim H^{0,1}= dim H^1(S, O_S)=0 $). On voit que les groupes de cohomologie de sheaf (faisceau) sont des donn\'ees et par la dualit\'e de Serre (dimension complexe) ,on les connait tous, il n'y a que le premier et le dernier, do\`u la carcateristique Sheaf egale a 2.

\smallskip

De meme pour la cohomologie de Hodge (dimension reelle 4) on les connait tous ainsi que les nombres de Betti.

\smallskip

Le theoreme de Max Noether donne $ X=S=K3 $:

\smallskip

$ \chi(X, O_X)=\chi(S,O_S)=2= \dfrac{ \int c_1\wedge c_1+ \int c_2}{12} $, avec $ c_1 $ la premiere classe de Chern et $ c_2 $ la deuxieme classe de Chern, definies dans le developpement de la courbure  $ \Omega $

\smallskip

$ det (\Omega-t I)= 1+ tr(\Omega) + \Omega\wedge. .. = 1+ c_1+c_2 $,

\smallskip

L'integrale $ \int c_2 $ est le 2\`eme nombre de Chern dans le cas d'une surface il est egale a la caracteristique d'Euler-Poincar\'e, dans le cas dune surface $ K3 $, cest donc  $ \chi (K3)= \int c_2=24 $.

\smallskip

Ici, quand on contracte on elimine des termes du type $ dx^j\wedge dx^k $ et $ c_2 $ est une quatre forme ("forme volume") car on a contracte une  8-forme.

\smallskip

Plus generalement, on definit les nombres de Chern comme:$ \int c_k $ avec $ c_k $ la $ k-$ieme classe de Chern.

\smallskip

Comme la surface $ K3 $ est de dimension 4 et de classe de Chern nul, il existe une metrique d'Einstein-Kahler de constante 0, voir le livre d'Aubin. Donc, cette vari\'et\'e possede une metrique d'Einstein de constante nulle, donc la courbure de Ricci est nulle et donc la courbure scalaire aussi.

\smallskip

(Remarquons d'abord qu'avec $ (K3, h, \nabla) $ on a une structure hermitienne et par Siu, $ K3 $ est Kahler, c'est a d ire que la structure Riemannienne $ (K3, g=Re(h), \nabla_{g}) $ est compatible avec la structure complexe, les connexions de Chern et Levi-Civita coincident.

\smallskip

On sait que la classe de Chern est nulle et donc, il existe une metrique d'Einstein-Kahler de constante 0).

\smallskip

On utilise la formule donnant la caracteristique d'Euleur-Poincar\'e en fonction du tenseur de Weyl et de la courbure de Ricci et Scalaire, pour enfin dire que le tenseur de Weyl est non nul.

$$ 0\not = \chi(K3)\equiv\int |Weyl|^2+ (Ricci_g, S_g)= \int |Weyl|^2 $$ 

Donc, une surface $ K3 $ possede une metrique non-loclamment conformement plate.

\smallskip

{\bf Remarque importante:} 

\smallskip

Pour une surface $ K3 $, il existe des points $ x_0 $ tels que: $ Ricci (x_0)=0 $ et $ Weyl(x_0)\not =0 $. Du point de vue de la physique et de l'infinitisimal, un tenseur $ T $ a une valeur donn\'ee en $ x_0 $: on fait tendre $ T(x) $, vers  $ T(x_0) $, $ dT \to 0 $, en infinitisimal en physique, il suffit de prendre $ dT(x_0)=0 $, il est constant en infinitisimal. Donc: $ dWeyl(x_0)=dRiemann(x_0)=0 $. On peut dire que du point de vue de la physique et de l'infinitisimal, il y a des points o\`u $ K3 $ est localement symetrique. Ceci pour donner un exemple de vari\'et\'e Ricci plate non localement conformement plate et localement symetrique. En physique et en infinitisimal, c'est possible \`a partir d'une surface $ K3 $. Du point de vue de la physique et de l'infinitisimal, une d\'eformation locale de la surface $ K3 $ \`a partir de la surface $ K3 $.

\smallskip

////////////////

\smallskip

a) Le Projectif complexe de dimension complexe $ n \geq 2 $ est une variete Kahlerienne, c'est dire qu'elle possede (3 structures) dont la Riemannienne et complexe, qui sont compatibles, les connexions de Levi-Civita et de Chern sont les memes, pour la metrique de Fubini-Study. Il est d'Einstein et de courbure sectionelle non constante, donc non localement conformement plat. Il est orientable.

\smallskip

Pourquoi la courbure sectionelle du projectif complexe de dimension $ n \geq 2 $ est non-constante: C'est ecrit dans le Gallot-Hulin-Lafontaine.

\smallskip

a) Ils utlisent l'equation des champs de Jacobi, pour determiner le tenseur de Riemann.

\smallskip

b) Ils utlisent la multiplication par le nombre complexe, $ i $, la structure complexe, $ J $. Et le fait que $ H_x $ sont les vecteurs orthogonaux a $ z $ et $ i z $, $ z\in {\mathbb S}^{2n+1} $.

\smallskip

Ils utilisent la variations $ H(t,s) $ et les champs de Jacobi.

\smallskip

La courbure sectionelle  pour $ u, v $ est calcul\'ee a partir du champs de Jacobi, $ Y(s)=cos s u+ sin s Jv $:

$$ K(u,v)=1+3 sin^2 s. $$

qui n'est pas constante.

\smallskip

b) Le Projectif complexe de dimension complexe 1, est isometrique a la sphere de dimension 2 reelle et donc de curbure sectionnelle constante. C'est une surface de Riemann reelle (donc localement plate).

\smallskip

Pourquoi, le Projectif complexe de dimension 2 est isometrique \`a la sphere de dimension 2. Voir le Gallot-Hulin-Lafontaine:

\smallskip

1-On utilise le quotient par un groupe de Lie. Comme le groupe de Lie agit proprement et transitivement et..., cela veut dire qu'il existe une structure de vari\'et\'e sur le quotient $ {\mathbb S}^3/{\mathbb S}^1 $ telle que l'application de passage au quotient soit une submersion.

\smallskip

2-Le fait qu'on ait une submersion, on utilise la definition, localement, on a des sous-vari\'et\'es, on construit, par des cartes, une base du noyau et apres un proc\'ed\'e d'orthogonalisation de Gramm-Schmidt, on construit un supplementaire orthogonal. l'application de passage au quotient est bijective sur ce supplementaire. ce qui permet de constuire localement un "lift".

\smallskip

3-On definit la metrique sur le quotient de telle maniere que l'application de passage au quotient soit une isometrie, on la construit a partir du lift et par projection sur $ H_x $.

\smallskip

(Remarque: soit $ p(x)=[x] $ l'application de passage au quotient et $ \gamma $ une isometrie, $ x=\gamma^{-1}o\gamma (x) \Rightarrow [x]=[\gamma (x)] \Leftrightarrow  p(x)=po\gamma(x) $, ceci est important dans la consturction de l'isometrie en prouvant que cela ne depend pas de point de la fibre).

\smallskip

On voit que pour le cas de $ {\mathbb S}^3 $ et $ {\mathbb P}^1({\mathbb C})\equiv {\mathbb S}^3/{\mathbb S}^1 $, que la metrique du projectif est celle de la 3-sphere pour des vecteurs particuliers $v, w $, par exemple orthogonaux \`a $ z $ et $ i z $ avec  $ z\in {\mathbb S}^3 $. ( le produit scalaire sur $ H_x $ est egal a celui  pour  l'espace tangent au quotient).

\smallskip

Par exemple si on prend $ z=(1,0)=[(1,0), (0,0)] $ en notation complexe, les $ v, w $ sont de la forme $ v=(0, v_1) $, $ w=(0, w_1) $. Puis en utlisant l'apllication $ H $ du Gallot-Hulin-Lafontaine $ H(u,v)=(2u\bar v, |u|^2-|v|^2) $ et les chemeins $ u(t)=cos t z+ sin t v $, $ w(t)= cos t z + sin t w $, (on se ramene \`a ce cas par une roation car la vari\'et\'e est invariante par $ U(n+1) $, les rotations complexes,(elle est construite \`a partir de $ S^{2n+1} $), on a:

$$ < H(u(t))', H(v(t))'>_{{{\mathbb S}^2}_{t=0}}=<v,w>_{{\mathbb S}^3} = <v,w>_{{\mathbb P}^1(C)}. $$

Cela veut dire que,

$$ H^* (g_{{\mathbb S}^2})=g_{{\mathbb P}^1({\mathbb C})}. $$

Remarque: Les revetements riemanniens, sont similaires, au lieu d'avoir une submersion, on a un diffeomeorphisme local. Voir, le Gallot-Hulin-Lafontaine, pour la contruction du Tore et des projectifs reels (leur metriques).

\smallskip

Dans le cas des revetements, puisqu'il y a un diffeomorphisme local, la courbure sectionelle se conserve (voir le Hebey). On voit alors, dans la construction du projectif r\'eel et le Tore, que la courbure sectionelle est constante, donc, elles sont conformement plates.

\smallskip

On a la meme chose avec les groupes de Lie et les espaces homogenes. On a une fibration qui est une submersion et on construit une metrique $ G- $ invariante par le meme proc\'ed\'e que les submersions et les revetements. De plus, si par exemple, on considere $ P_n(C) $, les groupes  unitaires (les rotations comlexes) $ U(n+1) $ et $ SU(n+1) $ agissent transitivement et proprement sur $ P_n(C) $ par l'intermediaire de la sphere $ {\mathbb S}^{2n+1} $.

\smallskip

1- $P_n(C) $ est \`a la fois identifiable \`a $ U(n+1)/U(1)\times U(n) $ et $ SU(n+1)/S(U(1)\times U(n)) $. De plus, comme on la vu dans le cas des submersions, pour connaitre la metrique de $ P_n(C) $, il suffit de se placer sur la sphere $ S^{2n+1} $, car on prend $ H_x $ l'orthogonal de $ z $ et $ iz $, avec $ z\in S^{2n+1} $. Pour ce qui concerne les groupes de Lie, la formule donnant la metrique $ G-$ invariante est egale \`a ($ C $ designe le produit scalaire dans $ {\mathbb C} $):

$$ <x|y>_{[e], P_n(C)} =\int_{H_1} <Adh_1.x|Adh_1.y>_{C} dv_{H_1} $$

Or, dans l'espace des matrices, la derivee de l'adjoint est connue et est egale, (on raisonne par rapport aux chemins et en particulier de l'exponentielle):

$$ Adh_1.x=Ad_{G/H}h_1.x=h_1^{-1}.x.h_1,\,\, Adh_1.y=h_1^{-1}.y.h_1, $$

avec l'identification, $ x=[x]=u $, $ [x]\in SU(n+1)/S(U(1)\times U(n)) $ et $ x=u\in H_x \subset {\mathbb C} \equiv P_n(C) $.

\smallskip

Car ici, $ G/H =P_n(C) $ et on prend le plan tangent et dans ce cas il faut prendre des vecteurs de $ H_x $ donc de $ S^{2n+1} $, (raisonner par rapport a l'exponentielle et les chemins). ($ t\to \exp(tx) $ est un bon chemin et en plus quand on derive on a l'adjoint). (on se ramene au produit scalaire de $ S^{2n+1} $.

\smallskip

$ G=SU(n+1), H=G_m=S(U(1)\times U(n)) $. 

\smallskip

{\bf Remarque:} La metrique sur $ G/H $ est definie de telle maniere que la correspondance: $ F: [x] \in G/H \to u \in P_n(C) $ soit une isometrie. (voir le livre de Gallot-Hulin-Lafontaine, sur la correspondance via $ S^{2n+1} $ entre $ SU(n+1)/S(U(1)\times U(n)) $ en tant que groupe de Lie et $ P_n(C) $ via $ {\mathbb S}^{2n+1} $. Donc, on retrouve le fait que les metriques sont proportionnelles. (Eux, ils notent $ U(n+1)/U(1)\times U(n) $).

\smallskip

Gallot-Hulin-Lafontaine disent qu'on peut pendre "any scalar product on $ T_e(G/H) $. Comme $ G/H $ s'identifie par $ F $ \`a $ P_n(C) $, on choisit le produit scalaire tel que $ F $ soit une isometrie vers le projectif complexe. Ce qui revient \`a obtenir le produit scalaire du projectif complexe.

\smallskip

Finalement, on ecrit, puisque $ \bar h_1^th_1=1 $ et $ (\bar h_1^{-1})^th_1^{-1}=1 $:

$$ <x|y>_{[e],P_n(C)} =\int_{H_1} <(h_1^{-1}.y.h_1)^t|h_1^{-1}.x.h_1>_{C} dv_{H_1}= $$

$$ = \int_{H_1} <(h_1^{-1})^t.y^t.(h_1^{-1})^t|h_1^{-1}.x.h_1>_{C} dv_{H_1}= $$

$$ = \int_{H_1} <y^t|x>_{C} dv_{H_1}=\mu <x|y>_{{\mathbb C}},  $$

 avec, $ \mu >0 $ et $ x, y \in H_x $, on reconnait le produit de deux vecteurs de $ S^{2n+1} $ qui donne le produit scalaire dans $ P_n(C) $.

\smallskip

Ce qui a \'et\'e prouv\'e est que la metrique $G-$invariante est egale \`a la metrique usuelle du projectif complexe.

\smallskip

Pour savoir directement, que le projectif complexe est d'Einstein et de courbure sectionelle non constante. Utiliser le livre de Kobayashi-Nomizu, tome 2, pour avoir l'expression de la metrique ou et le livre d'Aubin, pour avoir l'expression du tenseur de Ricci, et par les calculs, directement, on prouve qu'il est d'Einstein. Pour a courbure sectionelle, voir aussi dans, le livre de Kobayashi-Nomizu, $ K, K_0 $, dans le chapitre Complex manifolds. Avec la metrique tir\'ee du livre de Kobayashi-Nomizu, la courbure sectionelle holomorphe est constante, ce qui induit que la courbure sectionelle Riemannienne est de la forme:

$$ K(X,Y)=(1+3 (cos a)^2)/4,\,\, a, \,l'angle \, de\, g(X,JY). $$

\smallskip

Ici, par la fomulation de Gallot-Hulin-Lafontaine avec la phrase "any scalar product", on transporte la metrique du projectif, pour avoir une structure d'espace homogene, une structure de plus sur le projectif complexe par l'isometrie (par definition et choix) de $ SU(n+1)/S(U(1)\times U(n)) $.

\smallskip

(On a une isometrie $ F $ entre la vari\'et\'e  homog\`ene $ SU(n+1)/S(U(1)\times U(n)) $ et $ P_n(C) $, la courbure de la var\'et\'e homoge\`ene est alors $ RiccioF $ (voir le Hebey), elle constante).

\smallskip

En tant que submersion:

\smallskip

Pour voir que la metrique de $ P_n(C) $ est celle de $ S^{2n+1}/S^1$, du Gallot-Hulin-Lafontaine, on ecrit que la metrique $ ds $ du Kobayashi-Nomizu est a la fois invariante par "rotation", on se place alors dans la carte de $ z=(1,0,\ldots 0) $ avec la carte $ (t_i)=(z_i/z_0)_{(1\leq i \leq n)} $, on ecrit alors $ ds $ en coordonn\'ees $ z_i $ et la metrique en $ z $ est euclidienne sur  les vecteurs de $ {\mathbb C}^n=H_x $, or la metrique de $ S^{2n+1}/S^1 $ dans $ H_x $ est la metrique de $ S^{2n+1} $ qui est la metrique euclidienne de $ {\mathbb C}^n=H_x $. Donc la metrique $ ds $, dite de Fubini-Study dans le Kobayashi-Nomizu est la metrique dite de Fubini-Study dans le Gallot-Hulin-Lafontaine. 

\smallskip

{\bf Remarque sur les espaces homog\`enes:} Tout ce travail sur les groupes de Lie, sert en fait pour les ensembles de matrices. Cette Th\'eorie est construite pour les groupes matriciels, $ SU(n+1), SO(n) $... (pour ces cas particuliers, les vari\'et\'es sont de dimensions finies).

\smallskip

(Remarque sur les actions de groupes: le fait d'ecrire $ . $, le point, dans l'action de groupe, c'est une notation, pour dire que le groupe agit sur cet ensemble, comme c'est ecrit dans le debut du chapitre sur les groupes de Lie dans le Dieudonn\'e, c'est une notation d'ecrire la derivation, $ s\cdot k_x $ dans l'espace tangent $ T_e(G) $, pour ecrire la composition de la derivation, $ \partial \psi/\partial \xi_j (s, \xi_1,\ldots, \xi_n) k_j $. On voit bien par la formule de derivation de la composition, on a $ s\cdot (t\cdot k_x) $ remplace $ \partial_{\rho_i}(s, \rho)\partial^i_{\xi_j}(t, \ldots) $

\smallskip

(Cette notation, (d'action de groupe, $ \cdot $, le point), permet d'expliquer pourquoi l'espace tangent $ T_e(G/H) $ s'identifie \`a $ T_e(G)/T_e(H)= G_-/H_- $. Car en prenant des chemins, $ c: t\to [g(t)]$ dans $ G/H $ et en les derivant cela revient \`a prendre l'action de la classe dans le quotient:

$ \frac{dc(t)}{dt} =\{ \frac{d(g(t).h)}{dt} \cdot h, h\in H\}=\{ g'\cdot h, h\in H \}=[g'].T_e(H)=classe\, de \, g' $, car on derive un chemin suivant $ H $ donc dans  l'espace tangent \`a $ H $.)

\smallskip

Le paragraphe suivant, traite du cas du projectif complexe et l'exercice propos\'e dans le Gallot-Hulin-Lafontaine, on a vu comment la submersion construite permet d'identifier $ P_n(C) $ et $ S^{2n+1}/S^1 $, on regarde aussi comment le projectif est vu comme un espace homogene en comparant $ S^{2n+1}/S^1 $ et $ SU(n+1)/S(U(1)\times U(n)) $:

(***Voir l'exercice sur le projectif complexe dans le Gallot-Hulin-Lafontaine, ici, aussi, on utilise les chemins dans la determination des espaces tangents en $ e=I_n $ ou algebre de Lie, on derive par exemple $ (I_n+H)(I_n+H^t)=I_n $  en prenant $ dH/dt_{t=0}=A$, pour $ U(n)$, on voit alors que $ A $ est antihermitienne. Pour $ SU(n) $, on derive le determinant en $ I_n $ on a de plus $ Tr(A)=0 $ et l'application adjointe est l'application adjointe reellement, quand on calcule $ h x h^{-1} $ avec $ x= \left ( \begin{matrix}  0 & -\bar v \\ v & 0 \end{matrix} \right) $;  
$ h= \left ( \begin{matrix}  \lambda & 0 \\ 0 & A \end{matrix} \right) $; 
$ h^{-1}=\bar h^t= \left ( \begin{matrix}  \bar \lambda & 0 \\ 0 & \bar A \end{matrix} \right) $  

$ hxh^{-1}= \left ( \begin{matrix}  0 & -\lambda A\bar v \\ \bar \lambda Av & 0 \end{matrix} \right) $,

\smallskip

De plus pour voir que $ x \in T_e(G)/T_e(H)=G_-/H_-$ est de la forme precedente, on a $ T_e(G)$ est l'ensemble des matrices antihermitiennes de trace nulles et $ H=S(U(1)\times U(n))$, $ T_e(H) $ est l'ensemble des matrices antihermitinnes, comme on quotient par $ T_e(H)$ on elimine les termes centraux et $ x $ a la forme precedente. On voit alors que la condition d'irreductibilit\'e isotope est vraie ici,  la definition est:$ \forall \lambda, \forall A, \forall w\in W $, $ \bar \lambda A w \in W \Rightarrow W=\{0\} $ ou $ W=G_-/H_- $. Et aussi, pour le produit scalaire dans $ {\mathbb C}^n $ des matrices est invariant par l'adjoint, $ hxh^{-1} $, le produit scalaire choisit fait reference \`a la phrase du Gallot-Hulin-Lafontaine "any scalar product", ici c'est le plus simple, celui des matrices. En outre, ici, si on considere la matrice de la courbure de Ricci, $ Ricci $ c'est une matrice symetrique donc diagonalisable. Ici, pour le cas du projectif complexe, on prend le produit scalaire $ (M_1 \cdot M_2) $ avec $ M_1, M_2 $ deux matrices de la forme $ hxh^{-1} $ avec $ h, x, y $ trois matrices comme ci-dessus, cela revient \`a avoir le produit scalaire usuel $ <\bar \lambda \lambda \bar v^t \bar A^t Aw>_{{\mathbb C}} = <v|w> $, car $ A\in U(n) $ et $ \lambda \in U(1), |\lambda|=1 $. 

\smallskip

On regarde alors si l'application $ \tilde F $ du livre de Gallot-Hulin-Lafontaine est une isometrie avec les bons produits scalaire (ce qu'on a dit au debut de ce point avec $ [x] \in G/H=SU(n+1)/S(U(1)\times U(n)) $ et $ u \in P_n(C)=S^{2n+1}/S^1 $):

\smallskip

Soit $ \tilde F $, l'application du livre de Gallot-Hulin-Lafontaine; $ [x] \to \tilde [x.m]=q(x.m), [x] \in G/H=SU(n+1)/S(U(1)\times U(n)) $, $ \tilde [x.m] \in S^{2n+1}/S^1 $ et $ q: y \to \tilde [y] $, l'application du passage au quotient de $ S^{2n+1}\to S^{2n+1}/S^1 $. $ \tilde F $ est un diffeomorphisme. On prend alors $ m=z=(1,0,...,0) $ et on se place dans la carte de $ m $, $ (t_k)=(z_k/z_0) $, voir le Gallot-Hulin-Lafontaine dans la construction de la metrique de Fubini-Study de $ S^{2n+1}/S^1 $.

\smallskip

Alors $ T_m(S^{2n+1})=\{ (i\eta, \xi^1,..., \xi^n), \eta \in {\mathbb R}, \xi^j \in {\mathbb C}, 1\leq j \leq n \} $. Si on prend $ t\to x(t) $ un chemin dans $ SU(n+1) $ alors $ x'(t)=dx/dt \in T_e(SU(n+1)) $ est une matrice antisymetrique de trace nulle et $ (dx/dt) \cdot m= (0, v_1, v_2,...v_n)=(0,v) $ est orthogonal a $ m $ et $ i.m $, et donc, $ (dx/dt)\cdot m \in H_m $, le sous espace vectoriel qui permet de definir le produit scalaire dans $ S^{2n+1}/S^1 $, de meme pour $ (dy/dt) \cdot m = (0, w_1,..., w_n)=(0,w) \in H_m $. (Ne pas confondre $ H_m $ l'espace vectoriel ici dans la definition de la metrique de Fubini-Study et $ H_m $ l'espace isotrope dans la definition de $ G/H=G/G_m $). (On peut prendre $ x $ tel que $ dx/dt=x$ avec des termes centraux nuls, comme ci-dessus).

\smallskip

Comme par definition $ q $ est une isometrie sur $ H_m$, on a:

$$ <(dx/dt)\cdot m| (dy/dt)\cdot m>=<dq[(dx/dt)\cdot m]|dq[(dy/dt) \cdot m]>= $$

$$ = < (d/dt) [q(x.m)]|(d/dt)[q(y.m)]>=<(d/dt)[\tilde F(x)]|(d/dt)[\tilde F(y)]>= $$

$$ = <d\tilde F(dx/dt)|d \tilde F(dy/dt)>,$$

donc, $ \tilde F $ est une isometrie, ce qui permet de transporter la metrique de Fubini-Study.

\smallskip

De plus, le produit scalaire:

$$ <(dx/dt)\cdot m |(dy/dt)\cdot m >=<(dx/dt)\cdot m|(dy/dt)\cdot m>_{H_m}=<v|w>_{{\mathbb C}}= $$

$$ = [(dx/dt)\cdot (dy/dt)]=<(dx/dt)|(dy/dt)>, $$

est le produit scalaire sur le matrices et coincide avec $ (M_1 \cdot M_2) $ de ce qu'on a  dit precedemment et avec le produit scalaire dans $ {\mathbb C}^n=H_m $, de la metrique de Fubini-Study, si on prend par exemple $ dx/dt=x $ avec des termes contraux nuls. Le produit scalaire choisit sur l'espace homogene (de matrices) induit que, $ \tilde F $, une isometrie sur $ S^{2n+1}/S^1 $. C'est ce qu'on disait au debut de cette section sur les espaces homogenes, de maniere theorique, on choisit le produit scalaire de sorte que $ \tilde F $ soit une isometrie et donc on transporte la metrique. Ici, pour le cas du projectif complexe, le produit scalaire qu'on pris, fait de $ \tilde F $ une isometrie.***)

\smallskip

Concernant la courbure de Ricci et le fait qu'on a le projectif comme espace homogene d'Einstein, sans les calculs:

\smallskip

On a parl\'e de la notation $ . $ d'action de de groupe, comme le groupe de Lie agit par isometries, dans la derivation on obtient des termes du type, dans le tenseur de courbure :

$$ R(h'.X, h'.Y, h'.Z, h'.T)=R(X,Y,Z,T)$$

en composant \`a droite par l'inverse, on obtient l'application adjointe:

$$ R(h'.X.h'^{-1}, h'.Y.h'^{-1}, h'.Z.h'^{-1}, h'.T.h'^{-1}) = R(X,Y, Z,T) $$

ce qui veut dire que le tenseur de courbure est invariant par l'application adjointe $ h.x.h^{-1} $, donc la courbure de Ricci aussi et par invariance des sous-espaces propres, on a $ R=cg $ d'Einstein. Ceci est dit dans le livre de Besse, dans le cas des metriques $ G-$ invariantes, lorsque $ H, G$ agissent par isometries. Pour le cas particulier du projectif complexe, on savait par les calculs que c'est une vari\'et\'e d'Einstein, et par l'identification isometrique du Gallot-Hulin-Lafontaine, \`a $ SU(n+1)/S(U(1)\times U(n)) $, c'est une vari\'et\'e ou espace homog\`ene, une propri\'et\'e supplementaire du projectif complexe d'etre un espace homogene et d'Einstein.)

\smallskip

Le fait que la metrique soit invariante par l'adjoint, fait qu'elle invariante \`a droite et \`a gauche. C'est construit pour qu'elle soit bi-invariante. Ce qui induit des propri\'et\'es sur la courbure, car on a des isometries \`a droites et \`a gauche. C'est construit pour.

\smallskip

a) Les actions a droites et a gauches sont des isometries, soit $ s $ la metrique et $ L_g: x\to g.x $: $ L_g^*(s)=s \Rightarrow  L_g^*(R)=R $ voir le livre de Hebey pour l'expression de la courbure quand on a une isometrie Riemannienne. On obtient l'expression ci-dessus.

\smallskip

b) Ici pour le cas particulier du projectif, on a des matrices, l'application $ h\to hx$ avec $ h $ une matrice et $ x $ un vecteur, est lineaire, donc sa deriv\'ee est elle meme, d'o\`u le fait que ce qu'on a ecrit ci-dessus, $ h'=h $.

\smallskip

l'application: $ x(t)\to h.x(t)$, sa derivee est $ h.x'(t)=h.x$  avec $ x $ un vecteur.

\smallskip

Ce qui fait qu'ici pour le projectif complexe, on a $ h'.X.h'^{-1}=hxh^{-1} $ avec le fait que pour la courbure de Ricci, $ Ricci(hxh^{-1}, hyh^{-1})=Ricci(x,y) \, \forall h $, l'invariance de la courbure de Ricci par l'adjoint.

\smallskip

C'est une formule g\'en\'erale qu'on applique au cas particulier du projectif complexe et des matrices. $ SU(n+1)/S(U(1)\times U(n)) $.

\smallskip

( On a,  soit $ x $ un vecteur propre de $ R $ pour le valeur propre $ \lambda $, $ R(x)=\lambda x $, on ecrit:
$$ R(h.x.h^{-1},h.y.h^{-1})=R(x,y)=<y^t|Rx>=\lambda <y^t|x>=\lambda<(h.y.h^{-1})^t|h.x.h^{-1}>= $$
$$= <(h.y.h^{-1})^t|R(h.x.h^{-1})>, \, \forall h, \, \forall y $$
$$ \Rightarrow R(h.x.h^{-1})=\lambda (h.x.h^{-1}).$$

Donc, si on note $F_{\lambda} $ le sous-espace propre de $ R $ de valeur propre $ \lambda $, alors $ Ad(h)(F_{\lambda}) \subset F_{\lambda}, \, \forall h $, il v\'erifie la propri\'et\'e d'invariance. d'o\`u la conclusion.)

\smallskip

Donc, le projectif complexe vu du point de vue $ SU(n+1)/S(U(1)\times U(n)) $ comme espace homog\`ene et d'Einstein de constante $ s_0 $. Pour connaitre le signe de $ s_0 $, on utilise le Th\'eor\`eme du livre de Besse qui dit que si $ s_0 < 0$, il ne serait pas compact, ce qui n'est pas possible. Si $ s_0=0 $, il serait plat et donc Einstein et Conformement plat, par le r\'esultat du livre de Hebey, il serait de courbure sectionnelle constante , ce n'est pas possible, donc, $ s_0 >0 $.

\smallskip

Donc: concernant le caract\`ere localement conform\'ement plat et non localement conform\'ement plat de certaines vari\'et\'es:

\smallskip

1-les revetements Riemanniens dans la construction des metriques de courbures sectionelles constantes, comme les projectifs r\'eels, le Tore, les vari\'et\'es Hyperboliques.

\smallskip

2-les submersions riemanniennes, dans la construction des projectifs complexes.

\smallskip

3-les fibrations, et groupes de Lie, dans la constructions, des vari\'et\'es homog\`enes, puis par exemples, on compare, les metriques, comme pour le projectif complexe, qui est d'Einstein de constante $ s_0 >0 $. (La fibration est construite de telle maniere qu'on a une strucutre de vari\'et\'e lisse pour la vari\'et\'e quotient. La fibration sert dans la construction de la m\'etrique, les sections.).

\smallskip

(Par exemple, pour les revetements, il ais\'e d'avoir une structure de vari\'et\'e lisse sur le quotient. Un peu moins pour les submersions et les fibrations. De plus, dans la construction des metriques du quotient on s'assure que ca ne depend pas du point de la fibre).

\smallskip

4-On a les vari\'et\'es de courbures sectionnelles constantes, qui sont localement conformement plates. Le produit d'un cercle et d'une vari\'et\'e de courbure sectionelle constante, est plate.

\smallskip

5-Les sommes connexes, de vari\'et\'es conformement plates, peuvent etre plates.

\smallskip

6-Les sommes connexes, par le theoreme de D. Joyce, peuvent avoir une partie plate et une partie non plate et de courbure scalaire $ -1 $.

\smallskip

7-Pour construire une vari\'et\'e nonlocalement conformemtn plate, on construit des produit. On a l'a fait \`a partir des surfaces r\'eelles de dimension 2. Et en consid\'rant des produit  dont un n'est pas plat.

\smallskip

8-Les surfaces $ K3 $. La vari\'et\'e de W. Goldman en dimension 3 sur la quelle on peut mettre n'importe qu'elle metrique, elle sera non plate.

\smallskip

//////////////////////////////////////

\smallskip

8) La formule d'integration par parties est valable quand la vari\'et\'e est compacte sans bord et non orientable, voir le livre d'Hebey. Cela veut dire que dans les problemes variationels de degr\'e 2, il est possible de resoudre par la formulation variationelle un probleme elliptique comme dans le cas du Projectif en dimension 2, c'est fait dans le livre d'Aubin et ceci sur une variet\'e non orientable.

\smallskip

Par contre quand il sagit de fonction de Green, sur une variet\'e compacte, il est imperatif de supposer la veriet\'e orientable car, par exemple, pour la Parametrix, on a besoin d'utiliser l'integration par parties en dehors de petites boules, on se ramene a une variet\'e \`a bord, d'o\`u, il faut utiliser la formule de Stokes et donc une orientation.

\smallskip

///////////////////////////////////////

\smallskip

9) Pour les vari\'et\'es de dimension 3. On considere un Torus bundle $ M_{\phi} $ d'application $ \phi $ (voir le Hatcher pour la definition, $\phi $ est la monodromy map), il est orientable si et seulement si $ \phi \in SL_2({\mathbb Z}) $. William Goldman prouve quil existe un Torus bundle orientable ne possedant pas de structure conformement plate, $ \phi $ non periodique, (de valeur propres r\'eelles) (par exemple, $ \phi = \left ( \begin{matrix}  2 & 1 \\ 1 & 1 \end{matrix} \right) $  il n'y a  pas de metrique localemeent conformement plate, comme une variete lisse possede toujours une metrique Riemannienne $ g $, elle est dans ce cas non localement conformement plate, tenseur de Cotton-York, $ C \not \equiv 0 $. Ce theoreme dit  qu'il n'existe pas de metrique localemment conformement plate sur cette vari\'et\'e. 

\smallskip

1-Une structure est conformement plate s'il existe une metrique conforme \`a la metrique euclidienne. Cette definition est equivalente \`a celle qui dit qu'une structure conformement plate s'il existe un atlas dont les changements de cartes sont des applications de Mobius (compose\'ee de translation, rotation et inversion). Cette equivalence est demontr\'ee dans l'article de Kuiper (1949).

\smallskip

2-Cela veut dire que si on veut savoir si une vari\'et\'e a une structure conformement plate, il suffit d'etudier les tranformations conformes, donc le groupe conforme. C'est l'id\'ee de William Goldman, il caracterise le groupe fondamental et le groupe de transfomrmations conforme des varietes conformment plate (il y a une relation entre le groupe fondamental ete le groupe des transformations conformes). Si cette condition n'est pas satisfaite alors, la variete n'a pas de structure conformment plate. (le groupe fondamentale doit etre polycyclique...).

\smallskip

/////////////////////////////////////////////////////////

\smallskip

10) Sur le th\'eor\`eme de la masse positive, dans le cas plat et non plat, en dimension 3 et en dimensions $ \geq 3 $.

\smallskip

Voir aussi, dans le livre d'aubin ou la caracterisation de Kuiper de la structure conforme par une isometrie, sert \`a prouver le Theoreme de la masse positive dans le cas conformement plat. En effet:

\smallskip

a) Quand on a une vari\'et\'e complete, conform\'ement plate et simplement connexe, elle est conformement isometrique \`a une partie de la sphere dont on connait la masse $ m_0 =0 $, par la formule donn\'ee dans le livre d'Aubin, $ \alpha_0(P)=k\times \alpha_{{\mathbb R}^n}(P)=k\times 0=0 $.

\smallskip

b) On considere une vari\'et\'e $ M $ compacte, conformemnt plate. Alors le revetement universel $ \tilde M $, verfie les conditions du point a) precedent. Donc la masse est nulle.

\smallskip

c) Si de plus $ M $ n'est pas conform\'ement diffeomorphe \`a la sphere, le revetement universel contient plus de deux feuillets.(Par Kuiper une vari\'et\'e compacte connexe, conform\'ement plate plate et simplement connexe, est diffeomorphe \`a la sphere. Donc, ici, elle ne peut pas etre simplement connexe). Ce qui fait que l'ensemble not\'e $ W=\Pi^{-1}(\Pi(P)) $ dans le livre d'Aubin n'est pas reduit \`a $ \{P\} $. Par la construction de la fonction de Green minimale (not\'ee $ \tilde G_P $ de $ \tilde M $) et le principe du maximum, il existe une fonction $ \tilde H $, qu'on peut prendre ici pour $ G_Po\Pi $,($ \tilde g=\Pi^*(g) $ qui est conforme a $ \bar g = \Phi^*(g_0) $ donc de la forme $ u^{-1} Ho\Phi $...), telle que on ait $ \tilde H >\tilde G_{\tilde P} $ en $ P $. Or $ \tilde M =\Pi^*(M) $ et donc la fonction de Green de $ M $, verifie au moins l'in\'egalit\'e. $ \alpha_P (G_Po\Pi) > \tilde \alpha_{\tilde P}=0 $, car pour le revetement universel, d'apres le point a), la masse est nulle (et la distance se conserve localement, car on a un revetement donc une isometrie locale). Ainsi la masse de $ M $ en $ P $, $\alpha_P (G_P) >0 $, dans le cas o\'u $ M $ n'est pas conform\'ement diffeomorphe \`a la sphere. 

{\bf Remarques:} Sur la partie de la preuve du livre de Aubin (communication de M.Vaugon) concernant l'effet du changement conforme sur la masse $ \alpha $. On $ \tilde g =\phi^{4/(n-2)} g $ avec $ \tilde g $ et $ g $ des metriques conform\'ement plates $ \tilde d $ et $ d $. on se ramene au cas o\`u le facteur conforme est constant. La vari\'et\'e est localement conform\'ement plate. Donc la masse $ m $ ne depend pas des coordonn\'ees (Bartnik). Comme :

$$ \tilde \alpha(P)= m({\tilde G_P}^{4/(n-2)} \tilde g)= $$

$$= m(G_P\dfrac{g}{\phi(P)^{4/(n-2)}})= m({G_{P,\lambda}}^{4/(n-2)} g^{\lambda})=\alpha(P)/(\phi(P)^2). $$

Car le raisonnement dans le livre d'Aubin reste vrai so considere une metrique conforme du type $ g^{\lambda}=\lambda^{4/(n-2)} g $ avec $ \lambda $ un constante $ >0 $. Ici, $\lambda=\phi(P) $ et $ G_{P,\lambda}=G_P/\lambda^2 $.

\smallskip

Pour le Theoreme de la masse positive de Schoen-Yau en dimension 3, ils prouvent qu'il exsite une surface minimale complete ayant une propri\'et\'e (integrale de la courbure de Gauss positive) alors que par les calculs cette inetgrale est negative ou nulle. Pour cela ils utilisent le theoreme d'exsitence de surfaces minimales et leur proprietes (Formule de Gauss-Bonnet et Travail de Huber, inegalites isoperimetriques classiques et de Huber, la metrique est voisine de la metrique euclidienne ce qui fait qu'on peut approcher ces deux notions  (inegalites isoperimetrique classique)). En particulier, ils prouvent qu'elle est complete dans la premiere etape.(les ferm\'es born\'es sont compacts.

La positivit\'e de la courbure scalaire implique que la surface minimale construite a son integrales de Gauss Positive. Et puisque la metrique est "presque" eucllidienne, il y a conservation de l'inegalit\'e isoperimetrique et par la formule de Gauss-Bonnet, elle serait negative ou nulle. Ce qui n'st pas possible.

Il s'agit de construire une metrique conforme avec une propriet\'e de la courbure scalaire. 

Il s'agit de construire une surface minimale ayant une propri\'et\'e qui decoule de la courbure scalaire.

Cette surface minimale complete n'exsite pas , car cette propri\'et\'e n'est pas vraie, par le fait que la metrique est "presque" plate et que l'inegalit\'e isoperimetrique euclidienne se conserve et la fromule de Gauss-Bonnet.

{\bf Remarques:}

a) Ils prouvent que la masse $ m \geq 0 $, par l'absurde:

\smallskip

1-La masse $ <0 $ intervient dans la construction de la surface minimale complete. (raisonnement par l'absurde). (Ils prouvent qu'elle complete dans cette etape, car, en dehors des bouts, elle est dans un compact, et dans les bouts elle est entre deux plans et la metrique est equivalente \`a la metrique euclidienne.(les ferm\'es born\'es sont compacts)).

2-Une surface est dite de type fini, si $ \int |K| <+\infty $: c'est le cas ds la demonstration de la masse positive de Schoen Yau.

3-Quand on a une surface de type fini, si $ \chi(S')>0 $ (la caracteristique d'Euler-Poicar\'e $ >0 $, ce qui est le cas ici par Cohn-Vossen), alors $ S'$ est diffeomeorphe au plan. C'est ecrit dans Hulin-Troyanov(1992), premier theoreme, qui dit que, $ S'$ est diffeomorphe \`a une vari\'et\'e compacte $ S $ sans bord, sans un nombre fini de points. De plus il y a une formule de Gauss Bonnet, on a alors $ S'$ est la sphere moins un au plus un point d'apres cette formule de Gauss-Bonnet. Donc, c'est le plan. Il s'agit d'un diffeomorphisme entre $ S' $ et le plan.

On a d'abord une isometrie $ i $, puis une application conforme, $ F_1 $. la courbure de Gauss "se conserve". D'abord $ K_{i(S')}=Koi $, puis par l'application conforme on a un facteur $ |F_1'|^2 $, qui s'elimine dans le changement de variable. $ \int K=\int K_{F_1oi}=\int KoF_1oi |F_1'|^2 $, cette derniere inetgrale est la courbure de Gauss de $ F_1oi(S') $. Puis on applique le result d'Huber. (On considere $ i(S') $ au lieu de $ S'$, car on a une isometrie et la courbure de Gauss  "se conserve", $ \int_{S'} K = \int_{i(S')} Koi $ et $ ({F_1^{-1}})^*( g_{i(S')})=|F_1^{-1}|^2g_{{\mathbb C}} $, car 
$ F_1 $ est conforme, la metrique se transporte suivant un facteur de la metrique de $ {\mathbb C} $, on a des coordonnees isothermes globales, on peut appliquer ce que fait Huber (1957) ).

(Pour obtenir l'application conforme, on raisonne par rapport au revetement universel et la caracteristique d'Euler-Poincar\'e). 

4-Dans un article de Huber, il prouve que sur une surface ouverte, la formule de "Gauss-Bonnet", pour une suite de domaines:

$$ \int K=2\pi-\lim_{\sigma \to +\infty} L_{\sigma}^2/2A_{\sigma}. $$ 

5-La consturction de la surface est en fait une construction d'un courant (surface avec singualrit\'e, par une minimisation d'une fonctionelle relativement \'a une metrique), puis, voir le livre de Federer, la surface est reguliere en dehors du bord. $ \partial (N \cap S) =(\partial N \cap S) \cup (\partial S \cap N) $, par la formule de Stokes, $ S $ est loin du bord, car la courbure moyenne est positive (c'est une hypothese).

(Voir aussi , le livre de Sa Earp-Tubiana, sur le theoreme de Poincar\'e-Koebe et la calssification par le revetement universel).

Un theoreme d'uniformisation de Poincar\'e-Koebe, dit que une surface non compacte simplement connexe est conforme au disque unit\'e ou au plan. (voir l'article de Hulin-Troyanov(1992)). 

\smallskip

b) Ils prouvent que si la masse est nulle, $ m=0 $, alors l'espace est isometrique \`a $ {\mathbb R}^3 $. la preuve de ce point est assez claire dans l'article de Schoen-Yau, ils se ramenent au cas de masse $ <0 $ et de courbure scalaire $ \geq 0 $ par des changement de metriques conformes..., ce n'est pas possible, ou bien la courbure sclaire et de Ricci seraient nulles comme on est en dimension 3 $ Weyl \equiv 0 $, donc c'est la cas plat. Puis il faut utiliser l'argument sur le volume de Bishop: $ Ricci\geq 0 $ (ici, $ Ricci\equiv 0 $), la fonction $ f: r\to Vol_g(B_g(r))/Vol_{\delta}(B_{\delta}(r)) $ est decroissante sur la vari\'et\'e complete de depart $ M $, et tend vers $ 1 $ quand $ r\to 0 $ et $ r\to +\infty $ (par l'hypothese  que la vari\'et\'e est assymptotiquement plate, la metrique tend vers la metrique euclidienne quand $ r\to +\infty $), $ f $ est constante egale \`a 1 et le theoreme  de comparaison de Bishop dit dans ce que que $ M $ est isometrique a $ ({\mathbb R}^3, \delta) $. 

\smallskip

{\bf Ceci s'applique a une vari\'et\'e compacte sans bord:}

En effet, soit $ M $ cette vari\'et\'e compacte sans bord. $ P \in M $, en se palcant au voisinage de $ P $ et si considere les coordonnees geodesiques polaires alors la metrique est du type $ g_{ij}=\delta_{ij}+(termes)r^2+\ldots $, si on fait le changement de variable $ r=1/\rho $, alors :

1-on exclut le point $ P $ et la vari\'et\'e  devient asymptoique, car pour $ r=0 $ (en $ P $), ce n'est pas bien definit.

2-$ M-\{P\} $ ressemble a $ {\mathbb R}^n $.

3-on sait que les boules geodesiques sont de courbure moyenne $ <0 $ concave, voir la preuve dans un des livres sur la conjecture de Poincar\'e (Chow and al, par exemple). Mais le changement de variable $ r=1/\rho $ inverse le signe. Ces boules $B_P(r)-\{P\} $ sont les bouts du cas general de la masse positive les $ N_k $ qui sont convexes. Cette hypothese de courbure moyenne positives de la vari\'et\'e $ N= B_P(r_0)-\{P\} $ est verifiee dans le cas compact sans bord en se placant en un point $ P $ et changeant de variable.

4-C'est ecrit dans le livre d'Aubin, que la masse $ m=cte\times A =4(n-1) A $ avec $ A $ la partie reguliere de la fonction de Green. Donc elle est positive si $ m >0 $.

5-Reste a voir ce qu'est la m\'etrique $ g $ comme c dit dans le livre de Aubin, ils considerent $ G_P^{4/(n-2)} \times g \equiv r^{-4}g $ avec $ G_P $ la fonction de Green, et le changement de variable $ \rho=1/r $ ou $ r=1/\rho $ on a le comportement asymptotique, une vari\'et\'e asymptotique d'ordre 2 en coordonn\'ees geodesiques polaires.

(Dans le changement de variable $ \rho=1/r $, on ecrit, $\partial_{\rho}=-\dfrac{1}{r^2}\partial_r $, puisque c'est un changement de variable, on calcule alors $ g_{\rho \rho} $, de meme que pour $ g_{\rho \theta} $ et finalement $ g_{\theta \theta} $ qui se conserve, puiqu'il n'y a pas changement angulaire).

\smallskip

Donc, on voit que le cas compact sans bord se deduit du cas general. (Au voisinage du point $ P $ on est dans les bouts, mais ces vosinages sont des boules geodesiques dont on connait la courbure moyenne $ <0 $, et avec le changement de variables, pour avoir des coordonnes asymptotiques, ces boules deviennent de la vari\'et\'e $ N $ avec le bon signe sur la courbure moyenne $ >0 $.
 
On a: $ N=B_P(r_0)-\{P\} $, $ N_k=B_P(r)-\{P\}=B_P(1/\rho)-\{P\} $ et $ \partial N_k=S_P(r)=S_P(1/\rho) $ avec orientation invers\'ee. 

\smallskip

{\bf Remarque:} Sur la courbure moyenne de la sphere de rayon $ r $, $ S_P(r) $. On ecrit $ B $ est la seconde forme fondamentale:

$$ H= trace(B) $$

Mais, 

$$ B(e_i,e_j)=<\nabla_{e_i}(e_j)|N>, $$

avec $ e_i, e_j $ les elements de la base de l'espace tangent a la sphere et $ N $ la normale a la sphere.

On prend alors $ e_i, e_j \in \{\partial_{\theta } \} $ et on utilise la deriv\'ee covariante et le produit scalaire $ <e_i, N>=0 $, en permutant la deriv\'ee covariante, jusqu'a obtenir:

$$ trace(B)= trace [\nabla_N(<\partial_{\theta }|\partial_{\theta }>)]=\nabla_r (<\partial_{\theta }|\partial_{\theta }>)= $$

$$ = \dfrac{1}{r^2} \times r \times g^{\theta_i \theta_j}  g_{\theta_i \theta_j}+... = \dfrac{n-1}{r}+... $$

Donc (avec l'orientation positive), 

$$ H=\dfrac{n-1}{r}+...$$

(Aussi, en permuttant la derivation covariante dans $ <e_i|N>=0 $, on retrouve l'expression de la courbure moyenne (positive) ecrite dans l'article de Schoen-Yau, $ div(N)>0 $).

\smallskip

////////////////////////////////////////////////////////////////////////

\smallskip

{\bf Ceci concerne, les 3 articles de T.Aubin}:

\smallskip

1) Sur quelques probl\`emes de courbure scalaire. Journ. Functional.Analysis.2006.

\smallskip

2) D\'emonstration de la conjecture de la masse positive. Journ.Functional.Analysis. 2007.

\smallskip

3) The mass according to Arnowitt Deser Misner. Comptes rendus. Math. 2007.

\smallskip

-On est dans le cas positif $ \Rightarrow $ par un changement de metrique conforme on a la courbure scalaire $ R=1>0 $ et int\'egrable $\Rightarrow $ on est dans le cas o\`u, autour d'un point, par un changement de variable $ r=1/\rho, r >0, \rho >0 $, la vari\'et\'e est asymptotiquement plate d'orde $ \tau > (n-2)/2 $. Dans le developpement limit\'e de la metrique ($ g={\mathcal E}+ O(r^{2+\omega}) $), ceci correspond au cas $ \omega > (n-6)/2 $. Ceci correspond au cas ou la masse $ A $ en $ P $ existe et finie. La fonction de Green s'ecrit comme $ G(P,Q)=1/r^{n-2} + A+ O(r), r=d(P, Q)>0 $. L'article de Aubin sur la masse positive est clair. Il utilise une condition sur la fonction de Green, sur la courbure scalaire, $ \alpha >n-4 $, qu'il obtient dans le premier article, "Sur quelques problemes de courbure scalaire", l'idee est de considerer, les parties principles des fonctions de Green, il faut soustraire des produits de convolutions, il elimine, petit a petit, les parties singulieres de haut degr\'es li\'es \`a la partie principale de la courbure scalaire.

\smallskip

La condition sur la courbure scalaire, $ \alpha >n-4 $, permet de montrer que la masse $ A \geq 0 $, ceci est obtenu en considerant un probleme variationnel, dans l'article de Aubin.T. Sur la masse positive.

\smallskip

L'identit\'e sur la courbure de Ricci et un argument de R.Schoen, permet de dire que la masse $ A>0 $ sauf sur la sphere.

\smallskip

-Le premier article concerne la courbure scalaire, il evalue la partie principale de la courbure scalaire et l'associe \`a la fonction de Green. C'est en relation avec la determination des parties principales singulieres de la fonction de Green. Il faut soustraire des produits de convolutions successivement.

Par exemple:

$$ G=r^{2-n}-fr^{\mu+4-n}+H, $$

Alors,

$$ \Delta_{dist} G+RG=\delta_P, \,\, \Delta_{dist} (r^{2-n})=\delta_P,\,(\sqrt {|g|}=1), \,\, \Delta (-fr^{\mu+4-n})=\bar Rr^{2-n}=h_1r^{\mu+2-n}, $$

Donc,

$$ \Delta_{dist}H+RH=g_1=O(r^{\mu+3-n}), $$

Soit, $ H_1=G*g_1=\int_{B_P(r)} G(x,y) g_1(y)dV_g(y) $, alors: par Giraud, $ H_1=O(r^{\mu+5-n}) $, et,

$$ \Delta_{dist} H_1+RH_1=g_1, $$

Donc,

$$ \Delta_{dist} (H-H_1)+R(H-H_1)=0, $$

Donc, par les theoremes de regularit\'e: $ H-H_1\in C^{\infty} $, avec, $ H_1=O(r^{\mu+5-n}) $, donc: $ H=O(r^{\mu+5-n}) $.

\smallskip

Il d\'emontre une identit\'e de Pohozaev astucieuse et la relie \`a la fonction de Green. Comme la masse $ A >0 $, la condition pour que l'ensemble des solutions de l'equation de Yamabe soit compact, est verifi\'ee (premiere condition).

Concernant la compacit\'e de l'ensemble des solutions de l'eq. Yamabe: il y a un probleme, il faut prouver que les blow-up sont isol\'es simples, pour cela, il faut utiliser Pohozaev, or dans Marques.F.C, le terme $ P(r,h) \geq 0, M_iu_i\to h \in C^2(B_r(x_0)-\{x_0\}) $, $ P(r,h) $ est positif seulement si la dimension $ 4 \leq n\leq 7 $, on n'arrive pas \`a prouver que les blow-up sont isol\'es simples, car dans l'identit\'e de Pohozaev le terme $ P(r, h) \geq 0 $ seulement en dimensions $ 4\leq n\leq 7 $ et pour $ n\geq 8 $ on ne sait pas. Il se peut que la condition $ \omega >(n-6)/2 $ permette d'avoir $ P(r,h) \geq 0 $, c'est ce que Marques.F.C. fait quand $ Weyl(x_0)=0 $. Il est possible que le degr\'e d'annulation du tenseur de Weyl ($ \omega > (n-6)/2 $ et $ g ={\mathcal E}+O(r^{2+\omega}) $ et ordre de $ R =\mu \geq \omega $), implique que $ P(r,h)\geq 0 $.(et apres on a $ h $ singuliere en $ x_0 $ et l'appendice de YY.Li-M.Zhu, permet de dire que $ h $ est reelement une fonction de Green, c'est a dire que c'est une solution au sens des distributions. Au depart rien ne dit que la fonction singuliere $ h $ est une fonction de Green). Ceci pour la compacit\'e des solutions de l'eq. de Yamabe.

C'est possible en dimension $ n\geq 8 $:

1- Il faut utiliser la formule de Pohozaev dans l'article d'Aubin. Le terme principal est celui qui contient la courbure scalaire: $ \int (2R+\rho \partial_{\rho}R) u_i^2 $.

2-Ici les estimations pour prouver qu'on a des blow-up isol\'es simples, pour l'eq. de Yamabe est equivalent au probleme de la courbure scalaire prescrite sur la sphere, de YY.Li. 1995. Pour YY.Li, le terme principal est $ \int (x\cdot \nabla K) u_i^{2n/(n-2)} $. Pour le probleme de Yamabe et avec la formule de Pohozaev d'Aubin, c'est le terme de la courbure scalaire: le point 1) precedent.

3-Comme YY.Li sur la sphere, il faut estimer le terme: $\int r^{\mu} u_i^2 $, $ \mu $ est l'ordre de $ R $. Comme $ \omega >(n-6)/2 $ et la fonction de Green s'ecrit $ G=1/r^{n-2} + A + O(r) $, alors $ \mu >n-4 $. Donc, $ \mu \geq n-3 $. 

Le procede de YY.Li, pour passer de blow-up isol\'es aux blow-up isol\'es simples, revient \`a estimer la quantit\'e: $ \int r^{\mu} u_i^2 $. Or en utilisant l'estim\'ee $ C^0 $: $ M_i u_i\leq C r^{2-n}$, on obtient:

$ \int r^{\mu} u_i^2 \leq [\int_{r_i}^{\epsilon} (r^{\mu-n+3}) dr]/M_i^2= [o(1)+O(\epsilon)]/M_i^2$. Ce qui implique que la quantit\'e de la formule de Pohozaev $ M_i^2(\int (R+\rho\partial_{\rho} R) u_i^2) \to 0 $ quand $ i\to +\infty $ et $ \epsilon \to 0 $.

4-Ce qu'on a not\'e $ u_i $, peut etre la fonction $ u_i $ ou des fonctions $ v_i $, $ \xi_i $ obtenue \`a partir de $ u_i $ ou son rescaling, dans le procede de passage de blow-up isol\'es \`a blow-up isol\'es simples de YY.Li. L'eq. de maniere generale ne change pas \`a des facteurs, petits, pr\'es.

\smallskip

Le premier article, une partie du raisonnement est en relation avec une identit\'e sur la courbure scalaire: theoremes 5 et 6. Madani Farid a publi\'e un article pour justifier ces theoremes. L'idee est de faire des changement de metriques conformes de telle maniere qu'une certaine quantit\'e sur la courbure scalaire ait le bon signe. Puis il l'applique \`a la compacit\'e des solutions de l'eq. de Yamabe.(Le cas: $ \omega \leq (n-6)/2 $).

\smallskip

Aubin.T. a voulu traiter tous les cas possibles: $ \omega >(n-6)/2 $ et $ \omega \leq (n-6)/2 $. La note des compes rendus Math, le 3eme article, permet de dire qu'on peut considerer le cas o\`u $ \omega >(n-6)/2 $.

\smallskip

L'article sur la masse positive est assez detaill\'e et clair, $ \omega >(n-6)/2 $:

\smallskip

Dans l'article de la masse positive, Aubin. T. parle de changement de metrique conforme du type: $ \tilde g = g_{\epsilon}= \Psi_{\epsilon}^{4/(n-2)} g $ avec:

$ \Psi_{\epsilon} \equiv 1 $ dans $ B_{r_1}(P) $: la courbure scalaire est: $ R_{\epsilon}=R $ et $ g_{\epsilon}=g $.

\smallskip

L'op\'erateur, $ \Delta + R $ est coercif, on peut r\'esoudre des probl\`emes de Dirichlet et utiliser le principe du maximum.

\smallskip

On r\'esout: $ \Delta \Psi_{\epsilon} + R\Psi_{\epsilon} =0, \tilde R =\tilde R_{\epsilon}=0 $, dans $ B_{r_1+\epsilon}(P)-B_{r_1}(P) $, avec les conditions de Dirichlet, $ \Psi_{\epsilon}=1 $ sur $ \partial B_{r_1}(P) $ et $ \Psi_{\epsilon}=c \to +\infty $, sur $ \partial B_{r_1+\epsilon}(P) $.

\smallskip

Et, $ \Psi_{\epsilon}\equiv c\to +\infty $, sur $ V-B_{r_1+\epsilon}(P) $. La courbure scalaire est en $ \tilde R= R/c^{4/(n-2)} = o(1) r_1^{\alpha}, \alpha >n-4 $.

\smallskip

1) On obtient, une fonction Lipschitzeinne: $ 1\leq \Psi_{\epsilon} \leq c \to +\infty $ (principe du maximum et principe du minimum: pour avoir $ \Psi_{\epsilon} \leq c $, on utilise le principe du maximum. Pour avoir $ \Psi_{\epsilon} \geq 1 $, il faut utiliser le principe du minimum. $ Lw_1\leq 0 $ et $ \liminf_{x\to \partial \Omega} w_1 \geq 0 $, alors, $ w_1 $ v\'erifie le principe du minimum, si et seulement si, $ -w_1 $ v\'erifie le principe du maximum, ce qui est vrai ici. Le principe du minimum a \'et\'e utiliser par Li-Zhang, en dimensions 3 et 4). On a une metrique Lipschitzienne. On doit alors utiliser les notions de Th\'eorie de la mesure et de la Th\'eorie de la mesure g\'eom\'etrique, pour une metrique riemannienne Lipschitzienne et une connexion de Levi-Civita lipschitzienne.

\smallskip

2) Notion de derivation p.p: derivation presque partout. On peut d\'eriver 2 fois donc, on peut considerer les probl\`emes elliptiques d'ordre 2. Les fonctions sont au moins dans $ H^1(M) \cap C^0(M) $ ou $ H^1_{1+\epsilon}(M) \cap C^0(M-\{P\}),\epsilon >0 $. Les objets sont tels que, $ f $ et $ \partial^k f $ derivables presque partout et $ f\in L^{\infty}$, $ \partial^kf \in L^{\infty} $. On peut considerer alors, les objets, comme au sens regulier, sauf qu'il y a la notion de presque partout, on enleve $ \partial B_{r_1}(P) $ et $ \partial B_{r_1+\epsilon}(P) $. Par exemple, les exponentielles, sont definies avant d'atteindre ces 2 bords, on peut alors definir $ \tilde \partial_r $ et $ \tilde \partial_{\theta} $. 

\smallskip

3) La distance est Lipschtzienne. La parametrix est definie avant d'atteindre les 2 bords, et, est \`a support compact, ce qui permet de definir la fonction de Green: il faut considerer: $ \int_{M-\{\partial B_{r_1}(P)\cup \partial B_{r_1+\epsilon}(P)\}} \Gamma_i(P,R) H(R,Q)dV(R) $, ce qui permet de deriver par rapport \`a $ Q $ et faire l'integration par parties en $ Q $ dans $ W^{1,1+\epsilon} $, par Fubini, Fubini-Tonnelli, voir ce qu'on a dit sur la fonction de Green, point (24) ci-apres.

\smallskip

4) Alors que, pour une carte $ (\Omega, \phi) $ de $ M $, les $ \partial_k $, sont definis partout et sont Lipschitziens. On peut d\'eriver 2 fois, donc, on peut consid\'erer les probl\`emes elliptiques d'ordre 2. Voir le livre de Gilbarg-Trudinger, pour les problemes avec coefficients Lipschtziens, chapitre 8, les solutions sont au moins $ H^2_2 $ ou $ H^2_{1+\epsilon} $ ou $ W^{2,1+\epsilon}$, th\'eor\`eme 8.8. On a aussi le fait que la metrique est r\'eguli\`ere dans $ B_{r_1}(P) $.

\smallskip

5) Sur $ V-B_{r_1+\epsilon}(P) $:

\smallskip

En se placant, dans la carte exponentielle, en $ P $ de $ g $, suivant une g\'eodesique, qui rencontre les trois morceaux ($ P, Q_{r_1}, Q_{r_1+\epsilon}, Q $). On a dans $ V-B_{r_1+\epsilon}(P) $, par un raisonnement sur les chemins et chemins minimisants: $ \tilde r \geq r_1+\epsilon+ \tilde t $, avec: $ \tilde t= \tilde d(Q, Q_{r_1+\epsilon})=c^{2/(n-2)} d(Q,Q_{r_1+\epsilon})= c^{2/(n-2)} (r-r_1-\epsilon)=c^{2/(n-2)} t $, $ t= r-r_1-\epsilon=d(Q,Q_{r_1+\epsilon}) $, on ecrit alors: 

$$ \int_{\tilde t=0}^{+\infty} \tilde R (\delta^{2(n-2)}/\tilde r^{2n-4}) c^{2n/(n-2)} \sqrt {|g|} t^{n-1} dt= o(1) \delta ^{2(n-2)} r_1^{\alpha +4-n}. $$

avec, $ r=d(P,Q)$ la distance pour la metrique $ g $, et, $ \tilde r=\tilde d(P,Q) $, la distance pour la m\'etrique $ \tilde g=g_{\epsilon} $. On a aussi, $ r\leq \tilde r \leq c^{2/(n-2)} r=\tilde t+ c^{2/(n-2)}(r_1+\epsilon) $.

\smallskip

Les coordonn\'ees d'un point $ Q \in V-B_{r_1+\epsilon}(P) $,  sont: $ \exp_P(t\theta), t=d(Q,Q_{r_1+\epsilon}), \theta \in {\mathbb S}_{n-1}, y=t\theta $. Au lieu de considerer la distance $ r=d(P,Q) $, on considere $ t=d(Q,Q_{r_1+\epsilon}) $. On change de parametrisation du chemin ou de la g\'eod\'esique, on fait un d\'ecalage, l'origine se trouve sur la sph\`ere $ \partial B_{r_1+\epsilon}(P) $ au lieu de $ P $. On a la meme chose avec un d\'ecalage. On a: $ \partial_{\theta}=t \partial_y $ et $ \bar g_{\theta_i, \theta_j}= t^2 g_{\theta_i,\theta_j}$. On obtient la meme chose que les coordonn\'ees polaires usuelles sauf qu'on part de $ t=0 $, $ t=d(Q,Q_{r_1+\epsilon}) $.

\smallskip

On s'int\'eresse \`a l'int\'egrale: $ \frac{4(n-1)}{n-2} \int_{V-B_{r_1+\epsilon}(P)} |\nabla \tilde \gamma |^2 d\tilde V $.

\smallskip

La m\'etrique ne change pas dans $ B_{r_1}(P) $, la masse ou le $ A $ pour $ g $ est \'egale \`a la masse, ou la $ \tilde A $ pour $ \tilde g=g_{\epsilon} $. Le $ \tilde \gamma $ est le $ \gamma $ pour la metrique $ \tilde g= g_{\epsilon} $.

\smallskip

Aubin.T, fait un calcul precis de cette integrale: il fait un changement de metrique conforme de sorte que  l'integrale moyenne: $  \oint_{\partial B(r)} \sqrt{|g|} =1 $, pour estimer num\'eriquement cette int\'egrale, il trouve avec, les notations ci-dessus, que:

$$ \frac{4(n-1)}{n-2} \int_{V-B_{r_1+\epsilon}(P)} |\nabla \tilde \gamma |^2 d\tilde V <4(n-1)\omega_{n-1} \delta^{2(n-2)} (r_1+\epsilon)^{2-n} $$

Ce terme va etre inferieur au terme explicite de l'integrale sur $ B_{r_1}(P) $, voir ci-apres.

\smallskip

6) Sur $ B_{r_1+\epsilon}(P)-B_{r_1}(P) $: $ \tilde R=0 $, le $ \epsilon \to 0, \epsilon >0 $, sert \`a controler la quantit\'e: $ \int_{r_1}^{r_1+\epsilon} |\nabla \tilde \gamma |^2 $.  

\smallskip

On obtient: $ \frac{4(n-1)}{n-2}\int_{r_1}^{r_1+\epsilon} |\nabla \tilde \gamma |^2 =o(1)\delta^{2(n-2)} $. 

\smallskip

On a: $ r\leq \tilde r \leq c^{2/(n-2)} r $. Le $ \tilde \gamma $ est le $ \gamma $ pour la metrique $ \tilde g= g_{\epsilon} $.

\smallskip 

7) Sur $B_{r_1}(P) $: $ R=O(1)r^{\alpha}, \alpha >n-4 $: $ \sqrt{ |g|} =1 $,

$$\int_{\delta}^{r_1} R (\delta^{2(n-2)}/r^{2n-4}) r^{n-1} dr= -\delta^{2(n-2)}( \delta^{\alpha +4-n})+ O(1) \delta^{2(n-2)} r_1^{\alpha+4-n}. $$

L'int\'egrale: $ \frac{4(n-1)}{n-2} \int_{B_{r_1}(P)-B_{\delta}(P)} |\nabla \tilde \gamma|^2 d\tilde V $: (il fait une majoration pr\'ecise):

$$ \frac{4(n-1)}{n-2}\int_{B_{r_1}(P)-B_{\delta}(P)} |\nabla \tilde \gamma|^2 d\tilde V = 4(n-1)\omega_{n-1}\delta^{2(n-2)} \int_{\delta}^{r_1} \dfrac{1}{r^{n-1}} dr =$$

$$ = 4(n-1)\omega_{n-1}\delta^{2(n-2)} [-r^{2-n}]_{\delta}^{r_1}= 4(n-1)\omega_{n-1} \delta^{2(n-2)}(\delta^{2-n}-r_1^{2-n}), $$

Donc,

$$ \frac{4(n-1)}{n-2} \int_{B_{r_1}(P)-B_{\delta}(P)} |\nabla \tilde \gamma|^2d\tilde V \leq  4(n-1)\omega_{n-1} (\delta^{n-2}-\delta^{2(n-2)}r_1^{2-n}). $$

Finalement, en sommant les trois int\'egrales, on obtient:

$$ \frac{4(n-1)}{n-2}\int_{V-B_{\delta}(P)} |\nabla \tilde \gamma|^2d\tilde V \leq  4(n-1)\omega_{n-1} (\delta^{n-2}+o(1)\delta^{2(n-2)}). $$

\smallskip

////////////////

\smallskip

Concernant les estimations avec la fonctionnelle. On a:

\smallskip

$ f $ la fonction minimisante et $ \psi=k\delta^{n-2}G $, v\'erifient l'equation d'Euler, au moins au sens faible:

$$ \Delta f+Rf=\Delta \psi+R\psi=0 $$

et au bord: $ f=1, \psi=1+k(H(\delta,\theta)-A)\delta^{n-2}=1+O(\delta^{n-1}) $.

\smallskip

Alors, par le principe du maximum (l'operateur $ \Delta +R $ est coercif) et les estimations elliptiques:

$$ f-\psi=O(\delta^{n-1}), \,\, \nabla (f-\psi)=O(\delta^{n-2}), \,\, {\rm sur} \,\, \bar B_{2\delta}(P)-B_{\delta}(P).$$

Il se ramene au calcul avec $ \psi $. Le support de $\chi$ est dans $ \bar B_{2\delta}(P)-B_{\delta}(P)$.

$$ J(\tilde f)=J(f)+ f\cdot \chi+ \chi\cdot \chi =J(f)+(f-\psi)\cdot \chi+\psi\cdot \chi+\chi\cdot \chi, $$

Il se ramene au calcul avec $ \psi $ et \`a un calcul similaire \`a celui avec $ \psi $.

\smallskip

L'autre calcul de $ J(\psi) $, il utilise l'equation d'Euler de $ G $ et la formule de Stokes, un calcul explicite sur le bord, avec la d\'eriv\'ee normale du terme $ \frac{1}{r^{n-2}} $, $ \partial_r (\frac{1}{r^{n-2}})=-\frac{(n-2)}{r^{n-1}} $. (Il fait ce calcul pour eliminer le terme avec la courbure scalaire).

\smallskip

Ce calcul avec la formule de Stokes: $ \int_{\partial B_{\delta}(P)} \partial_{\nu} h \cdot h $, sert aussi pour ramener le calcul de $ J(\psi) $ au calcul de $ J(\phi) $ avec $ \phi = \psi+ \chi $, comme c'est ecrit dans l'article de Aubin.T. puis on fait la meme chose avec $ J(\psi-f) $ et $ J(f)$, en remarquant qu'au bord, les fonctions sont des perturbations de $ \chi $, comme dans l'article de Aubin.T: au bord, on a: $\psi $ ou $ f $, $ + $ une perturbation du type $ \chi $, fonction en $ O(\delta^{n-1}) $ et \`a support compact, car les fonctions prennent des valeures voisines de celles prises au bord.

\smallskip

On ecrit: equation d'Euler de $ \psi $:

\smallskip

$$ J(\psi)=\int_{\partial B_{\delta}(P)} \partial_{\nu} \psi \cdot \psi= \int_{\partial B_{\delta}(P)} \partial_{\nu} (\psi-f+f) \cdot (1+\chi) =$$

$$ = \int_{\partial B_{\delta}(P)} \partial_{\nu} (\psi-f) \cdot \psi + \int_{\partial B_{\delta}(P)} \partial_{\nu} f \cdot f + \int_{\partial B_{\delta}(P)} \partial_{\nu} f \cdot \chi$$

avec, un calcul direct sur la partie superficielle, $ |\partial B_{\delta}(P)|=O(\delta^{n-1}) $ et les estimations elliptiques $ \nabla (\psi-f)=O(\delta^{n-2}) $:

$$ \int_{\partial B_{\delta}(P)} \partial_{\nu} (\psi-f) \cdot \psi=O(\delta^{2(n-2)+1}), $$

et, equation d'Euler de $ f $

$$\int_{\partial B_{\delta}(P)} \partial_{\nu} f\cdot f=J(f) $$

et, equation d'Euler de $ f $,

$$ \int_{\partial B_{\delta}(P)} \partial_{\nu} f \cdot \chi= \int \nabla f \nabla \chi+ \int R f \cdot \chi = $$

$$ {\rm ici \,\,estimations \,\, elliptiques} \,\, = \int \nabla (f-\psi) \cdot \nabla \chi + \int R (f-\psi) \cdot \chi + $$

$$ + \int \nabla \psi \cdot \nabla \chi+ \int R\psi \cdot \chi = {\rm calcul \, de \, Aubin .T}=O(\delta^{2(n-2)+\beta}), \beta >0. $$

Donc,

$$ J(\psi)=J(f)+O(\delta^{2(n-2)+\beta}), \beta >0. $$

Le support de $ \chi $ est dans $ \bar B_{2\delta}(P)-B_{\delta}(P) $.

\smallskip

Apr\'es, on utilise le calcul precis de $ J(\psi) $ et $ J(f)=m\leq J(\gamma) $.

\smallskip

///////////////////////////////////////////////////////////////////////

\smallskip

11) Conernant les sommes connexes, $ M\sharp N $ de deux vari\'et\'es $ M, N $,  elles sont definies a partir d'une relation d'equivalence qui donne une vari\'t\'e quotient Topologique, c'est \`a dire $ C^0 $, \`a l'ordre $ 0 $ il n'y a pas de probleme, car par continuit\'e le recollement se fait naturellement, mais des qu'on cherche une structure  differentiable ou $ C^k, k \geq 1 $, c'es plus compliqu\'e, dans ce cas on utilise le Disc Theorem de Palais, qui dit qu'on peut trouver une isotopie, c'est \`a dire une homotopie lisse entre les bords et donc le recollement se fait de maniere lisse.

\smallskip

Remarquons que endehors de la sphere ou se fait le recollement, la vari\'et\'e est egale \`a $ M $ ou $ N $. 

\smallskip

Lorsque on met une metrique sur $ M $ et $ N $ alors par des functions cutoff on peut mettre une metrique, par recollement, sur $ M \sharp N $ et elle coincide avec celles de $ M $ et $ N $ en dehors des boules de recollement. voir par exemple les articles de Dominic Joyce (2001) et Mazzeo (1995).

\smallskip

Dans le cas Riemannien, en utilisant les cartes exponentielles en $ p \in B(p,2\epsilon) \subset M $ et $ q \in B'(q,2\epsilon) \subset N $, on a deux injections  naturelles $ i_1, i_2 $ et on doit faire le recollement entre $ Id_{B(p,2\epsilon)} $ et $ i_2oi_1^{-1} $. On definit une relation dequivalence en partant de couronnes $ B(p,2\epsilon)-B(p,\epsilon) $ , $ B(q, 2\epsilon)-B(q, \epsilon) $ vers les varietes $ M $ et $ N $. On choisit une carte pour $ N $ renversant l'orientation et puis on compose par une inversion et on obtient $ i_2 $.

\smallskip

alors on recolle les deux applications de cartes et le changement de cartes, $ i_2oi_1^{-1} $ est lisse, ce qui permet le passage de l'atlas de $ M $ a celui de $ N $ de maniere lisse. (Les changements de cartes sont lisses de $ M $, $ N $ et $ M-B(p,\epsilon)\cup N-B'(q,\epsilon) $, par $ i_2oi_1^{-1}$.)

\smallskip

Comme cette carte permet le passage au deux atlas de $ M $ et $ N $, alors $ M\sharp N $ est orientable. 

\smallskip

On a: $ X=M\setminus B(p,\epsilon) \cup N\setminus B'(q,\epsilon) $ est compact, le graphe de la relation $ Gr=\{ (x,x), x\} \cup \{(x, i_2oi_1^{-1}(x)), x\} $ est ferm\'e, donc, $ X\setminus \sim $ est s\'epar\'e. Soit $ O=B(p,2\epsilon)\setminus B(p,\epsilon) \cup B'(q,2\epsilon)\setminus  B'(q,\epsilon) $, et $\pi $ l'application de passage au quotient. Alors $ \pi (O) $ est un ouvert de $ X\setminus \sim $ car $ \pi^{-1}[\pi(O)]=O $. 

\smallskip

Soit : $ \Phi $ l'application de recollement de cartes, definit precedemment sur $ \pi(O) $, alors $ \Phi $ est bien definit et est un homeomorphisme sur son image qui est un ouvert de $ {\mathbb R}^n $, par le theoreme d 'invariance du domaine de Brower (car $ i_2 $ et $ i_1 $ se recollent au bord qui est une sphere, par definition). C'est une carte, compatible avec l'orientation car, les changement de cartes \`a droite et \`a gauche sont l'identit\'e. 

$$ \Phi([x])=i_1(x),\,\, x\in B(p,2\epsilon)\setminus B(p,\epsilon), \,\, \Phi([x])=\dfrac{i_2(x)}{|i_2(x)|^2}, \,\, x\in B'(q,2\epsilon)\setminus B'(q,\epsilon), $$

avec , $ [x]=\pi(x) $, l'application de passage au quotient. $ i_2 $ inverse l'orientation de $ N $. Comme on compose avec une inversion, au final, on a des jacobiens de determinant positifs, de plus, sur la sphere, les orientations sont invers\'ees, ce qui permet un recollement correct (car sur le premier bout on a $ i_1 $ avec la bonne orientation et sur le deuxieme bout on a $ i_2 $ qui inverse l'orientation, et sur la sphere $ |i_2|=1 $).

\smallskip

On voit que $ \Phi o \pi $ est continue par definition et donc $ \Phi $ est continue.

\smallskip

On a bien, $ (\Phi, \pi(O)) $ une carte de recollement.

\smallskip

Elle permet le passage au deux atlas. Comme $ M $ et $ N $ sont orientables, $ M\sharp N $ est orientable.

\smallskip

Avec $ (\Phi, \pi(O)) $ on a une carte de recollement, et,

$$ \Phi o[i_1]^{-1} =\Phi o[i_1^{-1}] =Id_{B(p,2\epsilon)\setminus \bar B(p,\epsilon)}, $$

$$ \Phi o[i_2/|i_2|^2]^{-1} =\Phi o[(i_2/|i_2|^2)^{-1}] =Id_{B'(q,2\epsilon)\setminus \bar B'(q,\epsilon)}, $$

On voit que les changements de cartes sont lisses et conservent l'orientation.

\smallskip

En considerant les metriques $ g_M, g_N $ et $\Phi^*(\delta_{{\mathbb R}^n}) $, par des fonctions cutoff, on a par recollement, une metrique sur $ M\sharp N $:

$$ g_{M\sharp N}= \chi_1 g_M+\chi_0 \Phi^*(\delta_{{\mathbb R}^n})+\chi_2 g_N, $$

par exemple.

\smallskip

Pour voir que sur les bords, les inversions sont invers\'ees, on considere un observateur $ a $ de $ M $ tel que le repere $ (a, T) $ est direct avec $ T $ une orientation du bord. on deplace par une rotation, ou selon un chemin circulaire, vers un observateur $ b $ de $ N $ alors le repere $ (b, T') $ est direct si et seulement si $ T'=-T $.

\smallskip

D'autre part, on voit que dans la carte de recollement, on a compos\'e avec une inversion. Si $ M $ et $ N $ ont des structures locallement conform\'ement plates, par Kuiper, il existe pour chaque vari\'et\'e un atlas dont les changment de cartes sont des applications de Mobius. En utlisant la carte de recollement, on obtient un atlas sur $ M\sharp N $ form\'e d'applications de Mobius. Donc, sur la somme sonnexe on a une structure lcoallement conform\'ement plate. (ceci est dit dans la papier de Kulkarni).

\smallskip

12) Concernant la classification des surfaces. Cela peut se faire grace aux fonctions de Morse. Par recurrence sur le nombre de point critiques. 

\smallskip

1- On verifie que sous certaines hypotheses d'orientabilit\'e des bords et quand on a deux homeomorphismes entre deux ensembles, alors on a un homeomorphismes entre les recoll\'es.

\smallskip

2- On utilise une recurrence sur le nombre de points critiques des fonctions de Morse.

\smallskip

On note $ g_1, g_2 $ les applications de recollement, dans le Gramain, $ g_1=f $, $g_2=f' $. Ici, dans le cas des recollements, on prend, $ g_1=id $ et $ g_2=id $, les recoll\'es sont les ensembles de depart ( en fait ils sont diff\'eomorphes aux ensembles de depart et $ f_*=id, f'_*=id $, on a inversion des orientations quand on recolle, et, $ h_*=id, h'_*=id $, $ h $ et $ h' $ conserve l'orientation et l'inversion de l'orientation, respectivement). Les ensembles $ A $ et $ A' $, sont des ensembles de niveau pour une fonction de Morse $ C^{\infty} $, donc ce sont des vari\'et\'es compactes $ C^{\infty} $, il y a un nombre fini de composantes connexes, sinon (en se placant sur le bord, vecteur tangent et vecteur rentrant, on aurait un point critique dans l'intersection de 2 lignes de niveaux), il y aurait un point critique.

\smallskip

 Pour l'orientabilit\'e des bords, $ A $ et $ A' $ sont (hom\'eomorphes implique diff\'eomorphes en dimension 1) diff\'eomorphes \'a des cerles et, \'a droite (pour les sommes connexes), les orientations sont invers\'ees, donc le diff\'eomrphisme conserve l'inversion d'orientation, son Jacobien est positif. D'o\`u, l'application $ k_0 $ conserve l'orrientation. 
Ce n'est pas toujours vrai, car on n'a pas un diffoemorphisme global $h $ et $ h'$, par contre il y a une definition par l'homologie de l'orientation. Voir ci-dessus.

\smallskip

(on pouvait le voir en remarquant qu'on inverse deux fois l'orientation, donc au final on conserve l'orienation pour $ k_0 $).

\smallskip

Par exemple, on suppose que toute surface orientable ayant $ p \leq q-1 $ point critiques avec un maximum et un minimum et $ p-2 $ points critiques d'indice $ 1 $ est homeomorphe  \`a $ {\mathbb T}_p $. On peut commencer la recurrence pour $ p=0 $, la fonction de Morse a un maximum et un minimum et les ensembles de niveaux sont homeomorphes a des demi-disques et par recollement on a une sphere. A l'ordre $ p=1 $, pour le Tore, c'est ecrit dans le Gramain, la fonction  "cote" a 4 points critiques, un maximum, un minimum et deux points critiques d'indice $ 1 $ (au milieu), dans le Gramain, il colle deux cylindres (les cylindres eux-memes sont obtenues par recollement d'un ensemble d'indice 0 ou 2 avec un ensemble d'indice 1) et il obtient un Tore . Ou bien si on commence la recurrence  a partir du Tore, il suffit de voir que la surface est homeomorphe par recollement, au recollement de deux cylindres (les cylindres sont eux-memes obtenus par recollement d'ensembles d'indice 0 ou 2 et d'un ensemble d'indice 1, 0-1 et 1-2, de haut(0) en bas(1) ou de bas (2) vers le haut (1)) avec l'application de recollement $ g_1=id, g_2 =id $, et on obtient le Tore, $ g_2 $ est la deuxieme application de recollement). Dans le cas general, on sait qu'il existe une fonction de Morse $ u $, on decompose en deux, la partie de $ u $ ayant $ p-1 $ points critiques, on lui recolle un disque et par recurrence, elle homeomorphe \`a $ {\mathbb T}_{q-1} $ et le reste possede 4 point critiques, un maximum, un minimum et deux points critiques d'indice 1, on lui recolle un disque, qui  est, par recurrence, homeomorphe  au Tore , $ {\mathbb T}_1 $. Donc, comme c'est ecrit dans le Gramain, en revenant  \`a $ u $, il y a une partie homeomoprhe \`a $ V_{p-1} $ et l'autre  \`a $ V_1 $ et sur $ V_1 $ on inverse l'orientation. Or aussi sur une partie de la surface consid\'er\'ee, on inverse l'orientation, donc , on inverse l'orientation deux fois (se placer \`a droite, $  Y-g_1(A) $ et $ V_1={\mathbb T}_1-B $, $ B $ un disque), ce qui revient \`a conserver l'orientation en considerant la compos\'ee des deux homeomrphismes (c'est l'homeomorphisme $ k_0=h'og_1og_2^{-1}oh^{-1} $ du Gramain). Il reste a voir  donc avec ces hypoth\`eses, qu'il y a homeomorphismes entre les recollements, $ S= (S-A)\sharp (S'-B)=(X-A)\sharp (Y-g_1(A))=(X-A)\cup (Y-A) $, $ g_1=id $, et $ V_{p-1} \sharp V_1 $, ce qui est fait dans le cas orientable ($ A $ homemorphe au cercle dans le Gramain).

\smallskip

On ecrit $ (X-A) \cup (Y-A) $ pour dire qu'on recolle $ X $ et $ Y $ via $ A $. Les applications du Gramain, $ f=g_1=id, f'=g_2=id $ et le recoll\'e (pour la relation d'equivalence, $ x\equiv f(x)=x, x\in A $, $ f_*=-id, f'_*=-id $, est diffeomorphes a l'ensemble de depart. On a $ f_*=-id $ et $ f'_*=-id $ car on inverse les orientations de depart quand on recolle $ (+,h, A\to A',$ conservation de l'orientation, $ f', A'\to f'(A'),- $ inversion de l'orientation, $ h'^{-1}, f'(A') \to f(A), +, $ conservation de l'orientation, $ f^{-1}, f(A) \to A, - $, inversion de l'orienation). Donc, on a, $ h_*=id $, $ f'_*=-id $, $ (h'^{-1})_*= id $ et $ (f^{-1})_*=-id $. Et on ecrit car $ f=id $, $ f'=id $, $ f_*=id $, $ f'_*=id $, pour dire qu'on envoie un generateur d'orientation vers un generateur d'orientation (meme si elle invers\'ee).

\smallskip

On peut raisonner par les generateurs, si $ e_1 $ est genrateur de $ H_{1,A,X}=H_1(A, A-x_0,Z) $ et $ e_2 $ est gnerateur de $ H_{1,A',X'}=H_1(A',A'-x'_0,Z) $, alors, comme sur la partie de $ Y $ et $ Y' $ on inverse les orientations, on a $ -e_1 $ est genrateur de $ H_{1,A,Y} $ et $ -e_2 $ est genrateur de $ H_{1,A',Y'} $, comme $ f=id  $ alors , $ f_*(e_1)=-e_1 $. On ecrite, $ f_* =id $ dans le sens qu'elle envoie un genrateur d'orientation vers un genrateur d'orientation. et on ecrit $ f_*=-id $, pour dire que $ f_* $ inverse l'orientation de depart sur $ A $. 

\smallskip

On a la meme chose pour $ A' $ et $ f' $.

\smallskip

Des le depart on a deux homeomorphismes $ h $ et $ h' $, supposons pas exemple que $ h_*=-id $, $ h_*(e_1)=-e_2 $, comme pour $ A' $ on est dans l'epsace euclidien, il suffit de composer avec une symetrie $ t $ pour avoir $ t_*(e_2)=-e_2 $ donc, $ ho t $ conserve l'orientation sur le bord. $ (hot)_*(e_1)=e_2 $ et on fait la meme chose avec $ h'$. Finalement, on a deux homemorphismes $ hot $ et $ h'ot'$ ou $ h' $ (si $h'$ conserve l'inversion de l'orientation); qui conservent l'orientation sur les bords.

\smallskip

On a utilis\'e l'orientation definie par l'homologie pour prouver, par recollement que la surface $ S $ est homemorphe \`a la somme connexe de $ q $ Tores, muni d'une orientation homologique. Or La somme connxe de Tore est une variet\'e differentiable et donc, cette orientation homologique est equivalente \`a l'orienation usuelle du Tore. (On pouvait prouver ceci \`a chaque etape, en utilisant, le fait que la surface de depart est orientable et donc s'injecte dans $ {\mathbb R}^3 $, utiliser des "doubles", (voir le Vick, pour l'orientation homologique des vari\'et\'es sans bord et des variet\'es a bord), pour orienter les surfaces a bord de chaque etape, puis, obtenir des surfaces munies d'orientation homologiques et diffenrentiables). L'equivalence entre orientation homologique et differentiable pour une vari\'et\'e lisse est ecrite dans le livre de Bredon (geometry and tpology). Remarquons qu'une vari\'et\'e a bord $ W $ est orientable si et seulement si son interieur $ \dot W $ est orientable, il y a induction de l'orientation sur le bord $ \partial W $ et Il y a equivalence aussi entre orientation homologique et differentiable.(il y a un isomorphisme entre $ H_n(W,\partial W,Z) $ et $ H_n(W,W-x,Z) $). 

\smallskip

Concernant la structure diffentiable sur le recoll\'e de $ X-A\cup_{f,A} Y-A $, dans Milnor (Somme connexe et structure differentiable et Differentiable manifolds which are homotopy spheres), ils prouvent que les voisinages des bords de $ X-A $ et $ Y-A $, sont diffeomorphes a des tubes, puis recollent les tubes, c'est une carte de recollement, pour avoir une structure differentiable. Or, pour la somme connexe de deux vari\'et\'es riemanniennes, les cartes exponentielles sont des voisinages tubulaires des bords, il n'y a pas besoin d'utiliser ce que dit Milnor. Quant \`a ici, le recoll\'e $ X-A\cup_{f,A} Y-A $ est diffeomorphe \`a $ X-A\cup Y-A $, on a deja une strucuture diffenrentiable sur le recoll\'e, via ce diff\'eomorphisme. (la projection $ p:x\to [x] $ est un diff\'eomorphisme).

\smallskip

Quant a l'ecriture $ {\mathbb T}_{q-1}\sharp{\mathbb T}_1 $, c'est la somme connexe, comme on est dans $ {\mathbb R}^3 $, somme connexe et recollement coincident.

\smallskip

Quand on a une application $ s $ qui conserve l'orientation, cela veut dire qu'on a, sur une variete $ M_1 $ compacte sans bord de dimension 1, $ s_*=id $. Pour le voir, on a  si l'orientation est invers\'ee sur $ M_1 $ et $ M_2 =s(M_1) $, cela veut dire que :
$ -1 $ est generateur de $ H_1(M_1, M_1-{x_0}) $ et de $ H_1(M_2, M_2-s(x_0), Z) = H_1(s(M_1), s(M_1)-s(x_0), Z)= s_*(H_1(M_1, M_1-x_0, Z)) $. Donc, $ s_*(-1)=-1 $, c'est a dire que $ s_*=id $. Lorsque $ \phi_1 $ inverse l'orientation alors $ (\phi_1)_*=-id $. 

\smallskip

En utilisant une suite de Mayer-Vietoris, et le fait que $ M_1-x_0 $ est d'homotopie un segment donc a un point, donc d'homologie nulle on a:

\smallskip

$ s_* $ est un isomorphisme de $ H_1(M_1) $ vers $ H_1(M_2) $, comme il conserve l'orientation $ s_*=id $.

\smallskip

Pour le voir il suffit de raisonner sur les genrateurs $ s_*(e_1)=ke_2 $ et $ (s^{-1})_*(e_2)=k'e_1 $, donc, $ kk'=1 $, comme ces des entiens$ k=1 $ ou $k=-1 $, comme $ s $ conserve l'orientation $ k=1 $. On fait la meme chose pour $ \phi_1 $ et on obtient $ k=-1 $. 

\smallskip

On utilise l'homomorphisme surjectif d'Hurwicz qui est bijectif sur le cercle. De $\pi_1 (S_1) $ dans $ H_1(S_1) $

$$ \bar h: \,\, [\gamma] \to \tilde [\gamma] $$

On applique cela aux applications du Gramain, $ s=h'^{-1} o h $ (conserve l'orientation) et $ \phi_1 =\phi $, $ A=\phi(S_1 \times \{0\}) $, ici, $ f=id, f'=id $ et $ M_1=M_2=A $, puis on raisonne en considerant $ S_1 \times \{0\} $ au lieu de $ S_1 $.

$$ k=\phi^{-1} o s o \phi, $$

$$ \pi(k)[\gamma]=[\phi^{-1} o s o \phi o \gamma]= $$

$$ \bar h^{-1} (\tilde [\phi^{-1} o s o \phi o \gamma])= $$

$$ =\bar h^{-1}( ((\phi_1)^{-1})_* o s_* o \phi_* (\tilde [\gamma])) = $$

$$ =\bar h^{-1}\tilde [\gamma] = \bar h^{-1} \bar h ([\gamma]) =[\gamma], $$

O\`u on a utilis\'e le fait que $ \phi_*=id  $ si on se place sur $ X $ (on conserve l'orientation) ou $\phi_*=-id  $ si on se place de cot\'e de $ Y $ (on inverse l'orientation).

\smallskip

On pouvait remplacer $ [\gamma] $ par un generateur $ e_1=[\gamma_1] $ et raisonner sur les generateurs.

\smallskip

Donc,

$$ \pi(k)=id $$

Ceci dans le cas ou $ k(a)= a $, $ k $ fixe un point. Sinon, ($ k $ ne fixe pas $ a $), on compose avec une rotation $ r(e^{i\theta})=e^{iu} e^{i\theta} $, or en considerant $ F(e^{i \theta},s)=e^{i s u} e^{i \theta} $, $ s \in [0,1] $, on a une homotopie entre $ r $ et l'identit\'e $ id $, Donc, $ r_*=id $ de $ H_1(S_1) \to H_1(S_1) $.

\smallskip

Donc, $ k'_*=id $, puis on utlise le meme groupe fondamental en $ a $, $ \pi_1(S_1, a) $ et l'application d'Hurwicz, comme precedemment pour prouver que $ \pi(k')=id $.

\smallskip

On peut utiliser le lift dans le livre de Gramain et prolonger $ k $.

\smallskip

L'autre cas, $ Y $ rectangle et $ g_1(A) $, deux segments, correspond au cas non orientable. Pour le voir, on choisit une orientation  (choisir un repere mobile) sur $ Y $, celle-ci change quand on part de $ X $ et on lui recolle $ Y $ via $ A $. (ici, $ g_1=f $ est l'application de recollement). 

\smallskip

Notons que par Moise-Rado, homeomoprhe implique diffeomorphe en dimension  $\leq 3 $.

\smallskip

13) Pour l'inegalit\'e de Kato: on a $ u\in L^1_{loc} $ et $ \Delta u \in L^1_{loc} $ alors: $ u \in W^{1,1}_{loc} $, pourquoi ?

\smallskip

Comme c'est local, on ecrit dans une boule $ u $ en fonction de la fonction de Green (de la boule par exemple et $ n\geq 3 $ par exemple, $ n=2 $ c'est la meme chose):

$$ u= \dfrac{1}{r^{n-2}}* \Delta u+ \int_{\partial B} P*u, $$

Par Fubini, Fubini-Tonelli,

$$ \partial u = \partial (\dfrac{1}{r^{n-2}})* \Delta u+ \int_{\partial B} (\partial P)u, $$

avec $ \partial P $ regulier des qu'on considere des points interieurs et donc on peut deriver sous le signe $\int $ pour cette partie.

On a,

$$ u\in L^1_{loc} \,\,\Rightarrow \,\, \exists \,\, r \,\, \int_{\partial B_r} u(\sigma) d\sigma <+\infty, $$

Au sens des distributions on a pour $ \phi \in D(B_{r/2}) $:

$$ \int (-\partial \phi ) u=\int (-\Delta u ) \int (-\partial \phi)(\dfrac{1}{r^{n-2}}) +\,\, {\rm termes \,\, reguliers \,\, de\,\, la \,\,fonction\,\, de\,\, Green} $$

Donc, apres integration par parties pour le noyau Newtonien,

$$ \int (\partial \phi) u = \int \int (\partial (\dfrac{1}{r^{n-2}})(-\Delta u) \phi+ \int \int \partial P u \phi dx d\sigma = \int a \phi, \,\, a\in L^1 ,$$

Donc,

$$ u\in W^{1,1}_{loc}, $$

Ceci est valable pour $ u $ reguliere. Soit $ \rho_n $ un molliffier, alors :

Si $ u\in L^1(B_r) $, alors:

$$ \rho_n* u(x)=\int_{B_r} \rho_n(x-y) u(y) dy, $$

et pour $ x\in B_{r/2} $:

$$\Delta (\rho_n*u)(x)=(\rho_n*\Delta u)(x), $$

On a aussi,

$$ \rho_n*u \to u, \,\, {\rm dans }\,\, L^1_{loc}(B_r), $$

et,

$$ \Delta (\rho_n*u)=\rho_n*(\Delta u) \to \Delta u,\,\, {\rm dans} \,\, L^1_{loc}(B_{r/2}), $$

On fait ce qu'on a fait avec $ u $, avec $ \rho_n*u $, comme on a des bornes uniformes en $ L^1_{loc} $, alors,
$$ \rho_n*u \in W^{1,1}_{loc}, \,\, {\rm et}\,\, \nabla (\rho_n*u) \to v \in L^1_{loc},$$

Donc,

$$ \int u \partial \phi = \lim_n \int (\rho_n*u) \partial \phi=-\lim_n \int \partial (\rho_n*u) \phi \to -\int v\phi, $$

Donc,

$$ u\in W^{1,1}_{loc}, $$

On a la meme chose si on considere un probleme variationel $ -\Delta u = f $, $ f\in L^1_{loc} $ et $ u \in L^1_{loc} $. On a meme mieux, $ u \in W^{1,q}_{loc} $. Une partie du travail de Brezis-Merle suppose ces hypotheses, par ce qu'on vient de voir, on utilise le th\'eor\`eme 1 de Brezis-Merle et les estimations elliptiques, pour poser le probleme quand on n'a pas de condition aux bord dans la formulation de Brezis et Merle. 

\smallskip

On peut retrouver ce resultat en se ramenant a une fonction harmonique en soustrayant une solution d'un probleme de Dirichlet.

Soit $ f\in L^1 $, $ \exists \, f_j \to f $ dans $ L^1 $ avec $ f_j \in C^{\infty}_c $, on resout:

$$ -\Delta u_j = f_j, \,\, {\rm avec} \,\, u_j=0, \,\, {\rm sur}\,\, \partial \Omega, $$

Alors, par Stampacchia ou Brezis-Strauss, on a:

$$ ||u_j||_{W^{1,q}_0} \leq C_q, $$

En passant \`a la limite en $ j $, on a l'existence de $ u_0 \in W^{1,q}_0 $ tel que:

$$ -\Delta u_0=f, $$

Alors,

$$ -\Delta (u-u_0) =0,\,\, {\rm et} \,\, u-u_0\in L^1, $$

On utilise le theoreme de Weyl pour avoir,

$$ u-u_0\in C^{\infty} $$

Donc,

$$ u\in W^{1,q}_{loc}, $$

On retrouve le Theoreme dans le livre de Dautray-Lions.

On fait la meme chose avec le Theoreme 1 de Brezis-Merle, sauf que les relations ne sont pas locales.

On a:

$$ -\Delta u=f \in L^1,\,\, u\in W^{1,1}_0(\Omega). $$

On resout pour $ f_j\to f $ dans $ L^1 $. $ f_j \in C^{\infty}_c(\Omega) $. (la regularit\'e du bord est au moins $ C^2 $ pour pouvoir appliquer le theoreme de dualit\'e de Stampacchia).

$$ -\Delta u_j=f_j, \,\, u_j\in H^1_0(\Omega), $$

On a:

$$ ||\nabla u_j||_q \leq C_q, \,\, 1\leq q <2. $$

On passe a la limte en $ j $.

On a:

$$ -\Delta u_0=f, \,\, u_0\in W^{1,q}_0(\Omega). $$

Finalement:

$$ -\Delta (u-u_0) =0, u-u_0\in W^{1,1}_0(\Omega). $$

Soit on utilise le principe du maximum dans $ W^{1,1}_0 $ pour avoir $ u=u_0\,\, p.p $.
\smallskip

Soit, on utilise le theoreme de regularit\'e d'Agmon pour avoir $ u-u_0 $ reguliere et appliquer le principe du maximum usuel.

\smallskip

Pour, $ u_j, f_j, |f_j| $ on applique le raisonnement de Brezis Merle (theoreme 1, on ne peut utiliser le potentiel Newtonien qu'avec des hypoth\`eses sur $ f_j $; integrables et born\'ees, par exemple, voir le livre de Gilbarg-Trudinger). Puis par le lemme de Fatou, l'inegalit\'e de Brezis Merle est valide pour $ u_j $ et $ u_0 $ donc pour $ u $. ($ u=u_0,\,\, p.p $).

\smallskip

14) Sur la fonction distance au bord: c'est bien ecrit dans le Gilbarg-Trudinger. Grace a la propri\'et\'e de la sphere interieure, il existe un voisnage du bord assez petit tel que la fonction distance est $ C^k $ et est egale \`a la distance (pour chaque point) \`a un point unique du bord, car on fait varier des spheres le long du bord. Voir, Gilbarg-Trudinger. Le domaine (et le bord) doit etre par exemple au moins $ C^2 $.

\smallskip

Par des cartes, on a un segment, on peut faire rouler une sphere. On prend son image par l'application de carte, qui est suffisemment reguliere pour que le domaine d'arriv\'ee contienne des ellipses, donc des cercles.

\smallskip

15) Pour le theoreme de Dualit\'e de Stampacchia:

\smallskip

On part de l'equation:

$$ -\Delta u_i =V_i e^{u_i}, $$

avec la condition de Dirichlet. ($ \Delta= \partial_{11}+\partial_{22} $).

On a:

$$ u_i \in W^{2,k}\cap C^{1,\epsilon}(\bar \Omega), \,\, u_i= 0 \,\, \text{on }\,\,\partial \Omega. $$ 

On a:

$$ \int_{\Omega} |\Delta u_i|dx \leq bC. $$

On considere un vecteur $ f=(f_1,f_2) \in L^{q'}(\Omega) $ avec $ q' >2 $. On utilise le Theoreme de Lax-Milgram dans les Hilbert, (voir Gilbarg-Trudinger), pour avoir une solution de $ z \in W_0^{1,2}(\Omega) $ de l'equation:

$$  -\Delta z= div (f), $$ 

avec la propri\'et\'e suivante (voir la preuve de l'in\'egalit\'e de Harnack dans le Gilbarg-Trudinger):

$$  ||z||_{L^{\infty}} \leq C||f||_{q'}. $$ 

On met $  u_i $ comme fonction test dans l'equation de $ z $ , on a, en valeur absolue a l'exterieur des integrales:

$$ \int_{\Omega} f \cdot \nabla u_i dx = \int_{\Omega} \nabla z \cdot \nabla u_i dx = -\int_{\Omega} z \Delta u_i dx \leq C'||f||_{q'}, $$

On prend le supremum dans $ L^{q'} $, et on a:

$$ ||\nabla u_i ||_q \leq C_q, \, \forall \, 1\leq q <2. $$

Dans le probleme de Brezis-Merle, on a:

$$ ||\nabla (u_i-u)||_{L^q(\Omega_{\epsilon})}=o(1),  $$

c'est la convergence interieure d\^ue au travail de Brezis et Merle. Apr\'es, par l'inegalit\'e de H\"older:

$$ ||\nabla (u_i-u)||_{L^q(\Omega-\Omega_{\epsilon})} \leq |\Omega -\Omega_{\epsilon}|^{1/q-1/r} ||\nabla (u_i-u)||_{L^r} \leq C_r |\Omega-\Omega_{\epsilon}|^{1/q-1/r} =o(1), \,\, 1\leq q <r <2. $$

On obtient finalement:

$$ ||\nabla (u_i-u)||_q =o(1). $$

\smallskip

16) Pourquoi dans le livre d'Aubin, il considere les solutions $ u $ pour un operateur du type:

$$ a_{ij}(M) \partial_{ij} u+b_j(M) \partial_j u+c(M) u =f(M) ? $$

Dans le livre d'aubin $ a_{ij}(M), b_j(M) $ sont des champs de tenseurs, pourquoi tenseur ?

\smallskip

C'est li\'e au cartes et au changements par cartes et \`a la formulation variationelle.

\smallskip

1) Dans le livre d'Hebey, il donne la definition d'un probleme pos\'e de maniere variationnelle:

\smallskip

Pour les vari\'et\'es:

$$ \int_M a_{ij} \partial u \partial \phi dV_g +b_j u \partial \phi dV_g+ cu \phi dV_g = \int f \phi dV_g $$

a) Il faut que le probleme soit pos\'e sur une vari\'et\'e de maniere globale, donc independemment des cartes. C pour cela que dans le livre d'Hebey, il dit pour toute carte $ h $.

\smallskip

b) Si on revient \`a un probleme pos\'e sur un ouvert de $ {\mathbb R}^n $, il y a deja une carte disponible, la carte euclidienne, $ (\Omega, id) $. Donc pour $ {\mathbb R}^n $, il suffit que le probleme soit pos\'e en coordonn\'ees cartesiennes.

\smallskip

Et apres, par la formule de changement de variable, on se ramene \`a toute carte de l'ouvert. Notons que la notation $ dV_g $ l'explique bien car par un changement de variable, on retrouve la Jacobien, qui donne l'element de volume pour la metrique euclidienne.

\smallskip

Mainteneant quand on fait le changement de variables, il faut que le probleme soit invariant, c'est \`a dire que les coefficients ne changent par or ceci est possible si on les suppose comme champs de tenseurs.

\smallskip

on ecrit par exemple:

$$ a_{ij}(M)=a_{ij}[x(M)] ={a_{kl}}[y(M)] \partial_x y_i^k \partial_x y_j^l, $$

dans le changement  de variables, cette loi est celle des tenseurs. 

\smallskip

C'est pour cela qu'on suppose les coefficients, des champs de tenseurs. Pour que par changement de variables, ils soit invariants et ceci est possible si ils obeissent \`a la loi des tenseurs.

\smallskip

Et pour les ouverts de $ {\mathbb R}^n $. Il suffit que le probleme soit pos\'e en coordonn\'es cartesiennes. Alors on la pour toute carte.

\smallskip

Une forme bilineaire est un tenseur. Un champ de formes bilneaires est un champ de tenseurs. Un champ de vecteurs est un champ de tenseurs.

\smallskip

17) Pour quoi, quand on veut appliquer le principe du maximum sur une vari\'et\'e, on a besoin d'un op\'erateur global:

\smallskip

a) C'est parce qu'on a consid\'er\'e une vari\'et\'e, qui est cosntitu\'ee de cartes. Il ne faut pas que ca depende de la carte.

\smallskip

Par exemple, si on veut appliquer la principe du maximum en $ P $,

\smallskip

On ecrit : $ Lu(P)=L(u o \phi^{-1}(x_1,x_2)) $, on voit que dans l'ecriture meme, il ne faut pas que ca depende de la carte et de l'application de carte $ \phi $.

\smallskip

b) Par recouvrement, dans l'intersection d'ouvert de carte, on ecrit:

$$ Lu(P)=L[u o\phi^{-1}]=L[u o\psi^{-1}] $$

dans cette ecriture, on a deux cartes, et ile ne faut pas que ca depende des deux cartes, pour cela il faut que l'operateur soit global.

\smallskip

c) On a la meme chose, si on raisonne que la deriv\'ee normale et le principe du maximum de Hopf.

\smallskip

18) Sur les applications conformes du bord en dimension 2.

\smallskip

On a parl\'e d'application conforme en dimension 2:

\smallskip

Pourquoi la definition d'un domaine analytique permet d'avoir une carte conforme locale du bord ?

\smallskip

a)-Sur le disque unit\'e, on a une transformation conforme directe du disque sur le demi-plan de Poincar\'e. Il n'y pas besoin de definir l'application $ \psi $ (voir ci-dessous).

\smallskip

b)-L'application $ \psi $, permet de definir des coordonn\'ees qui permettent d'utiliser une application conforme.

\smallskip

Ce qui suit est la traduction de ce qui a \'et\'e ecrit dans les preprints sur arXiv. Pour definir un carte "conforme" $ \gamma $, en dimension 2.

\smallskip

1-Il suffit de prouver que $ \gamma_1((-\epsilon,\epsilon))=
\partial \Omega \cap \tilde \gamma_1(B_{\epsilon})= \partial \Omega \cap \tilde \gamma_1(B_{\epsilon})\cap \{ |abscisse|< \epsilon \} $, pour $ \epsilon >0 $ assez petit. Avec $ \tilde \gamma_1 $ l'extension holomorphe de $ \gamma_1(t)=t+i\phi(t) $.

\smallskip

Pour voir ca, on raisonne par l'absurde. On a pour $ z_{\epsilon} \in B_{\epsilon} $, $\tilde \gamma_1(z_{\epsilon}) =(t_{\epsilon},\phi (t_{\epsilon})) $ avec $ |t_{\epsilon}|\geq \epsilon $. Comme $ \tilde \gamma_1 $ est injective sur une boule fix\'es au depart par le theoreme d'inversion locale, $ B_1 $  et  $ \tilde \gamma_1 =\gamma_1=t+i\phi(t) $ sur l'axe r\'eel, on a n\'ecessairement  $ |t_{\epsilon}|\geq 1 $. Mais par continuit\'e $ |\tilde \gamma_1(z_{\epsilon})|\to 0 $ parce que $ z_{\epsilon} \to 0 $. Et on utilise le fait que $ |\tilde \gamma_1(z_{\epsilon})|=|(t_{\epsilon},\phi(t_{\epsilon}))|\geq |t_{\epsilon}|\geq 1 $, pour avoir une contradiction.) ( Cela veut dire que pour un rayon assez petit quand le graphe sort de la boule, il n'y retourne plus). ( Ce fait implique, en raisonnant avec des chemins, quand un a chemin qui coupe $\partial \Omega $ in $ \tilde \gamma_1(B_{\epsilon}) $ en un point ce point a pour abscisse $ |abscisse|< \epsilon $.  Encore, une fois, ce fait (par un raisonnement par l'absurde avec le fait que $ \partial \bar \Omega = \partial \Omega $, il faut raisonner avec des chemins), implique que l'image de la partie superieure est entierement d'un cot\'e de la courbe et l'image de l'autre est de l'autre cot\'e de la courbe.

\smallskip

2- Posons: $ \psi : (\lambda_1,\lambda_2)\to M \in \Omega $ telle que: $ \overrightarrow{x_0M}=\lambda_1 i_1'+\lambda_2j_1' $ avec  $ (i_1',j_1') $ est une base telle que $ i_1'=e^{-i\theta} i_1, j_1'=e^{-i\theta} j_1 $. Et, $ \phi=id: (x_1,x_2)\to M $ telle que: $ \overrightarrow{OM}=x_1i_1+x_2i_2 $ la base canonique $ (i_1,j_1) $. Alors on a deux cartes de $ \Omega $, $ \phi $ et $ \psi $ et les affixes $ T_M=\lambda_1+i\lambda_2 $ et $ z_M=x_1+ix_2 $ sont tels que :(l'application de changement de cartes):

$$ T_M/e^{i\theta}+x_0=z_M=\phi^{-1} o \psi(\lambda_1,\lambda_2), $$

On a:

$$ \partial_{\lambda_1}=\cos \theta \partial_{x_1}+\sin \theta \partial_{x_2}, $$

$$ \partial_{\lambda_2}=-\sin \theta \partial_{x_1}+\cos \theta \partial_{x_2}, $$

Donc, la m\'etrique dans la carte $\psi $ ou en coordonn\'ees $ (\lambda_1,\lambda_2) $ est : $ g_{ij}^{\lambda}=\delta_{ij} $ et le Laplacien dans les deux cartes (coordonn\'ees), $ \psi $ et $\phi $ est le Laplacien usuel $ \partial_{\lambda_1\lambda_1}+\partial_{\lambda_2\lambda_2} $.

\smallskip

On ecrit:

$$ \Delta u(M)=\Delta_{\lambda} (u o\psi(\lambda_1,\lambda_2)) $$

Puis on applique l'application conforme $\tilde \gamma_1 $ qui envoie l'affixe $ T_M $, $ M $ dans un voisinage de of $ x_0\in \partial \Omega $ vers $ B_{\epsilon} $ et envoie $ T_M, M\in \partial \Omega $ sur l'axe r\'eel $ (-\epsilon,\epsilon) $ et l'autre partie vers l'autre partie de $ \Omega $ et $ \bar \Omega^c $.

\smallskip

On a: $ \psi:(\lambda_1,\lambda_2) \to M $ et $ \tilde \gamma_1:(\lambda_1,\lambda_2)\to (\mu_1,\mu_2) $. L'application consid\'er\'ee est: $ \psi o \tilde \gamma_1^{-1} $. C'est bien une carte. Pour le voir:

\smallskip

La premiere carte (carte usuelle) est: $ \psi o g^{-1} $ avec $ g $ est l'application $ g:(\lambda_1,\lambda_2)\to [\lambda_1,\lambda_2+\phi_0(\lambda_1)] $ et $ \phi_0 $ l'application qui permet de definir $ \partial \Omega $ comme un graphe.

\smallskip

L'application de changement de carte est:

$$ g o\tilde \gamma_1^{-1} $$

elle est reguliere d'un voisinage de $ 0 $ vers un autre voisinage de $ 0 $  et envoie la partie superieure du premier voisinage vers la partie superieure de l'autre voisinage (ainsi que les parties inferieures).

\smallskip

On voit bien que $ \psi o\tilde \gamma_1^{-1} $ est une carte avec la propri\'et\'e que $ \tilde \gamma_1 $ est conforme.

\smallskip

(l'application $ \psi o \tilde \gamma_1^{-1} $ est aussi la compos\'ee de, une rotation et une translation (donc conforme) et de $ \tilde \gamma_1^{-1} $ qui est conforme, si on se placait en coordonn\'ees cart\'esiennes au d\'epart. Mais, il faut voir que la carte $ \psi $ definit des coordonn\'es qui permettent la construction d'une application conforme. On n'a pas besoin de cela pour le disque unit\'e, on a directement l'application du demi-plan de Poincar\'e).

\smallskip

Remarque: $ \psi $ et $ \phi=id $ sont des cartes de $ \Omega $ et $ \psi $ est "presque" une carte de $ \partial \Omega $. Mais $ \psi o \tilde \gamma_1^{-1} $ et $ \psi o g^{-1} $ sont des cartes du bord.

\smallskip

Le probleme elliptique peut etre vu de deux manieres : donc, si on considere, au depart, les coordonn\'ees cartesiennes, on a une transformation conforme. Si, on considere le probleme pos\'e sur une vari\'et\'e a bord, on a une carte analytique.

\smallskip

3- "Caractersation des domaines $ C^k $": On peut dire que la definition d'un domaine $ C^k, k\geq 1 $ est equivalente \`a la definition d'une sous-vari\'et\'e $ + $ la condition : $ \partial \bar \Omega = \partial \Omega $ or $ {\dot{\bar \Omega }} = \Omega $.

\smallskip

19) Concernant l'element de volume et la mesure de Hausdorff sur le bord:

\smallskip

Avec le point 18), la carte usuelle du bord est: $ h =\psi o g^{-1} $ avec, 

\smallskip

$ \psi : (\lambda_1,\lambda_2)\to M \in \Omega $ telle que: $ \overrightarrow{x_0M}=\lambda_1 i_1'+\lambda_2j_1' $ avec  $ (i_1',j_1') $ est une base telle que $ i_1'=e^{-i\theta} i_1, j_1'=e^{-i\theta} j_1 $. Et, $ \phi=id: (x_1,x_2)\to M $ telle que: $ \overrightarrow{OM}=x_1i_1+x_2i_2 $ la base canonique $ (i_1,j_1) $. Alors on a deux cartes de $ \Omega $, $ \phi $ et $ \psi $ et les affixes $ T_M=\lambda_1+i\lambda_2 $ et $ z_M=x_1+ix_2 $ sont tels que :(l'application de changement de cartes):

$$ T_M/e^{i\theta}+x_0=z_M=\phi^{-1} o \psi(\lambda_1,\lambda_2), $$

et,

$ g $ est l'application $ g:(\lambda_1,\lambda_2)\to [\lambda_1,\lambda_2+\phi_0(\lambda_1)] $ et $ \phi_0 $ l'application qui permet de definir $ \partial \Omega $ comme un graphe.

La metrique sur le bord $ \partial \Omega $ s'ecrit localement:

$$ h^*(i^*(\delta))={\sqrt { 1+[\phi_0'(\lambda)]^2}} d\lambda, $$

avec $ \delta $ la metrique euclidienne de $ {\mathbb R}^2 $.(la metrique sur le bord est $ i^*(\delta) $ avec $ i:\partial \Omega \to {\mathbb R}^2 $, $ x\to x $).

\smallskip

On l'a explicitement, en calculant (le determinant contient un terme car on est en dimension 1, c la meme chose en dimension superieure),  $ h^*(i^*(\delta))=\delta (dh(\partial_{\lambda}), dh(\partial_{\lambda}))=[1+[\phi_0'(\lambda)]^2] $.

\smallskip

Avec $ i:\partial \Omega \to {\mathbb R}^2 $ , $ x \to x $ l'inclusion.

\smallskip

Et considerant l'application Lipschitzienne ( c'est aussi l'application de carte, locale, du bord $ \partial \Omega $ qui permet aussi de calculer $ g_{11} $, ce qu'on a fait precedemment):

$$ \lambda \to (\lambda,\phi_0(\lambda))=h(\lambda,0), $$

on a par la formule de l'aire, l'egalit\'e entre la mesure de Hausdorff et le Jacobien (element d'aire):

$$ \int_{[a,b]} {\sqrt { 1+[\phi_0'(\lambda)]^2}} d\lambda=\int_{h([a,b],0)\subset \partial \Omega} dH_1, $$

Finalement, on a l'egalit\'e suivante:

$$ d\sigma ={\sqrt { 1+[\phi_0'(\lambda)]^2}} d\lambda={dH_1}_{|\partial \Omega}. $$

avec cette remarque, on voit par exemple que $ H_1(\partial \Omega) <+\infty $, dans le cas $\partial \Omega $ compact.

\smallskip

20) Quand y a il une relation entre $ W^{1,\infty} $ et $ C^{0,1}(\bar \Omega) =Lipschitz(\bar \Omega) $.

\smallskip

a) Dans le cas d'un convexe ou une boule, c'est ecrit dans le livre de Brezis quand il approxime $ L^{\infty} $ par $ L^p $, dans le critere pour qu'une fonction soit $ W^{1,p} $.

\smallskip

b) Quand l'ouvert est Lipschitzien ou au moins $ C^1 $, c'est un exercice dans livre d'Otared Kavian. On utlise l'operateur de prolongement (Notez que la constante de prolongement ne depend pas de $ p $, on peut utliser le proceder de Brezis $ p \to +\infty $), puis on se ramene \'a un convexe, une boule, puis on revient aux fonctions de depart.

\smallskip

Dans les cas precedents, il y a egalit\'e entre les deux ensembles.

\smallskip

21)  pourquoi on ecrit dan la derivation (ici $ \nabla =\nabla_X $ pour un vecteur $ X $):

 $$ \nabla (<a|b>)=<\nabla a|b>+<a|\nabla b> $$

ceci est du \`a ( $ g $ est la metrique et $ C $ la contraction):

$$ <a|b>=C_1^3C_2^4 (g\otimes a \otimes b ) $$

La deriv\'ee covariantes commute avec la contarction et $\nabla g=0 $, car c la connexion de Levi-Civita:

$$ \nabla (<a|b>)=C_1^3C_2^4(\nabla g \otimes a \otimes b+g\otimes \nabla a \otimes b+ g\otimes a\otimes \nabla b) $$

Donc,

$$ \nabla (<a|b>)=C_1^3C_2^4(g\otimes \nabla a \otimes b+ g\otimes a\otimes \nabla b)=<\nabla a|b>+<a|\nabla b> . $$

\smallskip

22) Sur les coordonn\'ees geodesiques polaires : elles sont definies a partir de la carte exponentielle. Il existe un voisinage du point tel que pour $ y \in B(x,\epsilon) $, $ z \to exp_y(z) $ est un diffeomeorphisme. On a alors une carte. On definit les coordonnees geodesiques polaires. On s'assure que (lemme de Gauss):

\smallskip

a) 

$$ \dfrac{d}{dt}(<\nabla_t|\nabla_t>)= \nabla_t (<\nabla_t|\nabla_t>)=2<\nabla_t(\nabla_t)| \nabla_t>=0 $$

Car c'est une geodesique $\nabla_t(\nabla_t)=0 $. Donc, $<\nabla_t|\nabla_t>\equiv cte = \lim_{t\to 0}  <\nabla_t| \nabla_t> =1 $, car la carte exponentielle est normale en $ y $.

$$ \nabla_t(<\nabla_t|\nabla_{\theta}>)=\dfrac{1}{2} \nabla_{\theta}(<\nabla_t|\nabla_t>)=\dfrac{1}{2} \nabla_{\theta} (1)=0, $$
Ici, on utilise le fait que la connexion est sans torsion pour permuter la derivation en $ t $ et angulaire.

Donc,

$$ <\nabla_t| \nabla_{\theta}>\equiv cte= \lim_{t\to 0}(<\nabla_t|\nabla_{\theta}>)=0 $$

C'est le lemme de Gauss. 

\smallskip

On pouvait ecrire $\nabla_{\partial_i}(\partial_j) $ avec $ \partial_i,\partial_j \in \{ \partial_t,\partial_{\theta_j} \}. $

\smallskip

b) Les geodesiques sont minimisantes, c bien ecrit dans le livre d'Hebey et Gallot-Hulin-Lafontaine, remarquons, comme c ecrit dans le Do Carmo. Pour calculer la longueur d'un chmin passant par l'origine, on prend $ \int_{\epsilon}^1 $, puis on fait tendre $ \epsilon  $ vers $ 0 $. Car les coordonnees geodesiques polaires sont definit en dehors de $ 0 $, origine, puis on fait tendre le rayon $ r $ vers $ 0 $.

\smallskip

c) On voit comme c dit dans le Hebey, que les geodesiques sont minimisantes, l'image de boules ouvertes sont des boules ouvertes, de boules fermees, sont des boules fermees, l'image de spheres, sont des spheres.

Soit, $ \tilde z $ l'application de $ {\mathbb S}_{n-1} $ vers $ {\mathbb S}_{n-1}^r $ d\'efinie par $ \theta \to r \theta $.

Par d\'efinition, pour $ h $ une fonction d\'efinie au voisinage de $ x $:

$$ \partial_{j,\phi_0}(h)=\dfrac{\partial (ho\phi_0)}{\partial \bar \theta_j}=\dfrac{\partial (ho\exp_xozo\psi^{-1})}{\partial \bar \theta_j} =\dfrac{\partial (ho\exp_x)}{\partial z_k}\times r \dfrac{\partial (\psi^{-1})^k}{\partial \bar \theta_j}.$$ 

avec $ \phi_0= \exp_xozo\psi^{-1} $.

\smallskip

Ce qu'on peut ecrire:

$$ \partial_{j,\phi_0}(h)=\dfrac{\partial (ho\exp_x)}{\partial z_k}\times r \dfrac{\partial \theta^k}{\partial \bar \theta_j}.$$ 

et,

$$ \partial_{\theta^i,x}=r b_i^k(\theta)\partial_{z^k,x}. $$

On a donc, pour tout $ x\in B(x_0,\epsilon_0) $, $ [B(x,\epsilon_0), \exp_x^{-1}] $ est une carte normale en $ x $. Pour tout $ x\in B(x_0,\epsilon_0) $:

$$ g_{x,ij}(z)=g(z)(\partial_{z^i,x},\partial_{z^j,x}), $$

o\`u $ \partial_{z^i,x} $ est le champ de vecteur locale canonique li\'e \`a la carte exponentielle.

\smallskip

On note, $ a_i^k(z,x)=\dfrac{(\phi o \exp_x)^k}{\partial z^i}[\exp_x^{-1}(z)] $, alors, $ \partial_{z^i,x}=a_i^k(z,x) \partial_{u^k,\phi}, $

\smallskip

avec, $ \partial_{u^k,\phi} $, le champ de vecteur canonique relatif \`a la carte $ (\Omega, \phi) $, ce dernier ne d\'epend pas du point $ x $, contrairement \`a celui de la carte exponentielle en $ x $. Quant aux fonctions $ a_i^k $ elles sont $ C^{\infty} $ de $ z $ et $ x $. On obtient,

$$ g_{x,ij}(z)=g(z)[a_i^k(z,x) \partial_{u^k,\phi}; a_j^l(z,x)\partial_{u^l,\phi}]=a_i^k(z,x)a_j^l(z,x)g_{kl}(z), $$

o\`u les  $ g_{kl} $ sont les composantes de $ g $ dans la carte  $ ( \Omega, \phi ) $.

\smallskip

On a, $ z=\exp_x(y) $, $ y \in B(0,\epsilon_0) \subset {\mathbb R}^n $ et $ y=r\theta $ en coordonn\'ees polaires, ainsi, la fonction $ (x,r,\theta) \to g_{x,ij}[\exp_x(r\theta)] $ est $ C^{\infty} $ des variables $ x, r $ et $ \theta $.

\smallskip

D'autre part, par d\'efinition, $ g_{ij}^k(r\theta)=g_{[\exp_x(r\theta)]}(\partial_{\theta^i,x},\partial_{\theta^j,x}) $, (champs de vecteurs canoniques).

\smallskip

En utilisant le m\^eme type de calcul que pr\'ec\'edemment, on a,

$$ \partial_{\theta^i,x}=r b_i^k(\theta)\partial_{z^k,x} , $$

avec, $ b_i^j $ des fonction tr\'es r\'eguli\`eres. Donc,

$$ g_{ij}^k(r,\theta)=r^2g[\exp_x(r\theta)](b_i^k\partial_{z^k,x},b_j^l\partial_{z^l,x}), $$

d'o\`u,

$$ g_{ij}^k(r,\theta)=r^2b_i^k(\theta)b_j^l(\theta)g_{x,kl}[\exp_x(r\theta)]. $$

Donc, l'application, $ {\tilde g}_{ij}^k : (x,r,\theta) \to b_j^kb_j^lg_{x,kl}[\exp_x(r\theta)] $ est $ {C}^{\infty} $ des trois variables $ x,r $ et $ \theta $. On a,

$$ \partial_r {\tilde g}_{ij}^k(x,0,\theta)=b_i^k(\theta)b_j^k(\theta)c^m(\theta)\partial_mg_{x,kl}(x)=0 , $$

car, la carte exponentielle est normale en $ x $ et les $ g_{x,kl} $ sont les composantes de  $ g $ dans cette m\^eme carte. On a aussi,

$$ \partial_{\theta^m}{\tilde g}_{ij}^k(x,r,\theta)={\tilde b}_i^k(\theta) {\tilde b}_j^l(\theta)g_{x,kl}[\exp_x(r\theta)]+r{\bar b}_i^k(\theta){\bar b}_j^l(\theta){\bar c}_m^s(\theta)\partial_sg_{x,kl}[\exp_x(r\theta)], $$

En red\'erivant en $ r $, on obtient,

$$ \partial_r \partial_{\theta^m} {\tilde g}_{ij}^k(x,r,\theta)=u_{ijmr}^{klq}(\theta)\partial_q g_{x,kl}[\exp_x(r\theta)]+v_{ijmr}^{kl}(\theta)w^s(t)\partial_sg_{x,kl}[\exp_x(r\theta)]+ $$

$$ + rh_{ijrm}^{klst}(\theta)\partial_{st}g_{x,kl}[\exp_x(r\theta)]. $$

Donc,

$$ \partial_r \partial_{\theta^m} {\tilde g}_{ij}^k(x,0,\theta)=0, \,\, \forall \,_, x\in B(x_0,\epsilon_0), \,\, \forall \,\, \theta \in U^k . $$

Ainsi, on obtient:

$$ \partial_r {\tilde g}_{ij}^k(x,0,\theta)=\partial_r \partial_{\theta^m}{\tilde g}_{ij}^k(x,0,\theta)=0,\,\, \forall \,\, x\in B(x_0,\epsilon_0), \,\, \forall \,\, \theta \in U^k. \qquad \qquad (*) $$

\smallskip

On a la meme chose en considerant le determinant:

\smallskip

Comme, $ \sqrt { |{\tilde g}^k|}=\alpha^k(\theta) \sqrt {[det(g_{x,ij})]} $, on en d\'eduit,

$$ \partial_r (\log \sqrt {|{\tilde g}^k|})=\partial_r [\log (\sqrt { [det(g_{x,ij})]})]. $$

Par d\'efinition du d\'eterminant et en d\'eveloppant par rapport al ligne 1 par exemple, on obtient,

$$ det [g_{x,ij}][\exp_x(r\theta)]=\Sigma \,\, \Pi \,\, g_{x,kl}[\exp_x(r\theta)], $$

d'o\`u,

$$ \partial_r det[g_{x,ij}](x)=\Sigma\,\, \Pi \,\, [g_{x,kl}(x)] a^s(\theta) \partial_sg_{x,mn}(x) = 0, $$

Car, la carte exponentielle est normale en $ x $.

\smallskip

Finalement,

$$ \partial_r |{\tilde g}^k|(x,0,\theta)=0, \,\, \forall \,\, x\in B(x_0,\epsilon_0), \,\, \forall \,\, \theta \in U^k . $$

En utilisant le m\^eme raisonnement  que celui qu'on utilis\'e pour calculer, $ \partial_r\partial_{\theta^m} {\tilde g}_{ij}^k(x,0,\theta) $, on prouve que,

$$ \partial_r \partial_{\theta^m} |{\tilde g}^k|(x,0,\theta)=0. $$

En conclusion,

$$  \partial_r {\tilde g}_{ij}^k(x,0,\theta)=\partial_r \partial_{\theta^m}{\tilde g}_{ij}^k(x,0,\theta)=0\,\, \forall \,\, x\in B(x_0,\epsilon_0), \,\, \forall \,\, \theta \in U^k. \qquad \qquad (**) $$

$$ \partial_r |{\tilde g}^k|(x,0,\theta)=\partial_r \partial_{\theta^m} |{\tilde g}^k|(x,0,\theta)=0\,\, \forall \,\, x\in B(x_0,\epsilon_0), \,\, \forall \,\, \theta \in U^k. \qquad \qquad (***) $$

$$ \partial_r \partial_{\theta_m}\left [\dfrac{\sqrt {|{\tilde g}^k|} }{\alpha^k(\theta) } \right ] (x,0,\theta)=\partial_r \partial_{\theta_m}\partial_{\theta{m'}} \left [\dfrac{\sqrt {|{\tilde g}^k|} }{\alpha^k(\theta) } \right ] (x,0,\theta)=0 \,\, \forall \,\, x\in B(x_0,\epsilon_0), \,\, \forall \,\, \theta \in U^k. \,\, (****) $$

L'ouvert $ U^k $ est choisit de telle mani\`ere qu'en le r\'eduisant un peu, il reste un ouvert de carte.

\smallskip

On utilise alors, l'uniforme continuit\'e des d\'eriv\'ees succ\'essives pour obtenir les r\'esultats de la proposition 1.

\underbar {\bf Remarque:} 

\smallskip

$ \partial_r [ \log \sqrt { |{\tilde g}^k |}] $ qui est une fonction locale de $ \theta $, est la restriction d'une fonction globale, d\'efinie sur toute la sph\`ere $ {\mathbb S}^{n-1} $ \`a savoir $ \partial_r [\log \sqrt { det(g_{x, ij})}] $. On notera, $ J(x,r,\theta)=\sqrt { det(g_{x, ij})} $.

\smallskip

23) Le Laplacien en polaires :

\smallskip

Ecrivons le laplacien dans $ ]0,\epsilon_1[\times U^k $,

$$ -\Delta = \partial_{rr}+\dfrac{n-1}{r}\partial_r+ \partial_r [\log \sqrt { |{\tilde g^k|}] }\partial_r+\dfrac{1}{r^2 \sqrt {|{\tilde g}^k|}}\partial_{\theta^i}({\tilde g}^{\theta^i \theta^j}\sqrt { |{\tilde g}^k|}\partial_{\theta^j}) . $$

On a,

$$ -\Delta = \partial_{rr}+\dfrac{n-1}{r}\partial_r+ \partial_r \log J(x,r,\theta)\partial_r+ \dfrac{1}{r^2 \sqrt {|{\tilde g}^k|}}\partial_{\theta^i}({\tilde g}^{\theta^i \theta^j}\sqrt { |{\tilde g}^k|}\partial_{\theta^j}) . $$

Le  laplacien s'\'ecrit ( d\'ecomposition radiale et angulaire),

$$ -\Delta = \partial_{rr}+\dfrac{n-1}{r} \partial_r+\partial_r [\log J(x,r,\theta)] \partial_r-\Delta_{{S}_r(x)}, $$

o\`u $ \Delta_{{S}_r(x)} $ est le laplacien sur la sph\`ere $ {S}_r(x) $. 

\smallskip

$ \Delta_{{S}_r(x)} $, s'\'ecrit localement,

$$ -\Delta_{{S}_r(x)}=\dfrac{1}{r^2 \sqrt {|{\tilde g}^k|}}\partial_{\theta^i}({\tilde g}^{\theta^i \theta^j}\sqrt { |{\tilde g}^k|}\partial_{\theta^j}) . $$

On coclut que l'op\'erateur de droite en $ \theta $, dans la premi\`ere expression du laplacien, n'est que l'\'ecriture locale d'un op\'erateur global.

\smallskip

On pose, $ L_{\theta}(x,r)(...)=r^2\Delta_{{S}_r(x)}(...)[\exp_x(r\theta)] $.

\smallskip

L'op\'erateur $ L_{\theta}(x,r) $ est un op\'erateur d\'efinie sur les fonctions $ {C}^2({\mathbb S}_{n-1}) $ de mani\`ere globale (il ne d\'epend pas de la carte choisie sur $ {\mathbb S}_{n-1} $ et vaut localement,

$$ L_{\theta}(x,r) =-\dfrac{1}{\sqrt {|{\tilde g}^k|}}\partial_{\theta^i}[{\tilde g}^{\theta^i \theta^j}\sqrt { |{\tilde g}^k|}\partial_{\theta^j}] . $$

De plus cet op\'erateur est elliptique, on \'ecrit alors,

$$ \Delta = \partial_{rr}+\dfrac{n-1}{r} \partial_r+\partial_r [ \log J(x,r,\theta)] \partial_r - \dfrac{1}{r^2} L_{\theta}(x,r) . $$

Si, maintenant, $ u $ est une fonction d\'efinie sur $ M $, alors, $ \bar u(r,\theta)=uo\exp_x(r\theta) $ est la fonction qui lui correspond en coordonn\'es polaires centr\'ees en $ x $.

$$ -\Delta u=\partial_{rr} \bar u+\dfrac{n-1}{r} \partial_r \bar u+\partial_r [ \log J(x,r,\theta)] \partial_r \bar u-\Delta_{{ S}_r(x)}(u_{|{S}_r(x)})[\exp_x(r\theta)], $$

$$ r^2\Delta_{{S}_r(x)}(u_{|{S}_r(x)})[\exp_x(r\theta)]=-\dfrac{1}{ \sqrt {|{\tilde g}^k|}}\partial_{\theta^i} \left [ {\tilde g}^{\theta^i \theta^j}\sqrt { |{\tilde g}^k|}\partial_{\theta^j} [uo\exp_x(r\theta)] \right ] =L_{\theta}(x,r)\bar u . $$

Ainsi,

$$ -\Delta u =\partial_{rr} \bar u+\dfrac{n-1}{r} \partial_r \bar u+\partial_r [ \log J(x,r,\theta)] \partial_r \bar u-\dfrac{1}{r^2}L_{\theta}(x,r)\bar u . $$

D'autre part, il est clair qu'en se placent dans la carte exponentielle, on transporte la m\'etrique dans l'espace euclidien. Puis, en consid\'erant les coordonn\'ees g\'eod\'esiques polaires, on voit que la m\'etrique se decompose suivant les champs de vecteurs radial et angulaires. Comme l'espace se decompose suivant la partie radiale et angulaire sur la sphere, il est clair qu'en consid\'erant les champs de vecteur $ \partial_{r}, \partial_{\theta^j} $, on retrouve localement le laplacien radial et angulaire, qui est un laplacien global par la carte exponentielle. On retrouve par l'interm\'ediaire de la carte exponentielle un laplacien global qu'on decompose en coordonnnees g\'eod\'esiques polaires en partie radiale et angulaire par le lemme de Gauss.

\smallskip

L'op\'erateur $ L_{\theta}(x,r) $ est un laplacien sur $ {\mathbb S}_{n-1} $ pour une certaine m\'etrique qui d\'epend du param\`etre $ r $.

\smallskip

On a, $ \Delta_{{S}_r(x)}=\Delta_{i_{x,r}^{*}(g)} $ o\`u, $ i_{x,r} $ est l'application identit\'e de $ {S}_{r}(x) $ dans $ M $ et $ i_{x,r}^*(g)=\tilde g $ est la m\'etrique induite sur la sous-vari\'et\'e $ {S}_r(x) $ de $ M $.

\smallskip

Comme $ \exp_x $ induit un diff\'eomorphisme de $ {S}_r(x) $ vers $ {\mathbb S}_{n-1}^r $, on a:

$$ \Delta_{{S}_r(x)} u_{| {S}_r(x)}=\Delta_{i_{x,r}^*(g)} u_{| {S}_r(x)}=\Delta_{\tilde g} u_{| {S}_r(x)}=\Delta_{\exp_x^*(\tilde g), {\mathbb S}_{n-1}^r} uo\exp_x(v), $$

Soit, $ \tilde z $ l'application de $ {\mathbb S}_{n-1} $ vers $ {\mathbb S}_{n-1}^r $ d\'efinie par $ \theta \to r \theta $. Alors,

 $$ \Delta_{\exp_x^*(\tilde g)} uo\exp_x(v) =\Delta _{\tilde z^* [\exp_x^*(\tilde g)]} uo\exp_x(r\theta), $$

Pour le choix de la carte $ (\psi, U^k) $, sur $ {\mathbb S}_{n-1} $, la carte locale en coordonn\'ees polaires pour le point $ x $  est, $ (\phi_0, ]0,\epsilon_0[\times U^k) $, o\`u $ \phi_0= \exp_xo\tilde zo\psi^{-1} $. Les $ g_{jl}^k $ sont par d\'efinition:

$$ g_{jl}^k(r,\theta)=g_{\exp_x(r\theta)}(\partial_{j,\phi_0}, \partial_{l,\phi_0}), $$

avec, $ \partial_{j,\phi_0}, \partial_{l,\phi_0} $ les champs de vecteurs canoniques correspondant \`a la carte $ (\phi_0, ]0,\epsilon_0[\times U^k ) $.

\smallskip

Par d\'efinition, pour $ h $ une fonction d\'efinie au voisinage de $ x $:

$$ \partial_{j,\phi_0}(h)=\dfrac{\partial (ho\phi_0)}{\partial \bar \theta_j}=\dfrac{\partial (ho\exp_xo\tilde zo\psi^{-1})}{\partial \bar \theta_j}, $$ 

Rappelons, qu'on ecrit, $ \psi(\theta)=(\bar \theta_1,\ldots,\bar \theta_{n-1}) $ et $ \theta_0 = r $ (d\'erivations angulaires et radiale).

Si, on s'occupe des composantes angulaires, c'est \`a dire les champs de vecteurs, $ \partial_{j,\phi_0} $ avec $ 1 \leq j\leq n-1 $, $ j \in {\mathbb N} $, alors,

$$ \partial_{j,\phi_0}(h)=\dfrac{\partial (hoi_{x,r}o\exp_xo\tilde zo\psi^{-1})}{\partial \bar \theta_j}, $$

$$ \dfrac{\partial (ho\phi_0)}{\partial \bar \theta_j}=\dfrac{\partial (hoi_{x,r}o\phi_0)}{\partial \bar \theta_j}, $$

Donc, si on note $ \bar \partial_j, \bar \partial_l $ les champs de vecteurs canoniques de la sph\`ere par rapport \`a la carte $ (\psi, U^k) $, on obtient, en notant $ d $ d\'esigne la diff\'erentielle,

$$ \partial_{j,\phi_0}=d(i_{x,r}o\exp_xo\tilde z)(\bar \partial_j), $$

Donc,

$$ g_{jl}^k(r,\theta)=g_{\exp_x(r\theta)}[d(i_{x,r}o\exp_xo\tilde z)(\bar \partial_j), d(i_{x,r}o\exp_xo\tilde z)(\bar \partial_l)], $$

d'o\`u, en utilisant la d\'efinition du pull-back,

$$ g_{jl}^k(r,\theta)=\tilde z^*[\exp_x^*[i_{x,r}^*(g)]](\bar \partial_j,\bar \partial_l)=\tilde z^*[\exp_x^*(\tilde g)](\bar \partial_j,\bar \partial_l), $$

le dernier terme \`a droite dans l'expression pr\'ec\'edente, n'est que la composante $ j,l $, de la m\'etrique $ \tilde z^* [\exp_x^*(\tilde g)] $ d\'efinie sur la sph\'ere unit\'e $ {\mathbb S}_{n-1} $.

Finalement, on peut comparer les expressions locales et globales,

$$ -\Delta_{\tilde z^*[\exp_x^*(\tilde g)], {\mathbb S}_{n-1}} uo\exp_x(r\theta)=\dfrac{1}{r^2 \sqrt {|{\tilde g}^k|}}\partial_{\theta^i}[{\tilde g}^{\theta^i \theta^j}\sqrt { |{\tilde g}^k|}\partial_{\theta^j} uo\exp_x(r\psi^{-1})], $$

Maintenant, si on prend  la nouvelle m\'etrique sur $ {\mathbb S}_{n-1} $ d\'efinie par, 

$$ g_{{}_{x,r,{\mathbb S}_{n-1}}}=r^{-2} \tilde z^*[\exp_x^*(\tilde g)], $$

cette m\'etrique est bien d\'efinie et de plus,

$$ {g_{{}_{x,r,{\mathbb S}_{n-1}}}}_{jl}=r^{-2}g_{jl}^k=\tilde g_{jl}^k, $$

Nos raisonnements futurs se ferons avec $ r > 0 $. On obtient,

$$ \Delta_{g_{x,r, {}_{{\mathbb S}_{n-1}}}}=r^2 \Delta_{\tilde z^*[\exp_x^*(\tilde g)], {\mathbb S}_{n-1}}=L_{\theta}(x,r). $$

\smallskip

{\bf Remarque:} Dans le cas plat, un ouvert de $ {\mathbb R}^n $, on n'a pas besoin de tout cela. D\`es qu'on a une boule, on peut d\'efinir les coordonn\'ees polaires. On ecrit $ x=a_i+r\theta $ en polaires. Alors, on passe aux diff\'erentielles, on obtient: $ dx=dr\theta+rd\theta $ et quand on prend la norme euclidienne, $ dx^2=dr^2+r^2d\theta^2 $, on voit bien que $ d\theta^2$ est la metrique euclidienne restreinte \`a la sphere, qui est la metrique de la sphere. On obtient alors la formule du Laplacien en polaires dans le cas euclidien, un ouvert de $ {\mathbb R}^n $ (formule usuelle dans des cartes, ici, la carte polaire, $ (r,\theta) $).($ \Delta = \partial_{rr} +\frac{n-1}{r} \partial_r +\frac{1}{r^2} \Delta_{\theta} $, $ \Delta_{\theta} $ est la laplacien de la metrique de la sphere  $ g_{ij}(\theta)$, $ d\theta^2=g_{ij}(\theta)d\theta^i d\theta^j $).

\smallskip

///////////////////////////////////////////////////////////////

\smallskip

{\bf Sur l'application exponentielle:} Voir le livre de Hebey et Aubin. On alterne les notations du livre d'Aubin et de Hebey.

\smallskip

1) le flot exponentiel: Soit $ x \in M, X\in T_x(M) $, on considere une carte normale. Ceci permet d'avoir un repere orthonorm\'e plutard et identifier l'espace tangent \`a $ {\mathbb R}^n $. Le flot donne l'exponentielle: $ \exp_x(X)=\gamma_x(1,X) $.

\smallskip

2) Theoreme des eq.differentielles: on fait varier $ x $ dans $ M $. La dependance en $ x $ de $ \gamma_x $ est $ C^{\infty} $ ou reguliere.

\smallskip

3) Theoreme d'inversion locale, on determine un ouvert $ U $ tel que: $ (x, X) \to (x\exp_x(X)) $ est reguliere et inversible.

\smallskip

4) lemme 2.3.7 de Hebey: on utilise ce qu'on appelle, le parallelisme, on inscrit un $ n-$ cube contenant $ (0,0)=\Psi(x,0) $, en considerant la carte de trivalisation, carte locale, puis comme $ U $ est ouvert, son image par l'application de carte $  \Psi $, est un ouvert, dans le quel on inscrit un $n-$cube, puis on revient, par image reciproque, $\Psi^{-1} $, la nature de $ \Psi $, cylindrique, permet d'inscrire un ouvert dans $ U $ du type $ S_x(\epsilon) \times B_0(\epsilon) $:

$$ \Psi: \cup_y T_y(M) \to {\mathbb R}^n \times {\mathbb R}^n, \,\, U\subset \cup_y \{Y \in T_y(M), ||Y||\leq r_y, r_y >0 \} \subset \cup_y T_y(M), $$

parallelisation et topologie de l'union disjointe.

\smallskip

5) Donc, on a: $ \exists  \epsilon >0, S_x(\epsilon) \times B_0(\epsilon) \subset U $ tel qu'on ait le lemme 2.3.7 de Hebey: $ \forall y \in S_x(\epsilon) $, $ \exp_y $ realise un diffeomorphisme de $ B_0(\epsilon) $ sur $ S_y(\epsilon) $.

\smallskip

On ecrit, ce qu'on appelle, le lemme des voisinages geoedesiquement etoil\'es:

$$ \forall x\in M, \exists \epsilon >0, \forall \mu<\epsilon, \forall y \in S_x(\mu), S_y(\mu) \,\, {\rm est \,\, geodesique}. $$

Le lemme de Zorn, dit qu'il y a pour chaque $ x\in M $, un $\epsilon_x >0 $, $ \alpha_x/2 \geq \epsilon_x \geq \epsilon >0 $, maximal. (Il y a un sup pour les boules ferm\'ees et un autre sup pour les boules ouvertes. le fait de considerer les boules ferm\'ees, cela permet d'avoir de l'espace, car elles sont dans un compact: $ \alpha_x/2 $. Notations du livre d'Aubin). 

\smallskip

On prend, $ y $ dans un petit voisinage de $ x $ et vice versa, pour prouver, \`a l'aide de la propri\'et\'e precedente, avec la notions de limite inf et limite sup, pour des $ \mu $ particuliers (et la maximalit\'e de $ \epsilon_x $), $ \mu=(1-\epsilon') \epsilon_x >0 $ et $ \mu=(1-\epsilon') \epsilon_y >0 $, ($0 < \epsilon'\to 0 $), que la fonction: $ x\to \epsilon_x $ est continue.

\smallskip

////////////////////////////////////////////////////////////////////////////////

\smallskip

24) Sur les fonctions Green:

\smallskip

Le monograph de Druet-Hebey-Robert est assez clair et precis. On donne quelques explications sur les fonctions de Green et sur le fait que le dernier terme dans la decomposition de cette fonction est de classe $ C^1 $. C'est essentielement du au fait que le premier terme est donn\'e explicitement et on utilise une recurrence pour prouver que les termes dela decompostion sont Sobolev et continues en dehors de la digaonale. (Ou bien comme dans le monograph de Druet-Hebey-Robert, on fait la difference des fonctions en deux points et on remarque qu'on a des fonctions Lipschitziennes, Notons comme c'est ecrit dans Ambrosio-Fusco-Pallara,pour pouvoir savoir si une fonction est Sobolev, il suffit qu'on ait des derivees directionelles).

\smallskip

Les fonctions qu'on considere sont continues en dehors d'un ensemble de mesure nulle (un point ou la diagonale). Elles sont mesurables. Elles sont integrables par les theoremes de comparaisons.

\smallskip

Soit $ (\bar W, g) $ une vari\'et\'e Riemannienne compacte avec ou sans bord. Dans le cas sans bord, on la note $ (M, g) $.

\smallskip

Soit $ H $ la fonction definie comme dans le livre d'Aubin. Cette fonction ne depend pas de la fonction rayon d'injectiv\'e dans le cas d'une vari\'et\'e compacte sans bord. Elle depend (\`a gauche en P) du rayon d'injectiv\'e dans le cas d'une vari\'et\'e \`a bord. On sait que la fonction rayon d'injectivit\'e est continue sur $ W $ avec $ \delta(P)\leq d(P,\partial W) $, on peut alors la prolonger par $ 0 $ sur $ \partial W $. Alors $ \delta(P) $ est continue sur $ \bar W $.

\smallskip

On pose:

$$ \Gamma_1(P,Q)=-\Delta_Q H(P,Q),\,\, \Gamma_{i+1}(P,Q)=\int_W \Gamma_i(P,R)\Gamma_1(R,Q) dR. $$

Alors,

$$ |\Gamma_1| \leq C [d(P,Q)]^{2-n}. $$

Explicitement et par recurrence les $ \Gamma_i $ sont continues en dehors de la diagonale et pour $ P $ ou $ Q $ proche du bord, elles sont nulles (conditions au bord). On les prolonge par $ 0 $ au bord. Par le theoreme de la convergence dominee de Lebesgue, on prouve qu'elles sont $  W^{1,\infty}  $ en dehors de la diagonale en $ Q $ car elles sont Sobolev est leur derivees bornees. La fonction $ d(P,Q)=d_P(Q) $ est reguliere en $ Q $ en dehors de $ P $ et en fixant $ Q $ elle est Lipschitzienne en $ P $, on alors le fait que $ (R, Q) \to \partial_Q d_R(Q) $ est continue en dehors de la diagonale et $ L^{\infty}_{loc} $. Donc, par le theoreme de convergence dominee (enlever de petites boules) $ \Gamma_i* \partial_Q \Gamma_1 $ et $ \Gamma_i*\partial_Q H $ sont $ C^0 $ en dehors de la diagonale. Donc, les $ \Gamma_i $ et $ \Gamma_i*H $ sont $ C^1 $ en dehors de la diagonale.(Le point essentiel est la fonction distance).(Concernant la fonction distance, $ (R_i, Q) \to d(R_i,Q)=d(R_i, exp_{R_i}(r\theta)), Q=exp_{R_i}(r\theta),  \epsilon_0/2 \leq r\leq \delta(R_i)/3, \theta\in{\mathbb S}_{n-1}$ et la carte de reference est $ [B(R_0,\delta(R_0)/3),  exp_{R_0}] $ avec $ R_i \to R_0 $. On voit alors qu'en coordonnes geodesiques polaires centrees en $ R_i $ un calcul de $ \Delta_{g,Q} d(R_i, Q) \in L^{\infty} $ et $ d(R_i,.) \to d(R_0,.) $ et par les estimations elliptiques on a $ \partial_Q d(R_i, Q) \to \partial_Q d(R_0, Q) $ pour $ d(R_i,Q) \geq \epsilon_0 >0 $.)

\smallskip

Dans le cas o\`u $ W=M $ sans bord:

\smallskip

D'apres sa formule explicite $ H $ est $ C^1 $  en dehors de la diagonale. Elle est $ W^{1,1+\epsilon} $ de chaque variable pour un $ \epsilon >0 $ assez petit, par recurrence et par Fubini, Fubini-Tonnelli:

$$ \partial_P\Gamma_{i+1}(P,Q)=\int_M \partial_P\Gamma_i(P,R)\Gamma_1(R,Q)dR, $$

et,

$$ \partial_Q\Gamma_{i+1}(P,Q)=\int_M \Gamma_i(P,R)\partial_Q\Gamma_1(R,Q) dR.$$

Pour une variet\'e avec bord, on ne peut pas deriver par rapport $ P $, car l'expression de $ \Gamma_i $ et en paticulier de $ \Gamma_1 $ contient $ \delta_P $ qu'on sait seulement continue et pas forcement differentiable. Dans le cas avec bord, les fonctions $ \Gamma_i $ sont continues en $ P $ en dehors de la diagonales et $ L^{1+\epsilon}, \epsilon >0 $ seulement, mais ceci suffit pour prouver la symetrie par le theoreme d'Agmon) 

\smallskip

On a aussi, les estimees de Giraud. Ainsi de proche en proche, les $ \Gamma_i $ deviennent de plus en plus regulieres. on utilise le theoreme de Giraud pour prouver que $ \Gamma_k = r^{\lambda-n}* r^{2-n} $ avec $ \lambda + 2 >n $ (donc $\geq n+1 $ car $ \lambda $ est entier puisqu'on derive les $ \Gamma_i $). ( Si $\lambda >n-1 $ ($ \lambda \geq n $) c'est fini (on regarde $ \Gamma_{k-1} $), si $ \lambda =n-1 $ alors: $ \Gamma_k = r^{-1}* r^{2-n}$ et $ \partial_P \Gamma_k=r^{-2}*r^{2-n} $ et $ \partial_Q \Gamma_k=r^{-1}*r^{1-n} $ qui est en $ \log $ par Giraud et donc $\partial_P \Gamma_{k+1}$ et $ \partial_Q \Gamma_{k+1} $ sont $ C^0(M\times M) $ par Giraud et donc $ \Gamma_{k+1} $ est $ C^1(M\times M)$. (Utiliser le produit de convolution, car $\Gamma_{k+1} $ et sa deriv\'ee Sobolev, est continue pour ecrire qu'elle est $ \int \partial_P $, l'integrale d'une fonction continue donc $ C^1 $, il a convergence uniforme de la convolution, car la fonction et son Sobolev sont continues). Dans le cas avec bord, $ \Gamma_{k} $ continue implique que $ \Gamma_{k+1} $ est $ W^{1,\infty} $ en $ Q $.(On a mieux $ C^1 $ en $ Q $ car la fonction $ (R,Q) \to \partial_Q d(R, Q) $ est $ C^0 $ en dehors de la diagonale de $ W $  et $ H $ est a support compact (elle s'annule avant le bord) et donc ($ \Gamma_k* \partial_Q H $ est $ C^0 $ en dehors de la diagonale de $ \bar W $)).

\smallskip

Dans le cas d'un operateur $ -L=-\Delta+h $, la regualrite de $ \Gamma_{k+1} $ depend de celle de $ h $. Si par exemple $ h $ est $ C^1 $ alors comme precedemment $ \Gamma_{k+1} $ est $ C^1 $ de $ P $ et $ Q $. Bien sur le $ \Gamma_1 $ change:

$$ \Gamma_1(P,Q)=-\Delta_QH(P,Q)+h(Q)H(P,Q). $$

{\bf Remarque:} a) sur la borne inferieure de la fonction de Green $ G_{\epsilon} $ de $ -\Delta + \epsilon = -\nabla^i(\nabla_i)+\epsilon $, $ \epsilon >0 $.

\smallskip

1) On a $ G_{\epsilon} \geq 0 $ donc par la propri\'et\'e de la borne inferieure $ m_{\epsilon}= \inf_{M\times M} G_{\epsilon} $ existe.

\smallskip

2) Par definition $ m_{\epsilon}= \lim_{\delta \to 0} G_{\epsilon} (x_{\delta}, y_{\delta}) $. Si on suppose que $ d(x_{\delta},y_{\delta}) \to 0 $, par la propri\'et\'e du monograph de Druet-Hebey-Robert, $ m_{\epsilon}=+\infty $, ceci n'est pas possible et Comme $ G_{\epsilon} $ est continue en dehors de la diagonale. $ m_{\epsilon}=G_{\epsilon}(x_{\epsilon},y_{\epsilon}) $ avec $ x_{\epsilon} \not = y_{\epsilon} $.

\smallskip

Si on met un indice , $ i $, alors on peut supposer que, $ G_i=G_{\epsilon_i} $:

$$ \inf_{M\times M} G_i=m_i=G_i(x_i,y_i),\,\, x_i\not = y_i, $$

Apres, on peut raisonner sur la borne inferieure globale sur tout $ M\times M $, au lieu de prendre l'inf sur $ M $ en $ x $ puis en $ y $. 

\smallskip

b) sur les fonctions de Green des vari\'et\'es \`a bord. on remarque que par exemple que pour les varietes \`a bord, la parametrix est solution d'une EDP elliptique en dehors de la singularit\'e, on utilise les theoremes d'Agmon et de regularit\'e dans la construction complete de la fonction de Green dans le cas des varietes a bord, voir le monograph de F.Robert, c'est essentielement l'utilisation des noyaux, la formule sommatoire et les th\'eor\`emes d'Agmon.

\smallskip

///////////////////////////////////////////////////////////////////////////////////////////////////////////////////////////////////

\smallskip

25) Sur la formule de Gauss-Bonnet:

\smallskip

C'est bien ecrit dans le livre de J.Lafontaine, edition 2010.

\smallskip

1) relation entre l'indice d'un champ de vecteurs et la courbure de Gauss.

\smallskip

2) relation entre indice d'un champ de vecteur, fonction de Morse et caracteristique d'Euler-Poincar\'e.

\smallskip

On a une surface riemannienne compacte sans bord $ (M,g) $.

\smallskip

On part d'un champ de vecteur $ X $ ayant des zeros. On consid\`ere le normalis\'e $ X_1=X/||X|| $, il n'est pas definit aux zeros de $ X $. Puisqu'il y a produit scalaire, on considere $ X_2 $ son orthoganle de sorte qu'on ait un repere orthonorm\'e.

\smallskip

$ (X_1, X_2) $ est un repere mobile ou "moving-frame". On a une connexion, donc on a de nouveaux symboles de Christoffel et une connexion forme $ \omega $. de sorte que $ d\omega =\Omega $ la courbure. celle ci n'est pas definie aux zeros de $ X $, puisqu'on considere $ X_1 $ et on divise par $ ||X|| $.

\smallskip

(Comme on l'a dit pour la surface $ K3 $, on a une connexion $\Rightarrow $ symboles de Christoffel $\Rightarrow $ connexion forme (par dualit\'e). Puis on derive $ d\omega $ on retrouve la courbure ou forme volume par dualit\'e. C'est bien ecrit dans le livre de J. Lafontaine quand il considere $ \nabla_{X_i}(X_j) $. Dans certains livres on note $ \omega =\omega_1^2 $. On peut le voir comme: quand on calcule $ \nabla_i dx^j=-\Gamma_{ik}^j dx^k $ pour les champs de vecteurs et les formes differentielles usuelles, ici pour $ X $ la forme $\omega $, on calcule $ \nabla_X \omega_i =-\omega_i^k(X) \omega_k $, changements de bases, nouveaux symboles de  Christoffel. )

\smallskip

Mais au vosinage des zeros on a une carte (canonique si on veut, normalis\'e) et le passage des champs issus de cette carte not\'es $ (Z_1,Z_2) $ \`a $ (X_1, X_2)$ se fait via un champ de rotations $ R=e^{i\phi} $ (plus precisement de $ Z_1 $ \`a $ X_1 $, on voit bien que $ X_1=e^{i\phi} Z_1 $ se contracte vers un vecteur du type $ e_1 \times e^{i\phi} $ car $ Z_1 $ est definit sur un voisinage de la singularit\'e et qu'on peut le contracter en un vecteur fixe). (On voit ici que $ R $ est homotope a $ X_1 $). On note $ \omega_0 $ la connexion forme pour ce repere mobile ou "moving-frame" $ (Z_1, Z_2) $. Alors on a la relation:

$$ \omega=\omega_0-d\phi. $$

Avec, $ \omega_0 $ definie au voisinage de la singularit\'e contrairement \'a $ \omega $. (on peut passer \'a la limite dans la formule de Stokes, voir ci-dessous).

\smallskip

La rotation $ R $ est homotope au champ $ X_1 $ donc: on a le meme degre topologique.

$$ deg (R)=deg (X_1). $$

Or,

$$ deg (R)= deg (rotation)= \int d\phi, $$

et par defintion,

$$ deg(X_1)=indice (X). $$

(Pour une fonction de Morse $ f $, on considerera $ X=\nabla f $ et l'indice se calcul explicitement dans une carte particuliere, car il ne depend pas de la carte).

\smallskip

Maintenant la formule de Stokes donne, c'est bien ecrit dans le Lafontaine:

$$ \int \Omega= \int d\omega =\int \omega \to \int d\phi=indice (X). $$

Ceci pour la relation entre courbure et indice de champs de vecteurs. (c'est li\'e a la connexion forme, qui induit la courbure et qui est produite par le repere mobile issu de $ X $).

\smallskip

Maintenant un theoreme, voir le Gramain ou Hirsch (Differential Topology), dit que sur une surface compacte, il y a une fonction de Morse $ f $. On applique, ce qui pr\'ec\`ede \`a $ \nabla f $. Ceci revient \`a calculer un nombre li\'e aux differents indice le la fonction de Morse et qui est bien detaill\'e dans le livre de Hirsch (Differential Topology). Ceci revient \`a calculer l'homologie d'un complexe cellulaire et qui donne la caracteristique d'Euler-Poincar\'e.

\smallskip

Pour ce qui est du calcul de la caractersitique d'Euler d'une somme connexe. C'est bien connu. Voir aussi dans le livre de J. Lafontaine. 

\smallskip

/////////////////////////////////////////////////////////////////////////////////////////////////////////////////////////

\smallskip

26) Equation de Liouville: Pourquoi, lorsque la courbure de Gauss d'une surface est constante, $ K $ alors $ K>0$ la surface est localement isometrique a la sphere (un ouvert de la sphere), $ K =0 $ elle est isometrique au plan (un ouvert du plan), $ K <0 $ elle est isometrique a un ouvert de l'espace hyperpbolique ?

\smallskip

En effet soit $ S $ une surface de metrique $ g $, alors, on peut definir des coordonn\'ees geodesiques polaires, $ (r, \theta) $, dans ses coordonn\'es $ g=dr^2+G(r,\theta) d\theta^2 $ . Alors la courbure de Gauss verifie l'equation differentielle suivante $ \partial_{rr} (\sqrt G)+ K\sqrt G=0 $, avec les conditions a l'origine $ \sqrt G(0,\theta)=0 $ et $ \partial_r(\sqrt G(0,\theta))=1 $. Clairement lorsque $ K $ est une constante , cette equation differentielle s'integre facilement et donne la metrique de la sphere ou du plan ou du plan hyperbolique. ($ (\sin r)^2 $, ou $ r^2 $ ou $ (\sinh r)^2 $). 

\smallskip

Si on note $ \phi: (r, \theta) \to \exp_x(r\theta) $, la carte definissant les coordonn\'ees geodesiques polaires pour la surface $ S $, on ecrit alors,  $ \phi^*(g)=dr^2+G(r,\theta) d\theta^2 =\psi^*(g_0) $, avec $ \psi $ la carte definisssant les coordonnes geodesiques polaires, selon qu'on soit sur la sphere, le plan, ou le plan hyperbolique (a partir de la carte exponentielle, $ \psi: (r, \theta) \to \tilde \exp_y(r\theta) $).

\smallskip

Alors, $ \phi, \psi $, qui sont des cartes, sont des diffeomorphisme avec $ (\phi o \psi^{-1})^*(g)=g_0 $ et $ \phi o \psi^{-1} $ est l'isometrie cherch\'ee.(Par exemple si $ K >0 $, $ \phi: \Omega_1 \subset S \to ]0, r_0[ \times ]a,b[, ]a,b[\subset {\mathbb S}_1 $ et $ \psi : \Omega_2 \subset {\mathbb S}_2 \to ]0, r_0[ \times ]a,b[, ]a,b[\subset {\mathbb S}_1 $, sont deux cartes dont les images sont les memes dans $ ]0, r_0[ \times ]a,b[, ]a,b[\subset {\mathbb S}_1 $ et $ \psi^{-1} o \phi : \Omega_1 \subset S\to  \Omega_2 \subset {\mathbb S}_2 $).

\bigskip

27) Sur le flot de Ricci et la conjecture de Poincar\'e: l'id\'ee est qu'on d\'eforme une metrique, a partir de l'equation du flot de Ricci pour obtenir une metrique de courbure de Ricci strictement positive, puis on r\'eutilise le flot de Ricci pour deformer la metrique pour obtenir une metrique de courbure sectionelle constante positive, (la vari\'et\'e est simplement connexe et compacte), elle est alors diffeomorphe a la 3-sphere:

\smallskip

a) On part du Resultat de Gao et Yau ou Lohkamp, il existe toujours une metrique de courbure de Ricci strictement negative. Puis on utilise le flot de Ricci.

\smallskip

b) Par le blow-up et les estimations de Hamilton, on deforme la metrique pour obtenir une metrique de courbure de Ricci positive ou nulle (blow-up, a la limite). Ici, c 'est dit dans le monograph de John-W Morgan et Gang Tian (Ricci flow and The Poincar\'e conjecture), on \'etend petit a petit le flot de Ricci par chirurgie (le travail de Perlman). La conjecture de Geometrisation. On passe d'une vari\'et\'e de courbure negative strictement \`a une vari\'et\'e de courbure positive ou nulle. (on change de metrique). (heuristiquement, si on note $ R_k =\max R_t $ avec $ R_t $ la courbure scalaire, en divisant par $ R_k$, on aurait, $ Ricci < b(n) <0 \Rightarrow \dfrac{Ricci}{R_k}>\dfrac{b(n)}{R_k} \to 0 $, si $ k \to + \infty $. En faisant le blow-up, a) on change d'echelle, b) il faut "recoller" les morceaux a la limite, pour avoir une vari\'et\'e et une metrique reguliere qui sera de courbure sectionelle positive ou nulle, c) la vari\'et\'e en question est not\'ee $ {\mathbb M} $ de metrique $ g $, dependant du temps et de l'espace et obtenue par chirurgie (\`a la limte), telle que $ {\mathbb M}(0)={\mathbb M}_{\infty}(0) = M $ et $ g(0)=g_{\infty}(0) $, avec $ M $ la vari\'et\'e de d\'epart sur laquelle on a une nouvelle metrique, $ g(0)=g_{\infty}(0) $ de courbure positive ou nulle.).

\smallskip

(Dans le cas blow-up avec $ R_k \not \to +\infty $, voir l'article de Hamilton, r\'eference 36 dans Morgan-Tian, on ne peut pas avoir collapse par les th\'eor\'emes de comparaison, voir Morgan-Tian (derniers chapitres, 15-16-17), on ne peut pas avoir ($ 0\not = \pi_1(H_1)\subset ou =\pi_1(M)=0 $), $ H_1 $ non compacte  est diffeomorphe a $ M $ qui est compacte, ce n'est pas possible, voir l'article de Hamilton, r\'eference 36 de Morgan-Tian et on obtient alors la convergence. Notons le flot de Ricci ici est defini sur $(0, +\infty) $, sinon , on l'etend petit a petit, comme pour le Theoerme des equations differentielles, de Cauchy, ou le Theoreme des Geodesiques avec en plus le tenseur de $ Riemann $ uniformement born\'e (derniers chapitres de Morgan-Tian, pour la non collapse et donc convergence dans ce cas $ R_k \not \to +\infty $). C'est ce qu'on a ecrit avant, dans l'ecriture ci-dessus entre $ Ricci, R_k, b(n) $, c'est la convergence vers une metrique de courbure de $ Ricci >0 $ partout.)

\smallskip

c) D\`es qu'on a une metrique de courbure de Ricci positive ou nulle, elle doit etre strictement positive quelque part, car comme on est en dimension 3 ($ Weyl=0 $), si $ Ricci=0 $, la vari\'et\'e ne peut pas etre compacte, car elle serait de courbure sectionelle nulle et simplement connexe, elle serait diff\'eomorphe \`a $ {\mathbb R}^3 $ (On a $ Riemann \geq 0$. Si $Riemann=0 \Rightarrow Ricci=0$ la vari\'et\'e compacte $ M $ serait diffeomorphe a l'espace euclidien, ce n'est pas possible. Donc $ \exists \, x_0 $ tel que $ Rimann_{x_0}\not =0 $ et $ \geq 0 $, en utilisant les valeurs propres de $ Riemann  \geq 0 $ qui devient alors symetrique ( Voir Morgan-Tian), il y a au moins une valeur propre $\lambda_{x_0} >0 $, comme $ Ricci_{x_0}=trace(Riemann_{x_0}) $ alors, (la trace ne depend pas de la base), $ Ricci_{x_0} =\lambda_{x_0}+\mu_{x_0}+\nu_{x_0} \geq \lambda_{x_0} >0$ et donc, $ Ricci_{x_0}(\theta, \theta) >0,\forall \theta \in {\mathbb S}_2$). Par un argument d'Aubin, il existe une m\'etrique conforme  telle que la courbure de Ricci est partout strictement positive.

\smallskip

d) On conclut gr\^ace au travail de Hamilton, on d\'eforme la metrique par le flot de Ricci jusqu'a obtenir une metrique de courbure sectionnelle constante positive. la vari\'et\'e est simplement connexe et compacte, c'est la 3-sphere.

\bigskip

28) Sur l'article de Chen-Li:

\smallskip

On a aussi la compacit\'e globale de Chen-Li, cit\'e dans la section 4.1. Le resultat de Chen-Li utilise le fait qu'on a compacit\'e au voisinage du bord lorsque $ ||\nabla \log V||\leq A $, puis il etend ce resultat lorsque $ ||\nabla V||\leq A_1 $, pour cela il utilise l'extension des resultats de Brezis-Merle (Theorem 1 de Brezis Merle), qui reste vrai dans des domaines Lipschitzien d\`es qu'on a la regularit\'e des solutions dans $ W_0^{1,2} $, car le principe du maximum est valable dans ce cas (on utilise l'integration par parties, qui est vrai des que la regularit\'e du bord est Lipschitzien. Pour l'inegalit\'e de Sobolev aussi et la resolution d'un probleme variationnel dans $ L^2 $. Lipschitz suffit). Ceci pour la fonction $ u_2 $ du debut de preuve de Chen-Li. (et aussi l'extension des fonctions harmoniques et la formule de Poisson, pour $ u_1 $ qui necessite une application conforme). Voir la preuve du corollaire de l'article de Chen-Li.

\smallskip

(Voir l'article de Sweers-Nazarov. Journal.Diff.Equations. 2007. Pour les conditions de regularit\'e $ W^{2,p} $ des problemes sur des ouverts Lipschitziens.)

\smallskip

Remarque sur la preuve de Chen-Li : pour etendre la partie $ u_1 $ harmonique, il faut supposer le domaine analytique, pour pouvoir utiliser une transformation conforme qui reste invariante par le Laplacien. Principe de symetrisation de Schwarz. Donc, le resultat de compacit\'e reste vrai avec la regularit\'e smooth, lorsque on suppose  $ ||\nabla \log V||\leq A $. Mais la regularit\'e du domaine doit etre suppos\'ee analytique lorsqu'on passe a $ ||\nabla V||\leq A_1 $.

\smallskip

Dans leur preuve Chen-Li, utilisent le fait que l'operateur est invariant par appplication de carte, ceci est possible si cette application est conforme, elle preserve le Laplacien. Puis, symetrise la fonction en symetrisant un probleme de Dirichlet, puis soustrayent les valeurs aux bords et ils obtiennent l'image de $ v_1 $. Alors $ u_1 $ est l'image de $ v_1 $ par l'application de carte. Maintenant pour construire $ v_1 $ ils utilisent une symetrisation d'un probleme de Dirichlet, qui requiert les solutions dans $ W^{2,p}\cap C^2(B_{\epsilon}) \cap C^1(\bar B_{\epsilon}), p >2 $ (la formule de representation de Green reste valable, dans ce cas, voir la preuve dans Gilbarg-Trudinger). Puis, ils utilisent la formule integrale de Poisson (qui necessite d'avoir l'operateur Laplacien).

\smallskip

a) Pour utiliser la formule de Poisson, on conserve le Laplacien : transformation conforme $\phi $.

\smallskip

b) Ils symetrisent $ uo\phi $ ils obtiennet une fonction $ u_v \in C^1(\bar B_{\epsilon}(0)) \cap W^{2,p} $.

\smallskip

c) Ils resolvent : $ -\Delta v_1=-\Delta u_v $ avec condition de Dirichlet sur $ B_{\epsilon}(0) $.

\smallskip

d) Ils utilisent la formule integrale de Poisson pour $ v_1-u_v \in W^{2,p}\cap C^2(B_{\epsilon}) \cap C^1(\bar B_{\epsilon}), p >2 $.
Sur le bord, il n'y a que les valaurs de $ u $.  

\smallskip

(Ce travail revient \`a symetriser une fonction harmonique qui n\'ecessite le theoreme de symetrisation de Schwarz
, qui necessite une application conforme, donc un domaine de depart $ \Omega $ analytique).

\smallskip

{\bf Remarque:} sur le fait que la r\'egularit\'e du domaine dans (De Figueiredo-Lions-Nussbaum et Chen-Li) est au moins $ C^{2,\beta}, \beta >0 $. Soit $ \alpha_y=<\nu_1|\nu(y)> $, l'angle avec la normale qui doit etre $ >0 $, avec $ y\in \partial \Omega $, qu'on devrait obtenir pour appliquer la m\'ethode moving-plane. Ici, $ \nu(y) $ est la normale au bord en $ y $: Pour definir l'angle, il suffit d'avoir une tangente  a la courbe de $ y \to \nu(y) $, et donc une d\'eriv\'ee de $ \nu(y) $, (l'application de carte $ \phi $ est au moins deux fois d\'erivable). Pour avoir une branche continue de 'cones', il suffit que $ y\to (\nu(y))'$ soit continue, (l'application de carte $ \phi $ est $ C^2 $, $ \phi''$ continue). Pour que l'angle $\alpha_y >0 $ soit  strictement positif, il suffit que la tangente a la courbe de $ \nu(y) $ soit non nulle, il suffit alors d'avoir la regularit\'e $ C^{0,\beta}, \beta >0 $ de la tangente (forme arrondie: continue et arrondie: $ C^{0,\beta} $).(Donc, il suffit que l'application de carte soit $ C^{2,\beta}, \beta >0 $).

\bigskip

29) Sur l'article de C.C. Chen et C.S. Lin: (Comm. Pure. Appli. Math. 1997):

\smallskip

1) Application du principe du maximum et le lemme de Hopf

\smallskip

2) Fonction auxiliere pour attenuer la perturbation due a la variation de la courbure scalaire prescrite et les conditions aux bords.

\smallskip

3) Ils s'assurent de la positivit\'e de la fonction finale obtenue comme somme d'une fonction et la fonction auxili\`ere. La positivit\'e sert dans le raisonnement par l'absurde quand on prend le sup des $ \lambda $ tel que ..., le fait de prendre un sup signifie qu'au dela du sup les fonctions sont negatives ou nulles en des points, les conditions dans cet articles assurent que le point limite reste dans le domaine (car les fonctions sont strictement positives au voisinage du bord), ils appliquent le principe du maximum et (ou) le lemme de Hopf au point limite (point minimum, on peut parler de principe du minimum ou principe du maximum). (raisonnement par l'absurde ( $\sup \lambda < \lambda_1=-1/4 $) et application du principe du maximum et le lemme de Hopf au point limite qui est cens\'e etre \`a l'interieur du domaine).

\smallskip

De plus, les fonctions $ h_{\lambda} $ sont $ C^1 $, pour le voir, elles sont definies par la representation integrale de fonction de Green, elles verifient au sens des distributions et d'Agmon $ C^2_0 $, une EDP, ici aussi on utlise Fubini, Fubini-Tonelli et les definitions des fonctions de Green au sens $ C^2_0 $ pour les fonctions test, pour se ramener au Theoreme d'Agmon et conclure que $ h_{\lambda } $ sont $ C^1 $. (Ecrire $ \Delta (h^{\lambda}-h^{\mu})=\bar Q_{\lambda} - \bar Q_{\mu} $, pour la continuit\'e en $ \lambda $).  Aussi, $ h^{\lambda} \geq 0 $ sont nulles au bord $\Rightarrow \partial_{x_1} h^{\lambda} \geq 0 $ et par le lemme de Hopf $ \partial_{\nu} (w_{\lambda}-h^{\lambda}) <0 $, (principe du minimum), avec $ \partial_{\nu} =-\partial_e=-\partial_{x_1} $, donc, $ \partial_{x_1} w_{\lambda} > \partial_{x_1} h^{\lambda} \geq 0 $).

\smallskip

4) pour prouver que les courbures prescrites sont nulles au point "blow-up", ils utilisent un raisonnement du type $ K(y)=\int_0^1 <\nabla K(y_i+sM_i^{-2/(n-2)}I_{\delta}(y))|M_i^{-2/(n-2)} I_{\delta}(y)> ds $ puis font la difference:

$$ K(y)=\int_0^1 <[\nabla K(y_i+sM_i^{-2/(n-2)}I_{\delta}(y))-\nabla K(y_i)]|M_i^{2/(n-2)} I_{\delta}(y)> ds + $$

$$ + M_i^{-2/(n-2)} <\nabla K(y_i) | I_{\delta}(y)> $$

a) Dans le premier ensemble $ I_{\delta}(y)$ est born\'e car $ y $ est born\'e, et son symetris\'e aussi, car $ \lambda_0 \leq \lambda \leq -1/4 $.

\smallskip

b) Dans la preuve de l'estimation asymptotique $ M_i^{2/(n-2)} |\nabla K_i(y)|^{1/(\alpha -1)} \leq C, $, on utilise le raisonnement ci-dessus des que le gradient est holderien ou $ C^2 $. ( comme la difference ci-dessus et un raisonnement du type $ |\nabla K (y_i+M_i^{-2/(n-2)}I_{\delta}(y))-\nabla K (y_i) | \leq B M_i^{-2\alpha/(n-2)} |I_{\delta}(y)|^{\alpha} $ dans la formule de Taylor integrale (ci-dessus). Et apres les integrales du type Giraud. Dans la partie qu'ils considere, les $ M_i^{-2/(n-2)} I_{\delta} $ sont born\'ees.

\bigskip

30) Sur le probleme de Brezis-Merle: Pourquoi la condition sur le potentiel $ V $ est $ W^{1,\infty} $ pour obtenir la compacit\'e des solutions $ u >0 $ ? 

1) On a $ e^u $ et le potentiel $ V $ sont li\'es par dualit\'e $ \int V e^u <+\infty $, $ e^u\in L^1, V\in L^{\infty} $. on considere le potentiel $ V $ comme une forme lineaire, une distribution et donc de possibles deriv\'ees (entieres).

Par dualit\'e et le fait qu'on considere le probleme de maniere globale: on a: la condition maximale et globale $ e^u \in L^1=W^{0,1} $, par dualit\'e des Sobolev, le minimum requis pour le potentiel $ V $ est $ W^{1,\infty} $. Si on veut diminuer la r\'egularit\'e de $ V \in W^{1-s,\infty}, 1> s >0 $, qui est presque l'espace $(1-s)-$holderien, alors par dualit\'e, $ e^u $ doit etre $ e^u \in W^{s,1} > W^{0,1}=L^1 $, or la condition $ e^u \in L^1$ est maximale.

\smallskip

a) $ e^u \in L^1=W^{0,1} $: pas de deriv\'ee pour $ e^u $ et $ e^u \in L^1$ $\Rightarrow $ au moins une deriv\'ee pour $ V $ et $ V,\nabla V \in L^{\infty} $ $ \Leftrightarrow $ $ V \in W^{1,\infty} $ au minimum:

\smallskip

b) Si on veut diminuer la r\'egularit\'e de $ V $ et utiliser la dualit\'e et l'integration par parties, on part aussi du fait qu'on a la compacit\'e pour $ V \in W^{1,\infty} $:

$ \int \nabla V \cdot e^u <+\infty \Rightarrow \int \nabla^{1-s} V\cdot \nabla^s(e^u) <+\infty  \Rightarrow $ $ e^u \in W^{s,1} > W^{0,1}=L^1 $, or la condition $ e^u \in L^1 $ est maximale.

2) Le cas holderien depend du point et du nombre de points, il est local. On a la compacit\'e avec $ \int V e^u \leq bC <16\pi $ ou $ \int V e^u \leq bC < 24\pi$ et $ V $ holderien uniformement. Brezis et Merle donnent des exemples de $ (u,V) $ qui blow-up avec $ \int V e^u \leq bC <16\pi $ ou $ \int V e^u \leq bC < 24\pi $.

\bigskip

31) {\bf Remarques:} Du point de vue de la logique math\'ematique (donc du point de vue math\'ematique), sur le blow-up et la compacit\'e en dimension 2 et concernant l'article de Wei-Ye-Zhou (dans Annales de l'IHP, 2008):

\smallskip

i) On donne des conditions pour qu'un raisonnement blow-up au bord soit possible ( En plus, il y a des exemples comme celui de Brezis-Merle).

\smallskip

ii) On donne des conditions pour que le raisonnement mis en defaut par Wei-Ye-Zhou (voir point, iv, ci-dessous), de compacit\'e soit possible.

\smallskip

iii) On donne des conditions pour qu'un raisonnement par l'absurde soit possible.

\smallskip

iv) Dans Wei-Ye-Zhou (Annales de l'IHP, 2008), lorsque $ a= constante $ au voisinage du bord, les solutions sont born\'ees au voisinage du bord, par la m\'ethode directe, moving-plane, de, De Figueiredo-Lions-Nussbaum (sous-critique, 1982), Suzuki (1990), Chen-Li(1993) et Ma-Wei(2001), (Ma-Wei, en donnent une nouvelle preuve directe, dans convegence for Liouville Eq.). Il n'y a pas d'analyse blow-up a faire. Il n'y a pas de raisonnement par l'absurde \`a faire. Et dans la compacit\'e, on obtient \`a la limite $ 0=0 $, dans l'identit\'e de Pohozaev-Rellich, ce qui ne veut rien dire.

\smallskip

(Assertion, $ A: \exists U $ un voisinage du bord $ \partial \Omega $, $ \exists C >0 $, pour toute solution $u>0 $ de l'Eq avec condition de Dirichlet, $ ||u||_{L^{\infty}(U)} \leq C$. Wei-Ye-Zhou disent que $ A $ est vraie, ils ne peuvent pas supposer $ non A $ vraie, ils ne peuvent pas faire un raisonnement par l'absurde)

\smallskip

(Argument suppl\'ementaire \`a l'article de Wei-Ye-Zhou (Annales de l'IHP, 2008): Lorsqu'on met un potentiel $ V $: l'analyse blow-up est possible. Par contre la compacit\'e, elle est donn\'ee par Chen-Li, dans le cas $ V $ Lipschitzien, par la m\'ethode directe, moving-plane. Il n'y a pas de raisonnement par l'absurde \`a faire. Dans la compacit\'e, qui n\'ec\'essite $ V $ Lipschitzien, dans Wei-Ye-Zhou, on obtient \`a la limite $ 0=0 $, ce qui ne veut rien dire.) 

\smallskip

Ceci dans le cas $ a=constante $ au voisinage du bord et le potientiel constant. Wei-Ye-Zhou, mettent en defaut ce raisonnement par l'absurde. On donne des conditions pour que le raisonnement par l'absurde soit possible.

\smallskip

Dans le cas $ a\not=constante $, au voisinage du bord et le potentiel constant (le but de l'article de Wei-Ye-Zhou, Annales de l'IHP, 2008), la preuve n'est pas correcte, ce qui veut dire que le raisonnement par l'absurde n'est pas possible.

\smallskip

Concernant l'article: Note on the Chen-Lin result with the Li-Zhang method. L\`a aussi, {\underbar {on dit que}} le raisonnement par l'absurde est possible. C'est ce que prouvent Chen-Lin.

L'assertion: $ (A: \sup+\inf <+\infty)$. Alors: $(non\, A: \sup+\inf \to +\infty)$.

Chen-Lin font:

\smallskip

Etape 1: $ (non\, A: fausse)$ (Chen-Lin).

\smallskip

Etape 2: $ (non\, A: fausse) \Rightarrow (A: vraie)$.

\smallskip

Dans, Note on the Chen-Lin result with the Li-Zhang method: {\underbar {on dit que}} l'etape 1 est vraie ou possible; $ ((non \, A: fausse) : vraie)$.(hypoth\`ese $ non A: vraie $, implique contradiction). Le raiseonnement par l'absurde est possible dans Chen-Lin.

\smallskip

//////////////////////////////////////////////////////////////////////////////////////////////////////////////////////////////////

\smallskip
 
Pour revenir \`a l'article de Wei-Ye-Zhou:

\smallskip

1) Supposons qu'on puisse definir ou considrer, $ v $, la symetris\'ee de Schwarz par rapport \`a l'axe des ordonn\'ees de $ u $. Alors, $ v $ et ses deriv\'ees successives a droites et a gauche sont born\'ees. De plus $ v $ est solution au sens des distributions:

$$ -\Delta_{distr} v= We^{|v|}= Heaviside $$

On se place au voisinage de $ 0 $. Si on derive une fois par rapport a la premiere variable, $ \partial_1 v $ est solution au sens des distributions, un dirac sur l'axe des ordonn\'ees:

$$ -\Delta_{distr} (\partial_1 v)= -\partial_{1, distr} (\Delta v)=-\partial_1(\Delta v)+ \int_{-\epsilon}^{\epsilon} k \phi(0,x_2) dx_2, $$

avec $ k=-2\Delta u(0,x_2)>0, x_2\in(-\epsilon, \epsilon) $, au voisinage de $ 0 $,

$$ -\Delta_{distr} (\partial_1 v)=-\partial_1(\Delta v)+ \int_{-\epsilon}^{\epsilon} k\delta_{(0,x_2)} dx_2 $$

On note $ G((0,x_2),y)= -\frac{1}{2\pi} \log ||(0,x_2)-y||=-\frac{1}{2\pi}\log|(x_2-y_2)^2+y_1^2|$ le potentiel Newtonnien et $ h_0(y)=\int_{-\epsilon}^{\epsilon} k G((0,x_2), y) dx_2 $, $ h_0\in L^p, p>1 $:

$$ \int_{-\epsilon}^{\epsilon} \phi(0,x_2) dx_2= \int_{-\epsilon}^{\epsilon} \delta_{(0,x_2)} dx_2=\int_{B(0,\epsilon)} -h_0(y) \Delta \phi(y) dy, $$ 

C'est a dire que:

$$ -\Delta_{distr} (\partial_1 v-h_0) =-\partial_1(\Delta v) \in L^{\infty} $$

On a par les estimations elliptiques:

$$ \partial_1 v-h_0 \in C^1(B(0,\epsilon)) $$

Or, $ \partial_{11} v $ est reguliere et born\'ee en dehors de l'axe des ordonn\'ees  et $ h_0 $ est reguliere en dehors de l'axe des ordonn\'ees.

Donc, $ \partial_1 h_0 $ serait born\'ee en dehors de l'axe des ordonn\'ees. $ \exists C >0, \forall y\in B(0,\epsilon)-\{(0, x_2), x_2\in(-\epsilon,\epsilon)\}, \,\, |\partial_1 h_0(y)|\leq C $.

\smallskip

Or, $ h_0(y)=-\frac{1}{2\pi} \int_{-\epsilon}^{\epsilon} \log|(x_2-y_2)^2+y_1^2| dx_2 $, et si on derive par rapport a $ y_1 $ on a :

$$ \partial_1 h_0(y)|_{\{y_1 \to 0, y_2 \to 0, y_1\not =0\}} \equiv \frac{1}{y_1} \to +\infty $$

Ce qui est contradictoire avec le fait qu'en dehors de l'axe des ordonn\'ees, $ \partial_1 h_0 $ est born\'ee.

Donc, on ne peut pas considerer $ v $, $ v $ n'est pas definie, n'existe pas, avec autant d'informations.

2) De meme, si on considere, $ w $ la fonction reflexion par rapport \` a l'axe des ordonn\'ees, son laplacien serait comme precedemment une somme de Dirac + une fonction $ L^{\infty} $. Si on considere $ \partial_1 w-\partial_1 h_0 $ elle serait $ C^0(B(0,\epsilon))$, ce qui impliquerait que $ \partial_1 h_0 $ serait born\'ee en dehors de l'axe des ordonn\'ees, comme precedemment, c'est pas possible.

$$ -\Delta_{distr}  w=-\Delta w + \int_{-\epsilon}^{\epsilon} k_0 \phi(0,x_2) dx_2, k_0 >0, $$

$$ -\Delta_{distr} w =-\Delta w + \int_{-\epsilon}^{\epsilon} k_0\delta_{(0,x_2)} dx_2 $$

\smallskip

Donc pour 1) et 2) on ne peut pas considerer, $ v $ et $ w $, ces fonctions n'existent pas, avec autant d'informations.

\smallskip

////////////////////////////////////////////////////////////////////////////////

Concernant le "blow-up" au bord:

\smallskip

{\bf Il faut fixer une mesure:}

On veut qu'il existe une sous-suite $ (u_{i_j})_j$ et des points $ y_1, \ldots, y_m$ telle que cette sous suite converge sur tout compact de $ D_m=\bar \Omega-\{y_1,..., y_m\}$. Pour cela on fixe un nombre maximal $ m $. Qui dit que ce $ m $ existe ? et faut que toute sous-suite de la sous-suite converge aussi, cela veut dire qu'il faut que $ V_{i_j} e^{u_{i_j}}$ converge, au moins faiblement. Cela veut dire, que c'est possible si on a au moins, de la convergence, meme faible  $ V_{i_j} e^{u_{i_j}} \rightharpoonup \mu $. Cela est possible si il y a convergence faible vers une quantit\'e $ \mu $.

Par exemple, on fixe, les points de concentration et $ D_m $, soit $ (K_n) $ une suite exhaustive de compacts de $ D_m $. On doit avoir, $ \liminf_{\epsilon_{i_j}} \int_{B(y_k, \delta)} \epsilon_{i_j} e^{u_{\epsilon_{i_j}}} dx \geq \alpha >0 $, on prend le sup sur chaque compact $ K_n $, ces points convergent vers $ y_0 \in K_n $ pour une sous-suite, il faut alors que la sous-suite $ (u_{i_{j,k}})_k $ verifie:

$ \liminf_{\epsilon_{i_{j,k}}} \int_{B(y_0, \delta_0)} \epsilon_{i_{j,k}} e^{u_{\epsilon_{i_{j,k}}}} dx <\alpha $, qui dit que c'est vrai ?

c'est a dire que la nouvelle sous-suite, conserve, la notion de convergence, en dehors des points de concentrations de depart, $ y_1,..., y_m $. cela veut que, a chaque fois qu'on extrait une sous-suite, cette suite converge sur  les compacts de $ D_m $. Or ceci n'est que la defintion de la convergence vers une quantit\'e fix\'ee, de la suite de depart.(Sinon, on aurait un autre point de concentration, et on arrive pas a fixer une sous-suite qui converge, ca se casse a un certain moment, on n'aura pas construit de sous-suite concergente et l'arret dit qu'on n'a pas construit de sous-suite convergente, ce n'est pas le bon candidat a la convergence).

Quand on a une suite pour la quelle, toute sous-suite , converge et est stable, cela veut dire que la suite de depart converge vers une quantit\'e fix\'ee. ceci est possible si on a au moins la convergence faible. C'est ce qu'on a dit pr\'ecedemment.

\smallskip

///////////////////////////////////////////////////////////////////

\smallskip

{\bf On ne peut pas considrer une mesure sur $ \bar  \Omega $, comme dans le deuxieme livre de Brezis, Functional Analysis, Sobolev Spaces and PDEs:}

\smallskip

Supposons qu'on puisse considrer une mesure sur $ \bar \Omega $. L'espace total sur le quel on se place et sur lequel on considere les mesures (de Radon) est: $ X=\bar \Omega $. On se place donc sur $ X=\bar \Omega $. les fonctions $ u_{\epsilon} \in C^0(\bar \Omega) $, c'est a dire qu'il faut prendre en compte les points de $\partial \Omega $ et les voisinages $ V_{\partial \Omega} $ des points de $ \partial \Omega $. Il faut definir ce que c'est $ L^p(\bar \Omega) $ et en particulier $ L^1(\bar \Omega) $, pour cela il faut une mesure de Radon $ \mu_0 $ sur $\bar  \Omega $, c'est difficile \`a definir, comme dans le livre d'Aubin et definir l'integrale Riemannienne, par des cartes. Supposons qu'on puisse definir $ \mu_0 $. Alors, comme $ u_{\epsilon} \in C^0(\bar \Omega) $, cette fonction $ \epsilon e^{u_{\epsilon}} \in L^1(\bar \Omega) $ qui est un espace qui n'a pas beaucoup de sens, or on sait d'apres l'hypothese du probleme que: $ \epsilon e^{u_{\epsilon}} \in L^1(\Omega) \not = L^1(\bar \Omega)$. On voit alors que $\epsilon e^{u_{\epsilon}} $ est dans deux espaces qui n'ont rien avoir l'un avec l'autre, ce qui veut dire que c'est pas bien definit. De plus $L^1(\bar \Omega)$ n'a pas beaucoup de sens(on le definit comme dans le livre d'Aubin par des cartes, localement, on definit $ \mu_0 $ localement par des cartes, comme dans le livre d'Aubin. Cet espace $ L^1(\bar \Omega)$ n'est pas familier, c'est pour cela qu'on dit qu'il n'a pas de sens, mais si on se place sur l'espace total $ X=\bar \Omega $, il faut bien definir $ \mu_0 $). On a quelque chose du type: $ \epsilon e^{u_{\epsilon}} \in L^1(\bar \Omega)\not = L^1(\Omega) \ni \epsilon e^{u_{\epsilon}} $, avec deux espaces $ L^1(\bar \Omega) $ et $ L^1(\Omega) $, qui n'ont rien a avoir l'un avec l'autre, car, les tribus Boreliennes sont differentes et n'ont rien avoir l'une avec l'autre, et les mesures sont differentes, une definie sur l'espace localement compact $\Omega $, pour la mesure de Lebesgue, et l'autre, $ \mu_0 $ definie sur $ \bar \Omega $. Alors qu'il faut considrer les mesures de Radon sur $\bar \Omega $, ce qui n'est pas le cas de celle de $ L^1(\Omega) $.(+ classes de fonctions qui n'ont rien a voir les unes avec les autres, il faut voir aussi par rapport aux classes de fonctions).

\smallskip

Pour la construction d'une mesure de Radon par les cartes sur un domaine a bord regulier $ \bar \Omega $, comme dans le livre d'Aubin. On a dit que ce n'est pas naturel, ceci est du \`a la partition de l'unit\'e. Il se peut que le support d'une des fonctions de la partitions touche le bord $ \partial \Omega $. Il y a un probleme avec la partition de l'unit\'e.

\smallskip

///////////////////////////////////////////////////////////////////////////////////////////////////////////////

{\bf Sur le blow-up} concernant les articles: {\it "About Brezis-Merle Problem with Holderian condition"} et les autres articles en dimension 2, et en particulier {\it "A compactness result for an equation with Holderian condition"}.

\smallskip

Concernant l'article {\it "About Brezis-Merle Problem with Holderian condition"}:

\smallskip

1) {\it Rescaling:} ou {\it changement d'echelle }: on se ramene \`a la situation d'un espace localement compact, $  \Omega_0=B(0,\epsilon), \epsilon >0 $, plus pr\'ecis\'ement: $ B(0,\epsilon_1), \epsilon_1>0,\ldots, B(0, \epsilon_k), \epsilon_k >0 $ et ici il faut, \`a chaque \'etape, fixer une mesure, c'est implicite, on se ramene \`a Brezis-Merle \`a chaque \'etape. En faisant un {\it blow-up de premiere esp\`ece} ou un {\it \'eclatement de premiere esp\`ece}.

\smallskip

2) Nombres fini d'etapes. Et chaque etape, on se ramene \`a un espace localement compact sur le quel une mesure se concentre. Tendant vers le bord. Processus d\'escendant ou d\'ecroissant.

\smallskip

Concernant l'article {\it "A compactness result for an equation with Holderian condition"}:

\smallskip

1) {\it Hierarchie:} on ramene la concentration de la mesure  d'un espace $ \Omega $ de dimension 2 \`a une concentration de la mesure sur une sous-vari\'et\'e compacte sans bord de dimension 1(qui est une vari\'et\'e compacte sans bord de dimension 1), par la formule de Stokes.

On "transporte" la mesure 2-dimensionnelle en une mesure 1-dimensionnelle sur une vari\'et\'e sans bord par la formule de Stokes. La masse totale sur $ \Omega $ se "transporte" en la masse totale sur $\partial \Omega $. 

Dans ces articles, on resout le probleme de concnentration de la mesure ou "blow-up" sur le bord:  sur tout le domaine \`a bord, ou, sur la vari\'et\'e compacte \`a bord de dimension 2, $\bar \Omega $.

\smallskip

2) Il y a deux mecanismes: 

-soit on se ramene \`a un espace localement compact sur plusieurs etapes en nombre fini. Par blow-up ou \'eclatement. 

-soit on se ramene \`a une vari\'et\'e compacte sans bord de dimension strictement inferieure, par "transport" de la mesure par la formule de Stokes.

\smallskip

3) Dans le blow-up, on ne peut pas utiliser n'importe quelle fonction cutoff, $\eta_{\epsilon} $, il faut prendre une fonction solution d'un probleme de Dirichlet $ \tilde \eta_{\epsilon} $: car on met la valeur absolue quand on \'evalue $ \int_{\Omega} |-\Delta [(u_i-u)\tilde \eta_{\epsilon}]| dx $. Il y a un probleme de valeur absolue et quand on utilise la formule de Stokes. Il y a 2 possibilt\'es, soit on prend une solution d'un probleme de Dirichlet. Soit on prend directement $ \eta_{\epsilon} $, mais il faut utiliser: a) L'in\'egalit\'e de Poincar\'e, et  b) prouver que $ ||\nabla (u_i-u)||_q =o(1), 1\leq q<2 $. Dans tous les cas, il faut une \underbar {fonction test}, et ces deux proc\'ed\'es (soit $ \tilde \eta_{\epsilon} $ solution d'un probleme de Dirichlet, soit directement $ \eta_{\epsilon} $ avec in\'egalit\'e de Poincar\'e et $ ||\nabla (u_i-u)||_q=o(1), 1\leq q<2$), sont equivalents, du point de vue, nombre d'etapes et nombre de notions: 

Fonction test +: 

\smallskip

a) $ \tilde \eta_{\epsilon} $ + probleme de Dirichlet, ({\bf Donn\'ees effectives}: ne d\'ependent pas de la solution), ou, 

\smallskip

b) $ \eta_{\epsilon} $ + inegalit\'e de Poincar\'e ({\bf Donn\'ee fictive}: d\'epend de la solution) et $ ||\nabla (u_i-u)||_q=o(1), 1\leq q<2$.

Les points de vue a) et b) sont equivalents: en tant que: nombre d'etapes, nombre de notions.

\smallskip

Hierarchie dans la th\'eorie de la d\'emonstration:

\smallskip

c) preuve avec donn\'ees effectives: donn\'ees r\'eelles, est meilleure que:

\smallskip

d) preuve avec donn\'ees fictives.

\smallskip

Par exemple: dans la modelisation: donn\'ee fictive, on approxime quelque chose d'approch\'ee (la solution, qu'on modelise), le taux d'erreur est plus grand. Donn\'ee effective: on n'approxime pas quelque chose d'approch\'ee, le taux d'erreur est plus petit. Il faut aussi estimer l'erreur: estimation d'erreur.

\smallskip

4) Dans la {\bf Compacit\'e}: il y a 2 etapes.

\smallskip

Par exemple, pour ce qui concerne les articles {\it "A compactness result for an equation with Holderian condition"} et {\it "A uniform boundedness result for solutions to the Liouville type equation with boundary singularity" }, il y a deux etapes:

\smallskip

\underbar {Etape 1:}

\smallskip

a1) blow-up,

\smallskip

a2) Transformation conforme, la partie du bord devient plate, un segment.

\smallskip

a3) Formule de Rellich-Pohozaev. Int\'egration par parties pour un potentiel lipschitzien donc, dans $ W^{1,\infty} $.

\smallskip

a4) masse $ m=0$.(ici aussi, on utilise le blow-up pour avoir la masse nulle).

\smallskip

\underbar {Etape 2:}

\smallskip

b1) blow-up,

\smallskip

b2) Formule de Stokes, apparition de la mesure sur le bord.

\smallskip

b3) Concentration de la mesure.

\smallskip

b4) masse $ m=\alpha_1 >0$.(ici aussi, on utilise le blow-up pour avoir la masse, $ m=\alpha_1>0$).

\smallskip

Dans le blow-up, on a l'existence d'une \underbar {mesure de Radon} $ \mu \geq 0 $ sur $ \partial \Omega $, puis, on montre une \underbar {decomposition de Lebesgue}: (mesure Riemannienne de Lebesgue sur le bord + mesures singulieres): $ \mu=\mu_L+\mu_s $, avec $ \mu(\cdot)=\int_{\partial 
 \Omega} (\partial_{\nu} u)(\cdot) d\sigma + \sum_{j=1}^N \alpha_j \delta_{x_j} $.

\smallskip

Maintenant, on veut connaitre l'equation, singuliere verif\'ee par $ u \geq 0 $: on fait un raisonnement approximatif.

\smallskip

On a: $ \partial_{\nu} u_i \in L^1(\partial \Omega) $, uniform\'ement, donc, par dualit\'e, on obtient (approximativement), avec le fait suivant: $ 0=u_i \to u=0 $ sur $ \partial \Omega $, la convergence faible dans un $ L^{1+\epsilon}(\partial \Omega), \epsilon >0 $, on a la notion de convergence faible pour le dual, on prend des fonctions $ \phi \in C^{\infty}(\bar \Omega) $, car on sait que $ \tilde \phi \in C^{\infty}(\partial \Omega) $ se prolonge en une fonction $ \phi \in C^{\infty}(\bar \Omega)$ (par exemple, harmonique):

(Il faut d\'efinir la dualit\'e, on suppose qu'elle est bien d\'efinie):

$$ <u_i,\partial_{\nu} \phi>_{\partial \Omega}=-<\partial_{\nu} u_i, \phi>=-\int_{\partial \Omega} \partial_{\nu} u_i \phi d\sigma, \,\, \forall \phi\in C^{\infty}(\bar \Omega), $$

en passant \`a la limite, on obtient:

$$ <u,\partial_{\nu} \phi>_{\partial \Omega}=-\int_{\partial \Omega} \partial_{\nu} u \phi d\sigma - \sum_{j=1}^N \alpha_j \phi (x_j), \,\, \forall \phi \in C^{\infty}(\bar \Omega). $$

ceci est l'equation singuliere de $ u \geq 0 $. Les fonctions $ u_i \geq 0 $ convergent \`a l'interieur du domaine vers $ u \geq 0 $. On voulait connaitre le type d'equation avec singularit\'e verif\'ee par $ u \geq 0 $. Pour Brezis-Merle, on a le potentiel logarithmique quand il y a singularit\'e interieure. Ici, on a voulu connaitre les singularit\'es de $ u \geq 0 $, au bord.

\smallskip

//////////////////////////////////////////

\smallskip

{\bf Remarque sur la preuve dans le cas d'un domaine analytique $ C^{\omega} $ et la preuve dans le cas d'un domaine $ C^{\infty} $}:

\smallskip

Ceci concerne l'article, "A compactness result for an equation with Holderian condition" et le preprint "About Brezis-Merle problem with Lipschitz condition": 

\smallskip

($ \beta=0 $, aussi, dans ce cas, quand le domaine est \'etoil\'e on peut borner le volume par l'identit\'e de Pohozaev, et on a alors, la compacit\'e, sans condition sur le volume. C'est pour cela qu'on considere ce cas $ \beta=0 $, qui est important, car on peut borner le volume et utiliser le degr\'e topologique, voir le point (32) ci-apr\'es). 

\smallskip

Pour $ 0 <\beta <1/2 $, on est oblig\'e d'utiliser une transformation conforme. Pour $ \beta =0 $, il y a 2 preuves possibles. 

\smallskip

Ici, pour $ \beta=0 $, on dit pourquoi, la preuve de la compacit\'e avec volume born\'e, avec la transformation conforme, dans le cas d'un disque ou d'une ellipse ou d'une couronne, passe avant la preuve avec la structure $ C^{\infty} $.

\smallskip

Quand on considere un disque ou une ellipse ou une couronne, on sait, au meme temps, que ce sont des domaines analytiques et $ C^{\infty} $, soit, en utilisant la bijection conforme avec le demi-plan de Poincar\'e (donc une carte ananlytique et holomorphe, donc $ C^{\infty} $), soit en se placant, par exemple, en $ x=0 $ et on considere la fonction: $ y=\sqrt {1-x^2} $ (on sait, d\`es le d\'ebut qu'elle est developpable en serie entiere) . On sait qu'elles sont analytiques et holomorphes, au meme moment qu'elle sont $ C^{\infty} $. On a la meilleure regularit\'e possible d\`es le debut, au meme moment qu'une regularit\'e plus faible. On a la meilleure r\'egularit\'e possible: analytique et holomorphe.

\smallskip

Donc, la preuve dans le cas analytique, pour les cas particuliers, d'un disque ou d'une ellipse ou d'une couronne, avec une transformation conforme, passe avant la preuve avec un domaine $ C^{\infty} $. Car, on connait explicitement les objets, et pour caracteriser ces objets, on a des applications analytiques et holomorphes et conformes, d\`es les d\'ebut. 

\smallskip

Par coh\'erence de la preuve, dans le cas d'un disque ou d'une ellipse ou d'une couronne, on utilise la structure analytique et holomorphe et conforme.

\smallskip

Ceci est li\'e au degr\'e d'une preuve: pour un disque ou une ellipse ou une couronne:

\smallskip

Preuve 1: holomorphie: $ C^{\omega} \Rightarrow $ compacit\'e. On n'a pas besoin de diminuer la r\'egularit\'e.

\smallskip

Preuve 2: lisse: $ C^{\omega} \Rightarrow C^{\infty} \Rightarrow $ compacit\'e. On est oblig\'e de passer par une diminution de la regularit\'e, alors qu'on sait d\`es le d\'ebut que le domaine est analytique.

\smallskip

Concernant l'article: "Compactness of the set of solutions to elliptic equations in 2 dimensions": on a le parametre: $ \epsilon_i \to 0, 1\geq \epsilon_i \geq 0 $, qui peut etre nul ou non nul, qui fait qu'on utilise la preuve pour un domaine regulier $ C^{\infty} $.

\smallskip

////////////////////////////////////////////////////

\smallskip

Remarque sur le cas $ \alpha =0 $ de l'article de Journal of Math.Sci.Univ.Tokyo.2016:

\smallskip

Supposons par exemple qu'on puisse d\'efinir le probl\`eme pour $ \alpha =0 $, cela veut dire qu'on peut definir la quantit\'e: $ 1=|x|^{-2\alpha}= 0^0= \lim \epsilon_i^{-2\alpha_i} $, pour toutes suites $ 0 < \epsilon_i \to 0 $ et $ 0 < \alpha_i \to 0 $. Or en prenant des suites particulieres dont l'une domine l'autre, par exemple, $ \epsilon_i =e^{-1/\alpha_i^2} \to 0 $, on obtient: $ \epsilon_i^{-2\alpha_i}= e^{-2\alpha_i \log \epsilon_i}=e^{2/\alpha_i} \to +\infty $, ce n'est pas possible. Donc, on ne peut pas d\'efinir le probl\`eme pour $ \alpha=0 $, le cas r\'egulier ne peut pas etre inclus dans le cas g\'en\'eral $ \alpha \in (0,1/2) $. Ce qui veut dire qu'il faut ecrire un autre article traitant du cas regulier. Ce qui a \'et\'e fait dans l'article de Comm.Math.Analysis. 2018, ($ \beta=0 $). Ce qui veut dire aussi que $ \alpha=0 $ s\'epare les cas $ \alpha \in (0,1/2) $ et $ \alpha < 0$, ce qui veut dire qu'il faut ecrire un autre article traitant le cas $ \alpha <0 $, c'est l'objet d'un preprint.

\smallskip

On dit cela pour expliquer qu'on n'ecrit pas plusieurs fois la meme chose. On doit r\'ediger chaque cas.

\smallskip

///////////////////////////////////////////////////////////////////////////////////////////////////////////////////////////////////

\smallskip

(32) {\bf Sur le degr\'e topologique et l'equation avec non linearit\'e exponentielle:}

\smallskip

Voir l'article de Brezis-Turner. Comme on l'a dit concernant l'article de De Figueiredo-Lions-Nussbaum, ils considerent un operateur $ M $ tel que $ M (0)=0$, nulle en $ 0 $. On peut considerer:

$$ -L u= V(x)(e^u-1), \,\, 0 < a \leq V(x)\leq b <\lambda_1, \,\, ||\nabla V||_{\infty} \leq A $$

$ f(u)=e^u-1 $ et $ f(0)=0 $, $ (-L)^{-1} $ un operateur compact. voir l'article de De Figueiredo-Lions-Nussbaum: $ M= (-L)^{-1}(V(x)f(u))$. $ \lambda_1 $ premiere valeur propre de $ -L $ avec condition de Dirichlet.

\smallskip

Voir l'article de Brezis-Turner. Une des conditions pour prouver l'existence de solutions est : $ g(x,u,p)/u <\lambda_1 $ pour $ u $ voisin de 0. et la presence de la nonlinearit\'e exponentielle en $ +\infty $ assure que les $ t \geq 0 $ de l'article de Brezis-Turner sont uniform\'ement born\'es. Il reste \`a borner les solutions $ u $ positives ou nulles pour appliquer la th\'eorie du degr\'e topologique. Ceci se fait par l'identit\'e de Pohozaev en ramenant l'equation pr\'ec\'edente  a l'equation qu'on a consid\'er\'e dans des articles pr\'ec\'edents avec seulement $ e^u $ et non $ e^u -1 $. En resolvant $ -Lw_1= V, -Lw_2=t J $ avec condition de Dirichlet, ($ J=\phi_1 $ premiere fonction propre) ces fonctions $ w_1, w_2 $ sont uniform\'ement born\'ees dans $ C^2(\bar \Omega) $ car $ t $ est born\'e unifrom\'ement et $ V $ aussi. Donc, on s'est ramen\'e, dans un domaine r\'egulier etoil\'e par rapport a l'origine a notre equation considr\'ee dans des articles pr\'ec\'edents et l'identi\'e de Pohozaev donne le r\'esultat de borne uniforme de $ u $, comme voulu dans Brezis-Turner.

\smallskip

Voir l'article de De Figueiredo-Lions-Nussbaum, pour l'id\'ee sur le degr\'e topologique et son applications.

\smallskip

Voir l'article de: H.H.Zou dans Calculus of variation. Pour une id\'ee sur l'application du degr\'e topologique.

\smallskip

Voir aussi, l'article de Brezis-Turner.

\smallskip

Donc: 

\smallskip

Degr\'e \, topologique  $\Rightarrow  \forall V\in C^{\infty}, 0 < a \leq V \leq b < \lambda_1, \,\, ||\nabla V||_{\infty} \leq A, \,\,  \exists \, u \in C_0^2(\bar \Omega), \, -L u = V(x) (e^u -1), $ et $ u >0 $ dans $ \Omega $. Pour ce qui nous concerne, on prend $ V $ lipschitzien.

\smallskip

Soit $ w $ tel que $ -L w= V $ avec condition de Dirichlet. Par le principe du maximum, th\'eor\`eme de Dualit\'e et les estimations elliptiques, $ w \geq 0 $ et $ w>0 $ dans $ \Omega $ et $ ||w||_{C^2(\bar \Omega)} \leq c_1(a, b, A, \Omega) $. Soit $ W = Ve^{-w} $, alors $ w $ et $ u+w $ sont deux solutions distinctes de $ -Lv=We^v $ avec conditions de Dirichlet et $ 0< a_1 \leq W \leq b <\lambda_1, $ et $ ||\nabla W||_{\infty} \leq A_1$, avec  $ ||v||_{C^2(\bar \Omega)} \leq c(a,b, A, \Omega) $. La solution $ v=w+u$ est la translat\'e de $ u $. Donc, elle est aussi topologique. On d\'ecale $ u $ de $ w $ pour obtenir $ v $.

\smallskip

Donc: Pour tout $ V $ on a $ W $ dans le meme type de classe avec deux solutions distinctes par le degr\'e topologique et l'ensemble des solutions est born\'e dans $ C^2(\bar \Omega) $.

\smallskip

On a: Pour tout $ V $ il existe $ W $ dans la meme classe que $ V $ tel que l'equation $ -Lv=We^v $ (avec condition de Dirichlet), possede au moins 2 solutions distinctes dont une topologique et ces solutions sont born\'ees uniform\'ement dans $ C^2(\bar \Omega) $.

\smallskip

Donc: si on note $ C $ le cone $ C=\{ u \in C^1(\bar \Omega), u \geq 0 \} $.

\smallskip

Degr\'e topologique $ \Rightarrow deg=i_C=-1 \Rightarrow $ il y a un point fixe pour $ u \to M(u) $ $ \Rightarrow  \forall \, V\in C^{\infty}, \exists \, W \in C^{\infty}, \,\, W(x)=V(x)e^{-w}, $ tel que l'equation $ -Lv=We^v $ possede 2 solutions distinctes dans $ C_0^2(\bar \Omega) $ pour $ W $ avec $ 0 <a_1 \leq W\leq b_1 <\lambda_1 $ et $ ||\nabla W||_{\infty} \leq A_1 $ et $ ||v||_{C^2(\bar \Omega)} \leq c(a, b, A, \Omega) $. Pour ce qui nous cencerne, on prend $ V $ lipschitzien.

\smallskip

1) Donc on a: Pour tout $ V $ il existe $ W $ dans la meme classe que $ V $ tel que l'equation $ -Lv=We^v $ (avec condition de Dirichlet), possede au moins 2 solutions distinctes dont une topologique et ces solutions sont born\'ees uniform\'ement dans $ C^2(\bar \Omega) $. On a un r\'esultat d'existence de solutions topologiques pour l'equation de type Liouville (en dimension 2).

\smallskip

2) a)(Equation sinh-Poisson en mecanique statistique et en physique des plasmas). On peut prendre $ M=(-\Delta)^{-1}(V(x)g(u)) $ avec $ g(u)=\sinh (u)$. On obtient alors $ u >0 $ solution par le degr\'e topologique de l'equation sinh-Poisson, qui modelise la 2d-turbulance en mecanique statistique et apparait en physique des plasmas. Voir les prints de Pablo Figueroa et de Joel Spruck. Donc: si on note $ C $ le cone $ C=\{ u \in C^1(\bar \Omega), u \geq 0 \} $, on a $ deg=i_C=-1 $  et l'application $ u \to M(u) $ possede un point fixe.(Voir l'article de De Figueiredo-Lions-Nussbaum et ce qu'on a ecrit sur le degr\'e topologique, au debut de ce print). Et on a une solution topologique de: $ -\Delta u= V(x) \sinh(u) $ avec condition de Dirichlet et $ u\not \equiv 0 $, $ u \geq 0 $ (par le principe du maximum $ u >0 $ dans $ \Omega $), pour $ 0 < a \leq V(x) \leq b < \lambda_1 $ et $ ||\nabla V||_{\infty} \leq A $ ou $ V $ avec poids holderien $ (1+|x-x_0|^{2\beta}), \beta \in (0,1/2), x_0\in \partial \Omega $ avec $ \Omega $ particulier, une couronne et, $ A=\frac{a}{2(1+2^{2\beta})} $, la solution $ u $ est uniform\'ement born\'ee dans $ C^2(\bar \Omega) $, $ ||u||_{C^2(\bar \Omega)} \leq c(a,b, \beta,\Omega) $.

\smallskip

b) On peut faire la meme chose en considerant une singularit\'e au bord. On peut prendre $ M= {(-\Delta)}^{-1}(\frac{V(x)}{|x-x_0|^{2\alpha}} g(u))$ avec $ g(u)=\sinh(u) $, $ x_0 \in \partial \Omega $, $ \Omega $ une couronne et $ \alpha \in (0,1/2) $, (ca reste vrai avec $ \alpha \in (-\infty, 0) $), et $ \lambda_1 $ premiere valeur propre avec condition de Dirichlet et avec poids singulier, ($ -\Delta \phi_1=\lambda_1 \phi_1\frac{1}{|x-x_0|^{2\alpha}} $, qu'on resout par la methode variationnelle comme dans le cas regulier), et $ 0 < a \leq V(x) \leq b <\lambda_1 $, $ ||\nabla V||_{\infty} \leq A=\dfrac{(-\alpha+1) a}{2} $. Si on note $ C $ le cone $ C=\{ u\in C^1(\bar \Omega), u \geq 0\} $, on a $ deg=i_C=-1 $ et l'application $ u\to M(u) $ a un point fixe et on obtient une solution topologique de: $ -\Delta u= \frac{V(x)}{|x-x_0|^{2\alpha}} \sinh(u) $, $ u\not \equiv 0 $ et $ u\geq 0 $ (par le principe du maximum $ u >0 $ dans $ \Omega $) avec condition de Dirichlet et la solution $ u $ est uniform\'ement born\'ee dans $ C^1(\bar \Omega) $, $ ||u||_{C^1(\bar \Omega)} \leq c(a, b, \alpha, \Omega) $.

\smallskip

c) On obtient la meme chose avec $ \alpha \in (-\infty,0) $ avec une solution topologique avec singularit\'e au bord et born\'ee uniform\'ement dans $ C^2(\bar \Omega) $.

\smallskip

d) On peut faire la meme chose que dans le point a) avec l'operateur $ L $ du point 1), on prend $ M=(-L)^{-1}(V(x)g(u)) $ avec $ g(u)=\sinh (u) $. On obtient dans le cone $ C =\{ u\in C^1(\bar \Omega), u \geq 0 \} $, $ deg=i_C=-1 $ et $ u\to M(u) $ a un point fixe. On obtient une solution topologique de: $ -Lu=V(x) \sinh (u) $ avec condition de Dirichlet et $ u \not \equiv 0 $ et $ u\geq 0 $ (par le principe du maximum $ u >0 $ dans $ \Omega $), pour $ V $ tel que, $ 0 < a \leq V(x) \leq b <\lambda_1 $ et $ ||\nabla V||_{\infty} \leq A $, $ \Omega $ un domaine r\'egulier etoil\'e par rapport \`a l'origine, $ \lambda_1 $ premiere valeur propre de $ -L $ avec condition de Dirichlet, la solution $ u $ est uniform\'ement born\'ee dans $ C^2(\bar \Omega) $, $ ||u||_{C^2(\bar \Omega)} \leq c(a, b, A, \Omega) $.

\smallskip

{\bf Remarque:} a) Lorsque $ V\equiv a=b \geq \lambda_1$, par un argument d'integration et integration par parties de la solution l'Eq si elle existe; avec la fonction propre et de positivit\'e, il n'y a pas de solution \`a l'Eq. De meme pour $ b\geq V \geq a \geq \lambda_1 $, il n'y a pas de solutions. Par contre Il y a des solutions si $ b< \lambda_1 $. 

\smallskip

b) On peut utiliser le fait qu'il n'y a pas de solutions en $ \lambda_1 $, l'equation de Liouville $ -\Delta u= Ve^u $ avec conditions de Dirichlet et $ V $ lipschitzien:  pour calculer le degr\'e topologique $ d=0 $, en utilisant une homotopie (se placer autour de $ \lambda_1$, et on utilise, cette fois, le degr\'e topologique pour l'equation de Liouville, pour l'expoenentielle, $ e^u $, au lieu de $ \sinh $ ou $ e^u-1$: s'il n'y a pas de solution autour de $ \lambda_1 $, alors le degr\'e topologique vaut $ d=0 $, si on suppose qu'il y a des solutions autour de $ \lambda_1 $, on utilise une homotopie en partant de $ \lambda_1-\epsilon $ jusqu'\`a $ \lambda_1 $ et en $ \lambda_1 $ il n'y a pas de solutions donc le degr\'e vaut $ 0 $, pour avoir le degr\'e en $ \lambda_1-\epsilon $, $ d_{\lambda_1-\epsilon}=d_{\lambda_1}=d=0 $).  On consid\`ere alors: $ L=Id-\lambda (-\Delta)^{-1}(Ve^u) $, $ 0 < 1-\epsilon \leq \lambda \leq 1 $ et $ 0 <a=\lambda_1\leq V\leq b $, $ |\nabla V|\leq A $, les solutions de l'equation du type Liouville sont uniform.born\'ees par $ C_{\epsilon} >0 $. Alors, $ L=Id-K $ avec $ K $ compact, le degr\'e est bien d\'efinit dans $ B_{2C_{\epsilon}}(0) $: $ d_{\lambda}=d_{\lambda-\epsilon}= d_{\lambda_1}=d=0 $. On donne la valeur du degr\'e topologique de Leray-Schauder: $ d_{\lambda}=d=0 $ ($ \lambda >0 $ petit, donc $ 1>\epsilon >0 $ voisin de 1).

\smallskip

(On sait d'apr\`es ce qu'on a dit ci-dessus, que l'equation $ -\Delta v=We^v $ avec condition de Dirichlet, $ 0 <\tilde a \leq a_1 \leq W \leq b_1 <\lambda_1 \leq \tilde b $ et $ |\nabla W|\leq A$, a des solutions. On peut dire qu'il y a toujours des solutions avec un degr\'e nul: $ d_{\lambda}=d=0 $). 

\smallskip

Le degr\'e nul: signification: soit il n'y a pas de solutions, soit il y a au moins 2 solutions dont les indices se compensent, rotation \`a gauche et rotation \`a droite. En $\lambda=1 $, il n'y a pas de solutions, en $ \lambda >0 $ petit, il y a des solutions, donc, il y a au moins 2 solutions. 

\smallskip

////////////////////////////////////////////

\smallskip

Existence de solutions avec volume born\'e: exemple equation du type Liouville: on utilise la partie 4.55 et le theoreme 4.55 du livre d'Aubin.T:

\smallskip

Soient $ b, A, C >0 $, $ C>|D| >0 $, $ C $ grand (pour que la solution nulle appartienne \`a $ \theta $) et $ D $ un domaine born\'e analytique de $ {\mathbb R}^2 $:

\smallskip

On prend comme espace de Banach: $ B=C_0^1(\bar D) $, $ \theta=\{ u\in B, \,\, \int_{D} e^u < C \} $, $ \theta $ est un ouvert de $ B $. $ F_t=Id-(-\Delta)^{-1}(tV(x)e^u) $. Avec $ 0 \leq t\leq 1 $, $ 0\leq V(x)\leq b $ et $ ||\nabla V||_{\infty} \leq A $. (Dans $ \theta $, on a l'in\'egalit\'e stricte, pour avoir un ouvert).

\smallskip

D'apres les resultats de compacit\'e pour l'equation du type Liouville avec condition de Dirichlet : $ \exists c>0 $ tel que: $ ||u||_{C^1(\bar D)} \leq c $.

\smallskip

On prend: $ \Omega \subset \theta, \Omega=\{ u \in \theta, ||u||_{C^1(\bar D)} <2c \} $, $ \Omega $ est un ouvert aussi, born\'e dans $ B $. (In\'egalit\'e stricte pour avoir un ouvert). 

Alors, par le resultat de compacit\'e: $ 0\not \in F_t(\partial \Omega) $. Donc le degr\'e $ d_t=deg (F_t,\Omega, 0) $ est bien definit et pour $ t=0 $, $ F_0=Id $, donc: $ d_0=1 $. Donc: 

$$ d_t=d_1=d_0=+1 \not =0.$$

Donc, pour $ t=1 $, l'equation du type Liouville a une solution born\'ee uniform. avec volume born\'e uniform.

\smallskip

//////////////////////////////////////////////////////////////

\smallskip

On peut faire cela pour l'equation de du type Liouville avec singularit\'e au bord:

\smallskip

Soient $ b, A, C >0 $, $ C>\int_D|x|^{-2\alpha}dx >0 $, $ C $ grand (pour que la solution nulle appartienne \`a $ \theta $) et $ \alpha \in (0,1/2) $ et $ D $ un domaine born\'e analytique de $ {\mathbb R}^2 $, $ 0\in \partial D $:

\smallskip

On prend comme espace de Banach: $ B=C_0^1(\bar D) $, $ \theta=\{ u\in B, \,\, \int_{D} |x|^{-2\alpha} e^u < C \} $, $ \theta $ est un ouvert de $ B $. $ F_t=Id-(-\Delta)^{-1}(tV(x)|x|^{-2\alpha} e^u) $. Avec $ 0 \leq t\leq 1 $, $ 0\leq V(x)\leq b $ et $ ||\nabla V||_{\infty} \leq A $. (Dans $\theta $, on prend l'in\'egalit\'e stricte pour avoir un ouvert). 

\smallskip

D'apres les resultats de compacit\'e pour l'equation du type Liouville avec condition de Dirichlet : $ \exists c>0 $ tel que: $ ||u||_{C^1(\bar D)} \leq c $. 

\smallskip

On prend: $ \Omega \subset \theta, \Omega=\{ u \in \theta, ||u||_{C^1(\bar D)} <2c \} $, $ \Omega $ est un ouvert aussi, born\'e dans $ B $. (In\'egalit\'e stricte pour avoir un ouvert).

\smallskip

Alors, par le resultat de compacit\'e: $ 0\not \in F_t(\partial \Omega) $. Donc le degr\'e $ d_t=deg(F_t, \Omega, 0) $ est bien definit et pour $ t=0 $, $ F_0=Id $, donc: $ d_0=1 $. Donc: 

$$ d_t=d_1=d_0=+1 \not =0.$$

Donc, pour $ t=1 $, l'equation du type Liouville avec singularit\'e au bord, a une solution born\'ee uniform. avec volume born\'e uniform.

////////////////////////////

\smallskip

On a la meme chose pour l'equation avec operateur different du laplacien: $ \Delta+\epsilon (x_1\partial_1+x_2\partial_2), 0 \leq \epsilon \leq 1 $.

\smallskip

////////////////////////

\smallskip

On peut faire cela avec un potentiel holderien, $ V \in C^s, 1/2 <s <1 $, $ D=B_1(0) $ et $ bC<24\pi $. On a des solutions topologiques avec volume born\'e. Mais pour les contre-exemples explicites aux probl\`emes 1 et 2 de Brezis Merle, il faut utiliser la methode variationnelle, voir ci-apres.

\smallskip

/////////////////////////////////

\smallskip

On peut faire cela avec l'equation de la courbure scalaire prescrite en dimension 4, sur une vari\'et\'e compacte sans bord, en se ramenant par une homotopie \`a l'equation de Yamabe: on prend $ \theta =\{u\in C^{2,\alpha}(M), u > m >0\} \subset C^{2,\alpha}(M) $, $ K_t(v)=Id-(\Delta_g+\frac{1}{6}R_g)^{-1}(V_t(x)E(v)v^3) $, avec $ V_t(x)=a t+(1-t)V(x) $, $ 0 <a\leq V(x)\leq b<+\infty $, $ ||\nabla V||_{\infty} \leq \epsilon_1=\epsilon_1(a,b,m) \Rightarrow ||\nabla V_t(x)||_{\infty}\leq \epsilon_0=\epsilon_0(a,b,m) $. 

\smallskip

On a: $ E(v)=\int_M (|\nabla v|^2+\frac{1}{6}R_gv^2) $, l'op\'erateur conforme est coercif. 

\smallskip

Par le resultat de compacit\'e: $ ||u_t||_{C^{2,\alpha}(M)} \leq c $. On prend: $ \Omega =\{u\in \theta, ||u||_{C^{2,\alpha}(M)} \leq 2c \} $. Par les resultats sur l'equation de Yamabe, quand $ M $ n'est pas conform\'ement diffeomorphe \`a $ {\mathbb S}^4 $, par exemple, non localement conformement plate, alors: $ deg(K_0,\Omega,0)=-1 $. 

\smallskip

Pour le degr\'e topologique de l'equation de Yamabe:

\smallskip

$ \exists \Lambda_0 >0, \exists c_0>0, \forall 0<\Lambda <\Lambda_0, \forall c\geq c_0, deg (K_0,\Omega_{\Lambda, c}, 0)=-1 $, avec 

$ \Omega_{\Lambda,c}=\{u\in C^{2,\alpha}(M), ||u||_{C^{2,\alpha}(M)} < c, u>\Lambda >0\} $, pour l'equation de Yamabe. 

\smallskip

(D'apr\'es les travaux de Li-Zhu, Li-Zhang, Druet, Marques, sur l'equation de Yamabe: on a: pour $ \mu >0 $ tr\'es grand: $ deg(K_0,\Omega_{\mu,1/\mu}, 0)=-1 $).

\smallskip

On choisit $ m >0, m<\Lambda_0 $, on d\'etermine, $\epsilon_0>0$, on prend le potentiel $ V>0$, comme ci-dessus, puis, on d\'etermine $ c $ qui borne les solutions $ u >0 $, de l'equation de la courbure scalaire prescrite pour les $ V>0$ comme ci-dessus, avec $ c>c_0>0 $. Alors: $ \Omega=\Omega_{\Lambda,c} $. Et pour $ 0 < m <\frac{1}{\mu} $ et $ c>\mu>0 $, alors pour les solutions de l'eq.de Yamabe: $ deg(K_0,\Omega, 0)=deg(K_0,\Omega_{\Lambda,c}, 0)=deg(K_0,\Omega_{\mu,1/\mu}, 0)=-1 $

Donc:

$$ deg (K_1,\Omega, 0)=deg(K_t,\Omega, 0)=deg(K_0, \Omega, 0)=-1 \not =0. $$

Donc, on a une solution topologique de l'equation de la courbure scalaire prescrite avec $ u>m>0$.

\smallskip

///////////////////////////////////////////////////////////////////////////////

\smallskip

Le cas de systemes en dimension 2:

\smallskip

e) Ceci s'applique aux systemes du type $ -\Delta u=Ve^v, -\Delta v=We^u $ avec conditions de Dirichlet en considerant un operateur de vecteurs: 

\smallskip

$ M(u,v)=(-\Delta)^{-1}(V\sinh v, W\sinh u)+t(\phi_1, \psi_1) $, $ t\geq 0 $, avec singularit\'e interieure $ x_0 $, (pour le cas holderien d'odre $ \beta >0 $, ca marche avec une singualrit\'e pour une equation dans un systeme, reste le cas ou on a la singularit\'e interieure holderienne dans les deux equations du systeme.) ou singularit\'e au bord $ x_0 $ et $ \phi_1, \psi_1 $ premieres fonctions propres du Laplacien avec conditions de Dirichlet. avec singularit\'e ou pas, $ (-\Delta)^{-1} $ s'applique aux deux equations du systeme.

\smallskip

Dans le cas de singularit\'e interieure $ 0 = x_0 \in \Omega $, l'estimation a priori au voisinage du bord  est obtenue dans les articles de, De Figueiredo-Do O-Ruf au bord et interieure (Par la methode moving-plane, comme dans De Figueiredo-Lions-Nussbaum, De Figueiredo-Do O-Ruf le font et le disent dans leur article) et Bartolucci-Tarantello-Chen-Lin (eux utilisent Pohozaev-Rellich pour quantifier, calculer la masse) pour une equation et $ x_0 $  singularit\'e a l'interieure, c'est vrai pour une equation (reste le cas d'un systeme avec $0$ comme singularit\'e positive dans les deux equations du systeme).

\smallskip

Ici, on s'interesse au cas avec singularit\'e au bord, pour des systemes avec singularit\'e au bord ou poids Holderien avec singualrit\'e au bord, sur une couronne c'est possible.

\smallskip

{\bf Remarque sur les in\'egalit\'es du type Harnack $ \sup \,\inf $:} 

\smallskip

Relation avec la physique: (voir aussi, les points suivants): enroulement, distortion:

a) La th\'eorie des cordes (string theory) unifie les th\'eories de Jauge (des champs), des champs conformes, (th\'eorie des champs de Liouville), de Yang-Mills, et la gravitation quantique (Kaluza-Klein, $ n+1 $, $ n=3, 4, 5, 6, 9, 10 $ ou plus g\'en\'erale $ n\geq 3 $). Les in\'egalit\'es du type Harnack apparaissent dans ces trois th\'eories, on peut considerer alors que cette notion est une notion de la th\'eorie des cordes. Corde qui est un 'fil' ou un ensemble de particules.

\smallskip

b) Concernant la minoration du $ \sup \times \inf $, on a les exemples dans le livre d'Aubin (cas compact). Dans le cas plat, on utilise l'exemple de C.Chen-C.S.Lin dans Commun.Pure.Appl.math.1997, ici, $ n\geq 3$ et aussi l'article de C.C.Chen et C.S.Lin de 1999, "Blowing-up...". 

\smallskip

33) {\bf Equation de la courbure scalaire prescrite (et du type Yamabe) en dimension 4 et champs de Yang-Mills}:

\smallskip

a) exemple: 1-si on note $ SW $ la vari\'et\'e hyperbolique sans bord et orientable de dimension 3, de Seifert-Weber, alors $ SW\times S_1$ est une vari\'et\'e loc.conf.plate de courbure scalaire constante stric.negative de dimension 4, elle est li\'e aux champs de Yang-Mills. 2- En dimension 4, aussi on considere, la vari\'et\'e de Davis ou le produit de deux surfaces de courbure sectionelle non oppos\'ees et de courbure scalaire strc.negative (donc, non-loc-conf.plate $ S_{-1} \times S_{-1} $ (elle est d'Einstein) avec $ S_{-1} $ une surface de courbure scalaire constante $ <0 $), on obtient une vari\'et\'e Riemmannienne li\'e aux champs de Yang-Mills. On peut considerer sur ses vari\'et\'es l'Equation de la courbure scalaire prescrite (dimension 4).

\smallskip

Dans le cas positif (pour l'equation de Yamabe et du type Yamabe), $ S_2\times S_2 $, $ S_3\times S_1$.

\smallskip

b) En physique (description de la force nucl\'eaire responsable de la cohesion des protons-neutrons dans le noyau) le cas le plus important est port\'e aux vari\'et\'es de dimension 4 avec une metrique Riemannienne ou Lorentzienne.

\smallskip

c) Regarder l'article de, Andrzej Derdzinski, dont le titre est: "Riemannian manifolds with harmonic curvature", il y a la fonctionnelle de Yang-Mills et pour lui la deriv\'ee du tenseur de courbure donne les points critiques de la fonctionelle de Yang-Mills. En fait il s'agit de la vari\'et\'e, de la metrique et de la connection. C'est la courbure de Riemann qui est un champ de Yang-Mills.

\smallskip

d) Pour ce qui est des champs de Yang-Mills en dimension 4 et la th\'eorie de Seiberg-Witten, apparait la fonctionelle de Yamabe et l'invariant de Yamabe et la courbure scalaire et les metriques conformes, dans des articles de M.Gursky par exemple en dimension 4:

\smallskip

 Si on prend une vari\'et\'e Riemannienne compacte sans bord $ (M^4,g) $, de dimension 4 avec champ de Yang-Mills, cela veut dire que le tenseur de Riemann, $ Riem=F_{\nabla} $ est un champ de Yang-Mills. Comme la courbure scalaire $  S_g $ se d\'eduit du tenseur de Riemann (elle est li\'ee au tenseur de Riemann), alors, $ S_g $ est un objet de la physique et l'Equation de la courbure scalaire prescrite, en dimension 4, sur les vari\'et\'es de Yang-Mills, est li\'ee \`a la physique puisqu'elle contient $ S_g $. De meme pour la metrique conforme, qui est li\'ee \`a $ g $, la metrique Riemannienne $ g $ est appel\'e en dimension 4 un instanton.

\smallskip

La th\'eorie de Yang-Mills existe sur des vari\'et\'es non-compactes ou completes. (non compact, or complete Manifolds with harmonic curvature). Voir Taubes et les articles de Gabor Etesi et Yawei Chu. Exemples: vari\'et\'es d'Einstein (Einstein manifolds) $ n\geq 3 $ et les vari\'et\'es loc.conf.plates de courbure scalaire constante (pour $ n\geq 4 $, en particulier pour $ n=4 $) ou Ricci parallel, en particulier les vari\'et\'es de courbure sectionelle constante.

\smallskip

 Du point de vue de la physique, en dimension 4, on peut considerer les champs de Yang-Mills sur une vari\'et\'e compacte ou complete ou non compacte et localement on a une notion de champs de Yang-Mills. On peut considerer au depart la vari\'et\'e et puis on regarde ce qui se passe localement, on a encore des champs de Yang-Mills. Puis, on regarde ce qui se passe globalement ou localement (mesures, obsevation, \`a partir des donn\'ees et du champs de Yang-Mills).

\smallskip

Comme pour l'equation de Schrodinger: La donn\'ee est la courbure scalaire prescrite $ V $ (elle est prescrite, pulsion ou signal ou potentiel), la solution $ u $ est une "fonction d'onde" qui lie $ V $ \`a $ S_g $, $ S_g $ est un "champ scalaire de Yang-Mills". L'equation de la courbure scalaire prescrite en dimension 4 et sur une vari\'et\'e de Yang-Mills de dimension 4 est un objet de la physique. 

\smallskip

Aussi, c'est ecrit dans l'article de T.H. Parker (Gauge Theories on four dimensional Riemannian manifold, Comm.Math.Physics, 1982), si on considere un terme ("potentiel de Higgs") du type $ P(u)= au^4+mu^2 $ avec $ u >0$ la solution et $ a, m $ deux fonctions r\'eelles avec $ m \leq 0 $ possible, l'equation bosonique dans un champs de Yang-Mills devient du type Yamabe (ou de type courbure scalaire prescrite). ($ \Delta_g u + (\frac{s_g}{6}+2m) u=-4au^3 $, $ \Delta_g=-\nabla^i \nabla_i $ et $ s_g $ la courbure scalaire). On conclut, voir l'article de T.H. Parker (dans un champ de Yang-Mills, l'equation bosonique, induit (avec le fait qu'on a un champ de Yang-Mills, les lagrangiens s'annulent et la d\'eriv\'ee par rapport \`a la connexion est nulle), que le r\'esidu est nul et constitue une solution de couplage des champs).(le lagrangien Yang-Mills bosonique, $ Y=B_1+B_2 $ avec $ B_1 $ la fonctionnelle de Yang-Mills et $ B_2 $ la lagrangien bosonique, la loi de la particule: $ \dot Y= dY=0= dB_1+dB_2=0 $, si on se place dans un champ de Yang-Mills, alors, $ dB_1=0 $, ce qui implique que $ dB_2=0 $, si de plus, $ u >0 $ est solution de l'equation bosonique (cela veut dire que $ \frac{\partial B_2}{\partial u} =0 $), ceci implique que $ \frac{\partial B_2}{\partial \nabla}=0 $ et le r\'esidu est nul: $ \Sigma <\nabla_j u|\rho (e_i) u>=0 $ et donc, les deux equations du systeme de couplage sont v\'erifi\'ees. Donc l'equation bosonique (du type courbure prescrite) dans un champ de Yang-Mills suffit \`a determiner le systeme.Ces deux conditions, champs de Yang-Mills et equation bosonique constituent une solution particuliere du systeme.)

On peut aussi considerer chaque lagrangien seul dans n'importe quelle formulation, par exemple ici, celui de Yang-Mills seul et celui bosonique seul. Ici (dans l'article de Parker, T.H) on a le lagrangien de Yang-Mills bosonique.

\smallskip

//////////////////////////////////////////////////////////////////////

\smallskip

34) {\bf En dimension 2: equation de Liouville:} 

\smallskip

L'equation du type Liouville et sa relation \`a la physique et la chimie et l'astronomie:

\smallskip

a) Equation de la courbure scalaire prescrite sur la sphere de dimension 2 (vortex equation) ou un ouvert de $ {\mathbb R}^2 $ (vortex equation) ou surface de dimension 2 (vortex equation).  Mean-Field-equation. 

\smallskip

Equation de Liouville: equation de l'emballement thermique.

\smallskip

b)Voir les articles de Crandall-Rabinowitz, De Figueiredo-Lions-Nussbaum, Chen-Li, pour la provenence de ses equations: G\'eometrie (courbure de Gauss), Gazs, combustion, astronomie et astrophysique.

c) En dimension 2, aussi, c'est un cas particulier des champs de Yang-Mills, la theorie de Glashow-Weinberg-Salam, modelise les interactions (faibles) electromagnetiques des particules. C'est la theorie de jauge (theorie de champs) avec un groupe de jauge (groupe de symetries locales) $ U(1) \times SU(2) $ 

voir le livre de G. Tarantello Self-dual Gauge theories.

\smallskip

d) L'equation de Liouville ou de type Liouville apparait aussi dans le ph\'enom\`ene de "cordes cosmiques", un objet de l'univers "cosmic strings" et \`a ne pas confondre avec la theorie des cordes. Voir l'article de J. Spruck et Yisong Yang (cosmic strings).(On prend une metrique Lorentzienne $ ds^2=-dt^2 +dz^2+g_{ij}dx^idx^j $, la metrique $ g_{ij} $ est definie sur une surface, $ M $ de dimension 2 et apres on peut choisir $ g_{ij} $ conforme a une autre metrique, par exemple $ g_{ij}=e^u \delta_{ij} $ et on choisit le champ particulier). C'est l\`a, par exemple l'id\'ee g\'en\'erale.

e-En dimension 2: c'est aussi, la theorie des champs de Liouville en dimension 2. Qui a \'et\'e etudi\'ee pour des surfaces \`a courbure n\'egative et plus r\'ecemment sur la sphere de dimension 2. (LCFT, Liouville conformal field theory (2+1)).

\smallskip

De meme ici, du point de vue de la physique, de la chimie, astronomie ou la physique quantique, en dimension $ n=2 $, on peut considerer des sources, des champs sur un ouvert de $ {\mathbb R}^2 $ ou une surface compacte ou complete ou non compacte et localement on a une notion de champs local. On peut considerer au depart l'ouvert ou la surface et puis on regarde ce qui se passe localement, on a encore des champs. Puis, on regarde ce qui se passe globalement ou localement (mesures, obsevation, \`a partir des donn\'ees et du champs ou des parametres).

\smallskip

Ici aussi, comme pour l'equation de Schrodinger: La donn\'ee est la fonction $ V $, ($ V $ est la source ou pulsion ou le signal ou le potentiel), la solution $ u $ est une "fonction d'onde" ou emission, qui lie $ V $ \`a la donn\'ee de d\'epart, par exemple $ S_g $ (dans le cas d'un ouvert $ S_g=0 $), $ S_g $ est un "champ scalaire". Il se peut que ce ne soit pas $ S_g $ comme pour l'operateur $ \Delta + \epsilon (x_1\partial_1+x_2\partial_2) $, $ \Delta=\partial_{11}+\partial_{22} $.

\smallskip

/////////////////////////////////////////////////////////////////////////////////////////////////////////////////////////

\smallskip

35) {\bf En dimension $ n\geq 3 $: Relativit\'e g\'en\'erale ($ n=3 $) et Cosmologie quantique ($ n\geq 3 $):}

\smallskip

L'equation d'Einstein-Lichnerowicz, en cosmologie, cosmologie quantique: on associe la relativit\'e g\'enerale \`a la m\'ecanique quantique. La gravit\'e quantique.

\smallskip

Cette id\'ee de partir d'un espace Lorentzien $ n+1 $ et de choisir une metrique conforme, sur la vari\'et\'e $ M $, en partant de l'espace temps $ M\times (0,\delta) $,  se g\'eneralise aux dimensions superieures. (Il y a la correspondance, ADS-CFT(Anti-de-Sitter, Conformal field theory) propos\'ee par Juan Maldacena ($\Lambda <0 $, constante cosmologique strictement n\'egative dans l'equation d'Einstein). Les equations d'Einstein classiques et quand la dimension $ n\geq 3 $ c'est la "quantum-cosmology qui correspond au cas $ \Lambda =0 $, equation d'einstein du d\'ebut, (1915)). On a $ ds^2=-dt^2+g_{ij} dx^idx^j $ avec $ g_{ij} =u^{4/(n-2)} g_0 $, une metrique conforme a celle de d\'epart $ g_0 $ sur $ M $. 

\smallskip

L'id\'ee g\'en\'erale est la suivante: (On part de la formulation Lorentzienne $ (n+1), M\times (0, \delta) $, d'Einstein en relativit\'e g\'enerale et on introduit sur la variable espace des metriques conformes dans la fonctionelle du champ (ici au lieu de la fonctionelle du champ comme dans la th\'eorie de Yang-Mills, on a l'equation d'Einstein), on a des champs particuliers li\'es au changement de metriques conformes, tout cela dans le but d'etudier les interactions des particules et des astres (gravitation d\^ue \`a la courbure de l'espace temps ou \`a la formulation en $ n+1 $ et les ondes qui sont les effets de la gravitation ou des champs, en particulier gravitationnel, un autre champ est le champ $ T $ d'energie-impulsion, peut etre vu comme une donn\'ee, l'interaction est ecrite dans l'equation d'Einstein ci-dessous, gravitation, $ Ricci^{\gamma}, R^{\gamma} $, et les donn\'ees ou un autre champ, est le tenseur $ T $). C'est l'adpatation \`a $ (n+1) $ de la th\'eorie d'Einstein, en relativit\'e g\'en\'erale, qui est en $ (3+1) $. L'espace Anti-de-Sitter est une solution particuliere des equations d'Einstein, comme l'espace de Minkowski ou la metrique de Schwazrchild, l'inconnu ici est un espace Lorentzien $(V, \gamma) $ de metrique Lorentzienne $ \gamma $, qui verifie l'equation d'Einstein $ Ricci^{\gamma}-\frac{1}{2} R^{\gamma} \gamma = T $, avec $ Ricci^{\gamma}, R^{\gamma}, T $ respectivement le tenseur de Ricci et la courbure scalaire et le tenseur d'Energie-impulsion. Dans certains cas, cela revient \`a resoudre ce qu'on appelle le probleme de Cauchy en relativit\'e.par exemple $ T=0 \Rightarrow Ricci^{\gamma}=0 $. Un autre cas particulier, on en parl\'e ci-dessus, est de trouver des solutions du type $ V=M\times (0,\delta) $ et $ \gamma=-dt^2+g_{ij} dx^idx^j $ avec $ g=u^{4/(n-2)} g_0 $.  Resoudre les equations des contraintes pour $ (V, \gamma) $ particuliers avec la m\'ethode conforme de Lichnerowicz. On prend $ V=M\times (0,\delta) $ avec $ \gamma = -dt^2+g_{ij}dx^i dx^j $, on transfome les equations d'Einstein en 2 Equations des contraintes (Probleme de Cauchy, par les relations Gauss-Codazzi), puis on prend $ g=u^{4/(n-2)} g_0 $ sur $ M $, on transforme les 2 equations des contraintes en un systeme de deux equations elliptiques. Einstein scalar field Lichnerowicz equations. Voir, Hebey, Pollack, Pacard, Choquet-Bruhat, Bartnik.)

\smallskip

{\bf Remarque:} Pour ce qui est de la th\'eorie de Yang-Mills, on s'interesse \`a une fonctionelle du champ, on peut se ramener \`a la dimension 2 et utiliser la methode conforme comme dans Spruck-Yisong (cosmic-strings) ou d'autres th\'eories avec fonctionnelle de champs sur une vari\'et\'e Lorentzinne de dimension 4. ($ ds^2=-dt^2+dz^2+g_{ij} dx^idx^j $) (exemple, Glashow-Weinberg-Salam, combine fonctionnelle de Yang-Mills et theorie conforme en 2 dimensions). Ici dans ADS-CFT ($ \Lambda <0 $ ou bien $ \Lambda=0 $ et on retrouve les Equations d'Einstein classiques avec $ n=3 $ pour la relativit\'e g\'enerale et $ n\geq 3 $ pour la cosmologie quantique ou "quantum cosmology"), on remplace la fonctionnlle par l'equation d'Einstein. De meme, ici, dans la relativit\'e g\'enerale ($ n=3 $) et pour la cosmologie quantique ou "quantum cosmology" ($ n\geq 3 $), on remplace la fonctionnlle par l'equation d'Einstein.

\smallskip

a) il y a une application et r\'ealisation en physique, en $(3+1) $ et $ (4+1) $ et $ (6+1)$ de la correspondance de ADS-CFT, de Juan.Maldacena (on voit bien le probleme du bon espace-temps). Et aussi $ (2+1) $ (LCFT, Liouville conformal field theory).

La correspondance ADS-CFT correspond a l'equation d'Einstein avec une constante cosmologique $ \Lambda <0 $ (Anti-de-Sitter): $ Ricci^{\gamma}-\frac{1}{2}R^{\gamma} g+\Lambda g= T=0 $. Si $\Lambda =0$, on retrouve l'equation d'Einstein (1915) ecrite plus haut, et les th\'eories en question sont: $ n=3 $, la relativit\'e g\'en\'erale et $ n\geq 3 $ la cosmologie quantique ou "quantum cosmology". 

\smallskip

Il y a des auteurs qui font r\'ef\'erence \`a la correspondance ADS-CFT, qui essaie de trouver le bon espace-temps avec $ \Lambda <0 $. Dans la th\'eorie classique d'Eisntein($ \Lambda=0 $), c'est la "quantum cosmology", trouver le bon espace-temps pour $ n\geq 3 $.

\smallskip

b) Pour ce qui est des equations d'Einstein et methode conforme: ($ n\geq 3 $ ("quantum-cosmology") et $ n=3 $ qui correspond \`a la th\'eorie de la relativit\'e g\'en\'erale:) Voir les articles de Choquet-Bruhat-Isenberg-Pollack, dans general relativity and quantum cosmology. Cosmologie quantique pour $ n\geq 3 $. (il y a des auteurs qui font r\'eference \` la correspondance ADS-CFT, qui essaie de trouver le bon espace-temps pour $ n \geq 3 $ avec $ \Lambda <0$. Ici, on a la th\'eorie classique d'Einstein ($ n=3 $ relativit\'e g\'en\'erale et "quantum cosmology" pour $ n\geq 3$ avec  $ \Lambda =0 $)).

\smallskip

(Equation d'Einstein classique: $ \Lambda =0 $, (1915) cosmologie). Voir le monograph de: E.Hebey, D. Pollack, F. Pacard, Y. Choquet-Bruhat, R. Bartnik: dans la th\'eorie des champs de Klein-Gordon massive le potentiel $ V(\Psi)=\frac{1}{2} m^2 \Psi^2 $. (les donn\'ees $ (\Psi, \sigma, \tau, \pi) $ peuvent etre choisies librement ainsi que la metrique $ g_0 $ de la vari\'et\'e Riemannienne). 

\smallskip

Voir aussi les articles de Choquet-Bruhat-Isenberg-Pollack, la vari\'et\'e Riemanninne $ (M, g_0) $ peut etre compacte sans bord ou non compacte. Voir la These de C.Valcu.

c) On part des equations de contraintes qui derivent des equations d'Einstein et qui sont obtenues par les equations de Gauss-Codazzi (on contracte le tenseur de Riemann dans les equations de Gauss-Codazzi pour avoir le tenseur de Ricci). On contracte le tenseur de Ricci pour faire apparaitre la courbure scalaire et on separe les termes en variables temporelles et spaciales. Puis on utilise une metrique conforme.

d) Inversement, d\`es qu'on a une solution des equations de contraintes avec champ scalaire, il existe une solution des equations d'Einstein, un espace-temps maximal unique avec une metrique Lorentzienne, d'apres Choquet-Bruhat et Choquet-Bruhat-Geroch. C'est ce qu'on appelle le formalisme de Choquet-Bruhat-Geroch-Lichnerowicz.(Il existe un espace temps $(L,\gamma)$ maximal et unique \`a isometrie pr\'es avec $ \gamma $ Lorentzienne, et un plongement $ i:M\to L $, tel que $ i^*(\gamma)=g $, $(L,\gamma) $ 'prolonge' $ (M,g) $ et les autres parametres, comme le champ scalaire par exemple). Voir le print de B. Premoselli et l'article de A. Carlotto (The general relativistic constraints equations).

\smallskip

Pour chaque solution des equations des contraintes on a une solution des equation d'Einstein, on a une multitude d'univers ou plusieurs espace-temps. La multitude d'espace-temps, modelise aussi la multitude d'etats quantiques, en particulier c'est la combinaison entre la gravit\'e et les differents champs dont le champ electromagnetique et d'autres champs: Yang-Mills, Klein-Gordon, fluides, etc. Cela veut dire que du point de vue mathematique on a plusieurs espace-temps alors que, du point de vue de la physique, il y a plusieurs etats quantiques ou plusieurs fonctions d'ondes a l'instant $ t=0 $, dans le meme espace ou vari\'et\'e $ M $ avec des propri\'et\'es concernant les fonctions d'ondes ou etats quantiques (Plusieurs ondes au meme instant qui sont dans un espace commun, et, differents au meme temps). Ceci est un modele mathematique pour expliquer l'existence de plusieurs etats quantiques ou fonctions d'ondes, qui sont a la fois differents, dans plusieurs espaces differents, et dans un meme espace au meme temps: $ (M,g_0, \psi, \sigma, \tau, \pi) $, modele de base commun a tous les espaces, qui est lui meme un 'univers' de 'parametres'. Ph\'enomene physique $ \leftrightarrow $ Propri\'et\'e mathematique: i) Etat quantique $ \leftrightarrow $ fonction d'onde $ u >0 $, et, ii) extra-dimensions ou dimensions supplementaires $ \leftrightarrow $ propri\'et\'e des fonctions d'ondes: $ \sup u = f(\inf u) $, (notion d'enroulement, de torsion), et les 2 valeurs $ (\sup, \inf) $. 

\smallskip

Dans un article de Choquet-Bruhat (1968), on ne peut pas avoir existence globale d'un espace temps car d'apres Penrose et Hawking, les Eq de la relativit\'e d'Einstein developpent des singularit\'es. Mais on a existence locale d'un espace temps $ (-t_0, t_0)\times M, t_0 >0 $ et unicit\'e locale et globale.

\smallskip

Th\'eorie de la graviation quantique: du type Kaluza-Klein. La th\'eorie de Kaluza-Klein c'est en $ (4+1) $ comme espace-temps avec une dimension cach\'ee. Les th\'eories du 'type' Kaluza-Klein, c'est en $ (5+1), (6+1), (n+1), n=9,10 $ ou $ (n+1), n\geq 3 $ avec plusieurs dimensions cach\'ees.

\smallskip

Pour ce qui est des solutions des Equations d'Einstein par la m\'ethode conforme de Lichnerowicz-Choquet-Bruhat-York: on a:

i) Si on prend $ \Psi =0 $, $ \tau \equiv constante \not =0 $, $ \sigma =\pi = 0 $, alors $ W=0 $ est une solution et $ u >0 $ est solution de l'equation de Yamabe (dans "le cas n\'egatif", la vari\'et\'e n'est pas n\'ecessairement, compacte sans bord).

\smallskip

ii) Si on prend $ \Psi\equiv constante \not =0 $, $ \tau =0 $, $ \sigma = \pi = 0 $, alors $ W=0 $ est une solution et $ u >0 $ est solution de l'equation de Yamabe (dans "le cas positif", la vari\'et\'e n'est pas n\'ecessairement, compacte sans bord).

\smallskip

iii) Si on prend $ \Psi \not \equiv constante $, $ \sigma = \pi =0 $ et $ \tau =0 $ alors $ W=0 $ est une solution et $ u >0 $ est solution de l'equation du type courbure scalaire prescrite avec $ h=S_{g_0}-|\nabla \psi|_{g_0}^2 \leq S_{g_0} $ (dans "le cas positif", la vari\'et\'e n'est pas n\'ecessairement, compacte sans bord).

\smallskip
On choisit $ W $ champ de vecteur de Killing conforme ($ D W=0 $, conformal Killing vector field), pour avoir un champ conforme et etre compatible avec la th\'eorie des champs conformes. Ici, on a pris $ W=0 $ (pour dire qu'il y a au moins toujours une solution), mais on peut prendre $ DW=0$. on a une th\'eorie des champs, gravitation+Klein-Gordon, ou Yang-Mills, ou Electromagnetique, et une th\'eorie conforme des champs (on avait deja une metrique conforme introduite par la methode de Lichnerowicz-York-Choquet-Brhuat-Isenberg-Pollack, c'est le contexte conforme. Ici conforme pour dire que les objets sont invariants par transformations conformes, 'Scale invariants' ou 'Zoom invariants').

\smallskip

De meme ici, du point de vue de la physique, de l'astronomie ou de la cosmologie ou de la cosmologie quantique, en dimension $ n=3 $ ou $ n\geq 3 $, on peut considerer les champs sur une vari\'et\'e compacte ou complete ou non compacte et localement on a une notion de champs local. On peut considerer au depart la vari\'et\'e et puis on regarde ce qui se passe localement, on a encore des champs. Puis, on regarde ce qui se passe globalement ou localement (mesures, obsevation, \`a partir des donn\'ees et du champs ou des parametres, comme $ (\Psi, \sigma, \tau, \pi) $ et $ (M, g_0) $).

\smallskip

Ici aussi, comme pour l'equation de Schrodinger: La donn\'ee est la fonction $ V $ (du type courbure scalaire prescrite) qui peut etre constante (Yamabe), ($ V $ est prescrite (ou Yamabe, $ V\equiv constante $)), on peut prendre $ V=\frac{1}{2} m^2 \Psi^2 $, $ V $ est la pulsion ou le signal ou le potentiel), la solution $ u $ est une "fonction d'onde" qui lie $ V $ \`a $ S_{g_0} $ ou $ S_{g_0}-|\nabla \psi|_{g_0}^2 $, $ S_{g_0} $ est un "champ scalaire".

\smallskip

///////////////////////////////////////////////////////////////////////////////////////////////////////////////////////////////////

\smallskip

36) {\bf En dimension $ n\geq 3$: r\'esultat d'unicit\'e et de rigidit\'e: en Biologie Math\'ematique:}

\smallskip

Interpretation en Biologie Math\'ematique: la Chimiotaxie (Chemotaxis en anglais), est l'etude du mouvement de certaines cellules (les amibes, amoebae en anglais) dans leur environement, mouvement d\^u  \`a la liberation de substances chimiques par les amibes. les amibes se deplacent vers les endroits ou la concentration de cette substance chimique est la plus grande et s'agregent, forment des agregats. L'interpretation biologique du resultat d'unicit\'e est que les amibes ne peuvent pas former des agregats (ne se concentrent pas) quand on a une solution unique et constante au systeme de Keller-Segel. Voir l'article de Lin-Ni-Takagi. Le systeme de Keller-Segel peut se reduire \`a l'etude d'une seule equation, voir l'article de Lin-Ni-Takagi, Journal of Diff. Equations, 1988. Cette formulation existe sur les vari\'et\'es compactes sans bord de dimension $ n\geq 3 $, voir T. Hillen et K.Painter.

\smallskip

{\bf Remarque 1:} Il y a aussi les D-branes en dimension 2, cordes ouvertes avec condition de Dirichlet.(fait reference a l'equation du type Liouville (equation de mouvement des particules ou de la corde) avec condition de Dirichlet). Pour la theorie conforme des champs de Liouville. Il y a aussi les D-branes en dimension $ n\geq 3$, on a Kaluza-Klein avec conditions de Dirichlet.

\smallskip

{\bf Remarque 2:} En ce qui concerne l'existence de dimensions supplementaires (extra-dimensions) dans les th\'eories du type Kaluza-Klein: ceci est li\'e a l'existence  d'une nouvelle perticule qui est l'axion qui se deplace dans ces dimensions cach\'ees, on peut les detecter ou savoir comment les eliminer avec des outils, des telescopes, voir les articles ou preprint de Horvat, Krcmar, Lakic. Il y a le telescope du CERN, et aussi d'autres appareils, pour essayer de detecter ou eliminer ces particules que sont les axions. C'est dit dans le print de Horvat, Krcmar, Lakic, les dimensions supplementaires $ (5+1), (6+1) $ sont li\'es a l'astrophysique, elles ont des r\'ealisations. (4+1) c'est la formulation Kaluza-Klein initiale.

\smallskip
"Du type courbure scalaire prescrite"$ \leftrightarrow $ Schrodinger=Dynamique d'une particule non relativiste.

"Du type Yamabe"$\leftrightarrow $ Schrodinger=Dynamique d'une particule non relativiste.

\smallskip

//////////////////////////////////////////////////////////////////////////////////////////////////////////////////////////

\smallskip

37) Des Remarques sur le tenseur metrique en polaires \`a l'ordre 2: 

Pour ce point voir l'expression du tenseur metrique en coordonn\'ees g\'eodesiques polaires. On le calcul pour des vecteurs de la sphere particuliers en utilisant des transformations orthogonales. La carte choisit sur la sphere est celle des coorodonn\'ees spheriques sur la $ (n-1)-$sphere:

Ici, on utilise la formule suivante qui donne le coefficient de la metrique d'ordre $ r^2$:

$$\frac{1}{3} R_{kpql}(y_i)\theta^p\theta^q \partial_i\theta^k\partial_j\theta^l. $$

On a: $ \theta_1= \cos \tilde \theta_1, \theta_2=\sin \tilde \theta_1\cos \tilde \theta_2, ...., \theta_{n-1}=\sin\tilde \theta_1\sin\tilde \theta_2....\sin\tilde \theta_{n-1} \cos\tilde \theta_{n-1}, \theta_n=\sin\tilde \theta_1\sin\tilde \theta_2...\sin\tilde \theta_{n-2}\sin \tilde \theta_{n-1} $. le $ (n-1)-$uplet est $ (\pi/2,\pi/2,...,\pi/2) $, le vecteur $ e_n=(0,0,...,1) $. Alors le terme d'ordre $ r^2$ de la metrique est $ R_{innj}(t_0)$, $ t_0=P_0=x=y_i$.

\smallskip

On prend le $(n-1)-$uplet, $ (\pi/2,\pi/2,...,\pi/2)$.

\smallskip

Quand on remplace la composante $ i $ par celle de $ n $ et celle de $ n  $ par celle de $ 1 $, on a $ \tilde g_{11}=R_{niin}, 1\leq i\leq n-1 $ (et on sait que $ R_{nnnn}=0 $). On avait un terme contenant $ R_{1nn1}$.

\smallskip

Quand on conserve les composantes $ 1 $ et $ 2 $ et on remplace la composante $ i $ par celle de $ n $, on a $ \tilde g_{12}=R_{1ii2}, 3\leq i \leq n $ et on sait que $ R_{1112}=R_{1222}=0 $. 

\smallskip

De meme, pour $ \tilde g_{mm} $, on remplace la composante $ i $ par celle de $ n $ et celle de $ n $ par celle de $ m $. on obtient $ \tilde g_{mm}=R_{niin} $, pour $ i\not =m $ et on avait un terme $=R_{mnnm} $, et on sait que $  R_{mnnm}=R_{nmmn} $. Avec, $ R_{nnnn}=0 $.

\smallskip

De meme, pour $ \tilde g_{ij} $, on conserve les composantes, $ i $ et $ j $ et pour $ k \not = i, j $, on remplace la composante $ k $ par celle de $ n $. On a alors, $ \tilde g_{ij}=R_{ikkj} $ pour $ k\not = i, j $, mais on sait que $ R_{iiij}=R_{ijjj}=0 $. On avait un terme contenant $ R_{innj} $.

\smallskip

1) ce n'est pas vrai, que l'ordre soit forc\'ement superieur a $ r^2$, on ne sait rien de cela, dans le theoreme 1.53 du livre d'Aubin, on doit prendre $ b >0, a >0 $. Le tenseur metrique est de l'odre de $ r^2 $ et sa deriv\'ee est de l'ordre de $ r $. Il n'y a plus rien a faire.

\smallskip

2) Si on suppose que $ Ricci\not=0 $, il y a au moins un terme $ Ricci_{ij} \not =0$, cela veut dire qu'il y a au moins un terme de la courbure de Riemann $ R_{ikkj}\not =0$. Or on a vu par les transormations orthogonales qu'en se placant en un des vecteurs de base $ e_i $, on a le terme d'ordre $ r^2 $ de la metrique est un $ R_{ikkj} \not =0$. Cela veut dire qu'on ne peut pas l'eliminer.

Donc, on ne peut pas ameliorer l'ordre de la deriv\'ee du tenseur metrique, sur une vari\'et\'e quelconque en particulier $ Ricci \not =0$. Sur une vari\'et\'e quelconque, on ne peut pas depasser l'odre $ r^2 $ dans la metrique en coordonn\'ees geodesiques polaires. Il faut supposer au moins que $ Ricci\equiv 0$. Mais dans le cas general ce n'est pas possible.

\smallskip

3) Pour arriver par un changement de metrique conforme a avoir $ Ricci=0 $, voir le livre de Hebey, il faut avoir $ Weyl=0$ et resoudre une equation en utlisant les conditions de Frobinius, donc, il faut supposer une condition sur $ Weyl $ ou supposer que $ Weyl=0 $, ce qui veut dire que la vari\'et\'e est loc.conf.plate, ce qui ne sert a rien car on est dans le cas, non.loc.conf.plat.

\smallskip

4) Meme si on suppose par exemple tous les termes de la forme $ R_{ikkj}=0 $, avec un double indice quelque part, ce qui implique que $ Ricci=0 $, alors il y aura toujours un point ou on a un terme du type $ R_{ijkl} $, comme le tenseur de Riemann est non nul ceci implique que la metrique a un terme d'ordre $ r^2 $. Par exemple avec $ \tilde \theta_1\not= 0, \pi/2$ et $ \tilde \theta_2=...=\tilde \theta_{n-1}=\pi/2 $, on conserve le terme $ R_{2n13} $, en permutant les lignes, on conserve $ R_{2kl3} $, comme les indices $ 2 $ et $ 3 $ sont quelconques, on a au moins un point avec $ R_{1234} $ ou plus g\'eneralement $ R_{ijkl} $.

En fait, on a soit $ R_{2n13} $ soit $ R_{21n3} $(on a calcul\'e $ g_{23} $) et on regarde aussi $ R_{23n1} $ et $ R_{2n31}=-R_{2n13} $ (on echange certaines lignes, et on calcul $ g_{23} $), et on utlise l'identit\'e de Bianchi.

\smallskip

5) Quand on consid\`ere l'equation de Yamabe, cela veut dire qu'on considere une courbure prescrite $\equiv 1 $. Le gradient est nul dans toutes les directions. Cela veut dire que toutes les directions sont prises en compte ou bien aucune direction n'est previl\'egi\'ee, la symetrie de la sphere. La technique utilis\'ee dans la recherche d'estimation a priori est la technique moving-sphere, elle optimise la recherche d'estimations a priori et ne previligie aucune direction. Pour les dimensions $ n=3, n=4, n=5, n=6 $, l'ordre d'annulation du tenseur metrique en polaire, permet d'avoir l'estimation a priori, jusqu'a l'ordre $ r^2$, a partir de cet ordre, il y a une direction previligi\'ee, or ceci n'est pas possible car on a suppos\'e qu'aucune direction n'est previligi\'ee. Donc le cas limite est $ n=6 $, pour les dimensions superieures \`a 7, il faut supposer des hypoth\`eses supplementaires. Pour le cas general, le cas limite est $ n=6 $.

Pour le cas limite, $ n=6 $, on a utilis\'e cette invariance par symetrie de la sphere, car les solutions tendent vers une fonction radiale.

Pour la dimension $ n=4 $ et l'Eq. de la courbure prescrite ou du type courbure prescrite (presence d'un potentiel), on donne des conditions sur le potentiel (petit en norme Lipschitz, proche d'une constante et la constante $ >0 $ est quelconque, petite ou grande, oscille autour d'une constante $ >0 $), pour pouvoir appliquer cette methode moving-sphere en presence d'un potentiel (possibilt\'e d'une direction pr\'evil\'egi\'ee).

6) Si, on regarde la th\'eorie de Kaluza-Klein, ou des extra-dimensions, il faut que la propri\'et\'e soit vraie pour toute vari\'et\'e, pour que la Probabilt\'e que la formulation soit vraie, soit de $ 1 $. Ceci est possible en dimensions, $n=3, n=4, n=5, n=6 $, ou bien les conditions supplementaires soit compatibles avec la th\'eorie de Kaluza-Klein, par exemple, la compacit\'e de la vari\'et\'e et operateur conforme coercif (ce qui est prouv\'e dans ce cas, par YY.Li-L.Zhang, par contre Choquet-Bruhat-York-Isenberg-Pollack ont consid\'er\'e le cas ou la vari\'et\'e n'est pas necessairement compacte). Pour le cas g\'eneral, le cas limite est $ n=6 $.

-On a vu ci-dessus qu'il y avait une obstruction geometrique, de la courbure de Riemann (tenseur de Riemann): "obstruction de la gravit\'e" li\'ee \`a la courbure, car la gravit\'e est li\'ee \`a la courbure; obstruction d'astronomie, th\'eorique.

-On a aussi des conditions astrophysiques, obstructions numeriques, obstructions astrophysiques, concretes, exp\'erimentales, en, $ 5+1, 6+1 $, voir l'article de, G.F. Giudice et J.D. Wells; et le cas $ 4+1 $ aussi; Extra-dimensions, dans Review of Particle Physics, 2005-2006. Dans, Particle Data Group, nouvelle edition en 2021.

En ce qui nous concerne, la Th\'eorie de Kaluza-Klein, c'est en $ 4+1, 5+1, 6+1 $. Et le cas $ 3+1 $, c'est la relativit\'e g\'enerale. On a bien le cas limite est $ n=6 $.

Dans la Th\'eorie des Supercordes c'est en $ 9+1 =10 $ dimensions spatio-temporelles et les dimensions supplementaires forment un espace a $ 6 $ dimensions du type Aubin-Calabi-Yau.

Chaque Th\'eorie explique des ph\'enom\`enes de la physique, d'astronomie, d'astrophysique (comme le disent Horvat, Krcmar, Lakic, voir leur prints sur les axions et les extra-dimensions) ou(et) les interactions des particules.

\smallskip

//////////////////////////////////////////

\smallskip

On va mieux expliquer en utilisant les contre-exemples de, C.C.Chen-C.S.Lin, Journal of Diff.Geometry.1998:

\smallskip

On considere des metriques conformes (au sens des matrices):

$$ g_u=u^{4/(n-2)} g=u^{4/(n-2)} (dr^2+r^2(1+Weyl\cdot r^2+\ldots) d\theta^2), $$ 

la courbure scalaire de $ g_u $ est $ V $, les contre exemples de Chen-Lin, sont possibles si la platitude est d'ordre $ [(n-2)/2] $, avec un coefficient ne tendant pas vers 0. La courbure scalaire est deduite du tenseur de Riemann, qui est deduit de la metrique en derivant 2 fois. Donc, pour obtenir $ V $, il faut deriver 2 fois $ g_u $, donc il faut deriver 2 fois $ g $, or le developpement de $ g $ contient le tenseur de Weyl, a l'ordre 2, donc, quand on derive 2 fois $ g $, on obtient $ Weyl\not = 0 $. Donc la paltitude est possible jusqu'a l'ordre 2. Cela veut dire que du point de vue de la dimension, $ [(n-2)/2] =2 $, donc, la platitude est possible jusqu'a la dimension $ n=6 $, comme les coefficient devant le tenseur de Weyl, tendent vers $ 0 $, la dimension $ n=6 $ est incluse dans la platitude. cela veut dire que la dimension $ n=6 $ est le cas limite dans la notion de platitude. De plus l'exemple de Chen-Lin, dit que le meilleur qu'on puisse obtenir est $ \sup $ major\'e si $ \inf $ minor\'e, c'est a dire que l'in\'egalit\'e de Harnack explicite n'est pas facilement possible, c'est a dire que c'est difficile de l'obtenir.

\smallskip

Pour Chen-Lin, les contre exemples, commencent des que la platitude est d'ordre $ [(n-2)/2] $ pour $ V $, avec des coefficients ne tendant pas vers 0. Or il faut deriver $ g_u $  2 fois pour obtenir $ V $, mais $ g_u=u^{4/(n-2)} g=u^{4/(n-2)} (dr^2+r^2(1+Weyl\cdot r^2+\ldots)d\theta^2) $ (au sens des matrices), donc, il faut deriver 2 fois $ g $, or, deriver 2 fois $ g $ c'est obtenir $ Weyl\not =0$. Donc, les contre-exemples commencent des que $ [(n-2)/2]=2 $. Comme, pour le cas $ n=6 $, les coefficients devant $ Weyl $ tendant vers 0, les contre exemples, ne commencent que pour $ n\geq 7$. De plus, Chen-Lin, disent que pour la platitude $ [(n-2)/2] $, le meilleur qu'on puisse obtenir est $ \sup $ major\'e si $ \inf $ est minor\'e.(on n'a pas cela). Ce qui veut dire que le cas limite $ n=6$, ce cas est inclus dans les cas ou l'in\'egalit\'e est possible et que le meilleur qu'on puisse obtenir est $ \sup $ major\'e si $ \inf $ est minor\'e.

\smallskip

(La platitude de $ V $ est li\'ee \`a la platitude de $ \partial_{rr} (g) $, considerer la platitude de $ V $, c'est aussi considerer la platitude dans la classe: $ \partial_{rr}(g) $, il faut deriver 2 fois $ g $, c'est "dans la classe de derivation 2 fois". Si on prend $ V=1 $, il faut que cette platitude soit compatible avec la platitude de $ \partial_{rr}(g) $).

\smallskip

Donc, avec la notion de platitude et en particulier $ V=1 $ (equation de Yamabe), le cas limite est $ n=6 $ et le mieux qu'on puisse obtenir est $ \sup $ major\'e si $ \inf $ est minor\'e. On n'a pas forcement $ (\sup)^{\alpha} \times \inf \leq c $. On a (le mieux qu'on puisse obtenir): $ \sup $ major\'e si $ \inf $ est minor\'e pour le cas limite $ n=6 $, une in\'egalit\'e de Harnack implicite.

\smallskip

Donc, avec la notion de platitude, il faut considerer une classe de metriques plus petite. Pour $ n=6 $, en considerant n'importe quelle metrique $ g $ sur $ M $, le mieux qu'on puisse avoir c'est $ \sup $ major\'e si $ \inf $ est minor\'e, une in\'egalit\'e de Harnack implicite, on n'a pas forc\'ement une in\'egalit\'e explicite $ (\sup)^{\alpha} \times \inf \leq c $. Il faut diminuer la quantit\'e de metriques sur la vari\'et\'e $ M $. Par exemple pour les potentiels $ C1 $, en dimension 4, il faut prendre un espace sym\'etrique, Ricci plat. Quand on prend un espace symetrique, on inclut, les espaces euclidien pour lesquels, il y a des contre exemples de Chen-Lin \`a l'ordre: $ (n-2)/2, n\geq 5 $. Donc, il faut faire attention avec les espaces symetriques.

\smallskip

(Dans les contre-exemples de Chen-Lin (1998), il faut que: $ \int_{{\mathbb R}^n} Q(\xi+y) U_0^{2n/(n-2)} dy < 0 $, donc, les contre exemples sont valables, par exemples, qu'avec $ Q < 0 $, par exemple, le blow-up, ne peut pas etre un minimum ($ Q >0$, ca ne peut pas blow-up). D'o\`u, les resulats de Chen-Lin (1997) et L.Zhang (2007), avec $ Q=\Delta V >0 $ positive. Ca peut etre d'ordre $ (n-2)/2 $ et $ Q >0$. Les contre exemples, sont construits pour l'ordre $ \alpha \geq (n-2)/2 $, $ \alpha >1 $ et $ Q <0$).

\smallskip

Donc, avec la notion de platitude on n'a pas l'in\'egalit\'e de Harnack explicite pour $ n=6 $ et les contre-exemples commencent \`a partir de $ n\geq 7$.

\smallskip

Pour obtenir des resultats plus g\'en\'eraux concernant le $ \sup \times \inf $, voir la partie et la notion de: g\'eom\'etrisation.

\smallskip

///////////////////////////////////////////////////////////////////////////////////////////////////////////////////////

\smallskip

38) Sur,  $ \sup u_i \leq f(\inf u_i) $, en dimension 4 pour l'Eq. de courbure scalaire prescrite ou du type courbure scalaire prescrite en dimension 4.

\smallskip

On a une sorte de "compacit\'e" dans l'espace des solutions $ (u,V)$ avec $ 0 < a \leq V\leq b <+\infty $ et $ \nabla V $ petit, $ V $ proche d'une constante : La propri\'et\'e est vraie pour tout ensemble denombrable de solutions avec $ V $ proche d'une constante. Si, on prend un "ensemble continu" de solutions, on peut en extraire une suite pour la quelle $ \sup u_i \leq f(\inf u_i) $. Si $ (u_t)_{ (t\in I)} $ est un ensemble de solutions relativement \`a  $ (V_t)_{(t\in I)} $ proche d'une constante, alors "on peut avoir pour un sous ensemble denombrable" $ J\approx {\mathbb N}, J \subset I $, $ J $ d\'enombrable, avec $ V $ proche d'une constante, $ \sup u_i \leq c_J(\inf u_i), \, i\in J $. Dans l'ensemble des solutions $ (u,V) $, pour une suite denombrable de solutions $ J $, $ \sup u_i \leq c_J(\inf u_i) $, c'est une sorte de compacit\'e de la propri\'et\'e  $ P_{\{\sup,\inf\}}=\{ \sup u\leq c(\inf u)\} $, dans l'ensemble des solutions: $ (u, V) $ avec $ V $ proche d'une constante. Ceci est coherent avec la notion d'estimation a priori, car on a la compacit\'e apres extraction de sous-suite d'un ensemble donn\'e. 

\smallskip

-Le d\'enombrable et les ondes: dualit\'e onde-corpuscule. Pour les ondes on a un flot continu d'emissions, alors que pour les corpuscules on a un nombre d\'enombrable d'emissions distinctes. Les deux interpretations de la mati\`ere sont possibles en physique, donc, il suffit d'avoir la propri\'et\'e pour un nombre d\'enombrable ou une suite. Ce qui est le cas en dimension 4, pour l'Eq. de la courbure scalaire prescrite et l'Eq. du type courbure scalaire prescrite avec un potentiel proche d'une constante $ >0 $. On a encore les notions d'enroulement, de torsion (et extra-dimensions) et les valeurs $ (\sup,\inf)$, en d\'enombrable.

\smallskip
-Le continu et minorant-majorant. Dans le cas du continu, emissions continues, on utilise l'estimation a priori, qui, \`a partir du minorant $ m_0 >0 $ donne le majorant  $ M_0=c(m_0)=c(a,b, m_0, K, M) $ et cette fonction est d\'ecroissante de $ m_0 >0 $, donc, quand on prend $ m_0 >0 $ petit et voisin de $ 0 $, $ 0 < m_0 \to 0 $, $ c(m_0) $ croit, $ m_0 \searrow 0 \Rightarrow c(m_0) \nearrow $, on voit alors l'amplitude g\'en\'er\'ee par la variation du minorant $ m_0 >0$, ainsi que l'enroulement et la torsion et les valeurs, $ (majorant, minorant)=(M_0, m_0) $, $ M_0=c(m_0) $, (et extra-dimensions). Avec ici, $ u >0 $ solution relativement \`a $ V >0 $ proche d'une constante, dans le cas continu. On a une "in\'egalit\'e de Harnack" implicite entre minorant et majorant dans le cas continu: il existe une relation entre le minorant et le majorant, not\'ee $ c=c(m_0)=M_0 $, d\'ecroissante, avec l'in\'egalit\'e suivante  $ 0 < m_0\leq M_0=c(m_0) <+\infty $:

$$ 0 < m_0 \leq \inf_M u \leq \sup_K u \leq M_0=c(m_0) <+\infty.$$

\smallskip

Ici, on resout le probleme du cas $ n=4 $ de $ \sup $ major\'e si $ \inf $ est minor\'e de Chen-Lin, voir le point 37) precedent. Dans les contre exemples de Chen-Lin, la regularit\'e est $ C^{\alpha}, \alpha =(n-2)/2, \alpha >1 $, pour ne pas avoir $ \sup $ major\'e si $ \inf $ est minor\'e, qui correspond aux cas $ n\geq 5 $. Il reste le cas de la dimension $ n=4 $. Alors, l'article de Journal.Math.Anal.Appl.2012, repond a cette question: on a bien $ \sup $ major\'e si $ \inf $ minor\'e pour des solutions avec un potentiel $ V $ voisin d'une constante, sans etre une constante forc\'ement, Lipschitzien et petit en norme Lipschitz. On traite alors, le cas de la dimension $ n=4 $ (sur toute vari\'et\'e Riemannienne, non necessairement sans bord, on donne des conditions sur $ V $: Lipschitzien et petit en norme Lipschitz, qui correspond \`a une hypothese plus faible que $ C^{\alpha}, \alpha=(n-2)/2=1, n=4 $, mais la constante de Lipschitz est petite). Dans le print avec potentiel $ C^1 $, on traite aussi, du cas $ n=4 $ avec un potentiel avec amplitude non necessaierment petite, pour un espace localement symetrique Ricci plat.

\smallskip

Comme les contre-exemples de Chen-Lin sont faits avec des suites. Pour mieux expliquer le resultat de la dimension 4, il faut utiliser des suites: On prend $ (u_i,V_i) $ avec $ |\nabla V_i|\leq A_i \to 0 $, alors, pour chaque $ m >0 $, il existe un rang tel que $ A_i \leq \epsilon_0=\epsilon_0(a,b,m), \epsilon_0 >0 $, donc on a la compacit\'e locale \`a partir d'un certain rang: $ \sup_K u_i \leq c $ si $ \inf_M u_i \geq m >0 $, pour $ i\geq i_0 $.

\smallskip

Les contre-exemples de Chen-Lin, c'est pour $ n\geq 5 $. Pour $ n=4 $, dans le print de 2024; "Harnack inequalities for equations of type...", quand le potentiel $ V $ a une amplitude quelconque ($ |\nabla V|\leq A $, $ A $ quelconque), on ne prouve pas: $ \sup $ major\'e si $ \inf $ minor\'e, ce qui veut dire qu'il est peu probable de prouver $ \sup $ major\'e si $ \inf $ minor\'e pour n'importe quelle vari\'et\'e riemannienne $ (M,g) $ et $ n=4 $. On prouve $ \sup $ major\'e si $ \inf $ minor\'e si l'amplitude de $ V $, $ A \to 0 $. Mais si $ A >0 $ est quelonque ceci est peu probable que ca soit vrai, pour n'importe qu'elle vari\'et\'e. D'o\`u le resultat avec potentiel $ C^1 $, sur les espaces localement symetriques Ricci plats. Donc, pour $ n=4$, on n'a pas de contre exemples, et les contre exemples de Chen-Lin, sont vrais pour $ n\geq 5 $, mais, le print de 2024, "Harnack inequalities for equations of type..." peut etre considerer comme un exemple illustrant l'impossibilit\'e de prouver $ \sup $ major\'e si $ \inf $ minor\'e pour $ n=4 $ avec potentiel $ V $ d'amplitude quelconque $ A >0 $ et sur n'importe quelle vari\'et\'e Riemannienne. D'o\`u le print sur les potentiels $ C^1 $, sur les espaces localement symetriques Ricci plats.

\smallskip

 $ \sup $ major\'e si $ \inf $ minor\'e en dimension 4, par exemple: c'est un critere de compacit\'e et aussi un r\'esultat de compacit\'e, car, les fonctions: $ u_{\epsilon}(x)= \epsilon/(\epsilon^2+|x|^2) $ sont telles que: $ u_{\epsilon}(1)=\inf_{B_1(0)} u_{\epsilon} \to 0 $ et $ \max_{B_1(0)} u_{\epsilon}= u_{\epsilon}(0) \to +\infty $. Donc, pour avoir la compacit\'e locale, il faut supposer l'inf minor\'e.

\smallskip

///////////////////////

\smallskip

Dans le cas g\'en\'eral, d'une vari\'et\'e riemannienne de dimension $ n=4 $, il faut considerer ce resultat. Pour $ n=2$, une surface est loc.conf.plate, donc, on a le resultat d'un disque de $ {\mathbb R}^2 $ et potentiel $ V $ lipschitzien, on a alors l'in\'egalit\'e optimale $ \sup + \inf $. Dans le cas de la dimension $ n=3 $, on a $ (\sup u)^{1/3} \times \inf u <c $ avec  potentiel $ V $ lipschitzien, sur toute vari\'et\'e riemannienne de dimension $ n=3 $. Pour la dimension $ n=4 $, on $\sup $ major\'e si $ \inf $ est minor\'e, et $ \sup =f(\inf)$ de maniere implicite, pour le suites de solutions, avec potentiel $ V $ lipschitzien, sur toute vari\'et\'e riemannienne de dimension $ n=4 $. (Potentiel $ V $ lipschitzien: implique que le type d'in\'egalit\'e devient de plus en plus faible avec la dimension: sur toute vari\'et\'e riemannienne).

\smallskip

Le resultat de C.C.Chen-C.S.Lin, de 1997. Comm.pure.appl.math. (pour $ n=4 $ par exemple) est fait d'abord sur des boules de $ {\mathbb R}^4 $. Puis il est appliqu\'e aux vari\'et\'es localement conf.plates. Il n'est pas prouv\'e directement sur les vari\'et\'es localement.conf.plates (vari\'et\'es riemanniennes particulieres). Mais sur des boules de $ {\mathbb R}^4 $. D'abord sur un espace particulier: boules de $ {\mathbb R}^4 $, puis appliqu\'e aux vari\'et\'es localement conf.plates (vari\'et\'es riemanniennes particulieres). l'espace de reference est: boules de $ {\mathbb R}^4 $: in\'egalit\'e optimale $ \sup \times \inf $:

\smallskip

L'espace de r\'ef\'erence est : les boules de $ {\mathbb R}^4 $, pour s'y ramener: un changement de metrique conforme: outil de la g\'eom\'etrie conforme.

\smallskip

Le r\'esultat du pr\'eprint: "Inequality on manifolds of dimension 4 with $C^2$ potentials", en dimension $ n=4 $: est fait sur les espaces localement symetriques, directement. La vari\'et\'e peut etre localement conf.plate ou non localement conf. plate. L'espace de r\'ef\'erence: espaces localement symetriques (vari\'et\'es riemanniennes particulieres). Puis, appliqu\'e aux espaces localement symetriques, localement conf.plats ou non localement.conf. plats: in\'egalit\'e $ (\sup)^{1/3} \times \inf $:

\smallskip

Le changement de m\'etrique conforme n'est pas utile ici.(pour le cas non localement.conf.plat).

\smallskip

L'espace de r\'ef\'erence est: les espaces localement sym\'etriques: ne provient pas de la g\'eom\'etrie conforme.

\smallskip

////////////

\smallskip

1) Par exemple, avec potentiel $ V\equiv 1 $ (equation de Yamabe, $ n=4 $), on a: $ (\sup)^{1/3} \times \inf $, sur toute vari\'et\'e riemannienne.

\smallskip

2) Par exemple, avec potentiel $ V\in C^2 $ (equation de la courbure scalaire prescrite, $ n=4 $): $ (\sup)^{1/3} \times \inf $: espaces localement sym\'etriques.

\smallskip

3) Par exemple, avec potentiel $ V\in C^1 $ (equation de la courbure scalaire prescrite, $ n=4 $): $ \sup =(\inf) $: espaces localement sym\'etriques Ricci plats.

\smallskip

4) Par exemple, avec potentiel $ (V_i) $ lipschitziens, de constantes de lipschitz tendant vers 0, (platitude): equation de la courbure scalaire prescrite, $ n=4 $: $ \sup $ major\'e si $ \inf $ minor\'e et $ \sup=(\inf) $ pour des suites: toute vari\'et\'e riemannienne.(ici, on a utilis\'e la platitude).

\smallskip

////////////////////////////

\smallskip

{\bf Sur les exemples et contre-exemples:}

\smallskip

(Explication pour la dimension $ (n=4) $. On a:

\smallskip

preprint de 2024: "Harnack inequalities for equations of type prescribed scalar curvature": dans la platitude: $ V $ tels que: $ 0 <a \leq V \leq b <+\infty, \,\, |\nabla V|\leq A < +\infty $:

\smallskip

le cas $ (n=4) $ de "Harnack inequalities for equations of type...": $ (P):$ Pour toute vari\'et\'e Riemannienne $ (M,g) $: il y a une relation entre $ \sup $ et $ \inf $: $ \Rightarrow $: $ (Q):$ cela revient \`a considerer, par blow-up, des solutions relativement \`a des potentiels $ \tilde  V_i $ d'amplitude tendant vers $ 0, |\nabla \tilde V_i| \to 0 $.  Donc:

\smallskip

$ non (Q): $ si on considere les solutions relativement \`a des potentiels d'amplitude ne tendant pas vers $ 0 $, $ V > 0, \nabla V \not =0 $ (pas de rescaling ou changement d'echelle, avec $ \lambda \to +\infty $, ou pas de blow-up): $\Rightarrow $: $ non (P) $: $ \Rightarrow $: il existe une vari\'et\'e Riemannienne $ (M_0,g_0) $ pour la quelle le $ \sup $ et l'$ \inf $ ne sont pas li\'es: $ \Rightarrow $: pour cette vari\'et\'e Riemannienne $ (M_0,g_0) $, $ \sup $ major\'e si $\inf $ minor\'e, n'est pas vrai. Comme $ \sup $ major\'e si $ \inf $ minor\'e est vrai dans le cas plat, alors, la vari\'et\'e Riemannienne $ (M_0,g_0) $ est non.loc.conf.plate).

\smallskip

Le preprint de 2024, cit\'e ci-dessus, permet d'avoir des contre -exemples, en dimension $ n\geq 4 $, en particulier, $ n=4 $, comme on vient de le voir, et $ n\geq 5 $.(potentiel $ V>0, \nabla V\not =0 $, pas de rescaling ou changement d'echelle, avec $ \lambda \to +\infty $ ou pas de blow-up).

\smallskip

(Explication: on a un ensemble de fonctions: $ E=\{V\in C^1, V>0, \nabla V \not =0, \,\, ou, \,\, 0 < c_2\leq |\nabla V|\leq c_1<+\infty \}$, et un groupe de transformations: $ G=\{g: x\to x_0+x/\lambda \} $, agissant sur $ E $, alors on doit avoir $ g\cdot E \subset E $, alors, quand $\lambda \to +\infty $, il n'est plus possible que $ G $ agisse sur $ E $: $ g\cdot E\not \subset E $, homogen\'eit\'e de l'ensemble $ E $ ou invariance de $ E $. Donc, il n'est plus possible de faire le rescaling ou le changement d'echelle pour $ \lambda \to +\infty $).

\smallskip

////////////////////

\smallskip

On peut fixer chaque vari\'et\'e $ (M,g) $ et avoir des contre-exemples pour cette vari\'et\'e Riemannienne $ (M,g) $, pour $ n\geq 4 $,  avec $ E= \{ V\in C^1, V>0, \nabla V \not =0, \,\, ou, \,\, 0 < c_2\leq |\nabla V|\leq c_1<+\infty \} $: impossibilit\'e de faire le blow-up. 

\smallskip

Ici, comme, pour toute vari\'et\'e Riemannienne $ (M,g), n\geq 4 $, on a des contre-exemples,(on a: $ \sup $ et $ \inf $ ne sont pas li\'es, donc, on n'a pas: $ \sup $ major\'e si $ \inf $ est minor\'e), il faut determiner les vari\'et\'es sur les quelles, par exemple, en dimension $ n=4 $, $ \sup $ major\'e si $ \inf $ est minor\'e est vrai, c'est l'objet du preprint des potentiels $ C^1 $, espaces localement symetriques Ricci plat de dimension 4.(Pour ces espaces aussi, on a des contre-exemples: pour $ V >0, V\in E'=\{ V\in C^1,  V >0, 0 < c_2 \leq |\nabla V|\leq c_1 < +\infty \} $: impossibilit\'e de faire le blow-up, ici aussi).

\smallskip

Cas de la dimension $ n=4 $:

\smallskip

Dans le livre d'Aubin, pour $ n=4 $ et la courbure scalaire $ R=R_g >0 $, il y a des solutions au probleme de la courbure scalaire prescrite, c'est \`a dire en considerant un potentiel $ V >0, V \in C^1 $, avec des conditions sur $ V>0 $. On a, si considere $  E=\{R>0, 0 <c_2 \leq |\nabla R| \leq c_1 <+\infty \} $, on a: $ R=R_g \in E $, en utilisant le preprint de 2024: "Harnack inequalities for equations of type prescribed scalar curvature", on raisonne sur la courbure scalaire, on a des contre exemples: impossibilit\'e d'utiliser le blow-up: $ \Rightarrow \sup $ et $ \inf $ ne sont pas li\'es $ \Rightarrow \sup $ major\'e si $ \inf $ minor\'e n'est pas vrai. Ici, la condition $ 0 < c_2 \leq |\nabla R|\leq c_1 <+\infty $ implique qu'on est sur un espace non localement symetrique, car: $ \nabla R= 0 \Rightarrow \nabla Ricci =0 \Rightarrow \nabla Scalaire = 0 $.

\smallskip

Donc, on a des contre exemples dans les espaces non localement symetriques, d'o\`u, le preprint des potentiels $ C^1 $ sur les espaces localement symetriques Ricci plats.

\smallskip

////////////////////

\smallskip

On peut faire la meme chose pour l'eq. de Yamabe, sur une vari\'et\'e Riemannienne compacte sans bord, avec le preprint de 2024: "Harnack inequalities for equations of type prescribed scalar curvature", on raisonne sur la courbure scalaire, en considerant l'ensemble: $ E=\{R>0, 0 <c_2 \leq |\nabla R| \leq c_1 <+\infty \} $, avec, $ R=R_g \in E $: $ n\geq 4 $ (preprint 2024): impossiblit\'e de faire le blow-up: $ \Rightarrow $ $ \sup $ et $ \inf $ ne sont pas li\'es $ \Rightarrow $ $\sup $ major\'e si $ \inf $ minor\'e, n'est pas vrai $ \Rightarrow \sup \times \inf \to +\infty $. Comme sur une vari\'et\'e compacte sans bord, l'eq. de Yamabe a toujours des solutions, en particulier dans le cas positif, en particulier si $ R >0 $, on a alors des contre-exemples.

\smallskip

Pour avoir des contre exemples, on a impos\'e la condition supplementaire $ 0 <c_2 \leq |\nabla R|\leq c_1<+\infty $, il se peut que sans cette condition on n'ait pas ces contre-exemples.

\smallskip

C'est une r\'eponse \`a la question de YY.Li-L.Zhang.

\smallskip

/////////////////

\smallskip

On peut appliquer cela aussi \`a l'article:"About Brezis-Merle problem with Holderian condition", en considerant l'ensemble $ E= \{ V\in C^1(\bar B_1(0)), b\geq V\geq 0, \nabla V \not =0, \,\, ou, \,\, 0 < c_2\leq |\nabla V|\leq c_1<+\infty \} $ et $ \int_{B_1(0)} e^u dx \leq C $ avec $ bC <16\pi $ ou $ bC<24\pi$, on a des contre-exemples en dimension 2: impossibilit\'e de faire le blow-up. 

\smallskip

////////////////

\smallskip

On peut appliquer cela aussi \`a l'article: "Harnack inequalities for Yamabe type equations", en considerant l'ensemble, $ R=R_g $ est la courbure scalaire: $ E=\{R >0, 0 < c_2\leq |\nabla R|\leq c_1 <+\infty\} $: impossibilit\'e de faire le blow-up.

\smallskip

D'apres un Theoreme d'Aubin, on a existence de solutions quand: $ \epsilon < \frac{(n-2)R}{4(n-1)} =\frac{(n-2)R_g}{4(n-1)} $. On a aussi, les solutions constantes. Pour l'eq: $ -\Delta u_{\epsilon}+\epsilon u_{\epsilon} = u_{\epsilon}^{N-1}, N=\frac{2n}{n-2}, n\geq 3,  $ sur $ (M,g), n\geq 4 $, $ \Delta =\nabla^j \nabla_j $. Rien ne prouve que les solutions constantes et les solutions d'Aubin coincident. (On a consid\'er\'e l'ensemble $ E $).

\smallskip

Comme on a consid\'er\'e l'ensemble $ E $: on a alors, l'existence d'une sous-suite, quand $ \epsilon_i \to 0 $, de solutions non constantes de $ -\Delta u_{\epsilon_i}+\epsilon_i u_{\epsilon_i} = u_{\epsilon_i}^{N-1}, N=\frac{2n}{n-2}, n\geq 3,  $ sur $ (M,g), n\geq 4 $, compactes sans bord, relatives \`a l'ensemble $ E $, $ \Delta =\nabla^j \nabla_j $, car, on a : $ \sup_M u_{\epsilon_i} \times \inf_M u_{\epsilon_i} \to +\infty $. 

\smallskip

Donc, $ R=R_g >0 $ sur $ (M,g) $ compacte sans bord de dimension $ n\geq 4 $ est une condition n\'ecessaire et suffisante (si, seulement si) pour obtenir le r\'esultat d'unicit\'e. C'est une r\'eponse au probleme 2 de Brezis-Li.

\smallskip

//////////

\smallskip

On peut appliquer cela, au cas n\'egatif, en considerant l'ensemble: 

$ E=\{V\in C^1, V<0, \nabla V \not =0, \,\, ou, \,\, 0 < c_2\leq |\nabla V|\leq c_1<+\infty \}$: impossibilit\'e de faire le blow-up. Ici, on a la notion de potentiel hyperbolique ou de noyau hyperbolique, sans qu'on ait explicitement ce potentiel ou ce noyau.

\smallskip

(On a: $ \exists (u_i), \,\, u_i >0, \,\, -\Delta u_i+Ru_i=V_i {u_i}^{N-1}, \,\, V_i < 0, \,\, \exists K_0 \subset \subset M, \,\, \sup_{K_0} u_i \to +\infty $. On a la notion suivante: $ u_0>0, \,\,\Delta u_0=\delta_{x_0}, \Delta=\nabla^i(\nabla_i) $ avec $ u_0 >0 $ un objet math\'emtique et de la physique non explicit\'e).

\smallskip

//////////////////////////////////

\smallskip

Concernant le contre-exemple du probleme de Brezis-Merle et "About Brezis-Merle problem with Holderian condition":(on en a parl\'e dans le preprint "Cas d'existence de solutions d'edp", concernant l'existence de solutions au probleme variationnel en dimension 2, equation du type Liouville, avec potentiel)

\smallskip

On a: on prend $ \Omega=B_1(0) $ et $ 0 < a \leq V_{\epsilon} \leq b <+\infty $ avec $ b|\Omega|=b\pi< 1 $.

\smallskip

Soit $ \mu_{\epsilon} \geq 0 $:

$$ \mu_{\epsilon} =\inf \{||\nabla u||^2, u \in H_0^1(\Omega), \,\, \int_{\Omega} V_{\epsilon} e^u =1 \}. $$

En utilisant la condition $ b|\Omega|<1 $ et l'injection de Moser-Trudinger, on obtient (par l'absurde):

$$ \mu_{\epsilon} \geq \mu_0 >0, \,\, \forall \epsilon >0, $$

En utilisant les multiplicateur d'Euler-Lagrange, on obtient:

$$ \exists \lambda_{\epsilon}, \, \exists u_{\epsilon} \in H_0^1(\Omega), \,\, -\Delta u_{\epsilon}= \lambda_{\epsilon} V_{\epsilon} e^{u_{\epsilon}}, $$

avec condition de Dirichlet au bord.

\smallskip

On utilise la condition $ b|\Omega|<1 $ et  le principe du maximum, on obtient $ \lambda_{\epsilon} >0 $ et $ u_{\epsilon} >0 $.(par l'absurde).

\smallskip

En utilisant le Theoreme 1 de Brezis-Merle, on a : $ \lambda_{\epsilon} \not \to 0 $.(par l'absurde on aurait $ \mu_{\epsilon} \to 0 $, ce n'est pas possible). Donc: $ \lambda_{\epsilon} \geq  \lambda_0 > 0 $. En utilisant la premiere valeur propre du laplacien avec condition de Dirichlet (pour le disque unit\'e $ \lambda_1 \leq 6  $, on prend une fonction test dans le probleme variationnel $ u=1-r $.) et une integration par parties, on obtient:

$$ \lambda_{\epsilon} \times a \leq \lambda_ 1 \leq 6, $$

On pose $ W_{\epsilon} =\lambda_{\epsilon} V_{\epsilon} $:

\smallskip

On prend $ b=1/2\pi $ et $ a=2b/3 $ et $ 0 <c_2 \leq |\nabla V_{\epsilon}|\leq c_1 <+\infty $.(On a alors: $ \int_{\Omega} W_{\epsilon} e^{u_{\epsilon}} dx = \lambda_{\epsilon} \leq 18\pi <24\pi $).(on peut prendre $ b=1/2\pi $ et $ a=4b/5 $ pour avoir $ \int_{\Omega} W_{\epsilon} e^{u_{\epsilon}} dx = \lambda_{\epsilon} \leq 15\pi <16\pi $) 

\smallskip

On obtient: 

\smallskip

$ \int_{\Omega} W_{\epsilon} e^{u_{\epsilon}} dx = \lambda_{\epsilon} \leq 15\pi <16\pi $ ou

\smallskip

$ \int_{\Omega} W_{\epsilon} e^{u_{\epsilon}} dx = \lambda_{\epsilon} \leq 18 \pi <24\pi $ 

\smallskip

et: $ 0 <\tilde c_2 \leq  |\nabla W_{\epsilon}|\leq \tilde c_1 <+\infty $.

\smallskip

On peut alors appliquer ce qu'on a dit sur l'impossibilit\'e de faire le blow-up en utilisant l'article "About Brezis-Merle problem with holderian condition". Donc, pour une sous-suite: 

$$ \sup_{\Omega} u_{\epsilon} \to +\infty, $$

D'apr\'es Brezis-Merle, on a la compacit\'e locale, d'o\`u le blow-up au bord. Ceci est un exemple et un contre-exemple aux problemes 1 et 2 de Brezis-Merle. C'est une r\'eponse aux problemes 1 et 2 de Brezis-Merle.

\smallskip

////////////////////////////////////////////////////////////

\smallskip

{\bf Remarque:} On veut voir si, quand on a la condition au bord de Dirichlet, on peut quand meme, utiliser Brezis-Li-Shafrir: conernant l'impossibilit\'e d'utiliser le blow-up: est ce qu'on peut utiliser Brezis-Li-Shafrir, en ayant une condtion supplementaire, qu'on n'utilise pas ?

\smallskip

En prenant l'exemple de $ u_{\epsilon}(r)=\epsilon (1-r^2), \epsilon >0 $, dans la couronne $ B_1(0)-B_{\epsilon'}(0), \epsilon'>0 $, on a  $ V_{\epsilon}(r)= 4\epsilon e^{(r^2-1)\epsilon} $, avec  $ 0 < c_2 \leq |\nabla V_{\epsilon}(r)|\leq c_1<+\infty $, avec $ u_{\epsilon} $, une suite born\'ee ou compacte :

\smallskip

On ne peut pas appliquer Brezis-Li-Shafrir de l'impossibilit\'e d'utiliser le blow-up. Ici, on a de plus, la condition au bord $ u_{\epsilon}(1)=0 $. En particulier: quand on a une condition supplementaire qu'on n'a pas utilis\'e, comme la condition au bord, on ne peut pas appliquer Brezis-Li-Shafrir. 

\smallskip

Aussi, la contrapos\'ee de l'assertion $ (P): \forall (u_k)_k, ... $: est dans le sens: il existe une suite telle que  $ non (P) $ est vraie: $ \exists (u_k)_k,... $, et non, pour toute suite, il y a divergence.

\smallskip

////////) Plus pr\'ecis\'ement, on consid\`ere: Brezis-Li-Shafrir et le probl\`eme de Dirichlet: on explique mieux tout ceci:

\smallskip

 Soit $ E=\{ -\Delta u= Ve^u\} $ et $ (P): \sup +\inf \leq c $: alors considerer $ (E)$ : on a $ (P) $ et $ non (P) $.

\smallskip

 Soit $ \tilde E $: le probleme aux limites de Dirichlet: $ \tilde E=\{-\Delta \tilde u= \tilde V e^{\tilde u}, \,\, \tilde u = 0 \, {\rm \, sur \, le \, bord}\}$:

\smallskip

 alors on ne peut pas definir la propri\'et\'e $ (P'): \tilde u: \sup_K \tilde u \leq c $: la compacit\'e locale, \`a partir de Brezis-Li-Shafrir, car:

\smallskip
 
si on pose le probleme $ (P') $: rien ne dit qu'on ne tombe pas sur  une solution telle que: $ \tilde u =(u_i) $ verifiant $ non (P) : \sup_K u_i \to +\infty $.

\smallskip

Considerer l'equation, inclut le fait de considerer le probl\`eme de Dirichlet. Si on suppose le domaine regulier, alors considerer l'equation de Liouville, inclut le fait de considerer le probl\`eme de Dirichlet. 

\smallskip

Donc, on ne peut pas definir la propri\'et\'e $ (P') $: compacit\'e locale, du probleme aux limites $ \tilde E $, qu'on veut prouver.

\smallskip

Car considerer $ (P') $ \`a partir de $ E $, c'est avoir $ (P) $ et $ non (P)$.

\smallskip

Il faut definir le probleme aux limites $ \tilde E $ et prouver la compacit\'e locale, $ (P') $, en dehors de Brezis-Li-Shafrir.(avec blow-up).

\smallskip

Par exemple, Brezis-Merle, prouvent la compacit\'e locale, sans Brezis-Li-Shafrir, mais ici, c'est sans blow-up.

\smallskip

Ici, la propri\'et\'e $ (P')$: prend toutes les possibilit\'es. Si, on l'a lie \`a Brezis-Li-Shafrir, c'est prendre en compte $ (P) $ et $ non (P) $. Ceci, pour dire que parfois on utilise Brezis-Li-Shafrir, dans des preuves, il est possible de le faire, car, les suites qu'on considere, ne sont pas dans un enemble plus general, qui prend en compte $ (P) $ et $ non (P) $. On utilise Brezis-Li-Shafrir, sans definir un ensemble o\`u il y a toutes les possibilit\'es.

\smallskip

Par contre, quand on veut considerer, le probleme au limites $ \tilde E $, on doit definir la propri\'et\'e $ (P') $, de maniere generale, c'est a dire, prendre en compte toutes les possibilit\'es. Si de plus, on la lie \`a Brezis-Li-Shafrir, cela revient \`a considrer $ (P) $ et $ non (P) $. Ce qui veut dire qu'on peut tomber, d\`es le depart, sur des solutions blow-up, sans avoir pu definir la propri\'et\'e $ (P') $: donc, on ne peut pas definir $(P')$ et $ non (P') $, donc, on ne pas les prouver ou utiliser Brezis-Li-Shafrir pour les prouver.

\smallskip

//////////////) Ceci est un argument qui confirme ce qu'on a dit auparavant, sur le fait que le $ (\sup,\inf) $ n'a rien a voir avec la compacit\'e globale. 

\smallskip

//////////////////////////////////////////////////////

\smallskip

Algorithme reconnaissant le oui et le non avec la notion d'in\'egalit\'e de Harnack pour l'equation de la courbure scalaire prescrite:

Maintenant on regarde l'algorithme suivant: de la courbure scalaire prescrite: l'equation de la courbure scalaire prescrite commence \`a partir de $ n\geq 2 $: Ca depend de la dimension. Considerer l'equation de la courbure scalaire prescrite c'est considerer les dimensions $ n\geq 2 $: avec $ n=2 $ et $ n\geq 3 $:

\smallskip

Pour $ n=2 $, on a le cas plat uniquement, equation avec non-lin\'earit\'e exponentielle. 

Pour le cas $ n\geq 3 $, il y a le cas plat et non plat.(Parfois il n'est pas n\'ecessaire de distinguer ces 2 cas, plat et non plat. De meme, parfois il n'est pas necessaire de supposer l'operateur coercif ou la vari\'et\'e \`a bord, bord regulier).

//////////////////////////////

On a:

$$ Oui=\{ V \in C^1, \,\, 0 < a \leq V(x) \leq b <+\infty, \,\, |\nabla V(x)|\leq A \}: $$

$$  (P): \,\,\forall (u_k)_k, \forall x_0 \in M, \, \exists \, (u_{i_j})_j, \, telle \, que: \, (\sup u_{i_j},\inf u_{i_j}): \, en \, relation, \, autour \, de \, x_0, $$

/////// et,

$$ Non=\{V\in C^1,\, 0 < \tilde a \leq V(x)\leq \tilde b, \, 0 < c_2 \leq |\nabla V(x)|\leq c_1 \}: $$

$$ non (P): \, \exists \, (u_k)_k,\, \exists\, x_0 \in M, \,\, tels \,\, que: \, \forall (u_{i_j})_j, \,\, \,(\sup u_{i_j},\inf u_{i_j}): \,\, pas \,\, en \,\, relation. $$

////////////

On a:

$$ Si \,\, (n=2),\,\, Brezis-Li-Shafrir: cas \,\, plat \,\, uniquement \,\, car \,\, n=2: J.F.A.1993.$$

$$ Sinon, \,\, (n\geq 3), \, Cas \,\, plat, \,\, C.C.Chen-C.S.Lin: Comm.Pure.Appl.Math.1997. $$

$$ Sinon, \,\, (n=3),\, Cas \,\, coercif, \, Cas\, non \, plat,\,L.Zhang, \,\, Trans.Amer.Math.Soc. \,2009. $$

$$ Sinon, \,\, (n=3), \, Cas \,\, non \,\, plat, \,\, Mathematica \,\, Aeterna: 2011. $$

$$ Sinon, \,\, (n=4),\, Cas \,\, non \,\, plat, \,\, Harnack \,\,inequalities \,\, for \,\, equations \,\, of \,\, type ... \,\, Arxiv: 2024.$$

$$ Fin: \,\, (n\geq 5), \, Cas \,\, non \,\, plat, \,\, Harnack \,\ inequalities \,\, for \,\, equations \,\, of \,\, type....\,\, Arxiv: 2024.$$

\smallskip

//////) Cet algorithme reconnait le oui et le non  avec la notion d'in\'egalit\'e de Harnack.

\smallskip

////////////////

\smallskip

//////) Dans le cas n\'egatif, on a aussi, un algorithme qui reconnait le oui et le non avec la notion de $ \sup $ local major\'e.

\smallskip

Un algorithme est une succession de manipulation de nombres positifs par un ordinateur, on va supposer alors qu'on est dans le cas positif, pour l'ex\'ecution d'un algorithme. Et en sortie, temps d'ex\'ecution( $ >0 $). Dans ce cas, le premier des 2 algorithmes qu'on vient de pr\'esenter est l'algorithme recherch\'e: algorithme qui permet de reconnaitre le oui et le non avec la notion d'in\'egalit\'e de Harnack.

\smallskip

Chaque dimension et proposition, repr\'esente une porposition du calcul propositionnel. L'ensemble des propositions du calcul propositionnel est d\'enombrable. Il y a equipotence entre le nombre de dimensions et des propositions et le nombre de propositions du calcul propositionnel. C'est aussi le probl\`eme de satisfaisabilit\'e bool\'eenne: le probl\`eme SAT. Donc, c'est un algorithme du probl\`eme SAT.

\smallskip

Dans cet algorithme, il n'y a pas d'arborescence, ou de combinatoire. C'est un algorithme \`a l'echelle humaine et non machine. Il faut compter le nombre d'instructions. 

C'est un algorithme math\'ematique, voir sur le site internet: Geeks for Geeks: si $ N $ est le nombre d'entr\'ees, variables, hypoth\`eses, un algorithme mathematique n\'ecessite $ O(\log N)=C \cdot \log N $, en temps, ou, sa taille, en complexit\'e en temps est de $ C \cdot \log N = O(\log N) $. 

\smallskip

Dans l'algorithme qu'on pr\'esente, il y a une succession d'algorithmes rudimentaires, en entr\'ees: $ N_1, N_2,\ldots, N_k $, donc, il est en: $ k \cdot \log N_k = k\cdot C \cdot \log N = O(\log N) $.

\smallskip

Donc, c'est probable que ce soit un algorithme en temps polynomial, et d'apr\'es les notes de Richard Lassaigne (logique propositionnelle): ceci donne une r\'eponse positive au probl\`eme: $ p=np $. 

\smallskip

Ce qui veut dire que: $ p=np $ est vrai.

\smallskip

////////////////

\smallskip

{\bf Remarque:} on peut considerer l'assertion: $ (\tilde P): (1) ou (2) $, du preprint de 2024, "Harnack inequalities for equations of type prescribed scalar curvature". Donc, $ non (\tilde P): non [(1) ou (2)] $.

\smallskip

Est ce que cela  reconnait le oui ou le non ? 

\smallskip

On a: $ (\tilde P): (1) ou (2) $ reconnait le oui. Pour le non, on fait le raisonnement suivant: soit le non,  supposons que $ (1) ou (2) $ soit vraie, alors, si $ (1) $ est vrai, $ (1) ou (2) $ est vraie, si $ (2) $ est vraie $ \Rightarrow $ on a utilis\'e le blow-up, or ceci est impossible, d'o\` u, contradiction. Donc: $ (1) ou (2) $ vraie $ \Rightarrow $ contradiction: $ \Rightarrow $ $ [(1) ou (2)] $ n'est pas vraie $ \Rightarrow non (\tilde P) $. Donc, ca reconnait le non.

Donc, ca reconnait le oui et le non.

\smallskip

 On a le cas: $ \inf_M u_k=0 $, alors, on a: $ \sup_K u_k \times \inf_M u_k = 0 <1: (1') $. On rajoute une ligne: 

 \smallskip
 
 soit: $ J=\{k \in {\mathbb N}, \inf_M u_k >0 \} $, alors, $ J $ est fini $ \Rightarrow (1') : \sup_K u_k\times \inf_M u_k \leq c<+\infty$, et si $ J $ est infini $ \Rightarrow (1) ou (2) $:

 \smallskip
 
 $ (\tilde P) $: pour des sous suites: $ (\tilde P): (1') ou (1) ou (2) $. ici aussi ca reconnait le oui et le non avec le meme type de raisonnement qu'avant.

\smallskip

////////////////////////

\smallskip

Dans le preprint de 2024: "Harnack inequalities for equations of type prescribed scalar curvature": on le dit: car, on dit qu'on veut montrer que $ \sup_K u_k $ et $ \inf_M u_k $ sont li\'es. Si $\inf_M u_k=0 $, alors c'est trivial, leur produit est nul. donc, on s'interesse au cas o\`u $ \inf_M u_k >0 $. 

$ \inf_M u_k =0 \Rightarrow \sup_K u_k\times \inf_M u_k = 0 <1 $.

\smallskip

//////////////////////

\smallskip

On peut prendre $ (P): (1)+(2) $, du corollaire 1.2. du preprint de 2024: "Harnack inequalities for equations of type prescribed scalar curvature". (C'est ce qu'on a dit au d\'ebut).

Ici, on a bien $ \inf_{B_{r_0}(x_0)} u_{i_j} >0 $, car la boule ferm\'ee, $ \bar B_{r_0}(x_0) \subset M $, $ 0 < r_0 <\frac{inj(x_0)}{2} $, avec $ inj(x_0) $ le rayon d'injectivit\'e en $ x_0 $ et $ u_{i_j} $ est continue sur $ M $.

\smallskip

Dans la preuve du preprint de 2024:"Harnack inequalities of equations of type prescribed scalar curvature equations":

\smallskip

Ici, on prend, $ \delta_0(K)=\delta_0(\bar B_r(x_0))=\inf \{ \inf\{\delta_P/4, P\in K=\bar B_r(x_0)\}, \frac{(r_0-r)}{20}\} $, on a alors, $ K_{\delta_0} \subset \bar B_{r+\frac{3(r_0-r)}{20}}(x_0) \subset B_{r_0}(x_0) $ et $ r_k=\inf \{ [u_k(\cdot)]^{-\epsilon}, \inf \{ \delta_P/4, P\in K_{\delta_0}\}, (\frac{(r_0-r)}{20})^{(n-2)/2}\}$, on a alors, $ B_{2r_k^{2/(n-2)}}(t_k) \subset B_{r_0}(x_0)$ $ \Rightarrow \inf_{B_2(0)} v_k=r_k \inf_{B_{2 r_k^{2/(n-2)}}(t_k)} u_k \geq r_k\inf_{B_{r_0}(x_0)} u_k \geq r_k \inf_{\bar B_{r_0}(x_0)} u_k >0 $.

\smallskip

//////////////////////////////////////////////////

\smallskip

Explication de l'impossibilit\'e d'utiliser le blow-up:

\smallskip

a) raisonnement par l'absurde:

\smallskip

$ (P) $ est vraie $ \Leftrightarrow $ $ [ Non (P) \Rightarrow $ contradiction $ ] $: est vraie, $ \Rightarrow $: on a utilis\'e le blow-up: or ceci est impossible $ \Rightarrow $ $ [ Non (P) \Rightarrow $ contradiction $ ] $: n'est pas vraie, $ \Rightarrow (P) $ n'est pas vraie.

\smallskip

b) raisonnement direct:

\smallskip
$ (P) $ est vraie $ \Rightarrow $: on a utilis\'e le blow-up: or ceci est impossible $ \Rightarrow (P) $ n'est pas vraie.

\smallskip

////////////////////////////////

\smallskip

{\bf Remarque sur le rescaling ou changement d'echelle:} On considere la suite: $ u_{\epsilon}(x)=(\frac{\epsilon}{\epsilon^2+|x|^2})^{(n-2)/2} $, on veut voir si en faisant le rescaling, on n'obtient pas n'importe quoi, ou si on obtient une vraie in\'egalit\'e de Harnack: $ u_{\epsilon}(0)=\epsilon^{-(n-2)/2} $, donc, $ [u_{\epsilon}(0)]^{\alpha}=\epsilon^{-\alpha (n-2)/2}, 0 < \alpha < 1 $ et $ v_{\epsilon}(y)=\frac{u_{\epsilon}(u_{\epsilon}(0)^{-2\alpha /(n-2)} y)}{[u_{\epsilon}(0)]^{\alpha}} = \epsilon^{\alpha (n-2)/2} u_{\epsilon} (\epsilon^{\alpha} y)=\epsilon^{\alpha(n-2)/2} (\frac{\epsilon}{\epsilon^2+\epsilon^{2\alpha} |y|^2})^{(n-2)/2} $. Alors: $ v_{\epsilon}(0) = \epsilon^{-(1-\alpha)(n-2)/2} \to +\infty $ et  $ v_{\epsilon}(1)= \epsilon^{(1-\alpha)(n-2)/2}(\frac{1}{\epsilon^{2(1-\alpha)} +1})^{(n-2)/2} \to 0 $, quand $ \epsilon \to 0 $. 

\smallskip

Donc: $ v_{\epsilon}(0)\times v_{\epsilon}(1) \to 1 $, si $ \epsilon \to 0 $. On obtient bien une in\'egalit\'e de Harnack $ \sup \times \inf $. On n'obtient pas n'importe quoi, mais, une in\'egalit\'e de Harnack $ \sup\times \inf $.

\smallskip

Maintenant, on ecrit: $ [u_{\epsilon}(0)]^{1-\alpha} = v_{\epsilon}(0) \leq \frac{c}{v_{\epsilon}(1)}\leq \frac{c}{\frac{u_{\epsilon}(1)}{(u_{\epsilon}(0))^{\alpha}}}= c\times \epsilon^{-(1+\alpha)(n-2)/2}= \frac{c}{[u_{\epsilon}(1)]^{1+\alpha}} $. Donc, \`a partir de l'in\'egalit\'e $ v_{\epsilon}(0)\times v_{\epsilon}(1) \leq c $, on obtient l'in\'egalit\'e: $ [u_{\epsilon}(0)]^{1-\alpha}\leq \frac{c}{[u_{\epsilon}(1)]^{1+\alpha}} $, c'est une in\'egalit\'e de Harnack plus faible que celle qu'on avait au d\'epart: $ u_{\epsilon}(0)\times u_{\epsilon}(1) \leq c $, mais c'est quand meme une in\'egalit\'e de Haranck.

\smallskip

Donc, le rescaling, ou changement d'echelle, permet d'avoir une vraie in\'egalit\'e de Harnack pour les nouvelles fonctions $ v_{\epsilon} $. Et l'in\'egalit\'e de Harnack des nouvelles fonctions $ v_{\epsilon} $, permet d'avoir une in\'egalit\'e de Harnack pour les fonctions de d\'epart $ u_{\epsilon} $, plus faible que celle qui est connue, mais c'est quand meme une in\'egalit\'e de Harnack.

\smallskip

Tout ceci, c'est pour dire, que ce proc\'ed\'e de rescaling, ou changement d'echelle, permet, ici, d'obtenir quand meme, une in\'egalit\'e de Harnack. On n'obtient pas n'importe quoi.

\smallskip

Plus precis\'ement, si on a, par exemple, au d\'epart, une in\'egalit\'e: $ u_i(0) \times (\inf_M u_i)^k \approx 1 $, cela veut dire qu'on a: $ u_i(0) \approx 1/m, m >0 $ et $ \inf_M u_i \approx m^{1/k}, m >0 $, alors on ecrit: ($ c $ est la fonction du rescaling, ou changement d'echelle): $ c(\inf_M u_i/u_i(0)^{\alpha }) \approx c(m^{\alpha + 1/k}): \Rightarrow c(\inf_M u_i/u_i(0)^{\alpha}) \approx c((\inf_M u_i)^{k\alpha+1}) $, donc: $ u_i(0)^{1-\alpha} \approx c((\inf_M u_i)^{k\alpha+1}) $. Si le rescaling, ou changement d'echelle, conserve le type de fonction $ c(t) \approx t^{-\gamma}, \gamma >0 $ (invariance de l'equation), alors: $ u_i(0)^{1-\alpha} \approx (\inf_M u_i)^{-\gamma (k\alpha+1)} $: on obtient le meme type d'in\'egalit\'e qu'au d\'epart. 

\smallskip

Tout ceci pour dire qu'il y a conservation du type d'in\'egalit\'e, apr\'es avoir utilis\'e le proc\'ed\'e de rescaling, ou changement d'echelle. Et le rescaling, ou changement d'echelle, permet d'avoir une vraie in\'egalit\'e de Harnack, tout \`a la fin.(pas seulement pour les fonctions nouvelles, du rescaling, $ v_i $, mais pour les fonctions $ u_i $, apr\'es avoir effectu\'e ce proc\'ed\'e (le rescaling, ou changement d'echelle)).

\smallskip

Si au d\'epart on avait une relation explicite ou implicite entre sup et inf, alors, le rescaling ou changement d'echelle, conserve, ce type de relations entre le sup et l'inf.

\smallskip

Si on utilise le fait qu'il y a une relation entre sup et inf, apres rescaling, cela veut dire qu'on avait au depart, une relation entre sup et inf, pour les fonctions de depart (invariance de l'equation par rescaling), et plutot raisonner sur la dependance de $ u_i(0) $ en fonction de $ m >0, m = \inf_M u_i $.

\smallskip

/////////////////) Il ne faut pas faire le raisonnement: "Si on a $ u_i(0)^{1-\alpha} \leq c [u_i(0)^{\alpha} / (\inf_M u_i)]^k $, $ \alpha >1-\alpha, \, k>1 $" et obtenir une in\'egalit\'e evidente, qui ne donne rien sur le sup et inf: apres simplification par $ u_i(0) $:  $ 1\leq c/(\inf_M u_i) $, mais considerer $ u_i(0) $, comme une fonction decroissante, \underbar{explicite} ou \underbar{implicite}, de $ m >0 $ avec $ m=\inf_M u_i $ et dans cas on ecrit plutot: $ u_i(0)^{1-\alpha} \leq c m^{-\gamma_1}, \gamma_1 \geq 1 $, avec $ m= \inf_M u_i $ et \`a la fin, on obtient, une relation entre $ \sup $ et $ \inf $: $ u_i(0)^{1-\alpha} \leq c (\inf_M u_i)^{-\gamma_1}, \gamma_1 \geq 1 $. Ce raisonnement reste valable, si $ u_i(0) $ est fonction decoirssante \underbar{implicite} de $ m >0 $ avec $ m= \inf_M u_i $.

\smallskip

Pourquoi on fait ce raisonnement: parce qu'on utilis\'e $ [1/u_i(0)^{\alpha}] $ avec le inf pour controler le sup donc, la quantit\'e: $ [1/u_i(0)^{\alpha}] $, doit etre li\'ee \`a l'inf.

\smallskip

Ceci pour dire qu'on a bien une in\'egalit\'e de Harnack entre sup et inf avec le proced\'e de rescaling ou changment d'echelle, pour les fonctions de depart: $ u_i $. 

\smallskip

//////////////////////////////

\smallskip

On a exhib\'e un exemple explicite, pour lequel le rescaling permet d'avoir une vraie in\'egalit\'e de Harnack pour $ v_i $ et le proc\'ed\'e de rescaling ou changement d'echelle permet d'avoir une in\'egalit\'e de Harnack plus faible que celle, initialement connue pour $ u_i $.

\smallskip

Un autre argument important qui confirme qu'on a bien une in\'egalit\'e de Harnack avec le proc\'ed\'e de rescaling ou changement d'echelle, \`a la fin, est que: l'assertion $ non (P) $ implique que le $ \sup \times \inf \to +\infty $, on l'a dit dans les exemples. Donc, cette in\'egalit\'e de Harnack est bien li\'ee au $ \sup \times \inf $.(car $ \sup \times \inf \leq c <+\infty $ (sur des boules), $ \Rightarrow $ l'in\'egalit\'e de Harnack du preprint de 2024: "Harnack inequalities for equations of type prescribed scalar cruvature"). Donc, au lieu de le voir directement, on utlise la contrapos\'ee: $ non (P) \Rightarrow \sup \times \inf \to +\infty $.

\smallskip

On a les implications suivantes: $ [\sup \times \inf \leq c] \Rightarrow [(\sup)^{\alpha} \times \inf \leq c] \Rightarrow [\sup \leq c'(\inf)] \Rightarrow $: $ (P) $: l'in\'egalit\'e du preprint de 2024: "Harnack inequalities for equations of type prescribed scalar curvature".

avec $ \alpha \in ]0,1[ $ et $ c': m \to c'(m)>0, m >0 $, $ c' $ est une fonction decoiroissante de $ m >0 $.

\smallskip

Donc: $ non (P) \Rightarrow non [\sup \leq c'(\inf) ] \Rightarrow non [(\sup)^{\alpha} \times \inf \leq c] \Rightarrow non [\sup \times \inf \leq c] \Rightarrow \sup \times \inf \to +\infty $.

avec $ \alpha \in ]0,1[ $ et $ c': m \to c'(m)>0, m >0 $, $ c' $ est une fonction decoiroissante de $ m >0 $.

\smallskip

//////////////////

\smallskip

On a la relation: $ (\sup u_i)^{1-\alpha} = c(\inf u_i/(\sup u_i)^{\alpha}) $, si on pose: $ x= \sup u_i, m=\inf u_i $, alors: $ x^{1-\alpha}= c(m/x^{\alpha}), m >0, x >0 $, $ \alpha \in ]0,1[ $.

\smallskip

/////////////////

\smallskip

Maintenant, on s'int\'eresse \`a la fonction $ c >0 $, de cette in\'egalit\'e: si on note $ \tilde u=(u_k)_k $, la suite, alors la constante $ c $ depend de la suite choisie, mais pas de l'indice $ k \in {\mathbb N} $: $ \exists c >0, \,\, \forall k \in {\mathbb N} $. Soit on ne met pas la dependance de $ c $ en $ \tilde u $, soit on rajoute une nouvelle variable $ c = c_{\tilde u}= c(\tilde u, K,...,) >0 $, mais comme on a fix\'e la suite et on a essay\'e de voir s'il existe une constante qui lie le sup et l'inf, independemment de l'indice $ k $, alors, on a enlev\'e la dependence de $ c $ en $ \tilde u=(u_k)_k $.

Pour chaque suite $ \tilde u=(u_k)_k $ et pour chaque compact $ K $, on a determin\'e, une suite $ t=(t_k)_k $, de points telle que: $ \forall k\in {\mathbb N} $, $ u_k(t_k) $ et $ \inf_M u_k $, sont li\'es. On peut faire la meme chose avec $ t=(t_k)_k= t(\tilde u, K) $, $ t $ est fonction de la suite $ (u_k)_k=\tilde u $ et du compact $ K $. Comme on a fix\'e la suite $ \tilde u= (u_k)_k $, alors: il ne reste que la dependance en le compact  $ K $. Ou bien on rajoute une nouvelle variable:

$$ \forall \tilde u= (u_k)_k, \,\, \forall K,\,\, \exists t=(t_k)_k=t(\tilde u, K), \,\,\exists \,\, c= c(\tilde u, t(\tilde u, K), K,...,.)>0,...,\,\, \forall \,\, k \in {\mathbb N}: ... $$

On \'efface la dependance en la variable $ \tilde u=(u_k)_k $, pour ecrire:

$$ \forall \,\, K, \exists \,\, c=c(K,...,.)>0,..., \forall k\in {\mathbb N}: ...,$$

parce qu'on a fix\'e la suite $ (u_k)_k=\tilde u $. On a alors une relation entre le sup et l'inf.

\smallskip

//////////////////////////////////////

\smallskip

Dans le preprint de 2024:"Harnack inequalities for equations of type prescribed scalar curvature": dans le theoreme 1.1. on ne peut pas prendre n'importe quel compact $ K $, car, on a par l'argument de la premiere valeur propre et fonction propre et le lemme de Fatou: $ \int_ B u_k^{N-1} dx \leq c  \Rightarrow  \int_B \liminf_{k\to +\infty} u_k^{N-1} dx \leq \liminf_{k\to +\infty} \int_B u_k^{N-1} dx \leq c \Rightarrow $ preque partout $ y \in B: \liminf_{k\to +\infty } u_k(y) <+\infty \Rightarrow  $ presque.partout $ y \in B: \Rightarrow $ si on prend des singletons, $ K=\{y\} $, alors $ \liminf_{k\to +\infty} \sup_K u_k = \liminf_{k\to +\infty} u_k(y) <+\infty $. Donc pour toute suite $ (u_k)_k $ on a le point 1) seulement, du theorme 1.1 du preprint de 2024, on n'a pas le point 2) du theoreme 1.1 du preprint de 2024. Donc, il faut prendre un compact plus grand, par exemple, comme on l'a dit, des compact d'interieurs non vide, alors on ne peut pas faire ce qu'on vient de dire, avec la valeur propre et le lemme de Fatou.

\smallskip

On regarde le rapport \`a l'algorithmique et l'informatique:

le fait d'avoir le point 1) seulement du theoreme 1.1. du preprint de 2024: "Harnack inequalities for equations of type prescribed scalar curvature":

\smallskip

a) On met en evidence, ce qu'on appelle un "bug", le blocage au point 1) du preprint de 2024: c'est un "bug" math\'ematique: dans l'algorithmique et l'informatique.

\smallskip

b)

Si on veut le blocage au niveau du point 1) du preprint de 2024, alors il faut prendre les singletons comme compacts dans la preuve.

Si on veut eviter le blocage au niveau 1) du preprint de 2024, alors il faut prendre des compacts d'interieurs non vides.

On associe les notions de topologie (singletons comme compacts, compacts d'interieurs non vides) \`a l'algorithmique et l'informatique avec cette notion de "bug".

\smallskip

c) On associe les notions abstraites de la th\'eorie de la mesure et l'int\'egration \`a l'algorithmique et l'informatique.

\smallskip

(Pr\'ec\'edemment: si le principe du maximum n'est pas verifi\'e: il y a bug): Dans: Li-Zhang (dimensions 3,4: op\'erateur conforme coercif): c'est le principe du maximum: comparer le minimum sur le bord et l'interieur.(Si le principe du maximum n'est pas verifi\'e: il y a bug).

\smallskip

////////////////////////////////

\smallskip

{\bf Remarque importante:} pour pouvoir parler d'exemples et contre-exemples, il faut expliciter des solutions. Ceci est l'argument math\'ematique. On a donc, des preuves math\'ematiques. 

\smallskip

Mais du point vue de la physique (th\'eorie du tout) et de l'algorithmique ou de l'informatique, il suffit de prouver les assertions $ (P) $ et $ non (P) $ ou les prouver et les characteriser, pour pouvoir dire qu'il y a solution au probl\`eme. (Par exemple, pour les contre exemples, on a la condition sur les potentiels: $ V >0, 0 < c_2 \leq |\nabla V(x)|\leq c_1<+\infty $: impossibilit\'e de faire le blow-up: ceci permet de dire qu'il y a des exemples du point de vue de la physique ou algorithmique ou informatique, sans expliciter des solutions. Bien sur, il faut faire une demonstration ou une preuve avec le blow-up ou le rescaling ou le changement d'echelle).  

\smallskip

//////////////////////////////////////////////////////////////////////////////////////////////////////////////////////////////////////////

\smallskip

{\bf Sur les relations $ \sup $, $ \inf $, explicites et implicites:}

\smallskip

Relation explicite (il y a plus d'informations: notions d'enroulement+ de torsion). Relation implicite (il y a moins d'informations: notion d'enroulement). (On l'a deja dit dans la partie relativit\'e generale et cosmologie quantique, notion d'enroulement et de torsion et les valeurs $ (\sup, \inf)$).

\smallskip

a) {\bf Relation explicite} entre $\sup $ et $\inf$: $\sup $ et $ \inf $ li\'es. Enroulement+Torsion. Type liaison: $ \sup \updownarrow \inf $, ou $ \sup \leftrightarrow \inf $. Type de torsion en une fonction: $ x\to \frac{1}{x^t}, t>0$.
\smallskip

a1) Relation explicite. $ \sup $ et $\inf$ se controlent mutuellement. Stabilit\'e. Stabilt\'e des vibrations: en $\sup$ et $ \inf$. Premier type de stabilit\'e: stabilit\'e en normes essentielles, entre $ \sup $ et $ \inf $. Deuxieme type de stabilit\'e: en norme $ L^p_{loc}, p>N-1, N=\frac{2n}{n-2} $, stabilit\'e en norme $ W^{1,1+\epsilon}_{loc}$, pour ces normes, il y a un "tunnel" dans lequel, les ondes, leurs vibrations, leurs oscillations, sont stables, avec le fait que le $ \sup $ et $ \inf $ peuvent avoir de grandes amplitudes, mais restent confin\'ees  ou leur $p$-volume reste born\'e, stable, maitrisable: stabilit\'e des vibrations.

a2) {\bf relation ouverte} entre $\sup$ et $\inf $: de $ \sup $ vers $\inf $ et de $ \inf $ vers $ \sup $, $ \sup \leftrightarrow \inf $

a3) $\sup $ controle $ \inf $ et $ \inf $ controle $ \sup $. Intervalles de precisions de nombres pour $ \sup $ et pour $ \inf $.

a4) Bornes uniformes pour les normes, $ L^p_{loc}, p >N-1, N=\frac{2n}{n-2} $ et $ W^{1,1+\epsilon}_{loc} $.

a5) estimations a priori explicites: $\sup \leq c <+\infty$ si $ \inf \geq m >0$ et $ \inf \leq c <+\infty $ si $ \sup \geq M_0 $. On peut quantifier explictement $ c $ en fonction de $ m $ ou $ M_0 $.

\smallskip

b) {\bf Relation implicite} entre $ \sup $ et $ \inf $: $\sup $ et $ \inf $ li\'es. Enroulement. Type liaison: $\sup \uparrow \inf $ ou $ \sup \leftarrow \inf $.

\smallskip

b1) Relation implicite. $ \sup $ est control\'e par $\inf$. Dans un sens seulement. Stabilt\'e dans un sens, du $\inf $ vers $ \sup$.

b2) {\bf relation semi-ouverte} entre $\sup$ et $\inf $: de $ \inf $ vers $ \sup $, $ \sup \leftarrow \inf $

b3) $ \inf $ controle $ \sup $. Pas necessairement d'intervalles de precisions de nombres pour $ \sup $ et pour $ \inf $.

b4) estimation a priori implicite: $\sup \leq c <+\infty$ si $ \inf \geq m >0$, ca va dans un sens.

\smallskip

/////////////////////////////////////////////////////////////////////////////////////////////////////////////////////////////////////

\smallskip

39) Sur l'estimation asymptotique de l'article de D.Holcman: on ecrit: $ l_{\epsilon}=\epsilon^{(n-2)/2} \mu_{\epsilon} $, on s\'epere l'integrale en deux, la partie sur $ B_P(\delta) $ et celle sur $ M-B_P(\delta) $, puis on fait une changement de variable $ u=\epsilon r $ dans $ B_P(\delta) $:

1)

$$\epsilon^{-(n/2)+1}(\mu_{\epsilon})^{N-1} =\frac{l_{\epsilon}^{N-1}}{\epsilon^n}= $$

$$ =\int_{0}^{\delta/ \epsilon} r^{n-1} |\frac{1}{(1+r^2)^{(n-2)/2}}-\frac{\epsilon^{n-2}}{(\epsilon^2+\delta^2)^{(n-2)/2}}-\epsilon^{(n-2)/2}\mu_{\epsilon}|^{4/(n-2)} \times (\frac{1}{(1+r^2)^{(n-2)/2}}-\frac{\epsilon^{n-2}}{(\epsilon^2+\delta^2)^{(n-2)/2}}-\epsilon^{(n-2)/2}\mu_{\epsilon}) $$
 
 $$ \leq C+O(\delta^2)\frac{l_{\epsilon}^{N-2}}{\epsilon^2}+ O(\delta^n) \frac{l_{\epsilon}^{N-1}}{\epsilon^n} +O(\delta^n) \frac{l_{\epsilon}}{\epsilon^{n-4}}$$

Avec $ C >0 $, on a d\'evelopp\'e l'in\'egalit\'e pr\'ecedente et on a utilis\'e:$ |x+y|^{\alpha}\leq |x|^{\alpha}+|y|^{\alpha}, \alpha =\frac{4}{n-2}\leq 1 $ ($n\geq 6$). Si, $ l_{\epsilon} \to +\infty $, on multiplie par $ \epsilon^n >0$ on divise par $ l_{\epsilon}^{N-1} $ puis on fait tendre $ \epsilon $ vers $0 $, on obtient $ 0 < C'\leq 0 $ c'est pas possible. Donc $ l_{\epsilon} \to k \geq 0, k<+\infty $, on multiplie l'in\'egalit\'e precedente par $ \epsilon^n $ et on fait tendre $ \epsilon $ vers $ 0 $, on obtient $ l_{\epsilon} \to 0 $ quand $ \epsilon \to 0 $.

\bigskip

2) $ l_{\epsilon} \to 0 \Rightarrow \exists r_{\epsilon} \to +\infty $ tel que:

$$ \frac{1}{(1+r_{\epsilon}^2)^{(n-2)/2}}=l_{\epsilon}=\epsilon^{(n-2)/2}\mu_{\epsilon}, $$

2-1) Si $ 0 < r_{\epsilon} <\frac{\delta}{\epsilon} $: $ 0 \leq r \leq r_{\epsilon} \Rightarrow (r^2+1)^{-(n-2)/2} \geq (r_{\epsilon}^2+1)^{-(n-2)/2} $ et $ r\geq r_{\epsilon} \Rightarrow (r^2+1)^{-(n-2)/2} \leq (r_{\epsilon}^2+1)^{-(n-2)/2}=\epsilon^{(n-2)/2}\mu_{\epsilon}, $ apres on utilise le theoreme de convergence domin\'ee de Lebesgue:

$$\epsilon^{-(n/2)+1}(\mu_{\epsilon})^{N-1} =\frac{l_{\epsilon}^{N-1}}{\epsilon^n}= $$

$$ =\int_{0}^{r_{\epsilon}} \frac{r^{n-1}}{(1+r^2)^{(n+2)/2}}(1+ O(\frac{\epsilon^{(n-2)/2}\mu_{\epsilon}}{\frac{1}{(1+r_{\epsilon}^2)^{(n-2)}}}))+\int_{r_{\epsilon}}^{\delta/\epsilon} r^{n-1}(3(\epsilon^{(n-2)/2} \mu_{\epsilon})^{N-1})= $$

$$= C+O(\delta^n)\frac{(\epsilon^{(n-2)/2} \mu_{\epsilon})^{N-1}}{\epsilon^n}=C+O(\delta^n)\epsilon^{-(n-2)/2}\mu_{\epsilon}^{N-1} $$
 
avec $ C >0$. Donc,

$$ \mu_{\epsilon}\equiv \epsilon^{(n-2)^2/2(n+2)},$$

2-2) Si, $ r_{\epsilon} \geq \frac{\delta}{\epsilon} $, on utilise le meme type de calcul que dans le, 2-1) et le Theoreme de convergence domin\'ee de Lebesgue:

$$ \epsilon^{-(n/2)+1}(\mu_{\epsilon})^{N-1}  = \frac{l_{\epsilon}^{N-1}}{\epsilon^n}= $$

$$ =\int_{0}^{\delta/\epsilon} r^{n-1}(\frac{1}{(1+r^2)^{(n-2)/2}}+...)^{N-2}(\frac{1}{(1+r^2)^{(n-2)/2}}+...)= \int_{0}^{\delta/\epsilon}\frac{r^{n-1}}{(1+r^2)^{(n+2)/2}}(1+ O(\frac{\epsilon^{(n-2)/2}\mu_{\epsilon}}{\frac{1}{(1+r^2)^{(n-2)}}})) \to C >0, $$

avec, $ r\leq \frac{\delta}{\epsilon} \leq r_{\epsilon} \Rightarrow \epsilon^{(n-2)/2} \mu_{\epsilon} =(1+r_{\epsilon}^2)^{-(n-2)/2}\leq (1+r^2)^{-(n-2)/2}  $. Donc,

$$ \mu_{\epsilon}\equiv \epsilon^{(n-2)^2/2(n+2)}. $$

////////////////////////////////////////////////////////////////////////////////////////////////////////////////////////////////////////

\smallskip

40) D-Branes, cordes, supercordes, Yamabe, courbure scalaire, Kaluza-Klein:

Sur les D-Branes, c'est des cordes avec condition de Dirichlet au bord. C'est construit comme pour Kaluza-Klein, $ (n+1)$. D'apres Bennequin, dans Ast\'erisque, Seminaire Bourbaki, 2003, voir aussi, Collion-Vaugon, pour la fibration. On prend la vari\'et\'e totale, $ W $ en $ (n+1) $ de la forme $ W=X\times M $, avec $ M $ une vari\'et\'e de Kahler-Einstein ou Aubin-Calabi-Yau, $ K_6 $ de dimension complexe $ 3 $ et r\'eelle $ 6 $. Donc, Ricci plate donc scalaire plate. Et, $ X $ l'espace temps de dimension $ 4 $ de la forme $ (3+1) $, $ X=Y\times {\mathbb R} $ avec $ Y $ de courbure scalaire $ \geq 0 $ et du au th\'eor\`eme de la masse positive, la courbure scalaire est quelque part $ >0 $, car la courbure agit et courbe l'espace (aussi, mesures experiementales sur l'existence d'ondes gravitationnelles, donc, $ Weyl \not =0 $ ici en dimension 3 pour l'espace $ Y $ le tenseur de Cotton non nul quelque part, se rappeler aussi que la courbure scalaire est li\'ee \`a la courbure ou l'inverse du rayon de courbure et la masse est $ >0 $) et de masse $ m >0 $, voir le monographe de Jonathan.Rosenberg et David.Wraith, sur le fait que la courbure scalaire est $ >0 $ quelque part et est $ \geq 0$ (Titre: Positive Scalar curvature). Comme le formalisme des D-branes, est aussi li\'e \`a la reduction de Kaluza-Klein, on a l'equation de Yamabe ou de type courbure scalaire prescrite avec condtions de Dirichlet au bord, comme $ S_g \geq 0 $ et $ S_g >0 $ quelque part, on voit qu'on est dans le cas d'un operateur coercif pour les D-branes avec formalisme ou reduction de Kaluza-Klein, car on a $ Y\times M $ et $ S_g=Scalaire_Y+Scalaire_M=Scalaire_Y+0=Scalaire_Y \geq 0 $ et $ >0 $ quelque part. L'invariant de Yamabe est $ >0$ ou operateur conforme coercif.(on peut prendre aussi, $ M=K_6=T_6 $ le tore de dimension 6, Einstein-Kahler).

\smallskip

Voir aussi le print de, Mariana Grana et Hagen Triendl, Saclay. Kaluza-Klein mechanism and d-branes. Titre: String Theory compactifications.(ils disent que le monde ou espace est un brane, $ Y\times M$, c'est cens\'e etre une vari\'et\'e \`a bord.).

\smallskip

Voir aussi l'article de Schoen-Yau sur le th. de masse positive, ils disent que la vari\'et\'e de dimension 3 est sans bord ou avec bord, avec une metrique complete,  $ {\mathbb R}^n $ n'est ni avec bord ni sans bord, car "son bord" est $ \infty $, un point, mais cela ne peut pas ete un bord, car une vari\'et\'e \`a bord, son bord est une vari\'et\'e sans bord de dimension $ 2 $. On peut considerer par exemple la boule unit\'e diffeomorphe \`a $ {\mathbb R}^3 $ avec l'application $ x\mapsto x/(1-||x||) $, ou bien sans bord avec une vari\'et\'e proche de la sphere de dimension 3, ou une d\'eformation de la sphere de dimension 3.(par la projection stereographique, on en leve un point, le pole nord et l'hemisphere sud, et on applique la projc.steorographique de pole, le pole nord, c'est un bout, par exemple). Comme Schoen-Yau disent que la vari\'et\'e est soit sans bord ou avec bord, mais avec des coordonn\'ees asymptotiques dans les bouts, avec une metrique asymptotique qui la rend complete, et ici, c'est cens\'e etre un brane, avoir une membrane, ou un bord. On est dans le cas avec bord dans le th\'eor\`eme de la masse positve, mais on a toujours, la courbure scalaire $\geq 0$ et $ >0 $ quelque part. 

\smallskip

Pour l'Eq. de Yamabe, on a un potentiel constant $ >0 $ et l'operateur est l'operateur conforme, $ \frac{4(n-1)}{n-2}\Delta + S_g $, $ \Delta=-\nabla^i\nabla_i$, l'operateur conforme est coercif. mais quand le potentiel varie, l'operateur devient $ \frac{4(n-1)}{n-2} \Delta+h $ avec $ h=S_g-|\nabla \Psi|^2 \leq S_g $. Ce nouvel operateur reste encore coercif, car du point de vue de la physique, l'impulsion est proche d'une constante et diminue progressivement, un observateur exterieur, regarde l'effet de l'impulsion a l'instant $ t $ o\`u elle est maximale et au voisnage de l'impulsion maximale, c'est \`a dire voisine d'une constante. C'est \`a dire $ \Psi(x_0) +\epsilon $ avec $ \epsilon >0$ petit. Il y a une impulsion \`a l'instant $ t $, qui est constante et apres on regarde au voisinage de $ a= \Psi(x_0)=constante$ \`a \, instant $ t$ .

Ceci confirme et explique  ce qu'on a dit au point $(37)$, sur Kaluza-Klein classique en presence d'un potentiel en dimension $ n=4 $ ou plus g\'en\'eralement $ n\not = 4 $ . Du point de vue de la physique, un observateur exterieur regarde ce qui se passe pour une impulsion \`a l'instant $ t $ et au voisinage de l'impulsion. Oscillation autour d'une constante. Il y a impulsion, observateur, puis le potentiel diminue de $ \epsilon >0 $ petit. Quand on met un potentiel, un observateur exterieur regarde ce qui se passe autour d'une impulsion, c'est \`a dire autour d'une constante.

{\bf Remarque:} Ceci est une interpretation assez intuitive ( les espaces sont {\bf independants les uns des autres}, mais forment des espaces-temps dont la restriction a la vari\'et\'e voulue ou cherch\'ee, est Riemannienne). Concernant la formulation en $ (n+1) $, on ecrit: $ W=M\times X =M \times (Y\times {\mathbb R}) $ et par Jonathan. Rosenberg et David. Wraith (titre: positive scalar curvature), on a $ X=Y\times {\mathbb R} $ est un espace temps tel que la metrique Lorentzienne sur $ X $ est Riemannienne sur $ Y $ ou sa restriction. Puis on regroupe $ M $ et $ Y$ en un produit, car les espaces sont independants l'un de l'autre. Donc, on obtient un espace temps  total $ W $ tel qu'on ait une metrique produit sur $ M\times Y$, ainsi, les courbures scalaires s'additionnent et donnent une courbure scalaire  totale, $ S_g \geq 0 $ et $ >0 $ quelque part par la contribution de la courbure scalaire de $ Y $, donn\'ee par le th. de la masse positive.

Cas 1: $ Y $ compacte sans bord, alors la vari\'et\'e Riemannienne $ M\times Y $ est compacte sans bord. Ici les cordes sont ferm\'ees.

Cas 2: $ Y $ compacte \`a bord, alors la vari\'et\'e Riemannienne $ M\times Y $ est compacte \`a bord. Ici les cordes sont ouvertes. avec condition de Dirichlet. C'est le cas trait\'e dans ce point.

\smallskip

//////////////////////////////////////////////////////////////////////////////////////////////////////////////////////////////

\smallskip

41) Sur le cas n\'egatif de l'Eq. de courbure scalaire (et sous-critique) (solutions positives) et supercritique negatif (pour ce cas les solutions peuvent changer de signe): Ici on prouve la compacit\'e locale sans conditions sur la courbure ou la vari\'et\'e Riemannienne. Les estimations a priori locales dans le cas n\'egatif ne dependent pas de la courbure ou de la nature de la vari\'et\'e Riemannienne. Il n'y a pas de conditions sur la courbure, on a des estimations du type Keller-Osserman. Dans l'article de Marco Rigoli: il a suppos\'e la vari\'et\'e complete et a impos\'e une condition sur la courbure de Ricci, pour obtenir des estimations asymptotiques  a l'infini du type Keller-Osserman, \`a la distance \`a un point: $ |u|\leq C/[d(x,x_0)]^{\tau}, \tau >0$, par principe du maximum et des raisonnement du type theoremes de comparaisons, ici la notation $ |u| $ signifie qu'ils majorent et ils minorent les solutions, d'autres estimations sont donn\'ees pas L.V\'eron et M.F.Bidaut-V\'eron, pour d'autres types d'equations ou systemes. Estimations asymptotiques aussi pour la compl\'etude. Aussi, l'article de L.V\'eron, 1992, Journal d'Analyse Mathematique, ici, sur un domaine de $ {\mathbb R}^n $ et operateur "strongly elliptic" pour se ramener au laplacien dans $ {\mathbb R}^n $, operateur non en divergence forme. 

Il y a aussi l'article de Aviles-McOwen, Journ.Diff.Geometry.1988.(Examen de DEA). Ici, ils considerent le potentiel constant $ K\equiv -1$ sur une vari\'et\'e complete et considerent une suite exhaustive de compacts.

Il y a aussi l'article de Ratto.Rigoli.V\'eron. Math Zeitschrift. 1997. Ils prennent le cas ou la courbure prescrite $ 0 > K\in C^{\infty}(M) $  {\bf fix\'ee} et sur une {\bf vari\'et\'e compl\`ete} $ M $ non-compacte et disent qu'ils s'inspirent de la m\'ethode de Aviles-McOwen, donc aussi avec une suite exhaustive de compacts.

A aucun moment ces auteurs, ne posent le probleme de l'estimation a priori et de compacit\'e locale(ne dependant que des parmetres exterieurs, $a, b, A, M, B, g $, avec $ B $ le compact), avec ou sans conditions sur la courbure(courbure de Ricci, de Riemann...) ou la vari\'et\'e Riemannienne (localement conform\'ement plate ou non, compl\`ete, \`a bord, sans bord...). On a aussi, des estimations du type Keller-Osserman: pour toute carte normale en $ y_0 $,  $ \sup_{B_R(y_0)} uo\phi^{-1} \leq \frac{c}{R^{2/(q-2)}} $ par un la technique "blow-up" ou eclatement, on s\'epare les constantes, une estimation a priori du type in\'egalit\'e de Harnack (dans l'in\'egalit\'e de Harnack, il y a majoration par l'integrale, ici, il n'y en a pas, c'est une estimation a priori du type inegalit\'e de Harnack sans condition sur l'integrale. Autre chose, la constante $ c $ depend du point, car on est sur une vari\'et\'e, ca depend du point, alors que sur ${\mathbb R}^n $, en tout point on a une carte normale et c'est invariant par translation, ce qui fait que les constantes ne dependant pas du point, l'inegalit\'e de Harnack usuelle ne depend pas du point, ceci est li\'e au fait que l'espace est plat partout, c'est une hypothese forte, en tout point la carte est normale, c'est l'identit\'e. Mais  la constante de l'inegalit\'e de Harnack depend du rayon. La technique blow-up qu'on a utilis\'e fait que $ c $ et $ R $ sont s\'epr\'ees, la constante totale est plus explicite $ C=\frac{c}{R^{2/(q-2)}} $, on explicite la dependance en $ R $ de $ C $, et aussi en $ c $ de $ C $. Et l'estimation a priori est vraie pour tout point $ y_0 \in M $ et toute carte normale, il y a l'influence de la g\'eometrie, ce qui etait sous-entendu dans $ {\mathbb R}^n $. Pour tout $ y_0 $ et toute carte normale, on a l'estimation a priori au voisinage de $ y_0 $ c'est a dire qu'a l'interieur d'un voisinage de $ y_0 $ ca ne depend pas du point $ y \in V_{y_0}$, cela suffit comme estimation a priori, ca depend d'un rayon maximal $ R $ du point $ y_0 $ et la forme explicite de $ C $ fait que c'est une estimation a priori du type Harnack, comme on l'a dit ci-dessus, uniforme en $ y \in V_{y_0} $, comme si on etait dans le cas plat, localement uniforme avec dependance du type Harnack en $ R $ rayon maximal). Non seulement on a l'estimation locale globale dans une boule $ B(y_0,R)$, mais aussi, dans toute petite boule $ B_{R'}(y)$ de $ V_{y_0}=B(y_0, R) $:

$$ \sup_{B_{R'}(y)} u_io\phi^{-1} \leq \frac{c}{{R'}^{2/(q-2)}}, \qquad (**)$$

avec $c=c(y_0, a, b, A, \alpha, M, g) $, c'est une estimation a priori du type Harnack. On peut l'ecrire aussi pour la distance geodesique en comparant les boules pour les metriques $ g, \phi^*(g)=h $. ici, on a melang\'e boules de cartes, euclidiennes et boules geoedesiques, mais on a toujours l'estimation a priori du type Harnack:

a) Pour un point $ y \in V_{y_0}=B(y_0, R) $ et un r\'eel $ R' >0$ donn\'es: on a  l'estimation a priori globale par rapport \`a un rayon maximal $ R >0 $ li\'e \`a $ y_0 $. Estimation par rapport \`a un rayon maximal $ R >0 $ li\'e au point $ y_0 $. A la Bonnesen, du type in\'egalit\'e de Bonnesen, par rapport \`a un rayon maximal.

b) Soit par soucis d'homog\'en\'eit\'e, l'in\'egalit\'e $ (**)$ qui exprime comment les solutions $ u_i $ evoluent en fonction de $ y $ et $ R'>0 $  uniform\'ement dans $ V_{y_0}=B(y_0,R) $, ici $ R >0 $ est li\'e au point $ y_0 $, alors que $ y $ et $ R'>0 $ sont pris de maniere quelconque dans $ V_{y_0}=B(y_0,R) $. On a la formule $ (**)$ vraie pour tout $ y $ et $ R'>0 $ dans $ V_{y_0}=B(y_0,R) $. Estimation qui peut etre interieure et spherique. Par soucis d'homog\'eneit\'e, pour $ y $ et $ R' $, on veut voir comment evoluent les solutions $ u_i $ en fonctions de $ y $ et $ R'$: il n'y a pas de dependance en $ y $, mais ca depend d'un rayon maximal qui est li\'e \`a $ y_0 $, un autre point. $ (**) $ permet d'exprimer l'evolution de $ u_i $ en fonction de $ R'$, on considere une boule $ B_{R'}(y) $ et on regarde comment evoluent $ u_i $ en fonction des parametres de cette boule, donc, en fonction de $ y $ et $ R'$, ca ne depend par de $ y $ dans $ V_{y_0} $ mais ca depend de $ R'$. On a, $ y_0 \leftrightarrow R $ et $ y\leftrightarrow R'$.

c) Ces assertions sont equivalentes, l'estimation a priori locale globale dans $ V_{y_0}=B(y_0,R) $ est equivalente \`a l'estimation a priori du type Harnack:

$$ \sup_{B(y_0,R)} u_io\phi^{-1} \leq \frac{c}{R^{2/(q-2)}} \Leftrightarrow \sup_{B_{R'}(y)} u_io\phi^{-1} \leq \frac{c}{{R'}^{2/(q-2)}}, $$

avec, $ B_{R'}(y) \subset B(y_0,R), \forall R'\leq R, \forall y \in B(y_0,R) $.

\smallskip

Dans le cas super-critique negatif, on a comme consequence, un principe d'Harnack. On a un crit\`ere de compacit\'e si la vari\'et\'e est compacte sans bord et $ \int_M R <0$ avec $ R $ la fonction de l'Eq. $ \Delta u + R = Ve^u, \Delta = -\nabla^i \nabla_i, a\leq V\leq b <0 $.

\smallskip

Pour ce qui nous concerne, la compacit\'e locale ne necessite pas de conditions sur la courbure ou la vari\'et\'e Riemannienne. Les {\bf estimations a priori locales}, {\bf compacit\'e locale}, dans le cas n\'egatif, ne necessitent pas de conditions sur la courbure ou sur la vari\'et\'e Riemannienne.

\smallskip

(Ce qui suit,  a \'et\'e dit dans l'introduction du print "Cas d'existence de solutions d'EDP", quand on a parl\'e l'elimination des blow-up isol\'es simples par le th. de la masse positive. Mais cela reste vrai dans le cas strictement n\'egatif, sans th.de la masse positive. On l'a dit, compacit\'e: resultat d'existence par le degr\'e topologique et par consequent l'estimation a priori. Voir aussi pour les solutions variationnelles dans le livre d'Aubin chapitre 6, cas strictement n\'egatif.)

1) Consequence pour le cas strictement negatif avec potentiel negatif strictement et sur une vari\'et\'e compacte sans bord: ($ \Delta=-\nabla^i \nabla_i $, $ c_n $ la constante de l'operateur conforme, on se ramene a la courbure scalaire $ S=-1$):

On a une solution topologique pour $ c_n\Delta u-u=Vu^{N-1}, u >0, -\infty < a \leq V \leq b <0, ||V||_{C^{\alpha +\epsilon}} \leq A, \epsilon >0 $. Pour cela on considere l'homotopie d'operateurs compact: $ K_t(u)=(c_n\Delta +1)^{-1}(V_t u^{N-1}+2u)$, avec $ V_t=-(1-t)+tV, t\in [0,1] $. On a la compacit\'e, implique qu'on peut considerer le degr\'e de Leray-Schauder dans $ X =\{u\in C^{2,\alpha}(M), u \geq 0, ||u||_{C^{2,\alpha}} \leq 2c_0\} $, de sorte que $ (Id-K_t)(u)=0 $ n'a pas de solution $ u $ dans $ \partial X $, on peut alors considerer le degr\'e de Leray-Schauder dans $ X $. Par homotopie:

$$ deg(Id-K_1, X, 0)=deg (Id-K_0, X, 0), $$

Or, $ f=Id-K_0 $ a une deriv\'ee qui est un isomorphisme, car dans le cas negatif, la seul solution de $ (Id-K_0)(u)=0$ est $ u=u_0\equiv 1$ (l'eq. est: $ c_n\Delta u-u=-u^{N-1} $, soit on utlise le resultat d'unicit\'e dans le livre d'Aubin, Kadzdan-Warner, chapitre 6, soit, se placer au maximum et au minimum on a a la fois $ u\leq 1$ et $ u\geq 1$, voir le livre d'Aubin, pour l'unicit\'e dans le cas n\'egatif, chapitre 6). Et en $ u_0=1$, $ f'(h)=f'(u_0)(h), u_0=1 $, est un isomorphisme $ f'(h)=h-(c_n\Delta+1)^{-1}(-(N-3)h) $ et en appliqant $ c_n\Delta+1 $ le noyau est reduit a $ 0 $ et $ f'$ est surjectif. Car l'operateur $ c_n\Delta +(N-2) $ est inversible. Donc le degr\'e de $ f $ est l'index de $ f'=I-\tilde K_0 $, $ \tilde K_0 $ compact et $ f' $ est un isomorphisme et est egal \`a $ +1 $ ou $-1$, voir l'article de Leray-Schauder, 1934, titre: Topologie et \'equations fonctionnelles, Annales de l'E.N.S, voir aussi l'article resum\'e de Jean Mawhin, donc le degr\'e est non nul. Donc: le degr\'e de $ Id-K_1$ est non nul ce qui veut dire que l'eq. $ c_n\Delta u-u=Vu^{N-1}, u >0 $, a une solution topologique, et on a l'estimation a priori.

\smallskip

Si on divise l'eq. par $ u >0$ et on integre, on obient que l'integrale de $ u^{N-2} $ est uniformement minor\'ee par une constante $ >0$ et par le principe du maximum les solutions sont alors unifrom\'ement minor\'ees par une constante $ >0$, on utilise la compacit\'e et la convergence vers une fonction non identiquement nulle, puis le principe du maximum implique qu'elle est $ >0 $. On pouvait prendre l'ensemble $ X =\{u\in C^{2,\alpha}(M), u \geq \frac{c_0'}{2}>0, ||u||_{C^{2,\alpha}} \leq 2c_0\} $, c'est cet ensemble qu'il faut prendre. 

Ici, on n'a pas utilis\'e l'id\'ee du livre d'Aubin, dans le cas compact, se placer au maximum et au minimum, pour majorer et minorer les solutions, c'est une methode differente. 

2) Pour le cas super critique et negatif, sur une vari\'et\'e compacte sans bord, on a la meme chose avec l'eq. $ \Delta u-1=Ve^u, u \geq 0 $. Ici on a impos\'e un signe, mais le fait que la vari\'et\'e soit compacte, implique que les solutions sont minor\'ees, en se placant au point minimum, $ Q $, $ \Delta u (Q) \leq 0 $. Ici aussi on peut se passer de cette id\'ee du livre d'Aubin: en effet, on utilise la majoration uniforme des solutions, on resout un probleme de Dirichlet, $ w $, car $ Ve^u+1 $ est uniformement born\'ee (On se place dans des ouverts de cartes de bords r\'eguliers $ (U_x,\phi_x), x\in M$ et $ w=w_x $, onse ramene a $ {\mathbb R}^n $)et on soustrait $ w $, $ u-w$, on se ramene a des fonctions harmoniques, puis on utilise le principe de Harnack. Donc, soit $ u $ converge, soit, elle tend vers $ -\infty $, or l'intergrale de $ u $ est uniformement minor\'ee par une constante $ >0 $, ce qui n'est pas possible.(on est sur une vari\'et\'e compacte sans bord, on integre l'eq de $ u $). Donc, les $ u $ sont uniform\'ement minor\'ees.

L'application $ K_t(u)=(\Delta +1)^{-1}(V_t e^u+1+u) $ et $ X =\{u\in C^{2,\alpha}(M), u \geq \frac{k_0'}{2}, ||u||_{C^{2,\alpha}} \leq 2k_0\} $. Avec, $ V_t=-(1-t)+tV $ et $ -\infty < a\leq V\leq b <0$ et $ ||V||_{C^{\alpha +\epsilon }} \leq A, \epsilon >0 $. Avec, $ g(u)=K_0(u)=(\Delta+1)^{-1}(-e^u+1+u) $ et $ g'_{u_0=0}(h)\equiv 0$, donc, $ (Id-g)'_{u_0=0}(h)=h=Id $.

{\bf Remarques:} a) En ce qui concerne l'unicit\'e pour l'eq. $ \Delta u-1=-e^u$ sur une vari\'et\'e compacte. On veut prouver que $ u\equiv 0$. Alors on ecrit $ \Delta u=(1-e^u) $, on multiplie par $ u $ et on integre, on obtient: on a toujours $ (1-e^u)u \leq 0 $:

$$ 0\leq \int_M |\nabla u|^2=\int_M (1-e^u) u \leq 0.$$

Donc: $ u(1-e^u)\equiv 0$, donc, $ u\equiv 0$.

b) On a la meme chose avec l'eq. $ c_n \Delta u-u=-u^{N-1}, u >0$, on deplace $ u $ au second memebre et on multiplie par $ u-1 $ et on integre, on obtient un terme positif ou nul egal a un terme negatif ou nul, donc, nul et donc, $ u\equiv 1$.

$$ c_n \Delta u= u(1-u^{N-2}), $$
on multiplie par $ u-1 $ (on a toujours $ u (1-u^{N-2})(u-1) \leq 0 $)  et on integre:

$$ 0\leq c_n\int_M |\nabla u|^2=\int_M u(1-u^{N-2})(u-1) \leq 0.$$

Donc, $ u(1-u^{N-2})(u-1)=0$, d'ou, $ u\equiv 1$.

c) Tout ceci pour ne pas utiliser l'id\'ee du livre d'Aubin: se placer au minimum et au maximum des solutions. Dans le cas d'un potentiel, on a l'unicit\'e aussi, par la technique du point a) et b).

\smallskip

////////////////////////////////////////////////////////////////////////////////////////////////////////////////////////////////

\smallskip

42) Sur l'article de Brezis-Merle et le probleme 1 de Brezis-Merle. Dans leur introduction Brezis-Merle commencent par se poser le problemes des estimations a priori. $( p=+\infty) $

a) Ils trouvent une condition sur les bornes $ C_1 $ et $ C_2 $ pour avoir la compacit\'e jusqu'au bord. Ici, les potentiels peuvent changer de signe et $ C_1 $ et $ C_2 $ sont quelqconques, mais $ C_1\cdot C_2 <4\pi $, c'est une condition sur la "masse" ou l'invariant conforme qui est la masse. C'est une condition \`a la Bol ou \`a la Fiala ou \`a la Alexandrov pour les in\'egalit\'es isop\'erimetriques.

C'est une condition "masse-niveau d'energie" comme pour les in\'egalit\'es isoperimetriques. Et sur un invariant conforme. Pour avoir la compacit\'e globale.

b) Comme Brezis-Merle ont restreint  le champ des possibilit\'es pour $ C_1 $ et $ C_2 $. Ils ont trouv\'e une condition sur le potentiel $ V $ pour qu'il n'y ait pas de conditions sur $ C_1 $ et $ C_2 $ et l'invariant conforme $ C_1\cdot C_2 $. Cette condition est la positivit\'e de $ V $. Ils obtiennent la compacit\'e locale. Mais la compacit\'e globale n'est pas vraie, comme le montre leur exemple et contre-exemple de la fin de l'article.

C'est un exemple d'existence de solutions d'un probleme aux limites avec condtions de Dirichlet, avec les conditions de Brezis-Merle qui illustre bien la compacit\'e locale et met en defaut la compacit\'e globale. Ces exemples existent pour toute masse $ m=4\pi A, A>1 $.

\smallskip

{\bf Masse ou invariant conforme, conditions sur V (et ou $ C_1 $) et niveau d'energie ou masse. (Pour ce qui est de l'in\'egalit\'e isoperimetrique de Bol, Fiala, Alexandrov)}.

\smallskip

c) Donc, pour le probleme 1 de Brezis-Merle, se pose la question de la compacit\'e globale, jusqu'au bord. Pour cela, il faut des conditions:

c1) Soit une condition sur $ C_1$ et $ C_2 $ ? Comme la positivt\'e de $ V $ pour la compacit\'e locale, qui est une sorte de condition sur la borne inferieure de $ V $, comme pour les in\'egalit\'es, $ \sup, \inf $. 

Ou comme pour les espaces uniform\'ement convexes. Des valeurs explicites pour caracteriser un ensemble.

c2) Soit une condition sur la masse $ C_1\cdot C_2 $ ? Comme pour les in\'egalit\'es isoperimetriques de Bol, Fiala, Alexandrov. Ou Comme pour les espaces uniform\'ement convexes, des valeurs explicites pour caracteriser un ensemble. Ou les conditions sur l'invariant de Yamabe de Aubin en dimensions $ \geq 3 $.

Brezis-Merle dans leur exemple et contre-exemple insistent sur le fait que la masse $ m=4\pi A, A >1 $ est arbitraire.

c3) Soit une condition sur $ V $ ? C'est en partie le probleme 1 de Brezis-Merle, les potentiels convergent en norme $ C^0(\bar \Omega) $ ou comme dans la remarque qui le precede, des conditions fortes sur $ V $.

c4)Soit une condition sur les solutions $ u $ ? comme pour les theoremes dans les $ L^p $, Lebesgue, Beppo-Levi, Fatou, ou d'Ascoli-Arzela ou Dunford-Pettis, ou les in\'egalit\'es $ \sup,\inf $.

\smallskip

 d) Chen-Li, donnent une reponse a ce probleme (le probleme 1 de Brezsi-Merle) en imposant \`a $ V $ d'etre $ C^1 $ uniform\'ement. On voit alors qu'a partir de la difference $ |V(x)-V(y)|\leq A |x-y|$, on a au moins, aussi, la condition sur $ C_1\leq A \cdot diamtre(\Omega)$. Ici, on a trois conditions, $ V $, $ C^1$, uniform\'ement, et sur $ C_1 $. Et ils obtiennent la compacit\'e globale pour toute valeur de la masse, $ m= C_1\cdot C_2 >0 $ quelconque. Independemment de la masse. 
 
 On a aussi, les articles de, De Figueiredo-Lions-Nussbaum, cas sous critique et Suzuki, cas $ V $ constant et Ma-Wei, cas $ V $ constant. Compacit\'e globale par la m\'ethode moving-plane.

\smallskip

 e) Pour ce qui nous concerne, dans {\it "About Brezis-Merle Problem with Holderian condition"}. On fixe la regulrit\'e de $ V $, $s-$holderienne, $ \frac{1}{2} < s <1 $ et on a une condition sur la masse ou l'invariant conforme, $ C_1 $ et $ C_2 $ sont quelconques, mais $ C_1\cdot C_2 < 24\pi $, borne explicite de $ 24\pi $. Condition du type in\'egalit\'e isoperimetrique de Bol, Fiala, Alexandrov, ou de Brezis-Merle pour des potentiels changeant de signe.

\smallskip

 f) Maintenant, si on veut eliminer cette condition sur la masse ou l'invariant conforme. Il faut une condition sur le potentiel $ V $. C'est l'objet de l'article: {\it "A Compactness result for an equation with Holderian condition"}. On prend $ V=W\cdot (1+|x|^{2\beta}), \beta \in (0,1/2) $ avec $ W $, Lipschitz et uniform\'ement. Presence d'un poids et ce poids est holderien et Sobolev. On obtient la compacit\'e pour une masse $ m =C_1\cdot C_2 $ quelconque.

Cette condition du point $ f)$ est suffisante, car le poids engendre un cone holderien. Par exemple, si on prend $ V $ une fonction holderienne totale, elle verifie; $ |V(x)-V(y)|\leq B|x-y|^{2\beta}, \beta \in (0,1/2)$, si $ 0 \in \partial \Omega $ est sur le bord, on ecrit $ |V(x)-V(0)|\leq B|x|^{2\beta} $, c'est a dire que $ V(0)-B|x|^{2\beta} 
 \leq V(x)\leq V(0)+B|x|^{2\beta} $, on voit alors que $ V $ est equivalente au poids $ 1+k|x|^{2\beta} $ au voisinage de $ 0 $. Donc, il suffit de prendre directement le poids, pour avoir une bonne approximation d'une fonction holderienne.

 Le poids $ (1+k|x|^{2\beta}) $, force $ W $ a etre holderienne dans un cone centr\'e en $ 0 $.

Dans la pratique, par la methode des differences finies et des elements finis, on considere, des approximations, un maillage, il suffit de considerer les solutions en des points, pour cela il suffit d'avoir un porduit fini de poids en un nombre fini de points, qui peut etre grand. Ceci est possible et c'est une approximation. $ V=W\cdot \prod_{i=1}^k (1+a_i|x-x_i|^{2\beta}), a_i=\pm \epsilon_i, 1>>\epsilon_i >0 $, $ \epsilon_i >0 $ petit, $ W $ Lipschitz et uniforme. Avec $k$ tres grand si possible. 

Ceci concerne le probleme 1 de Brezis-Merle avec un produit de poids. On a une approximation par un produit de poids holderiens d'une fonction holderienne. Approximation des fonctions et des solutions en considerant un nombre fini de points du maillage. 

\smallskip

g) Maintenant, si on s'interesse au volume et une borne uniforme du volume lorsqu'on a un produit de $ k $ poids holderiens. On se place sur un domaine analytique et on utilise les id\'ees (en partie) du print {\it "Uniform bound for the volume of the solutions to Liouville type equations on the annulus"}. $ 0 <a \leq W \leq b $ et $ W $, $ C^1$ uniform\'ement, $ ||\nabla W||\leq A $.

On ne se place pas forc\'ement sur une couronne, mais sur un domaine analytique. Alors, pour tout $ x\in \partial \Omega $, il existe une tranformation conforme $ f_x $ d'un voisinage de $ x $ sur la demi boule unit\'e. Alors, avant de faire agir la transormation conforme autour d'un point $x_i$ du poids holderien, on utlise l'estimation locale par la m\'ethode moving-plane de Chen-Li et on se retrouve (en partie)dans la situation du print {\it "Uniform bound ... on the annulus"}. Tout est born\'e en dehors de la singularit\'e, sauf en la singularit\'e. Alors on utlise la tranformation conforme $ f_{x_i} $, le bord contenant la singularit\'e devient plat et on se retrouve (en partie) comme dans le print {\it "Uniform bound...on the annulus"}. On utlise la formule de Pohozaev avec le poids, $ 1+a_i|x-x_i|^{2\beta} $, on obtient: termes de bords $= O(1) + 0 = $ l'invariant conforme a un facteur pr\'es d\^u \`a $ f_{x_i}$ et on obtient en prenant le rayon du demi-disque assez petit que le volume est alors born\'e par une constante dependant de $ a, A, \beta, f_{x_i}, \Omega $, qui sont fix\'ees au depart. 

\smallskip

Comme pour le probleme de Brezis-Merle avec $ k $ poids holderiens, ou Lipschitziens, dans la pratique, par la methode des differences finies et elements finis, il suffit de considerer $ k $ grand, on a une approximation du volume et des fonctions et des solutions.

La compacit\'e au voisiange du bord de Chen-Li est obtenue par la m\'ethode moving-plane. Elle est locale (au voisiange de chaque pint du bord sous certaines conditions), il faut supposer le potentiel ou la courbure prescrite, $ V $, $ C^1$.

\smallskip

//////////////////////////////////////////////////////////////////////////////////////////////////////////////////////////////////

\smallskip

43) Pour revenir aux supercordes, dans l'article de 2006 dans Bulletin des Sciences Mathematiques, on a ecrit pour un point interieur $ x_0 $, $ x_0\in M $:

$$ u_i(x_0) \geq m >0 \Rightarrow \sup_M u_i \times \inf_K u_i \geq c=c(a,b, m, x_0, K, M, g) >0, $$

Si on prend par exemple $ x_0 \equiv K $ le compact $ K $, on peut ecrire (on a la meme preuve, car on fixe $ x_0 $ et le compact $ K$):

$$ \sup_K u_i\geq m >0 \Rightarrow \sup_M u_i \times \inf_K u_i \geq c, $$

On remplace $ m=\sup_K u $, on obtient l'in\'egalit\'e de Harnack "enroulement-torsion", suivante:

$$ \sup_M u \times \inf_K u \geq c(a,b, \sup_K u, K, M, g) >0,$$

la fonction $ m \to c $ est croissante de $ m >0 $. (dans le cas de d\'ependance lin\'eaire, on a l'inegalit\'e de Harnack classique d'un cot\'e).

\smallskip

On a la fois l'enroulement et la torsion et c'est une in\'egalit\'e "harmonique-Yamabe", car elle met en relation $ \inf_K u $ et $\sup_K u $ pour ce qui est harmonique, et $\sup \times \inf $ pour ce qui est Yamabe. ici, relation entre $ \sup_M u, \sup_K u, \inf_K u $. ou dans le sens $ \sup_M u, \inf_K u, \sup_K u $. On a les relations suivnates: ouverte+semi-ouverte: $ \sup_M u \leftrightarrow \inf_K u \leftarrow \sup_K u $. Ou si, on place tout ce qui est inf d'un cot\'e et tout ce qui sup de l'autre:

$$ \inf_K u \geq \frac{c(a,b,\sup_K u, K, M, g)}{\sup_M u}= C(a,b, \sup_K u, \sup_M u, K, M, g) >0, $$

ou Formelement, on ecrit:

$$ \inf \geq c(\sup)>0,$$

ou si on revient a l'ecriture initiale,

$$ \sup_M u \times \inf_K u \geq c(a,b, \sup_K u, K, M, g) >0.$$

On voit alors pour ce qui est $ n+1$ avec $ n\geq 3 $ (et en particulier $ n=9$), concernant les supercordes et Kaluza-Klein, l'enroulement, la torsion, qui sont m\'elang\'es, et, harmonique-Yamabe, les differentes relations d'Harnack m\'elang\'ees et ceci sans n\'ecessairement avec conditions au bord, Dirichlet ou Neumann ou mixtes.

{\bf Remarques: quelques exemples }

1) Pour les fonctions solutions d'edp du type:

$$ Lu=f,$$

On a des in\'egalit\'es de Harnack avec une contibution de la norme $ L^n $ de f:

$$ \sup_K u \leq C_1\inf_K u + C_2||f||_n, $$

Si $ f= Vu^{N-1} $, on voit qu'on doit avoir une contibution de $ ||u||_p, p=nN $, dans l'in\'egalit\'e de Harnack. Pour des equations plus generales que les fonctions harmoniques (presence d'un $ f$), il faut une information ou une contirbution de $ f $, si de plus on suppose $ f $ dependant de $ u $ alors, dans l'in\'egalit\'e de Harnack classique il faut une contribution d'une integrale de $ u $. Voir le Gilbarg-Trudinger pour les in\'egalit\'es de Harnack classiques avec une fonction $ f $.

On explique ainsi l'inegalit\'e de Harnack a trois parametres, ci-dessus, avec, $ \sup_M u, \inf_K u, \sup_K u $. La presence de 3 parametres.

2) Dans le Gilbarg-Trudinger, il y a l'in\'egalit\'e d'oscillation, dans la solution de, De Giorgi et Nash, de r\'egularit\'e holderienne. Le lemme d'oscillation: il fait apparaitre $ \sup_{B_R} u, \inf_{B_R} u $ et $\sup_{B_{R_0}}|u| $, avec $ R_0 $ fix\'e plus grand ou egal que $ R$. Voir livre de Gilbarg-Trudinger, les chapitres.8 et 9 et en particulier le theroreme.15.7.

On explique que dans l'in\'egalit\'e de Harnack a 3 parametres, ci-dessus, on prend $ M $ plus grand que $ K $ avec presence de $ \sup_K u$, $\inf_K u $ et $\sup_M u$. 3 parametres dont deux d'Harnack harmonique et deux de Yamabe.

3) on explique la presence de 3 parametres dans l'inegalit\'e de Harnack harmonique-Yamabe, qui est du a la non linearit\'e de l'equation.

4) Il y a l'exemple suivant:

$$ u_{\mu}(x)= \frac{\mu^{(n-2)/2}}{(\mu^2+|x|^2)^{(n-2)/2}},$$

a) quand $\mu\to +\infty $, on a l'inegalit\'e $ \sup_{B_k(0)} u_{\mu} \leq c \inf_{B_k(0)} u_{\mu} $, $ c>0 $,

b) quand $ \mu \to 0 $, on a l'in\'egalit\'e $ \sup_{B_1(0)} u_{\mu} \times \inf_{B_k(0)} u_{\mu} \geq c >0 $,

L'in\'egalit\'e ci-dessus, melange ce type d'in\'egalit\'es de Harnack. Fraction et produit.

\smallskip

Formelement, on ecrit:

$$ \inf \geq c(\sup)>0,$$

ou si on revient a l'ecriture initiale,

$$ \sup_M u \times \inf_K u \geq c(a,b, \sup_K u, K, M, g) >0.$$

On voit alors l'enroulement et la torsion en ce qui concerne les supercordes. On a alors les 2 relations suivantes \`a la fois: $ \leftarrow $ et $ \updownarrow $, en une relation. Comme ces deux relations sont melang\'ees, on a alors l'enroulement dans la torsion.

On a:

\be \sup_K u\times \inf_M u \leq c, \,{\rm ou}, \, \sup_M u\times \inf_K u \geq c >0:\,\,  inegalites \, tres \, fortes, \ee

\be (\sup_K u)^{\alpha} \times \inf_M u  \leq c, \,{\rm ou},\, (\sup_M u)^{1/\alpha} \times \inf_K u \geq c >0: \,\, inegalites \, fortes,\ee

\be \sup_K u \leq c(\inf_M u), \,{\rm ou}\, , \inf_K u \geq c(\sup_M u)>0: \,\,  inegalites \, faibles,\ee

\be \sup_M u \times \inf_K u \geq c(\sup_K u)>0:\,\,  inegalites \, tres \, faibles. \ee

a) Pour le $ (3) $ et $ (4)$, pour le (3), le sup controle l'inf. Pour le (4) le sup local controle l'$ \inf  \times \sup $, l'inf doit etre difficile a controler. Donc, $ (4) $ est plus faible que $ (3) $ qui est plus faible que $ (2) $ qui est plus faible que $ (1)$. 

L'in\'egalit\'e $ (4) $ est une in\'egalit\'e tres faible, on voit bien l'enroulement dans la torsion, m\'elange d'enroulement et de torsion. 

On l'a dit pr\'ec\'edemment, \`a l'origine, Kaluza-Klein c'est en $ 4+1 $ puis le model \`a evolu\'e et on considere $ 5+1$, $6+1$, puis il y a d'autres formulations en $ 9+1$. L'operateur en question est $ c_n\Delta +\tilde h $, $ c_n $ une constante universelle, celle de l'Eq de Yamabe, il est cens\'e etre coercif. On peut ne pas avoir de bord, cas d'une vari\'et\'e compacte sans bord, ou bien, on peut avoir un bord. Ceci est li\'e au Th. de la masse positive, qui implique la coercivit\'e de l'operateur et la presence ou non de bord. Puis,l'Eq. du modele est: $ c_n\Delta u+\tilde hu=Vu^{N-1}, h=S_g-|\nabla \Psi|^2, V=\frac{1}{2} m \Psi^2, N=\frac{2n}{n-2} $, avec $ S_g $, la courbure scalaire. Dans le cas $ \Psi\equiv constante $, on retrouve l'Eq. de Yamabe. On pose par exemple, pour etre coherent avec les Eq qu'on a consid\'er\'e: $ h=\frac{\tilde h}{c_n}$. Sur la vari\'et\'e Riemannienne $ (M,g) $, qui peut etre un produit, comme on l'a dit, sur le Th. de la masse positive.

\smallskip

b) On a:

$$ (1): 1-Torsion-Enroulement, $$

$$ (2): \alpha-Torsion-Enroulement, $$

$$ (3): Enroulement, $$

$$ (4): Enroulement-Torsion, \,enroulement \, dans\, la \,torsion $$

///////////////////////////////////////////////////////////////////////////////////////////////////////////////////////////////////////

L'in\'egalit\'e $ (4) $, enroulement dans la torsion, on a 3 parametres, peut etre vu comme liant 2 parametres, le (sup local) et l'($ \inf $ local $ \times \sup$ global). Le sup global est une barriere qui permet l'oscillation entre le sup local et l'inf local. L'in\'egalit\'e $(4)$ est une in\'egalit\'e tres faible, car le sup local controle l'inf $\times$ sup, l'inf doit etre difficile a controler. 

Si on a une id\'ee sur l'inf local et le sup global, il est possible d'avoir une id\'ee sur le sup local qui est intermediaire. In\'egalit\'e $ (4) $: in\'egalit\'e \`a 2 parametres indirecte. Intervalle de precisions de nombres. 

\smallskip

In\'egalit\'e $ (4) $: comme la fonction $ m\to c(m), m >0 $ est croissante de $ m >0 $ et que $ \sup_K u \geq u(z)>0, z\in K $, le point intermediaire peut etre $ z $, si $ C_{x,y}$ est une corde qui passe par $ z $(si $z\not \in K$, on agrandit $ K $), on ecrit: $ \sup_M u \times \inf_K u \geq c(\sup_K u)\geq c(u(z)) >0$, on voit alors que le point $ z $ est un point intermediaire et que $ u(z)$ est une valeur intermediaire entre $ \inf_K u$ et $ \sup_M u$. Le $ \sup_K u $ est alors, une valeur moyenne dans l'ensemble des valeures intermediaires. 

\smallskip
L'in\'egalit\'e $ (4)$ est une in\'egalit\'e \`a 2 parametres indirecte entre $ \inf_K u $ et $\sup_M u $. La valeur $\sup_K u$ est une valeur moyenne parmi les valeurs intermediaires. La fonction $ m\to c(m), m>0 $ mesure la dilatation de la corde. On a bien une {\bf in\'egalit\'e \`a 2 parametres indirecte ou dilat\'ee}.({\bf Dilatation } de la corde $ C_{x,y} $ en $ z $ $\forall z $ et si on prend les poids $ \inf_K u $ et $ \sup_M u $, on a dilatation de la corde $ C_{x\in \inf_K u, y\in \sup_M u} $, de poids $ \sup_K u $ et {\bf amplitutde de la dilatation: $ c(\sup_K u) $}). Ce ph\'enomene de dilatation ne se voit dans le cas d'une vari\'et\'e compacte sans bord, car la dilatation possede une borne inferieure, une limite ainsi que l'amplitude de la dilatation.

\smallskip

L'in\'egalit\'e $ (4)$: {\bf in\'egalit\'e \`a 2 parametres indirecte ou dilat\'ee}; avec amplitutde variable.

\smallskip
Dans un travail r\'ecent, si on suppose les D-branes, conditions de Dirichlet, on a l'in\'egalit\'e $ (3) $, moins tordue, qui met en relation 2 parametres directement: in\'egalit\'e \`a 2 parametres directe.

////////////////////////////////////////////////////////////////////////////////////////////////////////////////////////

{\bf Sur l'article au Bulletin des Sciences Math. 2006. et le print 2024, Harnack inequalities for equations of type prescribed scalar curvature. arXiv 2024.}

Voir l'article. Exotic origins of tensionless superstrings. Du point de vue de la physique, il y a 2 possibilit\'es: soit la convergence, soit l'eclatement, le blow-up. Le cas le plus interessant et important est le blow-up, \`a la limite, le blow-up, la corde, se tend et "eclate":

1) On voit bien que du point de vue de la physique l'eclatement \`a lieu \`a l'interieur de la vari\'et\'e $ M $: l'hypoth\`ese de la Physique est:

$$ \exists x_0 \in M, \exists (x_k)_k \subset M, x_k \to x_0, u_k(x_k)\to +\infty.$$

Or ceci, implique d'apres l'article du Bulletin. Sci. Math. 2006 que:

$$ \exists c>0, \, \forall K\subset \subset M, \, \forall k, \,\, \sup_M u_k \times \inf_K u_k \geq c >0, $$
Donc, du point de vue de la Physique, ou sous les hypotheses de la physique, l'article du Bull.Sci.Math. 2006. suffit, on a bien $ n+1, n=9 $. 

Le blow-up est une condition de la Physique. 

\underbar {Du point de vue Mathematique, il y a 3 parametres}, mais on voit comme on l'a dit dans le print, "Quelques remarques sur les vari\'et\'es, fonctions de Green, formule de Stokes". {\bf Il y a principalement 2 paramatres, le 3eme est un paramatre d'amplitude, mais qui est en fait une donn\'ee du point de vue de la Physique}, donc, {\bf il s'auto-\'elimine pour faire apparaitre 2 parametres}.

Ici, dans le point de vue de la physique, le $\sup $ local s'auto-elimine, car il est minor\'e, dans le sens de la donn\'ee.

\smallskip

1-1) Point de vue Mathematique: 

Hypoth\`ese: $( \sup_K u_k \geq m >0 )$: 1 parametre $ + \Rightarrow (\sup_M u_k \times \inf_K u_k \geq c(m) >0)$: 2 parametres = 3 parametres au total: $ (\sup_K u_k, \sup_M u_k, \inf_K u_k) $.

\smallskip

1-2) Point de vue de la Physique: 

Donn\'ee: $ (\sup_K u_k\to +\infty) \Rightarrow $: $ \sup_M u_k\times \inf_K u_k \geq c >0 $ : = 2 parametres $(\sup_M u_k,\inf_K u_k) $.

\smallskip

'Point de vue Mathematique' $ \not = $ 'Point de vue de la Physique'.

\bigskip

2) Dans le print, arXiv 2024. On a une alternative, qui confirme ce qui a \'et\'e dit dans l'introduction de l'article. Exotic origins of tensionless superstrings. Du point de vue de la physique, il y a 2 possibilit\'es: soit la convergence, soit une in\'egalit\'e de Harnack.

\smallskip

{\bf Ici, on a 2 parametres} et la condition blow-up ou eclatement est dans l'in\'egalit\'e de Harnack. On a au plus 2 parmatres et au moins 1 parametre.

\smallskip

Du point de vue de la physique, le blow-up est inclu dans le nombre de paramatres , soit, 2. Ici, il y a 2 parametres au maximum et la condition de la Physique, blow-up, est incluse directement dans l'in\'egalit\'e.

Ici, on ne peut pas eliminer la $ \sup $ local, car, on cherche \`a l'estimer, \`a le majorer (sens de l'in\'egalit\'e: la majoration). On estime et majore le $\sup $ local par une quantit\'e dependant de l'$\inf $ global.

\smallskip

2-1) Point de vue Mathematique:

\smallskip

\underbar {In\'egalit\'e de Harnack \`a 2 parametres: point de vue math\'ematique.}

\smallskip

Hypoth\`ese: $ (\sup_K u_k \to +\infty)  \Rightarrow $ relation entre $ (\sup_K u_k, \inf_M u_k) $ = 2 parametres au total = $(\sup_K u_k,\inf_M u_k) $.

\smallskip

2-2) Point de vue de la Physique:

\underbar {Le blow-up est dans l'in\'egalit\'e \`a 2 parametres: point de vue physique.}

\smallskip

Donn\'ee: $ (\sup_K u_k \to +\infty) \Rightarrow $ relation entre $ (\sup_K u_k, \inf_M u_k) $ = 2 parametres au total = $ (\sup_K u_k,\inf_M u_k) $.

Finalement:

a) Ici, {\bf 'le point de vue Mathematique' et 'le point de vue de la Physique' sont les memes}.

\smallskip

b) Point de vue de la Physique est \underbar {coh\'erent} avec le point de vue Mathematique.

\smallskip

c) 'Point de vue Mathematique'= 'Point de vue de la Physique'= in\'egalit\'e \`a 2 parametres.

\smallskip

3) 

3-1) On peut dire que le point de vue de la Physique etait deja resolu dans l'article au Bull.Sci.Math. 2006.

3-2) Dans le print, arXiv 2024. L'aletrnative, confirme, l'article de la physique, "Exotic origins of tensionless superstrings" et Les points de vue Mathematique et Physique sont les memes. (En plus de l'in\'egalit\'e de Harnack, nouvelle et de nouveau type et dans le sens de la majoration).

\smallskip

On peut dire que la coh\'erence entre le point de vue Mathematique et le point de vue de la Physique, est resolu dans le preprint 2024.

\smallskip

{\bf Concernant l'article aux Bull.Sci.Math.2006:}

\smallskip

Du point de vue EDP et in\'egalit\'es de Harnack: il y a la contrainte $ x_0\in M, u_i(x_0) \geq m >0 $: ceci est coherant avec ce qu'on appelle in\'egalit\'e de Harnack.

Pour les fonctions harmoniques, dans le Gilbarg-Trudinger, ils disent "non-negative harmonic functions": donc on a la condition $ u\geq 0$ implique cequ'on appelle in\'egalit\'e de Harnack pour les fonctions harmoniques positive ou nulles. Il y a bien une contrainte: la positivit\'e de la fonction harmonique. On a bien 2 parametres $ \sup_K u $ et $ \inf_K u $ avec la contrainte: $ u\geq 0 $.

Il y a aussi l'exemple de Brezis-Merle, corollaire 7, avec volume born\'e uniform. et potentiel $ C\geq V(x) \geq 0$. Ils ont une contrainte sur le volume et ils obtiennent une in\'egalit\'e de Harnack \`a 2 parametres.

Dans la'article du Bull.Sci.Math. 2006: on a une contrainte: $ u_i(x_0) \geq m >0 $, avec 2 parametres: $ \sup_M u_i $ et $ \inf_K u_i $. 

Donc, c'est coherent avec cequ'on appelle in\'egalit\'es de Harnack. On a une condition qui implique une in\'egalit\'e \`a 2 paramatres.(exemple, c'est ce qu'on a dit sur les fonctions harmoniques: positivit\'e des fonctions). Cette condition est compatible avec le phenomene physique de blow-up.

In\'egalit\'e de Harnack de, Bull.Sci.Math.2006: 

a) elle caracterise, l'enroulement et la torsion ou l'enroulment et la distortion.

b) Elle montre la stabilit\'e du systeme.

c) C'est vrai en toute dimension $ n\geq 3 $: qui inclut: la relativit\'e, Kaluza-Klein, la th\'eorie des cordes et la th\'eorie des supercordes. $ n=3, 4, 5, 6, 9 $. 

C'est bien une in\'egalit\'e de Harnack \`a 2 parametres. C'est coherent avec ce qu'on appelle in\'egalit\'es de Harnack.

Tout ceci pour dire que ces in\'egalit\'es, selon le point de vue, font apparaitre 2 parametres ou 3 parametres. L'exemple de fonctions harmoniques, positive ou nulles (contrainte), avec 2 parametres, est l'exemple typique. $ \Rightarrow $ l'in\'egalit\'e de Harnack \`a 2 parametres, avec une contrainte de l'article du, Bull.Sci.Math.2006. C'est coherent avec ce qui est appel\'e in\'egalit\'e de Harnack.

\smallskip

Dans la preuve de ce resultat (du theoreme 2 de cet article), on a utilis\'e aussi, comme dans le cas negatif(critique, dimension 2 et supercritique) de l'article de 2003, l'integration par parties, localement dans une boule geodesique petite: sous vari\'et\'e Riemannienne compacte \'a bord $ C^{\infty}$ plong\'ee dans une vari\'et\'e Riemannienne, plongement obtenu par la carte exponentielle.

\smallskip

Il y a aussi, le point de vue "metric geometry", de Cheeger, Anderson, Peters, Gromov: Eux ils donnent des conditions sur le volume, la courbure de Ricci, le diamtere, donc, la metrique, ...etc... pour montrer un resultat de compacit\'e. Ici, on donne une condition sur la metrique conforme: une condition de saut, ou de concentration en un point, pour avoir une in\'egalit\'e sur le $ \sup, \inf $, de la metrique conforme, donc, un resultat sur la metrique conforme, une formule de la moyenne, stabilit\'e en norme essentielles du $ \sup $ et $ \inf $: "mean value theorem", "oscillation", "conformal metrics", "stability". 

Ici, c'est une equation plus generale, qui inclut, les metriques conformes. Donc, ca inclut la "metric geometry", si on considere en particulier, l'eq. de la courbure scalaire prescrite ou de Yamabe.

\smallskip

{\bf Remarque:}

\smallskip

On a vu dans l'article. Bull.Sci.Math.2006, qu'on peut avoir, selon le point de vue une in\'egalit\'e de Harnack \`a 3 parametres ou 2 parametres. 

Si on suppose la condition de Brezis-Merle ou Brezis-Li-Shafrir: $ \inf_M u >0 $, on a, \`a partir de l'in\'egalit\'e \`a 3 parametres:

$$ \sup_M u\times \inf_K u \geq c(\sup_K u) >0, $$
 avec $ m \to c(m) >0 $ une fonction croissante de $ m >0 $, $ \sup_K u\geq \inf_K u \geq \inf_M u >0 $, (positivit\'e stricte):

$$ \sup_M u\times \inf_K u \geq c(\inf_K u) >0, $$

On a bien une in\'egalit\'e de Harnack \`a 2 parametres entre $ \sup_M u $ et $ \inf_K u $, en admetttant la condition de Brezis-Merle ou Brezis-Li-Shafrir: $\inf_M u > 0 $. Qui est possible en considerant le probleme sur $ (M, g) $ avec $ u >0 $ sur $ M $, puis on prend un sous-ensemble ouvert connexe, $ \tilde M $, alors on aura : $ \inf_{\tilde M} u >0 $. (Par exemple autour d'un point, comme le font Li-Zhang).

\smallskip

Comme le probleme se pose autour de chaque point $ x_0 \in M $, comme le font Li-Zhang. On prenant un ouvert relativement compact de $ M $, $ \tilde M \subset \subset M $. On a toujours $ \inf_{\tilde M} u >0 $. On a:

$$ \forall K \subset \subset \tilde M, \,\,\sup_{\tilde M} u\times \inf_K u \geq c(\inf_K u) >0, $$

1) Ici, comme l'operateur est coercif, par le principe du maximum, on a toujour: pour tout compact $ K $ de $ M $: $ \inf_K u >0 $. Il n'y a pas besoin de supposer la condition, de Brezis-Merle ou de Brezis-Li-Shafrir.

Par le principe du maximum, on a toujours: $ \inf_K u >0 $ (aussi, ici, on a directement $ \inf_K u >0 $, car $ u >0 $ est strict.positive sur $ M $ et continue sur $ \bar M $ et $ K $ est compact). Donc, on a l'inegalit\'e de Harnack sans la condition de Brezis-Merle ou de Brezis-Li-Shafrir.

\smallskip

2) Ce qui parait ambigu, est que la fonction $ m\to c(m)$ est croissante de $ m >0 $, on ne connait pas sa forme explicite. Par exemple, on ne peut pas avoir une fonction linaire, $ c(m) =c_0 m $, car l'exemple qu'on a donn\'e avant: $ x\to (\mu/\mu^2+|x|^2)^{(n-2)/2}, \mu \to +\infty  $, implique qu'on a : $ \sup_M u \geq c_0 \frac{\sup_K u}{\inf_K u} \to c_0 >0$, alors que $ \sup_M u \to 0, \mu \to +\infty $ (et aussi le cas $ \mu \to 0 $, on ne peut aps avoir $ c(m)=c_0 m^k, k >1 $, une puissance). Donc, parfois, on ne peut pas avoir $ c(m)$ lineaire. Ce qui veut dire aussi, que l'in\'egalit\'e qu'on a obtenue, est une vraie in\'egalit\'e de Harnack. Par exemple, on n'a pas une in\'egalit\'e triviale du type $ \sup_m u \geq c_0 >0 $.

\smallskip

3) on a:

$$ \forall K \subset \subset  M, \,\,\sup_M u\times \inf_K u \geq c(\inf_K u) >0, $$

ou bien:

$$ \forall K \subset \subset M, \,\,\sup_M u\geq \frac{c(\inf_K u)}{\inf_K u} >0, $$ 

avec, $ m\to c(m) >0 $ fonction croissante de $ m >0$.

Cette in\'egalit\'e \`a 2 parametres se deduit du th 2 de l'article du Bull.Sci.Math.2006. en prenant $ x_0 \in K $ et $ m=\inf_K u >0 $. Ce qui etait ambigu: il ne faut pas prouver l'in\'egalit\'e en prenant $ v >0$, $ \inf_K v >0$, mais en considerant un point $ x_0$ et on la prouve en prenant $ v(x_0) \geq m >0$. Si on suppose $ \inf_K v =m >0$, il n'y a rien a prouver. Mais en prenant $ v(x_0) \geq m >0$, on prouve l'in\'egalit\'e en prenant $ m=\inf_K u >0 $, puis on l'applique \`a $ v=u$, $ v $ particuliere, $ v=u$, comme $ x_0 \in K $, on a bien $ v(x_0)=u(x_0) \geq \inf_K u >0 $.

\smallskip

Le th 2 de l'article du Bull.Sci.Math.2006. ecrit, permet de deduire l'in\'egalit\'e \`a 2 parametres, en prenant, $ x_0 \in K $, $ m=\inf_K u >0 $ et $ u_i=u $.

\smallskip

Comme on l'a dit, la fonction $ c(m)$ ne peut pas etre triviale, en prenant, les fonctions particlieres ci-dessus. $ x\to (\mu/\mu^2+|x|^2)^{(n-2)/2}, \mu \to +\infty  $ ou $ \mu \to 0 $.

\smallskip

{\bf Important:} Cette in\'egalit\'e est obtenue, {\bf sans concentration de la mesure}. On n'a pas de notion de concentration de la mesure \`a l'interieur et au bord. On ne sait pas si la concentration de la mesure est vraie ou  pas.

\smallskip

Notions utilis\'ees:

\smallskip

Cette in\'egalit\'e est obtenue avec la notion de fonction de Green d'un op\'erateur coercif dans $ H^1(M) $ (ce qui inclut, la condtion de Dirichlet, de Neumann, mixte, conditions au bord lin\'eaires, non lin\'eaires (probl\`eme de Cherrier), c'est \`a dire, avec n'importe quelle condition au bord telle que l'op\'erateur soit coercif): ici ce qui compliqu\'e est la construction de la fonction de Green avec condition de Dirichlet. On a aussi la notion d'int\'egrale Riemannienne de Lebesgue, espace de Lebesgue $ L^p$, espaces de Sobolev et injections de Sobolev, cas de l'exposant critique et l'int\'egration pas parties; formule de Stokes-Green-Gauss-Riemann, proc\'ed\'e d'it\'eration. Raisonnement par l'absurde: \'energies tendant vers $ 0 $.

\smallskip

On a:

$$ \forall K {\rm \,\, compact\,\, de \,\,}  M, x_0\in K, \,\,\sup_M u\times \inf_K u \geq c(\inf_K u) >0, $$

avec $ m \to c(m) $ croissante de $ m >0 $.

\smallskip

On ecrit cela pour mettre en evidence l'enroulement et la distrotion, la rigidit\'e, la stabilit\'e, l'in\'egalit\'e $ (\sup, \inf) $ en dimension $ n\geq 3 $, {\bf in\'egalit\'e \`a 2 param\`etres}, avec la notion de pincement (pinching). 

\smallskip

L'operateur coercif consid\'er\'e est $ \Delta +h $, ce qui inclut l'op\'erateur conforme, l'eq. de Yamabe, de la courbure scalaire prescrite, l'eq. d'Einstein-Lichnerowicz, l'eq. de Schrodinger ($ n\geq 3 $, Schrodinger c'est $ n\geq 2 $, il est possible de completer le cas $ n=2 $). Ici, l'exposant est critique de Sobolev et la dimension de la vari\'et\'e est $ n\geq 3 $, ce qui r\'esout le probleme du $ (\sup, \inf) $ dans le cas coercif, dans tous les cas (pour ce qui concerne les equations): Yamabe, Courbure scalaire prescrite, Einstein-Lichnerowicz, Schrodinger ($n\geq 3$, il est possible de completer Schrodinger $ n=2 $ et exposant sous-critique pour $ n=2 $).

////////////////////////////////////////////////////////////////////////////////////////////////////////////////

\smallskip

44) les Equations, apparition et provenence: voir, Aubin,T. Lions, P-L. Rabinowitz, P-H.

\smallskip
a) Equation de Yamabe:  precisement elle s'appelle Eq de Yamabe. En G\'eometrie et Physique et Astronomie. Kaluza-Klein, cordes, supercordes et D-branes, en dimension 4, Yang-Mills. {\bf Dynamique relativiste}., energie, volume, potentiel...convexit\'e(minor\'e ou cas positif courbure scalaire prescrite positive+major\'e, par exemple Eq.Yamabe dim 4,5,6), concavit\'e (major\'e, ou cas n\'egatif, ou courbure scalaire prescrite negative+major\'e)...
\smallskip

b) Equation de la courbure scalaire  prescrite: on peut la nommer "Eq de Schrodinger particuliere". En G\'eometrie et en Physique, en dimension 4, Yang-Mills., energie, volume, potentiel...convexit\'e(minor\'e ou cas positif courbure scalaire positive prescrite+major\'e, Par exemple, Eq. dim 4, courbure prescrite), concavit\'e (major\'e)...

\smallskip

c) Equation de type courbure scalaire, comme on en a parl\'e. En cosmologie, Astronomie, Kaluza-Klein, cordes, supercordes, D-branes. En dimension 4, Yang-Mills., energie, potentiel...

\smallskip

d) Equation de Lane-Emden-Fowler: $ \Delta u= V(x) u^p; u >0, p >1 $, en astronomie.

\smallskip

e) Equation de Schrodinger:(stationnaire: du type Yamabe; du type courbure scalaire prescrite): $ \Delta u +b(x)u=f(x,u), u >0 $ avec $ f $ avec des combinaisons de termes avec exposant critique et sous-critique, $ n\geq 3 $, par exemple, $  f(x,u)=V(x)u^{N-1}+W(x)u^{\alpha}, 1<\alpha <N-1=\frac{n+2}{n-2} $. En physique, en optique non-lin\'eaire (du type Yamabe), {\bf Eq de Schrodinger: dynamique d'une particule non relativiste}. Energie, potentiel...

\smallskip

f) Eq de Liouville ou du type Liouville (courbure de Gauss): En G\'eometrie et Physique et Chimie et Astronomie. Gaz, Cordes, D-branes, en dimension 2., energie, volume, potentiel...convexit\'e(minor\'e ou cas positif courbure prescrite+major\'e), concavit\'e(major\'e)...

\smallskip

Par exemple:  Du point de vue de la g\'eometrie, par exemple, Dans le cas positif: courbure scalaire prescrite positive, notion de convexit\'e, et aussi le fait suivant: $ \sup <+\infty $ si $ \inf >0$: cela veut dire que si la vari\'et\'e est uniformement convexe alors elle loc. uniformement concave...

\smallskip

////////////////////////////////////////////////////////////////////////////////////////////////////////////////////////////////

\smallskip

45) Sur l'article de C.C.Chen. C.S.Lin. Communications on pure and applied math, 1997. Ils supposent les potentiels $ V $, $ C^1$:

Ceci est du au fait suivant: Ils cherchent a avoir comme donn\'ee, dans la convergence des fonctions blow-up, $ R_i^{n-2} ||v_i-U||_{C^2(B_{2R_i}(0))} \to 0 $. En faisant le rescaling, ou changement d'echelle, $ x=R_i y $, $ y\in B_3(0) $, il suffit et il faut que, la difference, $ w_i(y)=v_i(R_iy)-U(R_iy) $  tend vers $ 0 $ et appliquer les estimations de Schauder, aux fonctions, $ w_i $, $ \Delta (w_i-U) \to 0$, or, $ -\Delta (w_i-U)=R_i^2(V_i(y_i+R_iy)-V_i(y_i))(v_i)^{N-1} +... \to 0 $. C'est a dire qu'il faut que $ R_i^2(V_i(y_i+R_iy)-V_i(y_i)) \to 0 $, or cette difference n'est que la deriv\'ee de $ V_i $ au point blow-up: $ \frac{V_i(y_i+h)-V_i(y_i)}{h} $ avec $ h=\frac{1}{R_i^2} $. Donc, il faut que $ V_i $ soit derivable en $ y_i $ et $ V_i'(y_i)=0 $ et par le theoreme des accroissement fini, $ V_i'(c_h) \to 0 $ avec $ 0 < c_h< h=\frac{1}{R_i^2} $ et la deriv\'ee en $ 0 $ doit etre continue. Donc, $ V_i $ doit etre $ C^1$ en $ y_i $, comme le point $ y_i$ est quelconque, il faut que $ V_i $ et sa limite $ V $ soient $ C^1$. Il faut que pour toute suite $ (V_i) $ convergeant dans $ C^{\alpha} $ vers $ V $, il faut que $ V $ soit $ C^1 $, pour que ceci soit possible, il faut que $ \alpha =1 $, il faut que $ (V_i) $ soient $ C^1 $ et convergent dans $ C^1$ vers $ V $, $ C^1$. Donc, il faut que $ (V_i) $ soient $ C^1 $ et convergent vers $ V \in C^1 $ dans $ C^1$.

\smallskip

a) On a la condition de depart pour les $ V_i $.(Par exemple: $ 1+a_i|x|^{\alpha}\to 1, a_i \to 0 $).

\smallskip

b) On a la condition: \underbar {pour toute suite} $ V_i $ dans $ C^{\alpha}$ qui converge dans $ C^{\alpha}$ vers $ V $, $ V $ doit etre $ C^1$.

Pour que ce soit possible il faut que $\alpha=1$, $ V_i $ soit $ C^1 $ et converge dans $ C^1 $ vers $ V $. On a alors $ V $ est $ C^1 $.

///////////////////////////////////////////////////////////////////////////////////////////////////////////////////////

\smallskip

46) {\bf Volume=Energie en dimension 2. Eq. de la courbure scalaire prescrite:} Sur le volume et energie en dimension 2: dans le cas d'une surface riemennienne compacte sans bord, lorsque le potentiel est entre 2 nombres strict.positifs, on utlise la formule de Gauss-Bonnet. Lorsqu'on se place sur un domaine \`a bord (surface \`a bord particuliere), avec potentiel Lispchitzien(+poids Holderien ou domaine etoil\'e, par exemple dans l'article "A compactness result for an equation with Holderian condition") ou potentiel C1(+poids Holderien ou Lipschitzien) et avec des poids holderiens, on a la borne uniforme du volume=energie sans utliser la formule de Gauss-Bonnet.

Avec condition de Dirichlet. la condition de Dirichlet est aussi une deuxieme contrainte.

Volume global ou total=Energie globale ou totale.

\smallskip

//////////////////////////////////////////////////////////////////////////////////////////////////////////////////////

\smallskip

47) {\bf Volume local=Energie locale en dimension 4. Eq. de la courbure scalaire prescrite, cas plat:} En ce qui concerne l'eq. en dimension 4, $ -\Delta u=Vu^3, u >0, \Delta =\nabla^i(\nabla_i)$ sur un domaine $ \Omega \subset {\mathbb R}^4 $. On a une condition sur $ u_i: u_i\geq  k(a)\cdot A_i, k(a)=\dfrac{8e^2{\sqrt 2}}{3a {\sqrt a}} >0$, $ A_i$ la constante de Lipschitz de $ V=V_i $, pour avoir l'estimation optimale $ \sup_K u_i \times \inf_{\Omega} u_i \leq c $, en particulier le volume=energie est localement uniformement born\'e. Brezis-Merle, donnent le meme type de conditions pour borner localement le volume=energie en dimension 2 et en deduisent la compacit\'e locale.

Pour la dimension 4 (ce fait est illustr\'e par l'exemple suivant), il y a l'exemple usuel de volume localement born\'e sans compacit\'e locale: $ n=4, u_{\mu}(x)=\frac{\mu}{\mu^2+|x|^2}, x\in B_1(0), \mu \to 0 $.

a) Lorsque $ A_i \to 0 $: volume local=energie locale born\'e(e).

b) Lorsque $ A_i \to A >0 $: compacit\'e locale.

En dimension 2, on avait une condition de Dirichlet qui etait aussi une deuxieme contrainte. Ici, en dimension 4, on a une deuxieme contrainte qui permet d'avoir la borne uniforme locale du volume=energie. Brezis-Merle aussi, imposent une deuxieme contrainte pour majorer localement le volume=energie, et obtiennent la compacit\'e locale.

Volume local=Energie locale.

\smallskip

//////////////////////////////////////////////////////////////////////////////////////////////////////////////////////

\smallskip

48){\bf Energie et volume. Eq. du type Liouville. en dimension 2 avec condition de Dirichlet.:} il s'agit de l'Eq. $ -\Delta u+ \epsilon (x\cdot \nabla u)= V e^u $ avec condition de Dirichlet. $ \Delta =\partial_{11}+\partial_{22}, 0 \leq \epsilon \leq 1$, domaine etoil\'e par rapport \`a l'origine. On a une deuxieme contrainte, qui est la condition de Dirichlet.

/////////////////////////////////////////////////////////////////////////////////////////////////////////////////////////////

\smallskip

49) {\bf Energie locale. Eq. de Schrodinger:} En ce qui concerne l'eq. de Schrodinger $ -\Delta u=Vu^{N-1}+Wu^{\alpha}, u >0, n\geq 3, \frac{n}{n-2} \leq \alpha < N=\frac{2n}{n-2} $ ou l'Eq. de Schrodinger sur les vari\'et\'es, $ -\Delta u-\lambda u = n(n-2)u^{N-1}, n\geq 3, N=\frac{2n}{n-2} $. $ \Delta =\nabla^i(\nabla_i)$. On a les in\'egalit\'es optimales; $ \sup_K u \times \inf_{\Omega} u \leq c $ et $\sup_K u \times \inf_M u \leq c $. On en d\'eduit la borne uniforme locale des energies.

{\bf Eq. de Schrodinger:} dynamique d'une particule non relativiste.

Energie locale. 

\smallskip

La solution $ u >0 $ n'est pas solution de l'eq. de Yamabe, mais on peut considerer les metriques conformes $ \tilde g= u^{4/(n-2)} g $, alors on a, (si l'operateur $ \Delta -\lambda $ est coercif), le volume conforme est born\'e. On peut parler de volume conforme, en considerant les metriques conformes, sans que $ u>0 $ soit solution de l'equation de Yamabe. Mais $ u >0 $ est solution de l'equation de la courbure scalaire prescrite.

\smallskip

Ici, on g\'eom\'etrise l'equation ou le probleme. On prescrit une quantit\'e mathematique ou physique, qui donne une equation de Schrodinger. Puis, on la rend g\'eometrique, en considerant le potentiel $ V = [n(n-2)+(\frac{n-2}{4(n-1)} R_g + \lambda )u^{2-N} ]$, alors $ V \geq n(n-2) >0 $, on obtient l'equation de la courbure scalaire prescrite $ -\Delta u+ \frac{(n-2)}{4(n-1)} R_g u= Vu^{N-1}, u > 0 $, pour laquelle on a l'in\'egalit\'e de Harnack $ \sup \times \inf $ optimale et le volume local = energie locale, uniform\'ement born\'es: on ainsi g\'eometriser le probleme et la solution. Avec la metrique conforme $ \tilde g=u^{4/(n-2)} g, u >0, n\geq 3 $.

\smallskip

On part d'un probleme ou d'une solution, mathematique ou physique (calculs) et on obtient un probleme ou une solution g\'eom\'etrique (equation de la courbure scalaire prescrite, mathematique ou physique): g\'eom\'etrisation du probleme ou de la solution.

\smallskip

/////////////////////////////////////////////////////////

\smallskip

{\bf Remarque:} on a dit qu'on n'obtient pas forc\'ement l'eq.de Yamabe (astronomie), mais on a l'eq. de la courbure scalaire prescrite (physique, Schrodinger). Comme le potentiel s'ecrit: $ V=[n(n-2)+(\frac{n-2}{4(n-1)}R_g +\lambda) u^{2-N}] $, le terme $ u^{2-N}$ peut tendre vers 0 aux poins blow-up, et la fonction $ u $ converge faiblement vers $ u_0 $, il se peut que $ u_0 $ soit constante. Si on se place en des points blow-up $ x_0 $, il tend vers 0. (On se place autour de points blow-up), il se peut qu'il y ait un nombre infini d'autres points $ y_j $ autour de $ x_0 $, $ V $ oscille autour d'une constante. Donc, $ V \approx n(n-2)+k, k>0 $, autour de $ x_0 $, au moins, presque partout. C'est presque partout le potentiel constant et donc l'equation de Yamabe.

A la limite, on peut avoir comme solutions, les constantes, ou fonction constante par morceaux, donc, Yamabe par morceaux. Donc, on a bien l'eq. de la courbure scalaire prescrite (physique, Schrodinger), mais aussi, il se peut que ce soit l'eq. de Yamabe (astronomie) (\`a la limite, fonction constante, ou constante par morceaux). Obtenir l'eq.de Yamabe est un ph\'enomene possible. 

\smallskip

Hebey-Vaugon-Druet-Robert ont montr\'e que quand on a blow-up, il se peut que la fonction limite soit nulle, donc, constante nulle. ici $ u $ converge, au moins presque partout vers une fonction, il se peut que la fonction limite soit constante $ >0 $. Ceci est en relation avec la Th\'eorie physique du tout. On a au moins la convergence faible de $ u>0 $ vers une fonction, il se peut que la fonction limite soit constante $ >0 $ et donc le potentiel converge presque partout vers une constante $ V  \to n(n-2)+k, k >0 $, $ V $ s'ecrit comme $ V=n(n-2)+k+o(1), k >0 $ et $ o(1) $ est negligeable devant $ n(n-2)+k, k >0 $.(a la limite eq.Yamabe en astronomie, c'est possible). 

On pouvait dire cela des le debut, en invoquant la Th\'eorie physique du tout, le potentiel $ V $ d\'efinit une eq. de la physique g\'enerale, dite de Schrodinger, il peut etre constant (au moins \`a la limite). Ce qui inclut Yamabe (astronomie). Ici, on explique comment avec la notion de convergence presque partout, convergence faible, et points blow-up, on peut avoir $ V\to n(n-2)+k, k >0 $.

\smallskip

Ceci est dit pour faire remarquer qu'a la limite on obtient Yamabe.

\smallskip

Dans Einstein-Lichnerowicz, la seule maniere de combiner, l'eq. de la courbure scalaire prescrite, eq. de Schrodinger non lin\'eaire, et la relativit\'e g\'en\'erale, est de considerer l'eq. de Yamabe. Ici, \`a la limite, il est possible de de recuprer l'eq. de Yamabe. Pour ce qui nous concerne, on a l'eq. de la courbure scalaire prescrite.

\smallskip

Dans Einstein-Lichnerowicz, il y a eu separation de la variable temporelle et spaciale, puis il y a eu utilisation de metriques conformes. On a: \`a un instant $ t =0 $, par approximation, on obtient l'eq. de la courbure scalaire prescrite: car on part de $ -\Delta u+ (\frac{n-2}{4(n-1)}R_g-|\nabla \Psi|^2)u= \frac{1}{2} m \Psi^2 u^{N-1}, u >0, N=\frac{2n}{n-2}, n\geq 3 $, comme on consider\'e $ \Psi $ voisine d'une constante (impulsion \`a l'instant $ t=0 $, puis le potentiel se stabilise autour d'une constante), on peut supposer $ |\nabla \Psi | \approx 0 $, et on obtient l'eq. de la courbure scalaire prescrite, avec un potentiel $ V $ voisin d'une constante $ >0 $.(On en a parl\'e pour la dimension 4, d'o\`u les conditions de platitudes). L'eq, devient:$ -\Delta u+ \frac{n-2}{4(n-1)}R_g u= \frac{1}{2} m \Psi^2 u^{N-1}, u >0, n\geq 3, N=\frac{2n}{n-2} $, eq. de la courbure scalaire prescrite. On en a parl\'e pour la dimension 4. Ceci reste valable en toute dimension $ n\geq 3 $: $ \Psi $ voisine d'une $ constante >0  \Rightarrow |\nabla \Psi | \approx 0 \Rightarrow $ l'eq. de la courbure scalaire prescrite $ \Rightarrow $ eq. de l'invariance conforme $ \Rightarrow $ astronomie et g\'eom\'etrisation de l'astronomie.

\smallskip

On peut dire que l'eq. de la courbure scalaire prescrite est aussi l'eq. de la supersymetrie, comme l'eq. de Yamabe. L'eq. de Yamabe est un cas particulier de l'eq. de la courbure scalaire prescrite.

\smallskip

On peut dire qu'avec ce proc\'ed\'e d'approximation: que l'eq. de la courbure scalaire prescrite est une eq.de l'astronomie. Au d\'epart, on s'est focalis\'e  sur l'eq. de Yamabe (qui est obtenue aussi par approximation), on peut considrer de la meme maniere , l'eq. de la courbure scalaire prescrite. C'est une equation plus g\'en\'erale que l'eq. de Yamabe. C'est l'eq. de l'invariance conforme. (Laplacien conforme, courbure scalaire, tenseur de Weyl, eq. de la courbure scalaire prescrite: invariants conformes). Masse et Masse positive: courbure, scalaire, de Ricci, de Weyl.

\smallskip

Il y a le formalisme de Penrose: g\'eom\'etrie conforme+compactification+courbure+masse.

\smallskip

//////////////////////////////////////////////////////

\smallskip

On peut faire la meme chose avec l'equation qui contient $ V, W $ dans $ {\mathbb R}^n, n\geq 3 $, termes critique et sous critique: $ -\Delta u=Vu^{N-1}+W u^{\alpha}, u, V, W >0, n\geq 3, \frac{n}{n-2} \leq \alpha < N=\frac{2n}{n-2} $. On g\'eom\'etrise l'equation en l'equation de la courbure scalaire prescrite sur un ouvert de $ {\mathbb R}^n, n\geq 3 $, on obtient: $ -\Delta u= \tilde V u^{N-1}, u >0, \tilde V \geq a >0, n\geq 3 $, $ \tilde V= V+W u^{\alpha +1-N} $, la r\'egularit\'e de $ \tilde V $ est Lipschitz. La r\'egularit\'e de $ \tilde V $ est celle de $ V, W $: au moins Lipschitz. 

\smallskip

On a alors, l'in\'egalit\'e optimale $ \sup \times \inf $ et volume local = energie locale, uniform\'ement born\'es pour les metriques conformes $ \tilde g=u^{4/(n-2)} \delta $.

\smallskip

Exemple 1. Ce type de proc\'edure de g\'eometrisation, a \'et\'e evoqu\'e dans le livre d'Aubin: si on resout le probleme de Yamabe ou de courbure scalaire prescrite dans le cas sous-critique, alors on a r\'esolu le probleme de la courbure scalaire prescrite dans le cas critique. Avec une courbure scalaire prescrite particuliere.

\smallskip

Exemple 2. Quand on a r\'esolu par le degr\'e topologique, voir ce qu'on a ecrit avant, en dimension 2, l'equation: $ -\Delta u= V(e^u-1) $, avec condition de Dirichlet au bord, puis on a obtenu, une solution de l'equation de Liouville avec potentiel $ W=Ve^w $ o\`u $ w >0 $ est solution de $ -\Delta w= V $, avec condition de Dirichlet au bord. On a ainsi g\'eom\'etriser le probleme et la solution: Equation de Liouville. Voir, aussi, dans le livre d'Aubin, sur le probl\`eme de Nirenberg: on r\'esout une equation similaire \`a celle de Nirenberg, puis, on se ramene, par subsitution ou addition ou multiplication, (alg\`ebre de fonctions: on effectue des op\'erations alg\`ebriques sur les fonctions), \`a une solution du probl\`eme de Nirenberg.

\smallskip

////////////////////////////////////

\smallskip

Sur l'existence de solutions par le degr\'e topologique: ici, dans le cas positif:

\smallskip

(On peut faire la meme chose que pour l'eq.de la courbure scal.prescrite: dans le cas d'une vari\'et\'e compacte sans bord $ M $, loc.conf.plat. avec compacit\'e  pour $ V\equiv 1$, pour l'eq. $ -\Delta u+R_g u= u^{N-1}+W(x)u^{\alpha}, \dfrac{n+2}{n-2} <\alpha <N-1 $, l'energie born\'ee implique qu'il n'y a que des blow-up isol\'e et isol\'es simples, en nombre fini, en combinant avec l'article; in\'egalit\'e de Harnack et phenomene de concentration, on a la compacit\'e par le th. de la masse positive dans le cas plat et la formule de Pohozaev, voir les article de Y.Y.Li, et, le degr\'e topologique. Mais ici, on prefere parler d'energie locale, volume lcoal, pour cette eq.de Schrodinger sur une ouvert born\'e de $ {\mathbb R}^n $. Il y a des applications, dans la compacit\'e et le degr\'e topologique, voir les articles de Y.Y.Li, par le th.masse positive dans le cas plat et formule de Pohozaev, pour une vari\'et\'e compacte sans bord, loc.conf.plate. On a pr\'ef\'er\'e se limiter \`a la notion de $ \sup \times \inf $ et d'energie locale, volume local.)

//////////////////////////////////////////////////////////////////////////////////////

Si on veut etre plus precis et expliquer mieux: on commence par poser la question suivante: Trouver $ (u,V), u >0, V>0 $ solutions de:

$$ -\Delta u+\frac{(n-2)}{4(n-1)} R_g u= Vu^{N-1}, u >0, N=\frac{2n}{n-2}, n\geq 3, \quad (*) $$

tels que ($ K $ d\'esigne un compact quelconque de la vari\'et\'e Riemannienne $ (M,g) $ non necessairement compacte sans bord):

$$ \sup_K u\times \inf_M u \leq c,$$

et,

$$ \int_K u^{2n/(n-2)} dV_g \leq c, $$

Ici on est dans le cas positif. Ceci inclut le cas non loc.conf.plat.

\smallskip

////////////////////////////////////////////////

\smallskip

a) Dans le cas negatif on a:

$$ -\Delta u+\frac{(n-2)}{4(n-1)} R_g u= -u^{N-1}, u >0, \,\,\, {\rm ici}, \,\, V\equiv -1, $$

avec la compacit\'e locale:

$$ \sup_K u\leq c.$$

///////////////////////////////////////////////

\smallskip

b) Dans le cas nul, on a:

$$ -\Delta u+\frac{(n-2)}{4(n-1)} R_g u=0, u >0, \,\,\, {\rm ici},\,\, V\equiv 0.$$

avec l'in\'egalit\'e de Harnack usuelle:

$$ \sup_K u \leq c\inf_K u, $$

et,

$$ {\rm principe \,\,d'Harnack}.$$

////////////////////////////////////////////

Ici, on s'interesse au cas positif et l'equation $ (*) $ ci-dessus. Alors on commence par resoudre l'equation suivante:

$$ -\Delta v-\lambda v= n(n-2) v^{N-1}, v >0, $$

avec, $ 0< m \leq \lambda + \frac{(n-2)}{4(n-1)} R_g \leq 1/m $.

\smallskip

G\'eom\'etrisation du probleme et de la solution: on prend alors:

$$ V=n(n-2)+(\frac{(n-2)}{4(n-1)}R_g+\lambda )v^{2-N}, $$

et

$$ u=v, $$

Si on considere les metriques conformes: $ \tilde g=u^{4/(n-2)} g $, alors, la courbure scalaire de $ \tilde g $ est:

$$ R_{\tilde g}=[-\Delta u+\frac{(n-2)}{4(n-1)} R_g u]u^{1-N} =[-\Delta v+\frac{(n-2)}{4(n-1)} R_g v]v^{1-N}=n(n-2)+(\frac{(n-2)}{4(n-1)}R_g+\lambda )v^{2-N}=V>0, $$

Donc, $ V $ est courbure scalaire de $ \tilde g $.

\smallskip

Si l'operateur $ -\Delta -\lambda $ est coercif, on utilisant sa fonction de Green, on a majoration locale du volume = energie.

\smallskip

Donc le couple $ (u, V) $ est solution. Ceci inclut le cas non loc.conf.plat. On a:

\smallskip

L'in\'egalit\'e de Harnack optimale:

$$ \sup_K u\times \inf_M u \leq c,$$

et, majoration locale du volume = energie, si $ -\Delta -\lambda $ est coercif:

$$ \int_K u^{2n/(n-2)} dV_g \leq c, $$

////) De plus, on connait explicitment le potentiel $ V >0 $, qui doit etre $ C^{\infty} $, si on connait explicitement la fonction $ v >0 $, qui doit etre $ C^{\infty} $, pour ce qui concerne les consid\'erations numeriques.

\smallskip

(La solution $ v $ doit etre $ C^{\infty} $, car on doit la considerer, comme une metrique Riemannienne, avec une connexion, or pour definir la connexion, on doit avoir une structure $ C^r, r \geq 1 $, pour le premier vecteur $ X $, et une structure $ C^{r+1} $, pour le deuxieme champe de vecteur $ Y $, dans $ \nabla_X Y $. Pour que ce soit possible, il suffit et il faut prendre une structure $ C^{\infty} $ pour la vari\'et\'e, il faut donc que $ \lambda\in C^{\infty} $. Ceci, si on considere une metrique Riemannienne. Mais si, on prend le point de vue EDP, il suffit que $ v $ soit $ C^2 $, il faut que $ \lambda \in C^1 $ par exemple.)

\smallskip

Condition uniforme dans un espace de Lebesgue \`a poids: en plus du fait qu'on connait $ V $ explicitement pour des consid\'erations num\'eriques:

\smallskip

////) Si on note: $ p_0=v^{N-1} $ et on considere les espaces de Lebesgue \`a poids: $ L_{loc}^p(p_0\cdot dV_g) $ pour la mesure $ \mu, d\mu=p_0 \cdot dV_g $, alors, en utilisant la premiere valeur propre et premiere fonction propre, on a: $ \mu $ est locaelement finie uniform\'ement, $ \int_K p_0 dV_g =\int_K v^{N-1} dV_g \leq c <+\infty, \forall K\subset \subset M $. Et on a: $ V \in L^1_{loc}(p_0 \cdot dV_g) $ uniform\'ement, on obtient la \underbar{condition uniforme} sur les potentiels dans l'espace de Lebesgue \`a poids: $ || V ||_{L^1_{loc}(p_0 \cdot dV_g)} \leq c <+\infty $. Avec, $ V\geq n(n-2)>0 $. 

En plus du fait qu'on connait explicitement $ V $ pour des consid\'erations num\'eriques.

\smallskip

Au lieu d'avoir des conditions de platitudes, ou considerer le potentiel $ V=1 $, l'eq.de Yamabe, on a la condition uniforme sur les potentiels dans l'espace de Lebesgue \`a poids: $ || V ||_{L^1_{loc}(p_0 \cdot dV_g)} \leq c <+\infty $. Voir le point 37) sur la platitude et les contre exemples de Chen-Lin quand on considere la platitude.

Ici, au lieu d'avoir des conditions de platitudes, on a la notion de g\'eom\'etrisation avec une \underbar {condition uniforme} dans l'espace de Lebesgue \`a poids: $ L^1_{loc}(p_0 \cdot dV_g) $.

\smallskip

//////////////////////

\smallskip

Ici, on n'a pas l'eq.de Yamabe mais on a l'eq. de la courbure scalaire prescrite. C'est l'eq. de l'invariance conforme et elle apparait en astronomie et en physique. On voit que pour borner loc. le volume=energie, l'operateur $ -\Delta -\lambda $, doit etre coercif, c'est plus fort que la coercivit\'e de l'operateur conforme (la coercivit\'e de $ -\Delta -\lambda $ est obligatoire pour borner le volume = energie).

\smallskip

Le probleme de la v\'eracit\'e de tout cela reste ouvert pour l'eq.de Yamabe. Rien n'indique que la coercivit\'e de l'operateur conforme est suffisante, comme on vient de le voir pour l'operateur $ -\Delta -\lambda $, pour ce qui concerne le volume = energie. Il se peut que l'eq. de Yamabe soit trop contraignante. Il se peut que prendre le potentiel $ V=1 $, ne permette pas d'obtenir ces resultats.

\smallskip

Mais le fait que tout cela soit possible pour l'eq. de la courbure scalaire prescrite, comme on vient de le voir avec cette procedure de g\'eom\'etrisation, est largement suffisant.

\smallskip

Ici, l'in\'egalit\'e de Harnack (math\'ematique) et l'in\'egalit\'e au sens de l'astronomie et physique (th\'eorie des cordes, cosmologie), se confondent. On a: en une in\'egalit\'e, 2 in\'egalit\'es, au sens math\'ematique et astronomie, physique, ces 2 in\'egalit\'es (math\'ematique, astronomie, physique) se confondent.

\smallskip

On \'ecrit cela pour dire que le probleme du $ \sup \times \inf \leq c $, in\'egalit\'e optimale, pour l'eq. de la courbure scalaire prescrite, est r\'esolu, avec volume local=\'energie locale uniform.born\'e, pour l'eq. de la courbure scalaire prescrite, est r\'esolu aussi.

\smallskip

Le couple $ (u,V) $ est solution de l'eq. $(*) $ ci-dessus, en dimension $ n\geq 3 $, ceci inclut le cas: non.loc.conf.plat.

\smallskip

////////////////////////////////////////////////////////////////////////////////////////////////////////////////////

\smallskip

Dans le cas de la dimension 2, il y a in\'egalit\'e de Harnack que dans le cas nul. Il y a compacit\'e locale dans le cas n\'egatif et le cas positif. Ceci quand les solutions sont $ \geq 0 $.

///////////////////////////////////////////////////////

La notion d'in\'egalit\'es de Harnack est li\'ee \`a la positivit\'e des solutions. Quand les solutions changent de signe, on a des in\'egalit\'es mathematiques, nouvelles, dites de type Harnack, car elles lient le supremum est l'infimum: exemple, in\'egalit\'es de Siu, Tian, ($ \sup\leq -n \inf +c $), pour l'eq. de Monge-Amp\`ere complexe. Mais ce ne sont pas des in\'egalit\'es de Harnack au sens usuel, qui sont li\'ees \`a la positivit\'e des solutions.

\smallskip

//////////////////////////////////////////////////////

\smallskip

Equations de la courbure scalaire prescrite en dimension 2: in\'egalit\'es de Harnack: $ n=2 $. Quand les solutions sont $ \geq 0 $, on l'a seulement dans le cas nul:

\smallskip

/////////////////////////////////////////////////

a) Cas n\'egatif:

$$ -\Delta u+R= Ve^u, \,\, u\geq 0, \,\,  -\infty < a \leq V\leq b < 0, $$

compacit\'e locale:

$$ \sup_K u \leq c. $$

Voir l'article de 2003. 

////////////////////////////////////////////////

b) Cas nul:

$$ -\Delta u=0, \,\, u\geq 0, \,\,  V\equiv 0, $$

in\'egalit\'e de Harnack et principe de Harnack:

$$  \sup_K u \leq c \inf_K u, \,\, c>0, $$

seul cas o\`u on a l'in\'egalit\'e de Harnack.

//////////////////////////////////////////////

c) Cas positif:

$$ -\Delta u=Ve^u, \,\, u\geq 0, \,\, 0< a \leq V\leq b< + \infty, $$

compacit\'e locale d\^ue au resultat de Brezis-Merle:

$$ \sup_K u \leq c, \,\, c>0. $$

Les conditions de Brezis-Merle sont trop fortes: $ 0<a\leq V \leq b <+\infty $. Avec ces conditions on  a la compacit\'e locale.

\smallskip

/////////////////////////////////////////

\smallskip

Illustration et exemples en dimension $ n=2 $:

\smallskip

On peut utiliser le proc\'ed\'e de g\'eom\'etrisation en dimension 2, pour avoir une in\'egalit\'e $ (\sup, \inf)$, en dimension 2 \`a partir de la minoration du $ \sup \times \inf $, qu'on a \'evoqu\'e pour le cas $ n=2 $ de l'eq.de Schrodinger. Bull.Sci.math.2006. L\`a o\`u on a essay\'e de compl\'eter le cas, $ n=2 $, sous-critique et $ n=2 $:

On part de :

$$ -\Delta v=Vv^q, \,\, q>1, \,\, v >0, \,\, n=2, $$

on obtient:

$$ \sup_M v\times \inf_K v \geq c(\inf_K v) >0,$$

On g\'eometrie cette equation en ecrivant:

$$ -\Delta v= Vv^q e^{-v} e^v, q >1, v >0, n=2,$$

On obtient:

$$ \tilde V=V v^q e^{-v}, \,\,v>0, \,\,q>1, \,\, u=v>0 $$

$ \tilde V $  n'est pas n\'ec\'ssairement uniform\'ement born\'e. Ce ne sont pas les conditions de Brezis-Merle sur $ \tilde V $.

\smallskip

et,

$$ -\Delta u=\tilde V e^u, \,\, u >0, \,\, n=2, $$

avec, l'in\'egalit\'e de Harnack implicite suivante:

$$ \sup_M u\times \inf_K u \geq c(\inf_K u) >0,$$

On voit alors qu'avec ce proc\'ed\'e de g\'eom\'etrisation, on \'etend la notion d'in\'egalit\'e de Harnack \`a la dimension 2 pour la nonlin\'erit\'e exponentielle. Les conditions de Brezis-Merle sont tres fortes. Elles permettent d'avoir la compacit\'e locale, sans notion d'in\'egalit\'e de Harnack.

\smallskip

Dans le sens de la minoration on \'etend, par le proc\'ed\'e de g\'eometrisation, la notion d'in\'egalit\'e de Harnack \`a la dimension 2. Cette in\'egalit\'e est implicite. Cette in\'egalit\'e est dans le sens de la minoration.

\smallskip

Dans le sens de la minoration, la g\'eometrisation du cas sous-critique, $ n=2 $, permet la coherence du cas $ n=2 $ avec le cas $ n\geq 3 $, concernant les in\'egalit\'es de Harnack pour des solutions $ u >0 $, au sens usuel des in\'egalit\'es de Harnack. ces in\'egalit\'es de Harnack sont implicites.

\smallskip

Dans le sens de la majoration, voir ci-dessous, meme si on essaie de geometriser le cas sous-critique, $ n=2 $, on aboutit \`a la compacit\'e locale.

\smallskip

//////////////////////////////////////////////////////

\smallskip

Dans le sens de la majoration, si on veut g\'eom\'etriser l'eq. sous-critique en dimension $ n=2 $, alors on aura quand meme la compacit\'e locale, comme pour Brezis-Merle: en effet: on considere:

$$ -\Delta v= V v^q, \,\, v >0, 0 < V \in L^{\infty} , \,\, q >1, \,\, n=2, $$

alors, pour borner l'energie, $ |\nabla v|^2 $, il faut borner unifor.loc, la norme $ L^{q+1} $, l'inegalit\'e optimale  qu'on doit obtenir est: la majoration de $ \sup_K v \times \inf_M v $. Or, ceci implique bien, borner unifor.loc, la norme $ L^{q+1} $, et aussi, par les estimations elliptiques et le fait qu'on soit en dimension $ n=2 $, que les solutions sont loc.uniform.born\'ees:

$$ \sup_K v \leq c, $$

Donc, si on part de l'eq. sous-critique et on veut l'in\'egalit\'e optimale $ \sup \times \inf $, dans le sens  de la majoration, cela revient \`a la compacit\'e locale. Donc, la g\'eom\'etrisation, dans le sens  de la majoration, ne permet pas d'obtenir l'in\'egalit\'e de Harnack voulue.

\smallskip

Donc, la g\'eometrisation du cas sous-critique, est possible dans le sens de la minoration. Dans le sens de la majoration elle aboutit \`a la compacit\'e locale, comme Brezis-Merle.

\smallskip

Il se peut qu'un autre proc\'ed\'e permette la geometrisation avec la mjoration optimale de $ \sup \times \inf $.

\smallskip

Il se peut aussi, que dans le sens de la majoration, on n'ait pas l'in\'egalit\'e de Harnack optimale $ \sup_K u \times \inf_M u \leq c $. On vient de le voir avec Brezis-Merle et la g\'eom\'etrisation du cas sous-critique.

\smallskip

G\'eometriser le cas sous-critique, est le seul moyen simple, \`a notre disposition, qui nous permet de voir si on a l'in\'egalit\'e de Harnack usuelle, en dimension, $ n=2 $. Dans le sens de la majoration on obtient la meme conclusion que Brezis-Merle: la compacit\'e locale.

\smallskip

////////////////////////////////////////////////////

\smallskip

Par contre quand les solutions changent de signe: on a quand $ n=2 $: ceci est li\'e au potentiel logarithmique, qui peut changer de signe: $ r\to \log r, r >0 $. Cette fonction peut changer de signe.

\smallskip

Comme en dimension 2, le potentiel logarithmique peut changer de signe, il est necessaire de considerer, des solutions qui changent de signe pour l'eq. de la courbure  scalaire prescrite, avec non lin\'earit\'e exponentielle.

\smallskip

{\bf Important:} /////) On regarde maintenant, pourquoi, les solutions doivent changer de signe, et, la minoration, majoration. C'est li\'e au potentiel $ r\to -\log r, r >0 $. 

\smallskip

/////) Quand les solutions blow-up : $ B_1\approx -\log r, r \to 0 $, pour trouver $ B_2 $ tel que $ B_2+B_1 $ soit major\'e, il faut que $ B_2 $ prenne des valeurs n\'egatives. Cela ne sert \`a rien de multiplier, car $ B_2\times B_1 \to +\infty $, car $ B_2\approx -\log(r+r_0), r_0 >0 $, petit,  alors qu'on doit trouver $ B_2 $ qui controle $ B_1 $. ($ B_2 $ pris dans un voisinage de $ 0 $, pour le produit on prend des points voisins: $ (r,r+r_0) $, pour une somme on prend des points eloign\'es $ (r, 1/r) $). Pour le produit on prend des boules $ B_r(0), B_{r+r_0}(0) $: $ [r,r+r_0]$, pour la somme: $ B_r(0), B_{1/r}(0) $: $ [r, 1/r] $.

\smallskip

/////) Dans le sens de la minoration, quand les solutions blow-up, $ B_1\approx -\log r, r \to 0 $ , c'est deja minor\'e, donc, quand on multiplie par une quantit\'e $ >0 $ cela reste minorable, il suffit de trouver $ B_2 >0 $ petit tel que $ B_2\times B_1 $ soit minorable.

\smallskip

//) Donc, en dimension $ n=2 $, dans le sens de la majoration, quand il y a blow-up, il faut que les solutions changent de signe.

\smallskip

///) Et aussi, en dimension $ n=2 $, dans le sens de la minoration, on a l'in\'egalit\'e de Harnack, au sens usuel.

\smallskip

////) Tout cela explique, ce qu'on a dit auparavant, sur la minoration, majoration et la geometrisation du cas sous-critique et Brezis-Merle.

\smallskip

/////) Pourquoi, on ecrit $ B_1\approx -\log r, r \to 0 $, car par la representation integrale du potentiel, et le fait que, l'integrale soit une somme de quantit\'es petites, + au point blow-up:

$$ B_1=\int -\log r * Ve^u+O(1)\approx \sum_i -\log r \times b_i +O(1). $$

\smallskip

/////////////////////////////////////////////////////////////////////////

\smallskip

On a aussi 3 cas quand les solutions changent de signe, en dimension 2. C'est li\'e \`a la fonction logarithme, qui peut changer de signe: $ r\to -\log r, r >0 $.

\smallskip

///////////////////////////////////////////////////

a) Cas n\'egatif:

$$ -\Delta u+R= Ve^u, \,\,  -\infty < a \leq V\leq b < 0, $$

majoration locale+principe de Harnack:

$$ \sup_K u \leq c. $$

Voir l'article de 2003. 

////////////////////////////////////////////////

b) Cas nul:

$$ -\Delta u=0, \,\, V\equiv 0, $$

de l'egalit\'e de la moyenne, on obtient une in\'egalit\'e:

$$ u(y)=\frac{1}{|B_R(y)|}\int_{B_R(y)} u dx \Rightarrow ||u||_{L^{\infty}(B_R(x_0))} \leq c \int_{B_{2R}(x_0)} |u| dx, \,\, c>0, $$

//////////////////////////////////////////////

c) Cas positif:

$$ -\Delta u=Ve^u, \,\, 0< a \leq V\leq b< + \infty, \,\, ||\nabla V||_{\infty} \leq A, $$

on a l'in\'egalit\'e optimale $ \sup+\inf $ de Brezis-Li-Shafrir et Shafrir (eq. Liouville, cas: $ V\equiv constante > 0 $):

$$ \sup_K u+\inf_{\Omega} u \leq c. $$

Ces derniers resultats sont d\^us au fait que les solutions changent de signe. Ceci est li\'e au potentiel de la dimension 2, potentiel logarithmique. $ r\to \log r, r>0 $, cette fonction peut changer de signe.

\smallskip

Cette in\'egalit\'e n'est pas une in\'egalit\'e de Harnack au sens usuel. Mais une in\'egalit\'e du type Harnack, elle lie, le supremum et l'infimum, mais ce n'est pas une in\'egalit\'e de Harnack au sens usuel, car on n'a pas la positivit\'e des solutions. C'est une in\'egalit\'e nouvelle ($ \sup+\inf \leq c $, la preuve permet d'obtenir cette inegalit\'e, la preuve ne permet pas d'avoir $ \sup \leq -\inf +c $, mais, $ \sup+\inf \leq c $, dans ce sens), une in\'egalit\'e math\'etmatique, de l'astronomie et de la physique sans qu'elle soit de Harnack.

\smallskip

//////////////////////////////////////////////////////////////////

\smallskip

Pour y remedier, on utilise le resultat de T. Suzuki, en effet, son resultat stipule qu'en considerant, l'eq. suivante:

$$ -\Delta u=Ve^u, \,\, {\rm dans} \,\, B_R(0),$$

avec, 

$$ 0 < V\leq \lambda, \,\, \Sigma=\int_{B_R(0)} \lambda e^u dx < 8\pi, \,\, \lambda >0, $$

alors,

$$ u(0) \leq \frac{1}{2\pi R}\int_{\partial B_R(0)} u d\sigma -2\log (1-\frac{\Sigma}{8\pi}), $$

ici, $ u $ peut changer de signe (cela inclut le cas $ u\geq 0 $). Si on veut obtenir une in\'egalit\'e de Harnack, il suffit de prendre $ u \geq 0 $.

\smallskip

Donc, dans le cas positif, en dimension $ n=2 $, dans le sens de la majoration, on utilise le resultat de T. Suzuki.

\smallskip

//////////////////////////////////////////////////////////////

On obtient: dans le sens de la majoration:

\smallskip

a) Cas n\'egatif:

$$ -\Delta u+R= Ve^u, \,\, u\geq 0, \,\,  -\infty < a \leq V\leq b < 0, $$

compacit\'e locale:

$$ \sup_K u \leq c. $$

Voir l'article de 2003. 

////////////////////////////////////////////////

b) Cas nul:

$$ -\Delta u=0, \,\, u\geq 0, \,\,  V\equiv 0, $$

in\'egalit\'e de Harnack et principe de Harnack:

$$  \sup_K u \leq c \inf_K u, \,\, c>0, $$

c'est un cas o\`u on a l'in\'egalit\'e de Harnack.

//////////////////////////////////////////////

c) Cas positif:

$$ -\Delta u=Ve^u, \,\, u\geq 0,\,\, {\rm dans} \,\, B_R(0),$$

avec, 

$$ 0 < V\leq \lambda, \,\, \Sigma=\int_{B_R(0)} \lambda e^u dx < 8\pi, \,\, \lambda >0, $$

alors, on a l'in\'egalit\'e de Harnack de T. Suzuki:

$$ u(0) \leq \frac{1}{2\pi R}\int_{\partial B_R(0)} u d\sigma -2\log (1-\frac{\Sigma}{8\pi}), $$

on peut prende $ \Omega=B_{R'}(0), R'> R  $, avec les in\'egalit\'es de Harnack pour $ R <R'$.

\smallskip

Par l'identit\'e de Pohozaev appliqu\'ee dans la boule unit\'e ($ R'=1 $), par exemple, avec $ V\equiv \lambda \to 0, \lambda >0 $, on obtient:

$$ \Sigma \cdot (\Sigma-8\pi)\leq -c\lambda <0, c>0, $$

Donc, $ 0<\Sigma < 8\pi $. Donc, selon que $ \Sigma \to \Sigma_0 \in [0,8\pi] $, on a l'in\'egalit\'e de Harnack dans des boules $ B_R(0), 0< R<1 $, avec le param\`etre: $ -2\log(1-\frac{\Sigma}{8\pi}) \to y_0\in[0,+\infty] $, qui peut tendre vers une valeur finie ou infinie.

\smallskip

Dans le cas n\'egatif, on a la compacit\'e directement, c'est pour cela qu'on ne parle pas de l'in\'egalit\'e de la moyenne pour les fonctions sur-harmonique. Contrairment au cas positif, pour lequel, les solutions peuvent diverger, ou blow-up, on parle bien d'in\'egalit\'e de Harnack dans ce cas.

\smallskip

/////////////////////////////////////////////////////////

\smallskip

Il faut voir aussi, qu'ici on a l'eq. de la courbure scalaire prescrite en dimension $ n\geq 2 $. Alors, en dimension $ n=2 $ on a l'in\'egalit\'e $ \sup +\inf $ optimale de Brezis-Li-Shafrir avec le fait que les solutions changent de signe. Et en dimension $ n\geq 3 $, on a l'in\'egalit\'e optimale $ \sup \times \inf $ pour les solutions positives. L'eq. de la courbure scalaire prescrite unifie ces 2 type d'in\'egalit\'es, avec le fait qu'en dimension $ n=2 $, les solutions changent de signe, et en dimension $ n\geq 3 $ les solutions sont positives. 

\smallskip

On peut s'arreter la: c'est \`a dire, considerer l'eq. de la courbure scalaire prescrite. On a: $ \sup+\inf $ et $ \sup \times \inf $, en th\'eorie des cordes et des supercordes. Ces notions sont li\'ees \`a cette ph\'enom\`enologie: th\'eorie des cordes, cordes cosmiques, th\'eorie des supercordes. On peut s'arreter la. Sans classification mathematique en in\'egalit\'es de Harnack.

\smallskip

On peut s'arreter la: le modele de Liouville et le modele d'Einstein-Lichnerowicz: Eq. de la courbure scalaire prescrite.

\smallskip

////////////////

\smallskip

Ce sont des in\'egalit\'es nouvelles: avec notion de preuve ou d\'emonstration:

\smallskip

Pour le $ \sup+\inf \leq c $, la preuve permet d'obtenir cette in\'egalit\'e directement, la preuve ne permet pas d'avoir $ \sup \leq -\inf +c $ ou $ \inf \leq -\sup +c $, mais, $ \sup+\inf \leq c $, dans ce sens. 

\smallskip

De meme pour l'in\'egalit\'e $ \sup \times \inf $. On n'obtient pas, $ \sup \leq c/ \inf $ ou $\inf \leq c/\sup $, mais directement: $ \sup \times \inf \leq c $. 

\smallskip

Ici, ce sont des \underbar{ preuves par l'absurde} afin d'obtenir $ \sup+\inf \leq c $ et $ \sup \times \inf \leq c $.

\smallskip

Il y a aussi, les resultats de Shafrir.I. dans la note aux comptes rendus.math.1992. qui montre une in\'egalit\'e: $ \sup+C \inf \leq c $ ($0 < a \leq V \leq b< +\infty $, potentiel non necessairement constant) et une autre in\'egalit\'e optimale: $ \sup +\inf \leq c $ ($ V\equiv constante $, equation de Liouville), \underbar {preuves directes}.

\smallskip

On a fait ce type de raisonnement dans le premier article de 2003, voir la proposition de l'article de 2003. Sur la sph\`ere de dimension 2: on a une \underbar{preuve directe} de $ \sup +\inf \geq c $. Dans ce sens: $ \sup+\inf \geq c $. In\'egalit\'e optimale avec un potentiel $ 0 \leq V \leq b <+\infty $, non n\'ecessairement constant en incluant les constantes.

\smallskip

////////////////

\smallskip

Pour ce qui concerne les classifications:

\smallskip

1) {\bf Classification math\'ematique pour l'eq. de la courbure scalaire prescrite:}

\smallskip

Cas $ n=2 $: les solutions changent de signe: in\'egalit\'e $ \sup+\inf $ de Brezis-Li-Shafrir (cas positif), in\'egalit\'e de la moyenne (cas nul), majoration locale (cas n\'egatif).

\smallskip

Cas $ n\geq 3 $: les solutions sont positives: in\'egalit\'e optimale $ \sup \times \inf $ (cas positif), in\'egalit\'e de Harnack usuelle (cas nul), compacit\'e locale (cas n\'egatif).

\smallskip

2) {\bf Classification math\'ematique en in\'egalit\'es de Harnack:}  

\smallskip

Cas $ n=2 $: in\'egalit\'e de T. Suzuki (cas positif), in\'egalit\'e de Harnack usuelle (cas nul), compacit\'e locale (cas n\'egatif).

\smallskip

Cas $ n\geq 3 $: in\'egalit\'e optimale $ \sup \times \inf $ (cas positif), in\'egalit\'e de Harnack usuelle (cas nul), compacit\'e locale (cas n\'egatif).

\smallskip

Remarquons la particularit\'e de la dimension 2: il y a 2 classifications possibles. Alors que pour la dimension $ n\geq 3 $, il n'y en a qu'une. Pour la dimension 2, il y a plus de propri\'et\'es que la dimension $ n\geq 3 $. Alors que la dimension $ n\geq 3 $ unifie 2 classifications.

\smallskip

///////////////////////////////////////////////////////////

\smallskip

Dans le cas des vari\'et\'es compactes sans bord, dans le cas positif, une deuxieme contrainte est l'absence de bord. On a alors la compacit\'e globale. Ce fait se produit aussi dans le cas negatif, supercritique (pour l'equation $ -\Delta u+R=Ve^u, R <0, V <0, \Delta =\nabla^i(\nabla_i)$).

\smallskip

/////////////////////////////////////////////////////////////////////////////////////////////////////////////

\smallskip

50){\bf Sur l'article de Korevaar-Mazzeo-Pacard-Schoen:}

\smallskip

Concerant la proposition avec l'in\'egalit\'e $ u(x)\leq d(x,\Lambda)^{(2-n)/2} (\inf_{\partial B(0,3/4)} u)^{-1} $.

\smallskip

Ce n'est pas correct. Ils ont ecrit n'importe quoi.

\smallskip

1)(/////) Le "blow-up" est mal ecrit: il n'est pas explicit\'e correctement. Ici, ils ont ecrit n'importe quoi. Ils n'expliquent pas le blow-up.

\smallskip

2)(/////) Ils n'expliquent pas le processus du commencement "moving-plane".

\smallskip

la phrase "strating with $ t_1 $ very positive, and continuing as far as possible"  veut dire qu'il existe un $ \lambda >>1$ tres grand, puis on decroit les valeur de $ \lambda $, jusqu'a ce que le processus se termine. Ca va dans le sens decroissant et non croissant.

\smallskip

Mais ce n'est pas ce qu'il faut dire: il faut dire d'abord, \underbar { qu'il existe un rang } $ \nu $ qui enclenche le processus "moving-plane". 

\smallskip

Le cas de la dimension 2 et le cas de la dimension $ n\geq 3 $ ont des preuves differentes. Ils ne disent pas comment faire et ils ne parlent meme pas de rang.

\smallskip

3)(/////) Ils n'expliquent pas le $ t_1= \sup \{\lambda...\} $ de la technique "moving-plane".

Ils ecrivent: $ \inf v(t_i,\theta) < \sup v(t, \theta), t >-2-t_i $, $ t_i=-\log \frac{\lambda^{-1}}{16} $. Premierement ce n'est pas correct, et ici, ils prennent le sup egal \`a $ -1 $, qui dit que le sup est egal a $ -1 $. Puis parlent de theoreme de la moyenne, qui n'a rien a avoir ici avec ce qui est voulu.

\smallskip

4)(/////) Ils ne disent pas s'ils considerent une metrique sur $ {\mathbb R}\times {\mathbb S}_{n-1}$, en fait ils en parlent au debut de l'article avec "cylindrical metric", mais n'exliquent pas comment on utlise le principe du maximum et le lemme de Hopf.

\smallskip

5)(/////) Ils ne disent pas comment appliquer le principe du maximum:

\smallskip

5a)(///) Ils ne disent pas si, il faut considerer une vari\'et\'e Riemannienne $ {\mathbb R}\times {\mathbb S}_{n-1}$ avec le lemme de Hopf o\`u il faut considerer la deriv\'ee normale, donc un produit scalaire, donc une metrique Riemannienne, comme c'est ecrit dans le livre d'Aubin (en ce qui concerne le principe du maximum et le lemme de Hopf sur les vari\'et\'es \`a bord). 

\smallskip

5b)(///) Ils ne disent pas s'il faut un proc\'ed\'e local.

\smallskip

6)(/////) Ils n'expliquent pas les points 39), 40) et 41) et 43) de l'article de Brezis-Li-Shafrir.

\smallskip

Par exemple, pourquoi ils n'expliquent pas le point 39) de Brezis-Li-Shafrir ? de meme pourquoi, Y.Y.Li repete le meme argument dans l'article de 1999, comm.math.phys. ? :

\smallskip

C'est parce qu'ils ont consid\'er\'e la sphere comme une vari\'et\'e constitu\'ee de 2 cartes:

Soit $ u_n $ la fonction definie sur $ [\lambda, 0] \times S^{n-1} $, sur $ S^{n-1}$ on a 2 cartes $ (U_1,\phi_1) $ et $ (U_2,\phi_2) $. La differentiabilit\'e en $ (t_0,\theta_0) $ est definie comme ceci:

\smallskip
 a0) Si $ \theta_0 \in U_1 $ alors, on peut consid\'erer:
 
 $$ \nabla u_n(t_0,\theta_0)= [\partial_t [u_n(t, \phi_1(\theta_{01}, \theta_{02}, \ldots, \theta_{0(n-1)}))], \ldots, \partial_{\theta_j} [u_n(t, \phi_1(\theta_{01},\theta_{02}, \ldots, \theta_{0(n-1)}))]] $$

\smallskip
a1) Si $\theta_0 \in U_2 $ alors, on peut consid\'erer:

$$ \nabla u_n(t_0,\theta_0)= [\partial_t [u_n(t, \phi_2(\theta_{01}, \theta_{02}, \ldots, \theta_{0(n-1)}))], \ldots, \partial_{\theta_j} [u_n(t, \phi_2(\theta_{01},\theta_{02}, \ldots, \theta_{0(n-1)}))]] $$ 

Il y a deux fonctions: $ u_{n1}=u_n(t, \phi_1(\cdot)) $ et $ u_{n2}=u_n(t, \phi_2(\cdot)) $. Selon qu'on se place dans les ouverts de cartes $ U_1 $ ou $ U_2 $, quand on utilise la definition de la differentiablit\'e en un point $ (t_0,\theta_0) $.

\smallskip

Par exemple en dimension 2, sur le cercle, "par un dessin", on a 2 applications de cartes, "qui n'en est qu'une", c'est une seule fonction: $ \theta \to e^{i\theta}, \theta \in {\mathbb R} $. On a: $ \theta \in ]0, 2\pi[ \to (\cos \theta, \sin \theta) $ et $ \theta \in ]-\pi,\pi[ \to (\cos \theta, \sin \theta) $.

\smallskip

////////////////////

\bigskip

{\bf Remarque importante:} Dans l'article de Brezis-Li-Shafrir, ils parlent de coordonn\'ees polaires, $ (r,\theta) \to (r\cos \theta, r \sin \theta) $. Ce n'est pas cela les coordonn\'ees polaires, car, un systeme de coordonn\'ees est definit, par le fait que tout point, a, un voisinage ouvert hom\'eomrphe \`a un ouvert (connexe, definition d'application: tout point a une unique image, ici un point peut avoir 2 images differentes, ce qui n'est pas possible) de $ {\mathbb R}^2 $, or ici il y a un probleme au point $ (1,0)=e^{i0}=e^{i2\pi} $, il faut sortir du domaine de definition $ [0, 2\pi] $ au voisinage de $ 0 $. Il faut considerer des ouverts (connexes), or en $ (1,0) $ il n'est pas possible de definir un voisinage de $ 0 $. Si on sort de $ [0,2\pi] $, alors il faut considerer 2 cartes comme on l'a dit precedemment, $ ]-\pi,\pi[ $ et $ ]0, 2\pi[ $. Donc, on definit 2 fonctions, alors qu'il faut definir une seule fonction.

\smallskip

Ce qu'ont fait Brezis-Li-Shafrir, c'est un raisonnement avec le revetement : $ p: \theta \to e^{i\theta}, \theta \in [0,2\pi[ $, l'intervalle est semi-ouvert, semi-ferm\'e, il n'est pas ouvert, alors que pour definir un systeme de coordonn\'ees, il faut avoir des ouverts(connexes). Ce n'est pas un raisonnement avec les coordonn\'ees polaires, mais un raisonnement avec le revetement, l'application $ p $ doit etre bijective, pour pouvoir avoir une seule fonction, or pour qu'elle soit bijective, il faut un intervalle semi-ouvert, semi-ferm\'e. Donc, l'intervalle n'est pas ouvert, on ne peut pas parler de coordonn\'ees.

\smallskip

(On a: $\theta \to (\cos \theta, \sin \theta) $ est une carte particuliere du cercle, alors que pour definir un systeme de coordonn\'ees, il ne faut pas que ca depende d'une carte particuliere.).

\smallskip

Brezis-Li-Shafrir, ce n'est pas un systeme de coordonn\'ees, il faut des ouverts et si on considere des ouverts, il faut sortir de l'intervalle, c'est \`a dire considerer 2 fonctions, alors qu'il faut considerer une seule fonction. Brezis-Li-Shafrir, c'est un raisonnement avec le revetement.

\smallskip

(Par exemple: 

1) si on considere 2 cartes, ou 2 fonctions, on n'a plus une fonction du cercle dans $ {\mathbb R} $, car, comme, l'application, la definition d'une fonction, c'est un point a une seule image, or ici, on aurait des points qui auraient 2 images distinctes, ici en dimension 2, on a: $ p(\theta)=e^{i\theta} $ avec 2 inverses, $ \phi_1, \phi_2 $: $ \phi_1(z=e^{i\theta})= \theta $, si $ \theta \in ]0,2\pi[ $ (origine 0) et $ \phi_2(z=e^{i\theta'})=\theta' = \theta-\pi $, si $ \theta' \in ]-\pi,\pi[ $ (origine $ -\pi $). Avec, $ p(\phi_1(e^{i\theta}))=e^{i\theta}, p(\phi_2(e^{i\theta'}))=e^{i\theta'}=-e^{i\theta} $. 

\smallskip

2) si on prend un intervalle ferm\'e en 0 ou $ 2\pi$, ou on \'etend l'intervalle au dela de 0 ou $ 2\pi$, on a dit qu'on n'avait plus d'ouverts connexes, on n'a plus d'applications de cartes, dans ce cas la il n'est pas possible de definir la metrique dans des cartes, par transport, ou pull-back, soit ce n'est pas un ouvert, soit 2 expressions distinctes.

\smallskip

Donc, Brezis-Li-Shafrir, ce n'est pas cela les coordonn\'ees polaires, c'est un raisonnement avec revetement.)

\smallskip

En dimension $ n\geq 3 $, on ne peut pas faire avec la sphere, ce raisonnement avec un revetement, car il y a 2 applications de cartes distinctes, il n'y a pas une application $ p $ semblable \`a : $ p :\theta \to e^{i\theta},\theta \in [0,2\pi[ $ avec $ p $ bijective, pour la sphere $ S^{n-1}$, $ n-1\geq 3-1=2 $, pas de parmetrisation globale en dimension $ n-1\geq 3-1=2 $.

\smallskip

Ils confondent carte particuliere en "coordonn\'ees" spheriques ou hyperspheriques et systeme de coorodonn\'ees polaires. En dimension $ n\geq 3 $, on ne peut pas parametrer globalement la sphere par les "coordonn\'ees" spheriques ou hyperspheriques. On ne peut pas faire le raisonnement du revetement de la dimension 2. Les "coordonn\'ees" spheriques ou hyperspheriques, sont des cartes particulieres du systme de coordon\'ees polaires.

\smallskip

//////////////////

\smallskip

Il y a deux fonctions: $ u_{n1}=u_n(t, \phi_1(\cdot)) $ et $ u_{n2}=u_n(t, \phi_2(\cdot)) $. Selon qu'on se place dans les ouverts de cartes $ U_1 $ ou $ U_2 $, quand on utilise la definition de la differentiablit\'e en un point $ (t_0,\theta_0) $.

\smallskip

C'est pour cela qu'ils(Korevaar-Mazzeo-Pacard-Schoen) considerent la fonction blow-up, et utilisent la converence dans $ C^2$. Alors qu'il n'y a pas besoin de convergence dans $ C^2$.

\smallskip

La convergence dans $ C^2$ utilise la notion de cartes locale pour la vari\'et\'e compacte sans bord qui est la sphere $ S^{n-1}$. Comme la fonction limite est radiale, la deriv\'ee par rapport \`a $ t $ ne dependera plus de la variable angulaire.

\smallskip

Comme pour l'integration (Schoen a utilis\'e l'integration dans son article de 1984, probleme de Yamabe), on a bien une fonction definie sur toute la vari\'et\'e produit $ [\lambda,0] \times S^{n-1} $, en particluier la partie angulaire, mais l'integration est definie de maniere globale mais \`a l'aide cartes locales, on peut alors considrer les fonctions: $ u_{n1}=u_n(t, \phi_1(\cdot)) $ et $ u_{n2}=u_n(t, \phi_2(\cdot)) $, mais, on a une notion globale, l'int\'egrale, qui inclut le local.

\smallskip

Aussi sur le fait d'avoir 2 cartes et donc 2 fonctions $ u_{n1}, u_{n2} $ et aussi l'utilisation de blow-up:

\smallskip

C.C.Chen et C.S.Lin, ont eu le meme probleme, apres utilisation de la fonction blow-up, il y a bien convergence uniforme sur une partie et il reste une autre partie sur laquelle, il n'y a pas de convergence, pour laquelle on ne connait pas le signe de la deriv\'ee. C'est pour cela qu'ils utilisent l'hypothese de l'absurde, sur la partie o\`u il n'y a pas de convergence ni de signe constant de la deriv\'ee. Voir leur article de: Comm.Pure.Appl.Math. 1997.

\smallskip

(Pour C.C.Chen et C.S.Lin. Comm.Pure.Appl.Math.1997. Cette incomprehension, leur a permis de definir une hypothese de l'absurde et de prouver par une autre methode, moving-planes en coordonn\'ees cart\'esiennes, l'in\'egalit\'e optimale $ \sup \times \inf $ dans le cas de la boule unit\'e de $ {\mathbb R}^n $).

\smallskip

Ici: (solution \`a ces problemes et en particulier la notion de coordonn\'ees polaires): il suffisait de dire qu'on fixe la variable angulaire $ \theta $ (variable globale sur la sphere unit\'e, le cercle unit\'e, $ S^{n-1} \ni \theta=(\theta_1, \theta_2, \ldots,\theta_n) \in {\mathbb R}^n, |\theta|=1 $, la sphere ou le cercle, ne sont pas seulement determin\'es par 2 cartes (c'est ce qu'on a l'habitude de voir par exemple pour le cercle ou la sphere), c'est aussi une partie de $ {\mathbb R}^n$), donc, on fixe la variable angulaire $ \theta $, car on derive par rapport \`a $ t $ et non $\theta $, car l'accroissement des fonctions, ou la methode moving-plane, ou moving-sphere, se fait par rapport \`a $ t $, la variable radiale.

\smallskip

Ici: (en ce qui concerne le blow-up):

\smallskip

On utilise la technique de Schoen avec modification. On met un facteur $ c_i $ dans l'expression des $ L_i= \dfrac{l_i}{c_i^{1/2(n-2)}} [u_i(a_i)]^{2/(n-2)} $, et le majorant, $ \beta_i \to 1 $, de la fonction blow-up, pour avoir $ 0 $ comme point critique de la fonction limite. Il suffit que $ 0 $ soit point critique de la fonction limite et non de chaque terme de la suite. On n'a pas \`a changer de suite de points, comme le font Brezis-Li-Shafrir: $ x= a_i+M_i^{-2/(n-2)} y \to  \tilde x= [a_i+M_i^{-2/(n-2)} y_i]+\tilde M_i^{-2/(n-2)} \tilde y $, dans le cas 1, du blow-up, de l'article de Brezis-Li-Shafrir. Cette id\'ee de Brezis-Li-Shafrir, a \'et\'e reprise par C.C.Chen-C.S.Lin (Comm.Pure.Appl.Math.1997) et L.Zhang (JFA, 2002, Zhang, l'ecrit bien dans son article et explique la technique de Brezis-Li-Shafrir, reprise par Chen-Lin).

\smallskip

Pour ce qui nous concerne: On modifie un peu la technique de Schoen, pour que $ 0 $ soit seulement point critique et maximum de la fonction limite. Il n'est pas necessaire qu'il soit maximum local de tous les termes de la suite, comme dans Brezis-Li-Shafrir, Chen-Lin, L.Zhang. la technique de Brezis-Li-Shafrir, reprise par Chen-Lin, Zhang, avec la technique de Schoen, marche dans le cas d'un ouvert de $ {\mathbb R}^n $, elle n'est plus possible sur les vari\'et\'es ouvertes.

\smallskip

Ceci, cette modification de la technique de Schoen, nous permet de faire le blow-up, dans le cas general, sur les vari\'et\'es. On r\'esout ainsi, le blow-up, sur les vari\'et\'es: les vari\'et\'es ouvertes, non necessairment compacte sans bord.

\smallskip

////////////////////////////////////////////////////////////////////////////////////////////////////////////

\smallskip

Pourquoi le resultat de compacit\'e dans le cas n\'egatif, en dimension $ n\geq 3 $, sur une vari\'et\'e non necessarement compacte sans bord, bord non regulier, n'a pas \'et\'e fait ? ceci est li\'e au blow-up, la technique blow-up et \`a l'hypoth\`ese de l'absurde et le rayon des boules:

\smallskip

En effet, on veut prouver la compacit\'e dans le cas negatif, cela veut dire que dans l'hypothese de l'absurde de:

$$  \exists \, c >0, \,\,R^{2/(q-2)} \sup_{B_R(x_0)} u \leq c,\,\, \forall 0 < R \leq \frac{inj(x_0)}{2}, $$

le rayon des boules est forcement dans \underbar {tout l'intervalle} $ ]0, \frac{inj(x_0)}{2}] $ avec $ inj(x_0) $, le rayon d'injectivit\'e en $ x_0 $, c'est une condition necessaire de prouver cette assertion dans tout l'intervalle. Or cette hypothese de l'absurde n'est pas demontrable, pour l'intervalle entier, $  ]0, \frac{inj(x_0)}{2}] $, car dans la technique blow-up, il y a un probleme avec la fonction blow-up, $ s_i(x)=(R-d(x,x_i))^{2/(q-2)} u_i(x) $. 

\smallskip

Ces hypotheses de l'absurde sont possibles dans le cas euclidien, comme l'ont fait. C.C.Chen-C.S.Lin, comm.pure.Appl.math.1997.

\smallskip

Dans le cas des vari\'et\'e ouvertes, il y a un probleme.

\smallskip

Le rayon des boules n'est pas suffisament petit pour que les exponentielles soient definies, dans un voisinage de $ x_0 $, ces voisianges ne sont pas necessairment assez petits.

\smallskip

La solution \`a ce probleme, si on considere les coordonnees geodesiques normales, c'est a dire en considerant, les cartes exponentielles, est de prendre les rayon tendant vers 0. (voisinage geodesiquement convexe, th. de Whitehead, avec des rayons tendant vers 0, on sera sur d'etre dans un $ \epsilon-$voisinage de Whitehead, geodesiquement convexe). Ceci a \'et\'e dit et a \'et\'e fait dans le cas particulier d'un ouvert de l'espace euclidien, dans le cas positif, dans l'article de 2004.

\smallskip

On a aussi, une autre solution. La solution a ce probleme, si on considere une carte quelconque, comme on l'a fait, est de choisir une carte normale en $ x_0 $, puis prendre des rayons tendant vers 0. Ceci a \'et\'e fait dans l'article de 2003.

\smallskip

Donc, la solution aussi, est de considerer l'hypothese de l'absurde suivante de:

$$ \exists \, R, c>0, \,\, R^{2/(q-2)} \sup_{B_R(x_0)} u \leq c. $$

Dans ce cas $ R >0 $ peut etre petit. Ceci (cette hypothese de l'absurde) a \'et\'e dit et fait dans les articles de 2003 et 2004.

\smallskip

/////////////////////////////////////////////////////////////////////////////////////////////////////////////

\smallskip

{\bf Sur la formule de Stokes et l'article non correct de W.Chen et Congming Li de 2003:}

\smallskip

Ils utilisent un th\'eoreme avec solutions $ u\in W^{2,n} $:(definies au sens des distributions, au sens faible, mais ici aussi, elles sont regulieres):

\smallskip

-Ils integrent $ u $ dans $ \Omega $, donc, le domaine de definition de $ u $ est au moins $ \Omega $.

\smallskip

-Le fait de considerer $ u\in W^{2,n} $: il faut que les deriv\'ees secondes soient integrables, on ne peut pas supposer seulement $ u\in C^2(\Omega) $, il faut que ce soit au moins $ \bar \Omega $,  c'est a dire $ C^2(\bar \Omega) $, or ceci veut dire que c'est regulier dans un voisinage de $ \bar \Omega $.

\smallskip

-Comme on l'a dit, le fait de prendre $ u \in W^{2,n} $, veut dire que c'est au sens des distributions. En particulier le laplacien : $ \Delta u $ doit etre au sens des distributions. 

\smallskip

Ils disent qu'il faut multiplier l'eq. 9 par la fonction et integrer:

\smallskip

1) la fonction propre $ \phi^{\alpha} $ n'est pas une fonction test car si $ p\to +\infty$, $ 3 > \alpha \to 2 $ et elle n'est pas une fonction test. Le probleme doit etre pos\'e au sens des distributions, c'est dire prendre des fonction $ C^{\infty}_c $, smooth $ C^{\infty} $ a support compact. Ou bien au moins $ C^3 $, ce qui n'est pas vrai ici.

\smallskip

1) a) Pourquoi au moins $ C^3 $, car si le potentiel $ V \in C^s $, le potentiel newtonien est $ C^2 $, par soustration cela revient a chercher la regularit\'e d'une equation du type $ \Delta u = 0 $, or ici, il faut au moins qu'on derive l'equation pour avoir une regularit\'e $ C^2 $, par les Sobolev, c'est a dire $ H^3 $ en tant que Sobolev, c'est a dire que les fonction test doivent etre au moins $ C^3 $. Ce qui n'est pas le cas de la fonction propre qu'ils ont prise.

\smallskip

1) b) de maniere generale, quand on a un operateur d'orde $ k $, il faut une regularit\'e au moins d'orde $ k+1 $, pour les Sobolev pour avoir une regulrit\'e au final d'odre $ k $. En soustrayant un potentiel Newtonnien. Comme a l'ordre 2.

\smallskip

Pour la coherence de la th\'eorie. C'est pour cela qu'ils prennent $ C^{\infty}_c$ comme fonction tests.

\smallskip

-{\bf Donc: il faut prendre les fonctions $ C^{\infty}_c $. Ce qui n'est pas le cas de $ \phi^{\alpha} $}.

\smallskip

-Ici, il faut considerer l'article comme une illustration de la contradiction suivante: $ (P) $: ils utilisent l'integration par parties, car ils considerent le probl\`eme au sens des distributions, et, $ non (P) $: ils n'utilisent pas l'integration par parties. On en a parl\'e avant, dans les versions precedentes: on a $ (P) $ et $ non (P) $ dans la preuve, dans cet ordre dans la preuve: preuve erron\'ee.

\smallskip

///////////////////////////////

\smallskip

2) La formule d'integration par parties dans des domaines de $ {\mathbb R}^n $, ou la formule de Green-Riemann sur les vari\'et\'es. Ici, il y a un probleme de definition du gradient.

\smallskip

a) Soit: la formule de Green-Riemann, comme c'est prouv\'e dans le Az\'e. Le gradient est definit des le depart, par definition.

\smallskip

b) Soit: la formule de Green-Riemann, qu'on peut trouver dans la Gallot-Hulin-Lafontaine. A partir de la formule de Stokes. Il faut une notion de gradient, a partir de la metrique Riemannienne. 

\smallskip

3) Dans le livre de Az\'e Dominique ils prennent des fonctions $ f \in C^1(\bar \Omega) $, c'est une definition, c'est a dire qu'il existe une fonction reguliere dans une voisinage de $ \bar \Omega $, qui prolonge  $ f $, $ \tilde f \in C^1(V_{\bar \Omega}) $. Ils considerent les deriv\'es partielles qui sont bien definies sur $ \Omega $ puis prennent les primitives  et la limite ils ont les fonctions au bord.

\smallskip

Dans cette preuve, il n'y a pas besoin de definir un champ de gradient, car par les theoremes de Sobolev et les theoremes d'extensions, dans les problemes elliptiques, on a toujours ce cas la. C'est a dire des fonctions definies dans une voisinage de $ \bar \Omega $ qui prolongent les fonctions de depart.

\smallskip

4) La formule de Stokes, c'est bien ecrit dans le livre de Hebey. On se ramene localement a des demi espaces. Ici aussi la notion de primitive est utilis\'ee. Ils utilisent le th d'extension de Whitney, qui est valable dans les portions de demi-espaces. on se ramene a calculer les primitives de fonctions regulieres et on obtient les fonctions sur le bord qui coincdent avec les fonctions de depart.

\smallskip

La fomule de Green-Riemann, se deduit de la formule de Stokes, il faut avoir un champ de Gradient global, jusqu'au bord, car la vari\'et\'e n'est pas immerg\'ee dans un espace plus grand. Il faut definir un gradient partout, cela se fait a partir de la metrique et des cartes localement et on a gradient global.

\smallskip

On a: $ C^k(\bar \Omega)$, a plusieurs definitions: dans certains cas: c'est $ C^k(\Omega)$, avec des deriv\'ees qui ont une limite quand le point tend vers le bord. Cette definition est equivalente a celle qui dit que c'est la restriction de fonctions $ C^k(\Omega')$ avec $ \bar \Omega \subset \Omega'$, ouvert contenant $ \bar \Omega $, ceci est vrai lorsque $\Omega $ est un domaine regulier $ C^{k,\alpha} $. Voir, le gilbarg-trudinger, chapitre 6, sur les extensions.

\smallskip

En general, on prend la definition de l'extension a un ouvert plus grand, par exemple, quand $ \Omega $ n'est pas regulier.

\smallskip

Pour la formule de Stokes, sur un ouvert de $ {\mathbb R}^n $, on peut prendre les 2 definitions, il suffit qu'il y ai une limite quand les points tendent vers le bord. Il veut mieux prendre la definition avec les extensions, on est sur, qu'il y a une limite vers le bord.

\smallskip

Mais, en general, on prend la definition de l'extension. Comme dans l'article d'Agmon: "The $ Lp $ approach for the Dirichlet problem": car une fonction possedant des deriv\'ees sur $ \bar \Omega $, cette defintion est possible, si cette fonction possede des deriv\'ees sur un ouvert plus grand $ \Omega'$ contenant $ \Omega $. Il suffit de voir comment le th.6.2. est applique\'e, pour cela, il faut considrer le th.8.1.

\smallskip

La notion d'extension, est locale, par exemple pour les problemes de regularit\'e, par exemple au voisinage d'un point du bord, il suffit que le domaine soit $ C^{k,\alpha} $, au voisinage de ce point. Les 2 definitions sont equivalentes localement, des que, le domaine est regulier localement. D'ou le passage du local au global, en utilsant les 2 definitions. Par exemple, on peut, partir de la premiere definition pour, une demi-boule, un domaine Lipschitzien, et, avoir finalement la regularit\'e globalement avec la 2eme defintion, des extensions, globalement, sur l'ouvert regulier, ou la vari\'et\'e r\'eguliere.

\smallskip

Pour la definition 1, par exemple, pour la regularit\'e $ C^1 $ de la fonction, on utilise (comme on l'a dit pour la fonction de Green), qu'on a une fonction continue et Sobolev, si son Sobolev est continue, alors, elle est $ C^1 $. Pour $ C^k $, on it\`ere le proced\'e. Ceci est vrai, si le domaine de d\'epart est Lipschitzien.

\smallskip

Par exemple: pour prouver que la definition 1 est equivalente \`a la definition 2, ou pour utiliser la defintion 1, il faut une notion de carte reguliere, pour redresser le bord, et avoir des deriv\'ees directionnelles, par exemple: pour appliquer la formule de la moyenne \`a la fonction, et ses deriv\'ees, et passer \`a la limite, jusqu'au bord: ceci est possible, si on a des directions. Apres, si on suppose le domaine de depart $ C^k $, en se ramenant \`a une demi-boule, on peut peut faire l'extension, apres avoir prouver que la definition 1 implique que la fonction est $ C^k $, dans la demi-boule, suivant les directions existentes. Puis, on l'etant, comme dans le Gilbarg-Trudinger \`a la boule entiere.

\smallskip

Donc, on ne peut pas prouver qu'une fonction est $ C^k $, avec seulement, la definiton 1, car il faut des directions.

\smallskip

(Une fonction a des deriv\'ees continues sur $ \bar \Omega $, ssi, elle a une extension sur $ \Omega' \supset \bar \Omega $. Par exemple, concernant l'article d'Agmon: th.6.2. il faut voir comment il est appliqu\'e, quand on considere un domaine $ \Omega $, pour cela il suffit de se referer au th.8.1. on aura alors, la definition: avec extensions.)

\smallskip

{\it Solutions au sens faible: weak solutions:}

\smallskip

Soit $ L $ un operateur differentiel. On veut definir une solution $ u \in E $ faible d'une equation ou in\'equation. Pour cela on doit utiliser l'operateur $ L $, la seule possibilit\'e est de definir $ u \in E $, en utilisant l'adjoint $ L^* $, car $ L $ peut ne pas etre mis sous forme divergence. Il faut definir l'equation faiblement a partir de l'operateur $ L $, la seule donn\'ee est l'operateur $ L $. Dans ce cas, on utilise l'adjoint $ L^*$. Rien ne dit que l'operateur $ L $ peut etre mis sous forme de divergence. De maniere generale, pour un operateur $ L $.

Au sens des distributions:

$$ \int u L^* \phi dV = \int f \phi dV, \forall \phi \in C^{\infty}_c(\Omega), $$

C'est vrai si $ u \in E=L^1_{loc}(\Omega) $. Pour le laplacien, $ L^*=L =\Delta $, donc, on ecrit que $ u $ est solution faible si:

$$ \int u \Delta \phi dV = \int f \phi dV, \forall \phi \in C^{\infty}_c(\Omega), $$

Exemple: les fonctions de Green.

\smallskip

Lorsque $ u \in \tilde E=H_1^2(\Omega) $ et $ \partial \Omega \in C^{\infty} $, on peut integrer par parties et ecrire:

$$ \int \nabla u \cdot \nabla \phi dV = \int f \phi dV, \forall \phi \in C^{\infty}_c(\Omega), $$

-Il se peut que la definition de solutions faibles soit ce qui a \'et\'e dit dans cette derniere ligne, avec le gradient. Mais il faut prendre $ \phi \in C^{\infty}_c(\Omega) $. Or pour Aviles-McOwen, la fonction $ u\phi^q \not \in C^{\infty}_c(\Omega) $, pour cela il faut verifier qu'elle est $ \dot H_1^2(\Omega) $, or cet ensemble est definit par la notion de trace (au moins, si on n'arrive pas a prouver la densit\'e, il faut prouver que $ u\phi^q=\lim_{i} u_i, u_i \in C^{\infty}_c(\Omega) $), ce qui fait appel a la notion de trace donc de bord regulier.

\smallskip

-Il y a une difference entre les vari\'et\'es et les ouverts de $ {\mathbb R}^n $, par densit\'e et convolution, il suffit que les fonctions tests soient $ C^1_c(\Omega) $ pour $ \dot H_1^2 $ sans passer par les traces directement. Mais pour une vari\'et\'e quelconque de bord pas regulier, cette densit\'e n'est pas forc\'ement possible \`a prouver. 

\smallskip

-Il faut prendre aussi $ \phi\in H_1^2(\Omega) $, car l'espace total qu'on considere est $ \tilde E=H_1^2(\Omega) $. Donc, si on veut que la formulation soit vraie aevc $ \phi \in C^1_c(\Omega)$, il ne faut pas oublier de prendre aussi $ \phi \in H_1^2(\Omega)$. ( c'est pour cela, qu'ils considerent directement $ C^{\infty}_c(\Omega)$: \`a support compact en etant smooth $ C^{\infty} $. Or ici, la fonction "test" n'est pas $ C^{\infty} $. Il faut prouver qu'elle limite d'une suite $ C^{\infty} $, ce qui est pas tres possible, ou bien considrer des cartes, ce qui veut dire que la vari\'et\'e doit etre reguli\`ere.)

\smallskip

-Donc, les fonctions tests doivent etre dans $ C^{\infty}_c(\Omega) $ pour une vari\'et\'e quelconque et pour utiliser $ u\phi^q $ comme fonctions test, il faut passer par la notion de trace, donc, de bord regulier.

\smallskip

-Soit on doit appliquer la formule de Stokes, il faut un bord regulier. Soit, la densit\'e, prouver que la fonction est limite d'une suite reguliere smooth a support compact. Soit considerer, l'espace $ \dot H_1^2(\Omega)$, avec la notion de trace, ce qui n\'ecessite un bord regulier.

\smallskip

-Pour la formule $ \nabla (u\phi^{q})=u[\nabla (\phi^q)]+(\nabla u) \phi^q$: soit on la considere ainsi dans $ H_1^2(M)$ avec la notion de trace, donc de bord regulier.

\smallskip

Soit, on la prouve dans $ \dot H_1^2(M) $: c'est a dire pour la norme de $\dot H_1^2(M) $, par densit\'e c'est a dire qu'elle limite de fonctions $C_c^{\infty}(M)$, ce qui n'est pas possible, car $ u $ ne peut pas etre approch\'ee par des fonctions smooth \`a support compact.

\smallskip

Le seul moyen est de la considrer dans $ H_1^2(M) $ avec la notion de trace. Ce qui necessite un bord regulier.

\smallskip

Aviles-McOwen, dans leur theoreme d'estimation apriori, considerent le probleme au sens faible(avec des gradients, car ils considerent des problemes variationnels, $ H^1_2(\Omega)$ avec $ \partial \Omega $ r\'egulier) et non tres faible. Ils ne veulent pas utiliser la formule de Stokes.

\bigskip

{\bf Sur l'article de Aviles-McOwen et leur Theoreme de majoration locale et de compacit\'e locale:}

\smallskip

le Th. de majoration locale et de compacit\'e locale d'Aviles-McOwen:

\smallskip

1-Le probleme est pos\'e au sens des distributions: il faut faire une integration par parties pour se ramener au Th de Gilbarg-Trudinger. Pour cela il faut un bord lisse et la notion de trace pour appliquer la formule d'integration par parties ou formule de Stokes.

\smallskip

Aviles-McOwen, considerent le probleme au sens faible(avec les gradients, car ils considerent des problemes variationnels) et non tres faible. Ils ne veulent pas utliser la formule de Stokes.

\smallskip

a) $(\Omega,g) \in C^{\infty} $ avec $ \partial \Omega \in C^{\infty} $.

\smallskip

b) les solutions $ u \in H_1^2(\Omega) \Rightarrow $ trace, $ Tr(u) $.

\smallskip

2-Le resultat est local dans du global, sur tout $ \Omega $.(en particulier, ils obtiennent une sorte d'in\'egalit\'e de Sobolev \`a poids).

\bigskip

Pour ce qui nous concerne: on prend $ (M,g)=(\Omega,g) \in C^{\infty} $, avec un bord quelconque, qui n'est pas forc\'ement regulier. Que se passe t il ?

\smallskip

1-Dans le cas n\'egatif: sous-critique et critique et $ \partial \Omega \not \in C^{\infty} $ bord non regulier: 

\smallskip

On donne des conditions pour avoir un resultat local, pour qu'on conserve l'in\'egalit\'e locale.

\smallskip

a) $ u\in C^{2,\alpha}, \alpha >0$(pr\'efaisceau particulier), et $ \partial \Omega \not \in C^{\infty} $, bord non regulier: il n'y a plus de notion de trace et le bord n'est pas regulier. On ne peut pas utliser la formule de Stokes et conserver la nature globale du probleme.

\smallskip

Le probleme d'Aviles-McOwen(pos\'e au sens faible, au sens des distributions), est alors: pas bien pos\'e , n'est pas pos\'e correctement, il n'est pas defini.

\smallskip

b) On a donn\'e des conditions suffisantes pour avoir le resultat de compacit\'e locale.

\smallskip

c) Le resultat est local et obtenu de maniere locale: dans des cartes. Au debut de la proc\'edure on a: dessin + fl\`eche. cartes normales.

\smallskip

La th\'eorie des cat\'egories commence avec l'observation que de nombreuses propri\'et\'es des syst\`emes math\'ematiques peuvent \^etre unifi\'ees et simplifi\'ees par des dessins avec des fl\`eches.

\smallskip

d) La technique blow-up. On obtient une estimation a priori du type Harnack: un th de la moyenne a priori. Voir tout au debut du livre de Gilbarg-Trudinger, pour les fonctions harmoniques: {\it solid mean value inequalities}.

\bigskip

2-Pour ce qui concerne le cas n\'egatif supercritque et critique(non-linearit\'e exponenetielle), dimension $ n \geq 2 $, (aussi $ n=1 $): ici le bord peut etre regulier ou non: $ \partial \Omega \in C^{\infty} $ ou $ \partial \Omega \not \in C^{\infty} $, bord regulier ou non:

\smallskip

a) On donne des conditions pour avoir un resultat de majoration locale: $ u\in C^{2,\alpha}, \alpha >0 $ (pr\'efaisceau parcticulier).

\smallskip

Avec nos conditions, ce n'est pas le meme type d'equation. Pour nous c'est supercritique. De plus le Probleme d'Aviles-McOwen n'est pas bien definit ici pour l'Eq. supercritique et critique (non-linearit\'e exponentielle). 

\smallskip

Avec la condition $ \partial \Omega \not \in C^{\infty} $, le probleme d'Aviles-McOwen n'est pas bien pos\'e et n'est pas d\'efinit.

\smallskip

b) Dans ce cas (n\'egatif supercritque et critique, non-linearit\'e exponentielle), le probleme n'implique pas Gilbarg-Trudinger. Il faut faire un changment de fonction $ z=e^u $.

\smallskip

b) il faut faire l'integration par paries pour $ z $ et avoir la majoration locale de l'integrale de $ z^2 $, ce qui correspond a l'integrale locale de $ u^{\alpha} $ dans le probleme d'Aviles-McOwen.

\smallskip

c) Principe d'Harnack, les solutions $ u $ peuvent changer de signe.

\smallskip

d) Le resultat est local et obtenu de maniere locale: dans des cartes. ( On construit une sous-vari\'et\'e Riemmannienne compacte \`a bord $ C^{\infty} $ d'une vari\'et\'e Riemannienne, \`a partir de la carte: image reciproque d'une boule ferm\'ee de $ {\mathbb R}^n $, par exemple: smooth $ C^{\infty} $ compact Riemannian submanifold with boundary of a Riemannian manifold (of dimension $ n $ in dimension $ n $). Ce fait est connu pour les petites boules g\'eodesiques, dans l'existence de la premiere valeur propre des petites boules geoedesiques, voir par exemple: Karp-Pinsky, 1987. On a bien, par la carte exponentielle, une sous-vari\'et\'e Riemannienne compacte \`a bord plong\'ee dans une vari\'et\'e Riemannienne. Ici, c'est une carte quelconque. On admet que: l'existence d'une sous-vari\'et\'e compacte \`a bord plong\'ee dans une vari\'et\'e Riemannienne, est un fait connu, Karp-Pinsky.1987).

\smallskip

Dans ce point, on a: dessin + fl\`eche: La th\'eorie des cat\'egories commence avec l'observation que de nombreuses propri\'et\'es des syst\`emes math\'ematiques peuvent \^etre unifi\'ees et simplifi\'ees par des dessins avec des fl\`eches.

\smallskip

Pourquoi ce resultat n'a pas ete fait avant: il se peut que la raison soit la suivante: si on suppose que la solution est minor\'ee, $ u \geq -M >-\infty $, il suffit de considerer le cas o\`u $ u\geq 0 $, alors en multipliant par une bonne fonction test $ \eta ^4 $, on arrive directement a prouver que l'integrale de $ e^u $ est loc.uniform.born\'ee et par l'in\'egalit\'e de Harnack pour les fonctions sur-harmoniques, les solutions sont loc.uniform.major\'ees. Ce qu'on vient de voir est, qu'on a une relation du type $ \sup_K u \leq c(-M) $, ou bien $ \sup_K u \leq c(\inf u)$. C'est du Brezis-Merle. Comme Brezis-Merle ont trait\'e le cas positif, il se peut que le cas negatif eut ete evident et n'eut necessit\'e une preuve complete, car on le deduit des in\'egalit\'es de Harnack usuelles et de la methode de Brezis-Merle (integration et premiere fonction propre et premeire valeur propre). Ceci est peut etre la raison pour laquelle ce cas n'a pas ete trait\'e.

\smallskip

Mais ce qu'il faut voir est que, ce qu'on a prouv\'e, est plus precis, on n'a pas $ \sup_K u =c(\inf u)$, mais: $ \sup_K u \leq c $, le $ \sup $ est major\'e par une constante qui ne depend pas de l'inf. Effectivement, on l'a ecrit pour mettre en evidence cela, cette difference. Le $ \sup $ local est major\'e par une constante ne dependant pas de l'inf.

\smallskip

////////////////////////////////////////////////////////////////

\smallskip

{\bf Sur l'article de Aviles-McOwen:}

\smallskip

Pourquoi ils n'ont pas consider\'e, comme on l'a fait, dans une carte,(obtenue par un plongement), une sous-vari\'et\'e Riemannienne compacte a bord $ C^{\infty} $ de la vari\'et\'e Riemannienne $ (M,g) $ ?

\smallskip

1) Aviles-McOwen, ont consider\'e une varit\'et\'e compacte \`a bord $ (\bar \Omega, g) $: leur vari\'et\'e de depart $ M $ s'\'ecrit comme $ M=\cup_k \Omega_k $, une union denombrable de vari\'et\'es ouvertes (avec bords). Ceci est li\'e et est illustr\'e par leurs exemples \`a la fin de l'article o\`u ils considerent des cylindres: on voit bien que leur vari\'et\'e $ M $ est une union denombrable de cylindres, on prend un cylindre et on on le prolongent indefinement. Ils considerent le probleme, au sens des distributions, d'o\`u, ils ont besoin de la connexion \`a l'interieur des domaines, il n'y a pas besoin de bord. Sur les bord, il y a la notion de trace ou valeures limites.

\smallskip

Le fait d'ecrire $ M=\cup_k \Omega_k $: est une donn\'ee. Pour dire que ce qu'ils font sur  $ M $ est en fait ce qu'ils font sur chaque $ \Omega_k $. C'est une hypothese. Ceci est vrai pour un cylindre infini qu'ils decoupent en une suite exhausitve de compacts $ \Omega_k $. C'est aussi li\'e \`a leur exemples \`a la fin: des cylindres. 

\smallskip

C'est une hypothese d'ecrire: $ M=\cup_k  \Omega_k $ (les bords sont non compris, car le probleme se pose au sens des distributions, integration \`a l'interieur de $ \Omega_k $, avec la notion de trace ou valeures limites au bord): pour dire que ce qu'ils font sur $ M $ est en fait ce qu'ils font sur chaque $ \Omega_k $. La connexion de $ M $, est celle des ouverts $ \Omega_k $, pour tout $ k$, $ \nabla =\nabla^M=\nabla^{\Omega_k}, \forall k $ et donc: non pas la connexion sur tout $ \bar \Omega_k$, c'est \`a dire sans le bord. Car le probleme est pos\'e au sens des distributions, ce qui importe c'est \`a l'interieur du domaine, avec la notion de trace, ou valeures limites au bord. (Pour la trace, on a besoin de l'in\'egalit\'e de Sobolev a trace, qui s'obtient par densit\'e des fonctions $ C^{\infty}(\bar \Omega)$, localment dans des cartes, mais ce qui importe ce sont les valeures limites des fonctions au bord et la gradient $ L^p(\Omega)$ des fonctions, donc, le gradient \`a l'interieur de $ \Omega $, comme la formule d'integration par parties, l'inegalit\'e de Sobolev a trace, est une formule de trace, ce qui importe, sont les valeures \`a l'interieur, puis, on passe \`a la limite, pour obtenir les valeures au bord. $ \int_{\Omega} = \lim_{\epsilon\to 0} [\int_{\Omega_{\epsilon}} =_{integration \,\, reelle} \int_{\partial \Omega_{\epsilon}}] = \int_{\partial \Omega} $. Ceci veut dire qu'ils utilisent les valeurs \`a l'interieur, puis passent \`a la limite pour avoir les valeures au bord).

\smallskip

Ils(Aviles-McOwen) ne font pas d'integration par parties et ils considerent le probleme au sens faible (avec les gradients) et non tres faible. 

\smallskip

Ils n'utilisent pas le gradient jusqu'au bord (bord compris), la connexion de depart (sur l'espace total $ M $) est celle de chaque ouvert, sans tenir compte des bords.

\smallskip

Leur theoreme d'estimation a priori, ne n\'ecessite pas l'integration par parties et est fait sans notion de gradient jusqu'au bord (bord compris). Ici ils ont pris $ H^1_2(\Omega) $ avec $\partial \Omega $ regulier.

\smallskip

Cela veut dire:

$$ {\rm toute \, carte \, de \, } M \, {\rm en \, } x \, \Leftrightarrow \, {\rm est \, une \, carte \, de } \, \Omega_k \, {\rm en \,} x \in \Omega_k \, \forall k. $$
 
Ceci est vrai pour les cylindres. Et donc, tout coincide: la metrique sur $  \Omega_k $ est celle de $ M $ restreinte \`a  $ \Omega_k $. Il n'y a pas besoin de comparer les symboles de Christoffels en ce qui concenrne le bord de $ \Omega_k $, car ce qui est important, c'est la connexion \`a l'interieur du domaine $\Omega_k$ et elle coincide avec celle de $ M $, car c'est un ouvert de $ M $, il n'y a pas \`a regarder au bord pour ce qui est de la connexion(contrairement \`a la trace).

\smallskip

(les bords sont non compris, car le probleme se pose au sens des distributions, int\'egration \`a l'interieur de $ \Omega $ dans le theoreme de la premiere section, avec la notion de trace ou valeures limites au bord)

\smallskip

Donc, considerer $ M $, c'est equivalent \`a: considerer $ \Omega_k, \forall k $. C'est illustr\'e par les exemples: des cylindres.

\smallskip

2) Pour ce qui nous concerne: on a consider\'e une sous-vari\'et\'e Riemannienne compacte \`a bord $ \tilde  M $ (plong\'ee) dans une vari\'et\'e Riemannienne $(M,g)$: ici il y a la connexion $ \nabla^1=\nabla^M_{| \tilde M} $; de $ M $ restreinte \`a $
\tilde M $, (bord compris), et la connexion intrins\`eque de $  \tilde M $ en tant que vari\'et\'e Riemannienne: $ \nabla^2 $ (bord compris).(Il faut qu'il existe une et une seule connexion de Levi-Civita, elle doit etre unique, ici, il y a probleme de choix de connexion).

\smallskip

Ces deux connexions sont egales, en utlisant les symboles de Christoffels et le fait que les cartes de $ \tilde M $ (bord compris), sont construites a partir de la carte $(\Omega, \phi) $ du plongement.

\smallskip

Ici, la carte $ (\Omega, \phi) $ du plongement, est une carte de $ M $ mais pas n\'ecessairement celle de $ \tilde M $ (bord compris).

\smallskip

(En ce qui nous concerne), Ici, on fait l'integration par parties et n\'ec\'essit\'e du gradient jusqu'au bord (bord compris). Le raisonnement d'Aviles-McOwen, sans notion de gradient au bord, sans integration par parties, n'est plus possible, ici, $ \partial M $ n'est pas regulier $\Rightarrow $ les solutions doivent etre prises regulieres...etc...

\smallskip

/////////////////////////////////////

\smallskip

1) On a: {\bf on ne peut pas dire qu'on peut utiliser le theoreme d'Aviles-McOwen dans le th.2 de l'article de 2003} (apres avoir definit la sous-vari\'et\'e \`a bord compacte reguliere plong\'ee dans une vari\'et\'e Riemannienne, on admet que l'existence de cette sous-vari\'et\'e plong\'ee, est un fait connu, voir Karp-Pinsky. premiere valeur propre des petites boules g\'eodesiques.1987): car:

\smallskip

{\bf Theoreme d'Aviles-McOwen: $ \Leftrightarrow $ pas d'integration par parties}

\smallskip

Dans notre hypothese: les solutions $ u $ sont regulieres $ \Rightarrow $ pour se ramener \`a des solutions au sens $ H^1 $, il faut faire une integration par parties.

\smallskip

Donc: on a: ({\bf On utilise l'integration par parties}) et ({\bf on n'utilise pas l'integration par parties}). Ce qui veut dire qu'on fait: $ {\bf A} $ et $ {\bf non \, A} $: ce qui n'est pas possible. 

\smallskip

Dans cet ordre dans la preuve: $ A $ et $ non A $, ce n'est pas possible. 

\smallskip

Ce qui est parfois possible c'est de faire $ non A $ et $ A $: on n'utilise pas une proposition puis il y a un moment \`a partir duquel on peut l'utiliser. Mais ce qu'on ne peut pas faire, c'est: on utilise une proposition, puis apr\'es on n'utilise pas. Par exemple, il est important de voir comment les probl\`emes sont pos\'es: au sens des distributions, tr\'es faibelement, ou bien $ H^1 $, faiblement, ou bien au sens fort, ou bien dans un autre sens.

\smallskip

/////////////////////////////////

\smallskip

2) On regarde maintenant, si le th.1 de l'article de 2003 est compatible avec le th.2 de l'article de 2003:

\smallskip

2-1) Pour le th.1 de 2003, on a utilis\'e des potentiels $ C^{\alpha}, \alpha >0 $, uniform\'ement: condition sur les fonctions. Dans le th.2 de 2003, on n'a pas utlis\'e cette condition uniforme.

\smallskip

2-2) Pour le th.2 de 2003, on a utilis\'e la structure de vari\'et\'e de $ M $, on se ram\`ene \`a $ {\mathbb R}^n $, puis, on sait que les boules ferm\'ees de cet espace sont des sous-vari\'et\'es \`a bord, bord compris, smooth: condition topologique et geometrique. Puis on construit, la sous-vari\'et\'e Riemannienne plong\'ee dans la vari\'et\'e Riemannienne $ M $. Dans le th.1 de 2003, on n'a pas utlis\'e le bord de la vari\'et\'e \`a bord.

\smallskip

On a tendance \`a dire qu'il y a plus d'hypoth\`eses dans le th.1 que dans le th.2, de l'article de 2003. Ce n'est pas vrai, pour le th.1: donn\'ee de fonction, holderienne uniform., pour le th.2: donn\'ee topologique et geometrique, bord d'une vari\'et\'e \`a bord.

\smallskip

On voit que du point de vue des donn\'ees: les th.1 et th.2 de l'article de 2003, sont compatibles.

\smallskip

/////////////////////////////////////////////////////////////////////////////////////////////////////////////////////////////////////////

\smallskip

{\bf Remarques:}

\smallskip

Sur le premier article, 2003, "Differentes estimations de..."

\smallskip

Dans le print: "Quelques remarques dur les vari\'et\'es, fonction de Green et formule de Stokes" on a dit que lorsque $ \partial M $ n'est pas regulier, on ne peut plus utliser la formule de Stokes et le probleme d'Aviles-McOwen n'est pas definit:

\smallskip

Donc, il faut supposer les solutions $ u \in C^2 $.

\smallskip

Dans ce cas on plusieurs points de vues:

\smallskip

Le point de vue d'analyse fonctionnelle et EDP, et le point de vue groupes topologiques.

\smallskip

On suppose donc que les solutions $ u\in C^2 $: le but est de prouver des estimations a priori: dans $ L^{\infty}_{loc}, C^0_{loc}, C^1_{loc}, C^2_{loc} $.

\smallskip

a) Le point de vue groupe topologique: l'estimation dans $ C^0_{loc} $ suffit, pour avoir la notion de groupe topologique localement compact: c'est une bonne estimation.

\smallskip

b) Par contre, pour le point de vue analyse fonctionnelle, EDP, ou de fonctions, comme on a suppos\'e les solutions $ u \in C^2 $, il faut que l'estimation a priori cherch\'ee soit de la taille de l'hypothese, c'est \`a dire $ C^2_{loc}$. Or pour avoir $ C^2_{loc} $, il faut utiliser les estimations de Schauder, c'est a dire qu'il faut supposer $ C^{2, \alpha}, \alpha >0 $, or cela necessite que le potentiel $ V $ soit $ C^{\alpha} $ uniform\'ement.

\smallskip

c) Tout ceci pour dire que selon le point de vue, la regularit\'e forte uniforme du potentiel est n\'ecessaire.(Par exemple pour $ -\infty < a \leq V \leq b <0$, si on suppose par exmple qu'on peut avoir $ 0 > b \geq V_i \geq a_i \to -\infty $, il se peut qu'on ait uniform. sur une boule ouverte $ 0 >a_i/2 \geq V_i \geq a_i \to -\infty$, et en multipliant l'eq. par une fonction test et en integrant par parties, (il se peut aussi que les solutions $ u_i$ restent born\'ees), et on aurait, un terme qui tend vers $ -\infty $, et est egal a un terme born\'e uniform. Donc, pour qu'on soit sur qu'il n'y ait pas de contradiction, il faut supposer $ -\infty < a \leq V_i \leq b <0 $ uniform.)

Dans tous les cas $ V $ doit etre regulier $ C^{\alpha}, \alpha >0 $. Pour le point de vue groupe topologique, ce n'est pas n\'ecessaire de supposer la regularit\'e $ C^{\alpha} $ uniforme. Alors que pour le point de vue analyse fonctionnelle, EDP, et des fonctions, il est necessaire de supposer la regularit\'e $ C^{\alpha}, \alpha >0$, uniforme.

\smallskip

d) Quand on considerent les solutions $ C^2$ et on veut obtenir, des estimations $ C^2$, il faut considrer une eq et non une inequation. Ce n'est plus possible de considerer des inequations dans ce cas.

\smallskip

e) Quand on considere la regularit\'e $ C^2 $, une inequation peut s'ecrire comme une equation. En posant $ V= [(\Delta u+Ru)/(u^{q-1})] \in C^0$ ou $ V=[(\Delta u +R)e^{-u}] \in C^0, \Delta =-\nabla^i\nabla_i$.

\smallskip

En ce qui concerne le theoreme 1 de l'article de 2003. Dans le raisonnement par l'absurde: pour obtenir une contradiction: On utilise l'in\'egalit\'e de Holder au lieu de l'in\'egalit\'e de Jensen, car les fonctions peuvent toucher $ 0 $. On ne peut utiliser l'in\'egalit\'e de Jensen, car une fonction convexe, doit l'etre dans un intervalle ouvert, or le proc\'ed\'e diagonal, implique l'existence d'une fonction $ v \geq 0 $ sur tout $ {\mathbb R}^n $ solution de $ \Delta v= W(x_0) v^{q-1}, 2< q \leq N, \Delta = -\nabla^i\nabla_i $, apres on utilise les integrales superficielles, sur les spheres. On ne peut pas appliquer la formule de Jensen, car la fonction $ v $ peut parfois s'annuler, elle n'est pas strict.positive, et ses valeurs sont dans $ [0, +\infty[$ qui n'est pas ouvert, et la convexit\'e necessite des intervalles ouverts.

\smallskip

Les articles en question sont:

///////////////////////////////////////////////////////////////////////////////////////////////////////////////////////////////////////////

Quelques remarques sur les 2 premiers articles: "Differentes estimations...." de 2003 et "Majorations..." de 2004:

\smallskip

I) Remarques sur l'article "Differentes estimations": 2003: sur une vari\'et\'e Riemannienne:

\smallskip

1) th\'eor\`eme 1: cas n\'egatif sous-critique et critique, cas positif sous-critique:

a) Blow-up.

b) compacit\'e locale.

c) La vari\'et\'e peut poss\'eder un bord non n\'ec\'essairment regulier: smooth $C^{\infty} $ \`a l'interieur et pas forc\'ement smooth au bord.

d) Estimation a priori du type Harnack: th de la moyenne a priori.

\smallskip

2) th\'eor\`eme 2: cas n\'egatif critique (dimension 2) et supercritique:

a) majoration locale.

b) principe d'Harnack.

\smallskip

c) cas compact sans bord: existence par le degr\'e topologique de solutions topologiques: degr\'e =+1.(Pour le sous-critique positif aussi on a un resultat d'existence avec un potentiel $ V $, voir YY.Li pour le degr\'e =-1). {\bf C'est une cons\'equence.}

\smallskip

3) th\'eor\`eme 3: cas positif.

a) Minoration de $ \sup \times \inf $.

b) in\'egalit\'e optimale.

\smallskip

4) th\'eor\`eme 4 ou proposition: cas positif.

a) 2-sph\`ere.

b) Minoration de $ \sup +\inf $

c) in\'egalit\'e optimale.

\smallskip

II) Remarques sur l'article "Majorations...": 2004: sur un ouvert born\'e de $ {\mathbb R}^n $:

\smallskip

1) th\'eor\`eme 1 et corollaire 1: Cas positif, sous-critique tendant vers le critique. Eq.Lane-Emden-Fowler.

\smallskip

a) estimation asymptotique...Soit, le controle du sup: in\'egalit\'e asymptotique. Soit la constante de majoration: constante asymptotique.

b) corollaire 1: in\'egalit\'e du type $ (\sup)^{\alpha} \times \inf, \alpha \in (0,1) $.

c) crit\`ere pour avoir une in\'egalit\'e du type $ \sup \times \inf $: corollaire 1. Crit\`ere d'un nouveau type.(en general on a des critieres de compacit\'e).

\smallskip

2) th\'eor\`eme 2 et corllaire 2: cas critique:

a) in\'egalit\'e du type $ (\sup)^{1/3} \times \inf $ en dimension 3.

b) potentiel $ V $ lipschitzien non $ C^1 $ en dimensions 3 et 4.

c) avec une deuxieme contrainte:in\'egalit\'e optimale $ \sup \times \inf $ en dimension 4. 

d) Crit\`ere de compacit\'e: quand les constantes de Lipschitz $ A_i\to A >0 $. {\bf C'est une cons\'equence}.

e) Volume local = Energie locale born\'e(e). {\bf C'est une cons\'equence.}

f) $ \sup $ local major\'e  si $ \inf $ minor\'e par $ m >0 $ avec constante de Lipschitz $ A_i \to 0 $. Crit\`ere de compacit\'e. 

g) In\'egalit\'e de Harnack implicite pour une suite. {\bf C'est une cons\'equence}.

\smallskip

3) th\'eor\`eme 3: cas critique: cas radial:

a) in\'egalit\'e du type $ (\sup)^{\alpha} \times \inf $.

\smallskip

4) th\'eor\`eme 4: perturbation de l'equation par un terme sous-critique: Eq. de Schrodinger:

a) in\'egalit\'e optimale.

b) Energie locale born\'ee. {\bf C'est une cons\'equence.}

\smallskip

En ce qui concerne la minoration de $ \sup \times \inf$ de l'article de 2003: theoreme 3. D'apr\'es le probl\`eme de Yamabe, il y a 3 cas, selon que l'invariant de Yamabe $ \mu=0 $ ou $ \mu <0$ ou $ \mu >0$. Alors par un changement de metrique conforme, la condition $ 0 < a \leq V\leq b <+\infty $ et en int\'egrant l'equation apres changment de metrique conforme, on est forc\'ement dans le cas positif. On peut alors supposer que $ R=1$. 

\smallskip

Dans l'article de 2003, sur la minoration de $ \sup \times \inf $, on est forc\'ement dans le cas positif. Contrairement \`a l'article de 2006, "Estimations du type $\sup \times \inf $ sur une vari\'et\'e compacte": theoreme 1, o\`u on est oblig\'e de supposer que l'operateur est coercif.

\smallskip

On dit cela, pour dire qu'on n'ecrit pas la meme chose plusieurs fois.

\smallskip

Les articles en question sont:

\smallskip

///////////////////////////////////////////////////////////////////////////////////////////////////////////////////////////////////

\smallskip

{\bf Sur la dimension 4:} voir le print "Cas d'existence de solutions d'EDP". 

\smallskip

I) In\'egalit\'e de Harnack avec 2 contraintes liant $ u $ et $ V $, voir ci-dessus. Exemples: voir Brezis-Merle, corollaire 4, la contrainte de masse, $ \int_{\Omega} |V|e^u dx \leq \epsilon_0 < 4\pi $, il y a deux contraintes liant $ u $ et $ V $, l'Eq et la contrainte de masse, on a $ ||u^+||_{L_{loc}^{\infty}} \leq c <+\infty $, cela veut dire que si $ u\geq -M >-\infty $, on a la compacit\'e locale. Or ceci peut s'interpreter aussi comme $ \sup_K u\leq c(\inf_{\Omega} u) <+\infty $.

\smallskip

II) On a ecrit qu'il y a une in\'egalit\'e de Harnack pour \underbar {une suite} $ (u_i)$ relativement \`a une suite $ (V_i)$, telle que $ 0 < a \leq V_i \leq b<+\infty, ||\nabla V_i||_{\infty} \leq A_i \to 0 $:

$$ \forall K\subset \subset M, \exists +\infty > c(\cdot, K,M,g) >0, \,\, \sup_K u_i \leq c(\inf_M u_i, K, M, g), $$

la fonction $ m \to c(m,K,M,g)$ est d\'ecroissante de $ m >0$. Il se peut que la suite soit tr\'es longue, dense, comme $ {\mathbb Q} $ dans $ {\mathbb R} $.

\smallskip

Si on considere l'ensemble: $ E=\{u_i\in C^{2,\alpha}(M), u_i >0, \Delta u_i+\frac{1}{6} R_gu_i=V_i {u_i}^3, 0 <a \leq V_i \leq b <+\infty, ||\nabla V_i||_{\infty} \leq A_i \to 0\}$, ou tout simplement, $E=\{ u_i\} $, la suite, alors:

$$ \forall K\subset \subset M, \exists +\infty > c(\cdot, K,M,g) >0,\forall u\in E, \,\, \sup_K u\leq c(\inf_M u, K, M, g). $$

la fonction $ m \to c(m,K,M,g)$ est d\'ecroissante de $ m >0$.

Si on prend une autre suite $ v_i $ relativement \`a une suite $ W_i $ avec les memes hypotheses, en notant $ F=\{ v_i\} $, on a :

$$ \forall K\subset \subset M, \exists +\infty > \tilde c(\cdot, K,M,g) >0,\forall u\in F, \,\, \sup_K u\leq \tilde c(\inf_M u, K, M, g). $$

la fonction $ m \to \tilde c(m,K,M,g)$ est d\'ecroissante de $ m >0$.

\smallskip

Exemples:

\smallskip

1) Densit\'e, il se peut que l'ensemble $ E $ soit tr\'es "dense", comme $ {\mathbb Q} $ dans $ {\mathbb R}$. Et le compl\'et\'e de $ E $ soit consistant comme $ {\mathbb R} $ par  rapport \`a $ {\mathbb Q} $, ou les fonctions $ C^{\infty}_c$ \`a support compact dans les Sobolev. Il faut voir les compl\'et\'es pour les norme $ C^{2,\alpha}(M), C^{1,\alpha}(M), C^{0,\alpha}(M), C^0(M) $.

\smallskip

2)Les fonctions Harmoniques, on connait pas bien la taille de cet enesemble, il se peut qu'il soit de "mesure nulle" dans un espace plus grand. Comme $ {\mathbb Q} $ dans $ {\mathbb R}$. Bien qu'il soit consistant.

\smallskip

3) L'in\'egalit\'e de Harnack de Siu, pour l'Equation de Monge-Amp\`ere, portant sur les fonctions $ \phi_t$ index\'ee par un parametre $ t \in [\epsilon, t_0[$. Indexation et in\'egalit\'e de Harnack.

\smallskip

4) Cela depend de l'ensemble des solutions, ici, on a pris $ E $, on a bien une in\'egalit\'e de Harnack pour les \'el\'ements de $ E $.

\smallskip

\smallskip

{\bf Sur le print 2024: Harnack inequalities for equations of type prescribed scalar curvature:}

\bigskip

1) Blow-up, technique blow-up, ou technique eclatement et eclatement: plusieurs type d'eclatements.

\smallskip

2) Technique moving-sphere. symetrie sph\'erique. Enroulement.

\smallskip

3) Tenseur de Weyl: n'y apparait pas, n'y intervient pas, c'est du \`a l'eclatement et \`a l'enroulement. Pas de conditions sur le tenseur de Weyl. Eclatement et Enroulement.

\smallskip

4) Conditions de platitude du potenetiel $ V $: n'y apparait pas, n'y intervient pas, c'est du \`a l'eclatement et \`a l'enroulement. Pas de condition sur le potentiel $ V $. Eclatement et Enroulement.

\smallskip

5) Pas de changement de metrique conforme tel que $ Ricci(x_0)= 0 $, pas de changement de metrique conforme. c'est du \`a l'eclatement et l'enroulement.

\smallskip

6) Operateur $\Delta +h $ quelconque, non necessairement coercif + operateur principal: le laplacien.

\smallskip

7) On a bien une fonction $ m \to c(m) $ d\'ecroissante de $ m >0 $, mais cette fonction, qui, apres avoir fait l'eclatement, (apres avoir introduit le coefficient $ r_k $ de l'eclatement), $ c(r_k\inf_M u_k) $, n'est plus decroissante, on ne connait pas sa monotonie. Il y a rigidit\'e dans(le sens) la torsion (rigidit\'e \`a la torsion).

On a:

7)a) relation ouverte: $ \sup \times \inf $, pas d'enroulement dans la rigidit\'e. Enroulement et torsion.

\smallskip

7)b) relation semi-ouverte: $ \sup =c(\inf) $, \underbar {avec $ c $ decroissante}, enroulement dans la rigidit\'e. Enroulement.

\smallskip

7)c) relation semi-ouverte avec eclatement: $ (\sup)^{1-\epsilon}=c(\frac{\inf}{(\sup)^{\epsilon}})$, la fonction finale $ c $

\underbar{n'est plus decroissante}, rigidit\'e dans(le sens) la torsion et de l'enroulement. Pas de torsion, pas d'enroulement. Rigidit\'e.

\smallskip

7)b) est plus faible que  7)a) et 7)c) est encore plus faible que 7)b).

\smallskip

7)c): Rigidit\'e.

\smallskip

On voit l'effet de la courbure, scalaire, de Ricci, Weyl, de la non coercivit\'e de l'operateur: on n'a pas l'in\'egalit\'e optimale $ \sup \times \inf $. Mais on a une in\'egalit\'e de Harnack implicite interm\'ediaire.

\smallskip

Cas: orientable et non orientable:

\smallskip

Si $ M $ est orientable. Alors on a le theoreme 1.1.

\smallskip

Si $ M $ est non orientable, on considere le revetement universel $ \tilde M $ \`a 2 feuillets. Soit $ p: \tilde M \to M $ ce revetement \`a 2 feuillets et $ (u_k) $ une suite, solutions de l'equation de la courbure scalaire prescrite sur $ M $, alors $ u_k o p $ sont solutions de l'equation de la courbure scalaire prescrite sur $ \tilde M $. Si $ K $ est un compact d'interieur vide de $ M $ tel que $ \sup_K u_k \to +\infty $, alors: $ \tilde K=p^{-1}(K) $ est un compact (utiliser Bolzano Weirestrass et le fait que $ p $ est un revetement, par des suites extraites. Si $ (y_j) \subset \tilde K \Rightarrow (x_j)=(p(y_j)) \subset K $, on extrait de $ (x_j) $ une sous suite convergente dans $ K $ vers $ x_0 \in K $ et $ p^{-1}(x_0)=\{y_0, y_1\} \subset \tilde K $, puisque $ p$ est un diffeomorphisme local, on alors une sous suite de $ y_j $ qui converge vers $ y_0 $ ou $ y_1$). Ce compact est d'interieur non vide ($ p $ est localement un diffeomorphisme) et $ p(\tilde K)=K $ et $ \sup_{\tilde K} u_kop=\sup_K u_k \to +\infty $, donc on utilise le theoreme 1.1 pour avoir l'in\'egalit\'e sur $ \tilde M $ orientable et $ \inf_{\tilde M} u_kop=\inf_M u_k $, on remplace ces quantit\'es dans l'in\'egalit\'e. Donc, on a l'in\'egalit\'e sur $ M $ non orientable.

\smallskip

Pour le corollaire 1.2. Si $ M $ n'est pas orientable, on prend comme compacts: $ \tilde K=\bar B_r(y_0)\cup \bar B_r(y_1) $ pour le revetement \`a 2 feuillets $ \tilde M $ dans le theoreme 1.1. pour la suite $ u_kop $, sur $ \tilde M $. avec $ p^{-1}(x_0)=\{y_0,y_1\}, x_0\in M $. $ r >0 $ assez petit pour que le revetement $ p $ soit un diffeomorphisme de $ B_r(y_0)$ sur un voisinage de $ x_0 $ et aussi un diffeomorphisme de $ B_r(y_1) $ sur un voisinage de $ x_0 $.

\smallskip

Apr\'es: on a le corollaire 1.2. pour des voisinages:

$$ V_r^{x_0}=p(\tilde K)=p(\bar B_r(y_0) \cup \bar B_r(y_1)) \to_{r\to 0} \{x_0\}. $$

\smallskip

En fait, on utilise la preuve du theoreme 1.1 et du corollaire 1.2: \`a des compacts particuliers de $ \tilde M $ orientable, puis on revient \`a $ M $, en remplacant les quantit\'es.

\smallskip
 
Ceci est une premiere methode. Cette methode est la plus convenable.

\smallskip

Une deuxieme methode, consiste \`a prendre: 

\smallskip

Dans $ M $ non orientable: $ K=\bar B_r(x_0) $ avec $ \frac{inj(x_0)}{2}>r>0 $, assez petit pour que $ p $ r\'ealise un diff\'eomorphisme de voisinages de $ y_0, y_1 $ vers cette boule. On prend alors dans $ \tilde M $ orientable: $ \tilde K_r(y_0, y_1)= p^{-1}(\bar B_r(x_0)) $ et l'ensemble "maximum" est $ \tilde K_{r_0}(y_0,y_1) $. 

Le seul changement dans le corollaire 1.2. est de prendre le 

$$ 0 < \delta_0(\tilde K)\leq \frac{1}{8}\inf \{\tilde d[\partial (\tilde K_r(y_0,y_1)),\partial ( \tilde K_{r_0}(y_0, y_1))], \tilde d[\partial (\tilde K_r(y_0,y_1)), \partial (\tilde K_{r_0}(y_0, y_1))], ..., \}. $$

Afin d'avoir un infimum fini \`a la fin. Pour avoir de l'espace, pour pouvoir prendre l'infimum strictement \`a l'interieur de $ \tilde M $.

\smallskip

Par un dessin, cela ressemble \`a des deformations de boules.  On prend: la distance des bords de ces deformations de boules de rayons $ r>0 $ au bord de la  deformation de la boule "maximum", de rayon $ r_0 >0 $.

On a la meme preuve que dans le cas orientable, corollaire 1.2. avec des deformations de boules, par le diffeomorphisme local $ p $.

Ceci est un raisonnement un peu approximatif. Il faut prouver, cette derniere assertion, sur le minimum. Cette deuxieme methode est peu convenable.

\smallskip

///////////////////////////////////////////////////////////////////////////////////

\smallskip

En relation avec l'exemple de C.C.Chen-C.S.Lin: $ \sup $ major\'e si $ \inf $ minor\'e, cette assertion n'est pas toujours vraie, voir l'exemple dans l'article de C.C.Chen-C.S.Lin.1998. Journal.Diff.Geometry.

\smallskip

On va essayer d'expliquer, l'in\'egalit\'e:  $ (\sup_K u)^{1-\epsilon} \leq c[\frac{\inf_M u}{(\sup_K u)^{\epsilon}}]$ du print de 2024.

{\bf Critere de compacit\'e:} Si $ (u_k) $ est une suite qui verifie l'in\'egalit\'e de Harnack usuelle en $ x_0 $, alors, elle verifie le point 1) du print de 2024.

In\'egalit\'e de Harnack usuelle au sens suivant: 

$$ \forall (x_k), (y_k), x_k \to x_0, y_k \to x_0, \exists m >0:  u_k(x_k)\geq m \cdot u_k(y_k). $$

Alors en appliquant ceci aux points $ (t_k)$ et $ (\tilde t_k), \inf_{B(t_k, r_k^{2/(n-2)})} u_k= u_k(\tilde t_k) $, on ne peut pas avoir $ u_k(t_k) \to +\infty $. Donc, on a le point 1) du print de 2024. Alors, il existe un voisinage de $ x_0 $ relativement compact, $ K_0 $ et une sous-suite $ i_j $ et un reel $ C >0 $, tels que:

$$ \sup_{K_0} u_{i_j} \leq C. $$

Donc: on a le crit\`ere de compacit\'e suivant: Si la suite $ (u_k) $ verifie l'in\'egalit\'e de Harnack usuelle en $ x_0 $, alors, on a compacit\'e autour de $ x_0 $.

\smallskip

Donc, on n'a pas forc\'ement $ \sup $ major\'e si $ \inf $ minor\'e, comme on l'a dit et ce que disent Chen-Lin avec leur exemple. Mais, on a le critere de compacit\'e autour de $ x_0 $, si la suite verifie l'in\'egalit\'e de Harnack usuelle en $ x_0 $.

\smallskip

Donc: si $ (u_k) $ est une suite qui verifie en $ x_0 $ l'in\'egalit\'e de Harnack usuelle, alors, on a compacit\'e autour de $ x_0 $.

\smallskip

En conclusion: Donc: on a l'in\'egalit\'e de Harnack, $ \sup^{1-\epsilon} = c(\frac{\inf}{\sup^{\epsilon}}) $. De plus, on a le crit\`ere de compacit\'e suivant, et le resultat de compacit\'e suivant: si en $ x_0 $, la suite $ (u_k) $ verifie l'in\'egalit\'e de Harnack usuelle, alors, on a compacit\'e autour de $ x_0 $.

\smallskip

Ce n'est pas seulement un crit\`ere de compacit\'e, mais un resultat de compacit\'e, car on dispose de l'exemple suivant: $ x\to [\epsilon/(\epsilon^2+|x|^2)]^{(n-2)/2} $, cette suite de fonctions ne verifie pas l'in\'egalit\'e de Harnack usuelle en $ 0 $. Elle diverge en $ 0 $.

\smallskip

{\bf Important:}

\smallskip

1) L'in\'egalit\'e de Harnack est obtenue sans concentration de la mesure. Pas de notion de concentration de la mesure \`a l'interieur. Pas de notion de concentration de la mesure au bord. On ne sait pas si il y a concentration de la mesure ou pas.

\smallskip

2) L'in\'egalit\'e de Harnack est obtenue sans la notion de fonction de Green.

\smallskip

3) On reste dans le cadre de la "metric Geometry", car le potentiel est Lipschitzien pour la distance g\'eod\'esique.

\smallskip

4) On voit l'effet des courbures, scalaires, de Ricci, de Weyl, du potentiel $ V $, de la non coercivit\'e: on n'a pas l'in\'egalit\'e optimale $ \sup \times \inf $. L'\'energie ou le volume ne sont pas n\'ecessairement born\'es. (Ceci est une simple remarque: car l'in\'egalit\'e optimale n'est pas necessairement possible avec autant de contraintes: courbures, scalaire, de Ricci, de Weyl, et le potentiel $ V $. Meme si on met des conditions fortes sur le potentiel $ V $, il y aura les contraintes de depart: les courbures, scalaires, de Ricci, de Weyl. Aussi, meme si on suppose que la courbure de Ricci est nulle, il reste la contrainte, de la courbure de Weyl et aussi la contrainte de coercivit\'e de l'operateur. L'equation de Yamabe le montre bien, ici, le potentiel est $ V\equiv 1 $, il reste les contraintes de courbures, scalaire, de Ricci, de Weyl, et la coercivit\'e de l'operateur conforme). (Pour l'in\'egalit\'e $ \sup \times \inf $ optimale et volume local=energie locale born\'es uniform. voir, la partie g\'eom\'etrisation).

\smallskip

Mais ici, on a une in\'egalit\'e de Harnack et le crit\`ere de compacit\'e associ\'e: in\'egalit\'e de Harnack usuelle.

\smallskip

A cette in\'egalit\'e de Harnack implicite, on associ\'e un crit\`ere de compacit\'e. In\'egalit\'e dans le sens de la majoration. Controle du supremum local.(en ayant \`a l'esprit l'exemple: $ x \to [\epsilon/(\epsilon^2+|x|^2)]^{(n-2)/2} $).

\smallskip

5) Ceci inclut les eq.Yamabe, de courbure scalaire prescrite, d'Einstein-Lichnerowicz, de Schrodinger.

\smallskip

Notions utilis\'ees: 

\smallskip

Blow-up: eclatement math\'ematique. Changement d'echelle: rescaling: eclatement de premi\`ere esp\`ece. technique blow-up de Schoen modifi\'ee.

\smallskip

Alternative. Supremum local divergeant.(permet de definir l'hypothese de l'absurde apr\'es). 

\smallskip

Ici, il y a un probl\`eme de choix: il faut definir une hypoth\`ese de l'absurde, qu'on ne connait pas, et, conserver le $ c_k\to +\infty $, qui permet d'appliquer la technique blow-up de Schoen modifi\'ee et Caffarelli-Gidas-Spruck.

\smallskip

Reduction des supremums \`a des valeurs ponctuelles: rescaling.(permet de definir les objets: points et exponentielles: pour definir l'hypothese de l'absurde apr\'es rescaling). Et conserve le $ c_k \to +\infty $ de la technique blow-up de Schoen modifi\'ee. 

\smallskip

Rayon d'injectivit\'e. Fonction rayon d'injectivit\'e, carte exponenetielle.

\smallskip

Raisonnement pas l'absurde: hypoth\`ese de l'absurde: apr\'es rescaling: infinit\'e d'infimum minor\'es. Quantit\'e major\'ee si $ \inf $ minor\'e.

\smallskip

R\'esultat de classification de Caffarelli-Gidas-Spruck.

\smallskip

Le rescaling: le blow-up, l'eclatement, permet de reduire l'effet, d'annuler l'effet, des courbures, scalaire, de Ricci, de Weyl, et du potentiel $ V $.

\smallskip

"moving-sphere" method, deplacement de sphere: principe du maximum sur une vari\'et\'e diff\'erentielle compacte \`a bord, cylindrique, Lemme de Hopf. Op\'erateur elliptique sur une vari\'et\'e diff\'erentielle. Fonction auxili\`ere. Exposant critique de Sobolev.

\smallskip

Op\'erateur de Laplace-Beltrami sur une vari\'et\'e Riemannienne. Coordonn\'ees g\'eodesiques polaires: ces coordonn\'ees permettent de definir la fonction blow-up ($ w = e^{(n-2)t/2} u(e^t\theta) $) li\'ee \`a l'exposant critique de Sobolev.

\smallskip

/////////////////////////////////////////////////////////////////////////////////////////////////////////////////////////

\smallskip

{\bf Sur les Pbs du type Brezis-Merle et la compacit\'e (eq  de courbure scalaire)et le cas negatif pour l'Eq. de la courbure scalaire prescrite: groupes topologiques}

\smallskip

{\bf On regarde les resultats de compacit\'e d'une autre maniere: un autre point de vue}

\smallskip

1) a) On peut considerer le groupe conforme $ C(\Omega)=\{ f \in C^2(\Omega, \Omega), f^*(g)=e^u g \}$. Ici, $ C^2(\Omega,\Omega) $ est l'ensemble des transformations conformes  de $ \Omega $ vers $ \Omega $, avec condition de Dirichlet. Sans nuire \`a la generalit\'e: se donner $ f $ est equivalent \`a la donn\'ee de $ u $ avec condition de Dirichlet.

\smallskip

Dans le cas $ g=\delta $, {\bf Cet ensemble ($C(\Omega)$) n'est pas bien defini, car $Scalaire(f^*(\delta))= Scalaire(\delta)of\equiv 0 $, alors que $ Scalaire (e^u \delta)\not =0$.} voir le livre de Hebey, isometrie riemannienne.

\smallskip

On suppose par exemple, que $ C(\Omega)$ est bien defini(ce qui n'est toujours pas le cas), on a alors: resultat de compacit\'e = injection compacte dans ce groupe topologique.

\smallskip

Cet ensemble est un groupe. Muni de la metrique ou norme essentielle locale ($ sup $ local), comme pour les disributions, il suffit de regarder au voisinage de l'element neutre, qui est l'identit\'e. Ce groupe devient topologique.(en r\'ealit\'e la topologie r\'eelle est la topologie compacte-ouverte, voir le livre de Hebey).

\smallskip

La compacit\'e pour $ u >0 $ avec condition de Dirichlet, implique la compacit\'e pour $ f $. Ceci implique que ce groupe topologique devient compact. C'est un groupe de Lie. On a alors une mesure de Haar. 

\smallskip

b) Dans le cas de l'Eq. de la courbure scalaire prescrite (cas n\'egatif), on a la compacit\'e locale des $ u >0$ dans $ C(M)=\{f \in C^2(M,M), f^*(g)=u^{4/(n-2)} g, u >0\}$, implique la compacit\'e locale des $ f $. Le groupe est topologique (comme pour les distributions, au voisinage de l'element neutre, l'identit\'e, puis par multiplication), ce groupe topologique devient localement compact. d'ou l'existence d'une mesure de Haar. Il est alors possible de faire de l'analyse sur ce groupe de Lie localement compact.

Le cas n\'egatif de l'eq. de la courbure scalaire prescrite: On note $ R_g $ la courbure scalaire

Le groupe topologique:

$$ C(M) \,\,{\rm est \,\, {\bf un \,\,  groupe \,\, de \,\, Lie}\,\, {\bf localement}\,\, {\bf compact}\,\,} $$

$ R_{f^*(g)}=R_gof $, si $ R_gof\in [a,b]\subset ]-\infty, 0[, ||\nabla (R_gof)||\leq A $, $ C(M) $ en tant groupe topologique et la compacit\'e locale implique que c'est un groupe localement compact.

$$ f\in C(M) \Rightarrow \exists u >0, \,\, f^*(g)=u^{4/(n-2)} g $$

On suppose que les $ u >0 $ sont tels que: $ R_{u^{4/(n-2)}g} = V $ des hypotheses de compacit\'e locale: $ V\in [\tilde a,\tilde b]\subset ]-\infty, 0[, ||\nabla^{\alpha} V||_{\infty}\leq \tilde A $. Pour chaque diffeomorphisme $ f $, il existe un $ u >0 $ tel que $ V $ soit comme ci-dessus. Il se peut qu'il existe plusieurs $ u >0 $, mais on suppose que parmi les $ u >0 $ il y en a un tel que les hypotheses sur $ V $ ci-dessus soient vraies. Alors, il y a compacit\'e locale et le groupe topologique $ C(M) $ est alors localement compact.(Il se peut que $ C(M) $ soit reduit a un seul element $ \{id \}$, qui est compact, donc on peut supposer qu'il a d'autres elements).

\smallskip

Pour voir si il y a compacit\'e locale, on se place au voisinage de chaque element $ f_0 \in C(M) $, ce qui veut dire les conditions sur $ R_gof$ sont verifi\'ees pour $ f $ dans le voisinage de $ f_0 $, alors les $ u >0 $ qui proviennent des $ f $ dans un voisnage de $ f_0 $ fix\'e, sont localement compact. On a la compacit\'e avec les conditions suivantes sur $ R_g $: $ R_g \in [a',b']\subset ]-\infty, 0[, ||\nabla R_g||\leq A' $. Alors au voisinage de $ f_0 \in C(M) $, les conditions sur $ R_gof $ ci-dessus sont v\'erifi\'ees. 

\smallskip

D'apr\`es le livre de Hebey. La topologie sur $ C(M) $ est la topologie compacte-ouverte. Or, ce qu'on a dit precedemment concerne la topologie $ C^1_{loc}(M)$, or d'apres Ascoli-Arzela, on voit qu'on a directement la compacit\'e si on met la topologie $ C^1_{loc}(M)$. Il n'y a pas besoin du resultat de comapcit\'e locale des metriques conformes.

La topologie r\'eelle est la topologie compacte-ouverte. La compacit\'e locale est vraie si on arrive a prouv\'e la compacit\'e dans $ C^0_{loc}(M) $. Alors, les conditions sur $ R_g $ sont:$ R_g \in [a',b']\subset ]-\infty, 0[ $. (Il ne faut pas supposer $ C^1$, car on a la compacit\'e locale directement, sans passer par les metriques conformes).(voir, cas negatif, Bahoura).C'est une consequence.

\smallskip

On a aussi: sa s'applique aussi \`a:

\smallskip

a) En dimension 2, dans le cas n\'egatif, sur une surface Riemannienne de courbure $ R_g \in [a',b'] \subset ]-\infty,0[ $. L'eq est : $ \Delta u+R_g=Ve^u, \Delta=-\nabla^i(\nabla_i), V=R_gof\in [a',b']\subset ]-\infty,0[ $, ici, la surface n'est pas forc\'ement sans bord, de bord non necessairement regulier. Ici, $ R_g $ n'est pas constante negative forcement, elle peut l'etre si la surface n'est pas ferm\'ee ou de bord non regulier.(voir,cas negatif, Bahoura). C'est une consequence.

\smallskip

b)Ca ne marche pas forc\'ement pour un potentiel quelconque, car $ V=R_gof $ et il faut controler $ f $, ce qui n'est pas possible. Sauf le cas constant: $ R_g\equiv constante < 0 $. En particulier c'est vrai pour le resultat de YY.Li, sur une surface Riemannienne compacte sans bord de courbure scalaire constante stric.negative.(voir, cas negatif constant, surface ferm\'ee, YY.Li, Harnack inequality, the method of moving-plane 1999). Ou bien voir le livre d'Aubin, dans ce cas, (surface ferm\'ee et 
 $ R_g=constante <0 $, la solution est unique, il y a une seule solution, c'est $ u=0 $, dans ce cas c'est le groupe d'isometries, la surface est compacte sans bord, il est compact). C'est une consequence.

\smallskip

c) Pour ce qui est du resultat de compacit\'e locale de Brezis-Merle. Il faut consider une surface Riemmanienne, une metrique $ g $ dans le cas positif (toute surface Riemannienne est loc.conf.plate, on prend un voisinage dans lequel on peut appliquer le cas du changement de metrique conforme, conf.plat):

On a:

$$ \phi^*(g)=e^v \delta, \,\, C(\Omega)=\{f \in Diff (\bar \Omega), f^*(g)=e^u g,\}$$

avec condition de Dirichlet. Avec $ R_g \in [a,b] \subset ]0,+\infty[ $.
c'est une consequence.
\smallskip

d) Avec ce proc\'ed\'e, c'est vrai pour l'Eq. de Yamabe.avec condition de Dirichlet $ u=1 $ sur le bord. Ou vari\'et\'e compacte sans bord, (voir cas positif, YY.Li, M.Zhu, L.Zhang, C.C.Chen-C-S.Lin, T.Aubin, Bahoura).C'est une consequence.

\smallskip

Ce proc\'ed\'e de comapci\'e de diffeomorphismes, fonctionne si on considere une vari\'et\'e  ferm\'ee ou une vari\'et\'e ouverte. Par contre pour une vari\'et\'e \`a bord, il y a un probleme, soit on a la cmapcit\'e sur toute la vari\'et\'e y compris au voisinage du bord, soit c'est pas possible. Car, on ne peut pas considrer le bord, puis se limiter \`a l'interieur, car une base de voisinage qui engnedre les voinsiange de $ f_0\in Diff(\bar \Omega) $, peuvent contenir une partie du bord, donc, il faut avoir la compaict\'e au voisinage du bord (pour la metrique de la vari\'et\'e, meme en dimension 2, meme si loc.conf.plat, d\`es qu'on se place au voisnange du bord, on n'a pas forcement une metrique du voisinage du bord plate), et si on considere le bord, il faut avoir la compacit\'e au voisnage du bord.

Pour une vari\'et\'e ferm\'ee $ M $(sans bord), la topologie de $ Diff(M) $ est connue, il suffit d'avoir la compacit\'e dans $ C^0(M,M) $.

Pour une vari\'et\'e ouverte, de meme, on $ C^0_{loc}(M,M)$.

Par contre pour une vari\'et\'e \`a bord, on ne peut pas considerer le bord inclus, et faire de la convergence dans $ C^0_{loc} $. Car on considere $ f_0 \in Diff(\bar \Omega) $, le bord est pris en compte, donc, pour la convergence, on a besoin de l'information au voisinage du bord.

e) En dimension 2, une surface compacte ferm\'ee c'est la sphere, il se peut que ca blow-up, alors que, par exemple, le resultat d'Aubin, la compacit\'e, c'est precisement, pour une vari\'et\'e ferm\'ee non conf.diffeo. a la sphere (en dimension $ n\geq 3 $). En dimension 2, il n'y a qu'une possibilit\'e, c'est precisement la sphere de dimension 2.

De plus, dans le cas positif, on le sup depend de l'inf. Et c'est la sphere. le probleme ne se pose pas directement en prenant $ R_g \in [a,b] \subset ]0,+\infty[ $.

f) On a dit que le probleme ne se pose pas aussi, avec des conditions sur le gradient de $ R_g $, car, pour un diffeomorphism, $ f $, il faut controler $ f $ en norme $ C^1 $, c'est pas possible. Or les conditions, doivent etre prise sur $ V=R_gof $, pour qu'on ait des conditions sur $ V $ d'ordre $ C^1 $, il faut des conditions d'ordre $ C^1$ sur $ f $, or ceci n'est pas forcement possible. La seule possibilit\'e d'utiliser $ f $ est que $ R_g\equiv constante \not =0 $, alors $ R_gof=constante\not =0 $, alors $ \nabla (R_gof)=0 $.

g) Les conditions les plus naturelles sont $ R_g \in[a',b'] $, dans le cas positif en dimension 2, ca ne suffit pas et ca ne se pose pas, car $ \sup $ depend de $ \inf $, et c'est sur la sphere, il se peut que ca blow-up.

Par contre en dimension $ n\geq 3 $, c'est possible sur une vari\'et\'e ferm\'ee non conf.diffeomorph. \`a la sphere, $ R_g\equiv 1 $, Eq. de Yamabe, resultat de, YY.Li-M.Zhu, Druet, Marques, Schoen, Aubin.(loc.conf.plat et non loc.conf.plat, non conf.diffeo.\`a la sphere).

\smallskip

h) Il y a le cas negatif o\`u c'est possible avec la condition $ R_g\in [a',b'] \subset ]-\infty,0 [ $ et(ou) en particulier, $ R_g\equiv constante= a'=b' < 0 $, en dimension 2 et en dimension $ n\geq 3 $. (voir, cas n\'egatif, vari\'et\'e non necessaitement \`a bord regulier en dimension $ n=2 $ et $ n\geq 3 $, Bahoura). C'est une consequence.

\smallskip

i) Brezis-Merle et C.Li-W.Chen, en prenant une metrique de la forme: $ \phi^*(g)=e^v \delta $ avec $ R_g\equiv k=constante >0 $. avec condition de Dirichlet. L'equation devient, $ f^*(g)=e^u g $:

$$ \Delta_g u+R_g=R_gof e^u, $$
avec condition de Dirichlet. Ici, $ R_g=k>0 $. Il y a au moins une solution $ f=id $ et $ u=0 $.

$$ \Delta_{\phi*(g)} (uo\phi)+(R_go\phi)=(R_gofo\phi)e^{uo\phi}, $$

avec condition de Dirichlet. On peut ecrire:

$$ \Delta_{\phi*(g)} (uo\phi)+k=ke^{uo\phi}, $$

avec condition de Dirichlet. Avec $ \Omega $ domaine $ C^2 $ au moins. Comme $ \phi^*(g)=e^v \delta $.

L'equation devient:

$$ \Delta_{\delta} (uo\phi)+ke^v=ke^v e^{uo\phi}, $$

avec condition de Dirichlet. On elimine le terme $ ke^v $ en resolvant un probleme de Dirichlet.

Cette nouvelle equation sur un domaine au moins $ C^2 $ de $ {\mathbb R}^2 $, permet d'avoir la compacit\'e jusqu'au bord.

Donc le groupe conforme $ C(\Omega) $ est compact, si le th. de Myers-Steenrod est vrai pour la vari\'et\'e \`a bord, y compris le bord.

\smallskip

Dans l'article, "Compactness of the set of solutions to elliptic equations in 2 dimensions", le domaine $ \Omega $ est \'etoil\'e par rapport \`a un point et de classe au moins $ C^2 $. On a la compacit\'e dans ce cas.

\smallskip

{\bf Remarque:} ici aussi il y a un probleme avec le bord. On ne sait pas si le theoreme de Myers-Steenrod est vrai pour une vari\'et\'e \`a bord, bord compris. Donc, ce resultat, la compacit\'e de $ C(\Omega) $, n'est pas vrai(e) pour les domaines \`a bord, bord compris.

\smallskip

2) Pour ce qui concerne les resultats sur les problemes du type Brezis-Merle, d'habitude on considere des espaces vectoriels norm\'es complets ou Banach, ici, pour ce type d'eq. Il faut considerer des groupes topologiques. On a alors l'injection compacte de $ e^u\ni L^1 \hookrightarrow L^{\infty} $, par abus de language, considerer $ u >0 $ c'est considerer $ f $ la transformation conforme.

Donc, pour ce qui concerne le Probleme du type Brezis-Merle on a: 

Dans le groupe topologique $ C(\Omega) $, qu'on suppose bien defini (voir point i) ci-dessus, ce qui pr\'ec\`ede):

En termes de $ u $:

$$  L^1(e^{(\cdot)}) \hookrightarrow L^{\infty}\,\, {\rm est \,\, une \,\, injection \,\, compacte} $$

ou bien en termes de $ f $:

$$ 
  L^1(Jac(\cdot)) \hookrightarrow L^{\infty}\,\, {\rm est \,\, une \,\, injection \,\, compacte} $$

\smallskip

Dans le print "About Brezis-Merle problem with Lipschitz condition". Le domaine est au moins de classe $ C^2 $, s'il n'est pas etoil\'e par rapport \`a un point, on ne borne pas forc\'ement le volume, mais on a l'injection compacte. En termes de $ u $ il y a compaci\'e. Mais en termes de $ f $ il n'y a pas de comapcit\'e. 

\smallskip

Pour la structure de groupe des $ u $:

{\bf Remarques:} 

1) Chaque $ u $ est associ\'e \`a un $ f $ un diffeomorphisme, $ 0 $ est associ\'e \`a $ id $ ($ n=2 $), $ 1 $ est associ\'e \`a $ id $ ($ n\geq 3 $):

$$ u_1\star u_2=u_2of_1+u_1,\,\, n=2,\,\, {\rm et}\,\, u_1\star u_2=u_2of \cdot u_1,\,\, n\geq 3. $$

La loi de composition $ \star $ n'est pas forc\'ement commutative.

On a: $ u\star v=vof+u $, si on essaie de prouver que le groupe est topologique, on doit voir si $ (u,v)\to u\star v $ est continue, or il y a le diffeomorphisme $ f $ qu'on ne peut pas forc\'ement controler. Soit en $ L^{\infty}_{loc} $, donc c'est li\'e a la compacit\'e au voisinage du bord. Soit en $ L^{\infty} $, c'est li\'e au th. de Myers-Steenrod. De meme pour l'application inverse, $ u\to u^{-1}=-uof $.

Aussi, on ne peut pas forc\'ement le prouver localement, (quand on a $ ||df||\leq c <+\infty $), car quand on ecrit $ f(B_1) \subset B_2 $, pour des voisnages, $ B_1 $ de $ x_0 $ et $ B_2 $ de $ f(x_0) $, on a pas forc\'ement des voisinages geodesiquement convexes, car:

a) les metriques changent, quand on les restreint \`a $ B_1 $ et $ B_2 $, $ i_1^*(g), i_2^*(g) $ avec, $ i_1: B_1\to \bar \Omega, x\to x, i_2: B_2 \to \bar \Omega, x\to x $, ce n'est plus la metrique de depart $ g $, alors qu'il faut considerer la metrique $ g $. 

\smallskip

b) au voisnage du bord, on n'a pas le th. des voisinages geodesiquement convexes.

\smallskip

2) On regarde le groupe comme "g\'en\'er\'e" par les elements $ u $. Les metriques conformes. C'est un abus d'ecriture: on considere les $ u $ associ\'ees aux $ f \in C(\Omega) $, ou bien les metriques confromes.

\smallskip

3) Comme on l'a deja dit, en dimension 2, une surface compacte est isomorphe \`a la sphere. Il n'y a qu'une possibilt\'e la topologie et la g\'eometrie de la sphere, pour laquelle il y a l'in\'egalit\'e $ \sup, \inf $. On peut dire que dans le cas positif, en dimension 2, ca ne se pose pas, la compacit\'e, il se peut que ca blow-up. En dimension $ n\geq 3 $, il y a les cas non-conf diffeom \`a la sphere et non.loc.conf.plat. pour lesquelles il y a possibilt\'e de compacit\'e.

\smallskip

////////////////////////////////////////////////////////////////////////////////////////////////////////////////////////////

\smallskip

{\bf Sur la symetrie: pourquoi ?}

\smallskip

a) On considere un systeme: il est r\'egit par des equations. Par exemple, Kaluza-Klein ou supercordes: eq.de Yamabe et eq.du type courbure scalaire prescrite (cas potentiel non constant).

\smallskip

b) Maintenant, on considere plusieurs points de vue du systeme: cela veut dire: point de vue=symetrie en physique ou invariance=en math faire agir une transformation= d'autres points de vue du systeme.

\smallskip

Le point a) precedent veut dire qu'on considere un systeme en lui-meme sans actions.

\smallskip

b-1) Une theorie physique est dite symetrique si l'eq. du systeme prend en compte la particule et son symetrique. Par exemple, ce qu'on vient de voir ci-dessus. l'eq. de la courbure scalaire prescrite avec $ R_g $ et $ V=R_gof $, $ f $ la symetrie physique ou transformation conforme. On a bien une equation qui prend en compte $ x $ et $ f(x) $.

\smallskip

b-2) Par exemple, en theorie des cordes ou supercodes, la supersymetrie, etablit un lien entre des particules de spin entier et des particules de spin demi-entier. C'est un exemple concret d'apparition de la symetrie physique.

\smallskip

b-3) Le but est de determiner les caracteristiques physique du systeme: propri\'et\'es du systeme: en d'autres termes, les caracteristiques du groupe de symetrie puisqu'il concerne un meme systeme. Par exemple, le cas de la sphere: le groupe de symetrie ici est le groupe conforme: il est non compact.

\smallskip

b-3-1) La dimension 2, une surface ferm\'ee  c'est isomorphe a la sphere: groupe conforme "de dimension infinie"= groupe conforme non-compact. Aussi, la relation ou in\'egalit\'e de Harnack $ \sup, \inf $.

\smallskip

b-3-2) Par exemple, ce qu'on a vu ci-dessus, dans le cas negatif. groupe conforme est localement compact: le systeme est localement compact, rigide, epais, consistant, homogene, localement.

\smallskip

b-3-3) Par exemple, le cas positif, eq. Yamabe, cas non loc.conf.diff.\`a la sphere, le groupe conforme est compact: le systeme est compact, rigide, epais, consistant, homogene, globalement.

\smallskip

b-3-4) La compacit\'e locale (cas negatif et existence de mesure de Haar, le systeme dans cas est mesurable) et la compacit\'e (cas positif, non-conf.diffeomorph a la sphere): sont des caracteristiques du systeme physique avec symetrie physique ou invariant ou quand fait agir des transoformations conformes ou plusieurs poins de vues.

\smallskip

b-3-5) De meme pour les in\'egalit\'es de Harnack $ \sup, \inf $ quand on n'a pas la compacit\'e. En particulier en absence ou presence de symetrie. (par exemple l'eq. Yamabe avec $ R_g=1 $ en presence de symetrie on a encore $ V=R_gof =1 $, c'est cach\'e. Mais on peut consider l'eq. de Yamabe en absence de symetrie aussi).

Si on considere Kaluza-Klein ou cordes ou supercordes en presence de symetrie (c'est a dire considerer le groupe conforme, c'est a dire l'eq. de la courbure scalaire prescrite), la seule possibilit\'e de combiner le tout est de considerer l'eq. de Yamabe.

c) 
\smallskip

c-1) Dans Kaluza-Klein ou cordes ou supercordes: il y a 2 types de symetries. Les symetries dues \`a la gravitation, elle sont incorpor\'ees dans l'eq. puisqu'on a consider\'e la gravitation dans les equations: rotations, translations, groupe de Lorentz. Et la symetrie quantique, due aux champs exterieurs, comme l'electromagnetisme, ou Klein-Gordon, Yang-Mills. La symetrie quantique est dite symetrie de jauge, ici, c'est la symetrie conforme. Cela se passe comme pour la theorie conforme des champs de Liouville en dimension 2, le groupe conforme est le groupe de jauge ou groupe de symetrie quantique.

On voit alors, que faire agir une symetrie dans la theorie de Kaluza-Klein, c'est considerer, l'eq. du probleme, comme l'eq. de Kaluza-Klein ou cordes ou supercordes  et l'eq. de la courbure scalaire prescrite qui est l'eq. de l'invariant conforme ou de jauge.

Ce qu'on a dit avant, cela revient a consider l'eq. de Yamabe.

La supersymetrie est l'unifcation des 2 types de symetries ou ce qui revient \`a dire, considerer l'eq. de Yamabe. Comme on l'a dit precedemment, cette symetrie lie les particules de spin entier (les bosons) et les particules de spin demi-entier (les fermions).

\smallskip

c-2) Quand il n'y a pas de symetrie, on dit qu'il y a brisure de symetrie. Cela revient a considerer le systeme en lui-meme, sans actions de groupe conforme ou de symetrie conforme. Ou si on veut, l'eq. de Yamabe ou l'eq. du type courbure scalaire prescrite.

\bigskip

//////////////////////////////////////////////////////////////////////////////////////////////////////////////////////

\smallskip

{\bf Sur la dimension 3 et l'article YY.Li-L.Zhang. dans calculus of variations 2004.}

\smallskip

Il faut supposer l'operateur $ L_g $ coercif, pour $ n=3,4$:

Pour $ n=3 $ c'est du \`a :

[Preuve]: \{ Hypothese de l'absurde: (Etape 1: blow-up analysis)\ldots (Etape n) \ldots (Fin) \}

a) Il y a un probleme par rapport \`a la preuve: Toutes les etapes sont cens\'ees etre vraies, jusqu'a obtenir une contradiction.

b) Ils utilisent comme \underbar {donn\'ee}: comparer le minimum sur la boule ferm\'ee et le minimum sur le bord.(Il faut voir si ce n'est pas en contradiction avec le principe du maximum).

c) Si l'inf peut etre atteint a l'interieur, alors (il se peut que l'operateur soit coercif et les solutions seraient constantes et born\'ees), on ne peut pas commencer la preuve: l'etape 1, le blow-up analysis n'est pas possible. "il n'ya rien \`a prouver". Ce qui n'est pas vrai.

c-1) Si l'inf est atteint \`a l'interieur, les solutions sont constantes et ceci bloque l'etape 1: blow-up analsyis.

c-2) Supposons que la preuve est vraie: alors toutes les etapes sont vraies, jusqu'a la fin et on obtient la contradiction. Pour que ce soit vrai: il faut chaque etape soit vraie, en particulier l'etape 1: le blow-up analysis. Or en supposant l'inf atteint \`a l'interieur, (l'operateur peut etre coercif), les solutions sont alors constantes et born\'ees, ce qui fait que l'etape 1: le blow-up analysis n'est pas vrai.

\smallskip

d) Pour $ n=4 $ c'est du au meme fait que $ n=3 $ + \`a la fonction de Green. Pour la fonction de Green, ils utilisent le theoreme 1.36 du livre d'Aubin (il y a la distance geodesique de l'espace ambiant de depart $ M $ et la distance de l'espace, ou ouvert de carte $ d_{\Omega} $, le theoreme, 1.36 du livre d'Aubin dit qu'on peut toujours se ramener \`a cette distance de l'espace ambiant de depart).(De plus, ici, on a besoin de comparer des quantit\'es en fixant 2 points, la carte exponentielle sert quand on compare des quantit\'es en fixant un point, l'autre c'est l'origine. Alors que pour, Y.Y.Li-L.Zhang, ils veulent comparer des quantit\'es en fixant 2 points, pour cela, le theoreme 1.36 du livre d'Aubin est necessaire). On ne peut pas utiliser les valeurs propres des petites boules geodesiques, car, 1) on ne se retrouve plus dans l'espace ambiant de depart qui est $ M $, mais dans des petites boules, alors qu'ils utlisent $ \exp_{x_k}(y) $, donc c'est par rapport \`a l'espace ambiant de depart. 2) En comparant les rayons des petites boules geodesiques, pour les quelles, les valeurs propres sont grandes, il se peut qu'ils soient (ces rayons), plus ou moins grands que les rayons des boules $ B(x_k,2\epsilon_k) $, pour les quels on veut appliquer le principe du maximum. (Il faut comparer les rayons des differentes boules).

\smallskip

e) Donc, pour $ n=3,4$ l'operateur $L_g$ doit verifier le principe du maximum, donc doit etre coercif. Si on considere un autre operateur elliptique pas forc\'ement l'operateur conforme, le nouvel operateur doit verifier le principe du maximum, donc, doit etre aussi coercif, pour la dimension 3, 4, c'est a dire dans la preuve de YY.Li-L.Zhang. l'operateur doit etre coercif.

\smallskip

{\bf Remarques: sur la dimension 3:}

\smallskip

Cas plat: dans l'article "Majorations du type..." de 2004: le potentiel $ V $ est Lipschitzien et le bord du domaine born\'e est quelconque (le bord n'est pas n\'ecessairement regulier).

\smallskip

Alors que dans la methode de C.C.Chen-C.S.Lin: le potentiel doit etre $ C^1 $ et uniformement $ C^1$. On l'a dit deja.

\smallskip

Alors que dans la methode de Li-Zhang: le bord du domaine doit etre regulier. C'est du \`a: 

\smallskip

1) In\'egalit\'e de Poincar\'e: le domaine doit etre au moins Lipschitzien. Quand on considere un probleme elliptique ou le principe du maximum.

\smallskip

2) Le principe du maximum de Hopf: les solutions doivent etre $ C^2(\Omega)\cap C^1(\bar \Omega) $, pour avoir la regularit\'e des solutions jusqu'au bord, il faut qu'il y ait une notion de bord regulier.

\smallskip

3) L'operateur doit etre coercif: on l'a deja dit, mais ceci inclus le point 1) et 2).

Dans la methode Li-Zhang, dans le cas plat aussi: il y a une contrainte: il faut que le bord soit regulier. Operateur coercif. Il faut que le principe du maximum soit possible.(D'o\`u la citation de l'article de Berestycki-Varadhan-Nirenberg).

Dans la methode de Li-Zhang, dans le cas plat aussi: ils ont cherch\'e a comparer le minimum sur le bord et \`a l'interieur (des boules):

a) il se peut que le minimum soit atteint \`a l'interieur des boules, et aussi: 

b) il se peut que le minimum soit un minimum global (d'o\`u par exemple, la necessit\'e du principe du maximum sur tout l'ouvert $ \Omega $, ou espace ambiant).

c) le principe du maximum doit etre vrai localement partout, d'o\`u il faut une condition globale, il doit etre vrai sur l'espace ambiant.(et apres on l'applique localement. On considere des conditions sur les fonctions du type: $ C^2(\Omega') \cap C^1(\bar \Omega') , \forall \Omega'\subset \subset \Omega $: alors soit on a une infinit\'e de conditions, soit on considere une condition simple qui est $ C^2(\Omega) \cap C^1(\bar \Omega)$, qui englobe tout, et puisque on a pris $ C^1(\bar \Omega)$, il faut que le bord soit regulier, car une fonction reguliere, bord compris, implique, que le bord est regulier).

Tout cela necessite un bord regulier pour l'espce ambiant $ \Omega $.

\smallskip

Ce raisonnement se fait quand on considere un ouvert de  ${\mathbb R}^n $, muni d'une metrique Riemannienne $ g $, ou une vari\'et\'e $ (M,g) $, ou  dans le cas plat, pour un operateur $ L $ avec un terme lin\'eaire. 

\smallskip

Par contre, pour le laplacien euclidien, $ -\Delta $, il est coercif dans des boules, il n'y a pas besoin de supposer le domaine de d\'epart regulier. (On le dit, il faut que le principe du maximum soit vrai localement partout, ce qui est possible sur un ouvert de $ {\mathbb R}^n $).

\smallskip

Donc, dans le cas plat, avec l'op\'erateur laplacien euclidien, $ -\Delta $, en dimension 3, la technique de Li-Zhang, permet d'obtenir l'in\'egalit\'e optimale $ \sup \times \inf $ (et donc de borner l'energie maximale =volume maximal localement). Elle am\'eliore l'exposant du $\sup $ dans $ \sup^{1/3} \times \inf $ de l'article de 2004 (cette in\'egalit\'e permet de borner un volume conforme au sens de la convexit\'e). (Voir aussi, par exemple, le resultat de Siu et celui de Tian, concernant l'eq. de Monge-Amp\`ere complexe). 

\smallskip

Aussi (toujours en dimension 3, sur un ouvert euclidien), pour obtenir l'in\'egalit\'e optimale $ \sup \times \inf $ par la methode de Li-Zhang, on doit utiliser une donn\'ee en plus, au d\'ebut, le principe du maximum: le minimum au bord est le minimum \`a l'interieur, en plus de la technique moving-plane. On utilise 2 fois le principe du maximum dans la m\'ethode de Li-Zhang, la coercivit\'e et la technique moving-plane.(Th\'eorie de la d\'emonstration ou preuve). Ce qui n'est pas le cas dans l'article de 2004, on utilise, la technique moving-plane. (In\'egalit\'e $ \sup^{1/3} \times \inf $ de l'article de 2004: technique moving-plane. Alors que, in\'egalit\'e optimale $ \sup \times \inf $: principe du maximum + technique moving-plane, une donn\'ee en plus utilis\'ee par la m\'ethode de Li-Zhang (Th\'eorie de la preuve).):

\smallskip

1-En dimension 3, (potentiel $ V $ Lipschitzien), quand on n'utilise pas le principe du minimum, alors on a l'in\'egalit\'e: $ \sup^{1/3} \times \inf $ et on borne localement \underbar{un volume conforme} au sens de la convexit\'e. 

On n'a pas les memes constantes $ c>0, R>0 $, dans l'in\'egalit\'e, cela depend de la preuve. Ce n'est pas la meme inegalit\'e, cela depend de la preuve. Soient, $ c_1 >0, R_1 >0$ les constantes, pour ce cas, et $ (P_1) $ la proposition pour ce cas.

2-En dimension 3, (potentiel $ V $ Lipschitzien), quand on utilise le principe du minimum, alors on a l'in\'egalit\'e optimale: $ \sup \times \inf $ et on borne le volume conforme maximal, localement. 

On n'a pas les memes constantes $ c>0, R>0 $, dans l'in\'egalit\'e, cela depend de la preuve. Ce n'est pas la meme in\'egalit\'e, cela depend de la preuve. Soient, $ c_2 >0, R_ 2>0 $, les constantes pour ce cas, et $ (P_2) $ la propostion pour ce cas.

3-Alors, si $ R_1 \leq R_2, \,\, (P_2) \Rightarrow (P_1) $, et, si $ R_1 >R_2, \,\, (P_2) \not \Rightarrow (P_1) $. Donc, $ (P_1) $ n'a rien \`a avoir avec $ (P_2) $.

\smallskip

Maintenant, on s'interesse au cas non plat: Concernant l'article de Math.Aeterna. 2011. Tout d'abord, c'est le cas non plat.

\smallskip

L'operateur n'est pas necessairement coercif et le bord n'est pas necessairment regulier.

\smallskip

Difference avec l'article de 2004: "Majorations du type..."

\smallskip

Dans le cas non plat, il y a moins de possibilit\'es. Dans le cas plat, il y a plus de compacts, les ferm\'es born\'es  sont compacts. Alors que sur une vari\'et\'e Riemannienne, il y a moins de possibilit\'es de caracteriser les compacts: D\'enombrabilit\'e. Possibilit\'es. Aussi, potentiel Lipschitz pour une metrique Riemannienne n'implique pas potentiel Lipschitz pour la distance de $ {\mathbb R}^3 $. Dans l'article de Math.Aeterna 2011, en dimension 3, le potentiel est Lipschitz pour la distance geodesique, (les solutions) dependent de la metrique Riemannienne, aussi pour le cas o\`u la metrique Riemannienne est euclidienne. Alors que pour le cas d'un ouvert de ${\mathbb R}^3$, dans l'article de 2004, le potentiel est Lipschitz pour la distance usuelle de $ {\mathbb R}^3 $, (les solutions) dependent de la structure d'espace metrique pour la distance usuelle. Ce ne sont pas les memes ensembles de solutions.

\smallskip

YY.Li-L.Zhang ont consid\'er\'e un ouvert de $ {\mathbb R}^n, n=3,4 $, munit d'une metrique riemannienne $ g $. Ces hypoth\`eses sont moins g\'en\'erales que le fait de consid\'erer une vari\'et\'e riemannienne $ (M,g) $ de dimension $ n=3,4 $.

Consid\'erer une vari\'et\'e riemannienne $ (M,g) $, est le cas le plus g\'en\'eral possible en g\'eom\'etrie Riemannienne.

\smallskip

Ici, YY.Li-L.Zhang, ont eu un doute, car dans un pr\'ec\'edent article, celui de 2004 et dans le manuscrit de la th\`ese, 2003. On parle de translations ou decalage et la n\'ec\'essit\'e d'avoir une marge de manoeuvre. Cela peut paraitre comme une condition topologique. Pour YY.Li-L.Zhang, cela veut dire aussi, qu'il y a une condition topologique avant la condition metrique et un ouvert de $ {\mathbb R}^n $ saitsfait cette condition topologique, en prenant la carte $(\Omega, id_{\Omega}) $. Et apres, on munit cet ouvert d'une metrique Riemannienne qui satisfait la condition m\'etrique. 

En fait, ces consid\'erations ne se posent pas (topologique et metrique s\'epar\'ement). Il n'y a que la condition m\'etrique.(Cette condition est li\'ee au blow-up, qui utilise la distance, les 2 consid\'erations dont on a parl\'e avant ne se posent pas, il n'y a que la condition metrique. C'est li\'e au blow-up).

\smallskip

Il y a un article de Y.Naito et T.Suzuki, dans Pacific Journal of Math. 1999. dans lequel, ils utilisent la methode "moving-sphere" en dimension $ n\geq 2 $, pour des EDP  g\'enr\'erales, qui peuvent se ramener \`a des equations sur un espace de Riemann, une couronne (avec une metrique Riemannienne), de $ ({\mathbb R}^n, dx^2/|x|^2) $. Ceci, est un argument qui confirme que , YY.Li-L.Zhang, ont eu un doute et considerent un ouvert de $ {\mathbb R}^n $ munit d'une metrique Riemannienne.

\smallskip

(Supposons par exemple, que l'intention de YY.Li-L.Zhang, etait, de se ramener \`a un ouvert de $ {\mathbb R}^n $, par une application de carte $ \phi $: cela impliquerait qu'il faille considerer la metrique induite par la carte $ \phi^*(g) $ et donc, considerer, les boules geodesiques associ\'ees. D'autre part, ils utilisent les coordonn\'ees geodesiques de la metrique $ \phi^*(g) $, alors qu'on peut le faire pour la m\'etrique de d\'epart $ g $, sans passer par la carte).

(Supposons par exemple, que YY.Li-L.Zhang, soient concern\'es par le fait suivant: decalage et marge de manoeuve, c'est \`a dire quand les objets varient (en particulier): ils ecrivent, par exemple, qu'ils considerent une suite de metriques $ g_k $, alors qu'il faut considerer une metrique $ g $ fix\'ee. Il ne faut pas prendre une suite de metriques Riemanniennes $ g_k $, il faut considerer une metrique fix\'ee. (Il se peut que la suite de rayons d'injectivit\'es des $ g_k $ tende vers $ 0 $). On voit bien qu'il y a un doute).

Donc, on voit bien qu'il y a un doute chez YY.Li-L.Zhang. Ils considerent un ouvert de $ {\mathbb R}^n $, munit d'une metrique Riemannienne.

\smallskip

Pour T. Aubin. le seul exemple a disposition est l'article de Brezis-Li-Shafrir, qui traite du cas $ (\Omega\subset {\mathbb R}^2, {\mathcal E}) $. (Avec la notion de parametrisation et revetement). Car, en g\'eneral, quand on considere l'eq. de la courbure scalaire prescrite, on regarde la dimension $ n=2 $, sur les surfaces Riemanniennes, et apr\`es, le cas $ n\geq 3 $. Or, ici, on a seulement le cas $ n=2 $, sur un ouvert $ \Omega $ de $ {\mathbb R}^2 $, d'o\`u la question de savoir ce qui se passe d'abord sur un ouvert muni d'une metrique Riemannienne. Ce qui a \'et\'e fait en dimension $ n\geq 3 $ pour le cas particulier, $ (\Omega \subset {\mathbb R}^n, {\mathcal E}), n\geq 3 $, dans la th\`ese de 2003. (Avec la notion de coorodnn\'ees polaires).

\smallskip

On s'interesse au proc\'ed\'e suivant:

a) carte exponentielle: $ \exp_x(y) $. Sur un ouvert de $ {\mathbb R}^n $: $ \exp_x(y)= x+y $. (g\'eod\'esiques).

\smallskip

b) systeme de coorodnn\'ees.

\smallskip

c) technique "moving-plane".

\smallskip

On a alors:

-Dans Brezis-Li-Shafrir, c'est: parametrisation et revetement.

\smallskip

-Dans C.C.Chen-C.S.Lin, c'est: coordonn\'ees cartesiennes. Lin\'earit\'e et bilin\'earit\'e.

\smallskip

-Dans Y.Naito-T.Suzuki, 1999, c'est: une metrique Riemannienne particuli\`ere et coordonn\'ees cart\'esiennes.

\smallskip

-Dans la th\`ese de 2003, c'est: coordonn\'ees polaires. Non-lin\'earit\'e de la m\'etrique en polaires. Distance usuelle de $ {\mathbb R}^n $.

\smallskip

-Dans l'article de YY.Li-L.Zhang, c'est: coordonn\'ees g\'eodesiques normales. Non-lin\'earit\'e.

\smallskip

-Dans l'article de 2007 dimensions 5 et 6, c'est: coordonn\'ees g\'eodesiques polaires. M\'etrique compl\`etement non-lin\'eaire. Distance g\'eod\'esique.

\smallskip

-Dans les preprints de 2025, sur la dimension 4, avec potentiels $ C1 $ et $ C2 $, c'est: coordonn\'ees g\'eodesiques polaires et revetement.

\smallskip

(Dans la th\`ese de 2003: cas n\'egatif: cartes normales et vari\'et\'es Riemanniennes $ (M,g) $).

\smallskip

Il se peut que, YY.Li-L.Zhang aient eu un doute sur cette condition, (du point de vue topologique): translations ou decalage, et, marge de manoeuvre (et aussi avec le cas n\'egatif). c'est pour cela qu'ils considerent un ouvert de $ {\mathbb R}^n $, puis ils le munissent d'une m\'etrique Riemannienne $ g $.

\smallskip

On r\'epete que: consid\'erer une vari\'et\'e riemannienne $ (M,g) $, est le cas le plus g\'en\'eral possible en g\'eom\'etrie Riemannienne.

\smallskip

//////////////////////////////////////////////////////

\smallskip

{\bf Remarque:} pourquoi, on cite l'article de R.M.Schoen. Journ.Diff.Geometry.(18) 1983. et dans un cas particulier le livre d'Aubin. c'est li\'e \`a la partie: principe du maximum et le lemme de Hopf, sur une vari\'et\'e compacte \`a bord. 

\smallskip

Dans le livre d'Aubin, c'est ecrit pour la valeur propre du laplacien pour un probl\`eme aux limites avec condition de Dirichlet. quand il veut prouver que l'espace des fonctions propres est un espace vectoriel de dimension 1.

\smallskip

Dans l'article de R.M.Schoen. le lemme 1, theoreme 1 et theoreme 2. application du principe du maximum et lemme de Hopf, sur les vari\'et\'e compactes \`a bord, avec bord, une sous-vari\'et\'e compactes. Avec la notion de sym\'etrie.

\smallskip

C'est pour les lier au domaine de l'analyse sur les vari\'et\'es et l'etape 5 de l'article de Brezis-Li-Shafrir. appliquer et consolider cette machinerie. (technique moving-plane, principe du maximum sur une vari\'et\'e compacte \`a bord). 

\smallskip

(Concernant le revetement du cercle, ceci est ecrit dans le Dautray-Lions: le laplacien sur le cercle est la deriv\'ee seconde de $ u(\cos \theta, \sin \theta) $: $ \Delta_{\sigma}(\cos \theta,\sin \theta)=\frac{\partial^2}{\partial\theta^2} (u(\cos \theta, \sin \theta))$).

\smallskip
Un des points les plus importants de l'etape 5 de l'article de Brezis-Li-Shafrir est le point (41), l'in\'egalit\'e: $ \min, \max $, qui implique l'application et la consolidation de cette machinerie (cette technique est la technique moving-plane, propre \`a Brezis-Li-Shafrir).

\smallskip

////////////////////////////////////////////////////////////

On ne retrouve pas le resultat de l'article de 2004, de la dimension 4, en uitlisant la technique de Li-Zhang:
\smallskip

Dans la technique blow-up, la fonction, $ s_i(x)=(R_i-d(x,x_i))^{(n-2)/2} u_i(x) $: on ne peut pas considerer n'importe quel exposant $ \alpha >0 $ de $ R_i-d(x,x_i) $, il faut prendre un exposant $ \alpha $ li\'e \`a l'exposant de l'equation. Si on conidere $ \Delta u= Vu^q, q>2, \Delta=-\sum \partial_{jj} $, alors $ \alpha = 2/(q-2) $, pour l'exposant critique, $ N=2n/(n-2)$, il faut prendre, $ \alpha=(n-2)/2$. Tout ceci est li\'e au changement d'echelle, ou rescaling ou invariance: la distance $ d_R \cdot M_i^{-2/(n-2)} \approx R \leftrightarrow |x| \cdot M_i^{-2/(n-2)} \approx |y| \leftrightarrow \tilde d \cdot M_i^{-2/(n-2)} \approx d $, dans le blow-up.

\smallskip

Ce qui veut dire que la technique de Li-Zhang, n'est pas possible en dimension 3, avec des potentiels $ 1/2 -$ holderien et aussi n'est pas possible en dimension 4 avec condition de Lipschitz. Car par exemple, en dimension 4: la condtion au bord s'ecrit, on n'a pas necessairement: $ M_i^{-(1/2)} \not \leq 2 R_i $. La boule de rayon, $ M_i^{-1/2}$ n'est pas forc\'ement incluse dans la boule de rayon $ 2R_i$.

\smallskip

De meme, en dimension 3 aussi, la condition au bord s'ecrit: on n'a pas necessairement: $ M_i^{-1} \leq 2R_i $. D'ou par exemple la modification dans le th 4 de l'article de 2008, dans acta.math.sci. 2008. on a pris la sphere en dimension 3, pour le th 4.

\smallskip

Cela veut dire qu'on ne retrouve pas le resultat de la dimension 4, dans, l'article de 2004, avec la methode de Li-Zhang. La technique de Li-Zhang, ne permet pas de reprouver ou de retrouver le resultat de la dimension 4 de l'article de 2004, avec des potentiels Lipschitziens tendant vers 0.

\smallskip

////////////////////////////////////////////////////////////////////////////

\smallskip

De meme, pour le th.4 de l'article de 2019 dans Pac.Journ.Math. On retrouve pas ce theoreme a partir de l'article de C.C.Chen-C.S.Lin.1997. Car:

\smallskip

-le theoreme 1 de Chen-Lin.1997: le domaine de definition de la fonction blow-up est li\'e au domaine de definition de la fonction dans l'hypothese de l'absurde. Or, dans le th.1 de C.C.Chen-C.S.Lin.1997: il n'y a pas de restriction dans le domaine de definition de la fonction: la boule $ B_2(0)=2B_1(0) $.

\smallskip

Ce qui veut dire qu'on ne peut appliquer le th.1 de Chen-Lin.1997. au th.4 de l'article de 2019. dans Pac.Journ.Math.

\smallskip

-le theoreme 2 de Chen-Lin.1997: il y a une restriction dans le domaine de Def. dans l'hypothese de l'absurde pour la fonction de depart $ u_i $ : $ B_R(0) $. Donc, le domaine de def, de la fonction blow-up, doit etre restreint apres rescaling \`a: $ B_{2R}(0)=2B_R(0)$. Ici, on ne peut pas appliquer ce th.2 de Chen-Lin.1997 au th.4. de Pac.journ.Math.2019. pour la meme raison que ci-dessus($ M_i^{-(1/2)} \not \leq 2R_i $, c'est d\^u a la fonction $ s_i $, de la technique blow-up).

\smallskip

Ce qui veut dire qu'on ne retrouve pas le th.4. de l'article de 2019. Pac.Journ.Math. \underbar {entierement} avec la technique de Chen-Lin.1997. On utilise le resultat de platitude d'ordre 1 de Chen-Lin dans ce th.4 de l'article de 2019.

\smallskip

////////////////////////////////////////////////////////////////////////////////////////////////////////////////////////////////////

\smallskip

{\bf Sur les differentes topologies des vari\'et\'es \`a bord:}

\smallskip

1) Sur l'equivalence des topologies. Sur un ouvert de ${\mathbb R}^n $ avec bord regulier: la topologie des cartes est equivalente \`a la topologie des distances locales qui est equivalentes \`a la topologie de la distance Riemannienne qui est equivalente \`a la topologie de la distance de $ {\mathbb R}^n $ sur le domaine regulier. Et ceci jusqu'au bord, en incluant les cartes du bord.

\smallskip

2) Sur une vari\'et\'e Riemannienne \`a bord, intrinseque (non necessairement plong\'ee dans un espace plus grand). La topologie des cartes est equivalente \`a la topologie locale des distances qui est equivalente \`a la topoologie de la distance Riemannienne y compris le bord et la distance Riemannienne jusqu'au bord, bord compris. Le theoreme du livre de Hebey, fonctionne, avec les cartes du bord, quand on cherche \`a comparer les topologies des cartes du bord et la topologie des distances, bord compris.

\smallskip

/////////////////////////////////////////////////////////////////////////////////////////////////////////////

\smallskip

{\bf Sur le resultat de compacit\'e en dimension 2: application de: compacit\'e avec energie ou volume born\'es:}

Compacit\'e avec contrainte: 

\smallskip

1) Topologie des espaces de domaines \`a bord: en plus du fait qu'on a convergence de domaines avec bord et metriques conformes, suites de domaines \`a bord. geometrie metrique. Espaces d'espace ou espaces de domaines a bord muni de metriques Riemanniennes. Volume fix\'e: metrique blow-up, eclatement, ou effondrement d'espaces. ou compacit\'e de metriques, non effondrement d'espaces.

\smallskip

2) Resultat de compacit\'e de fonctions dans les EDP. et les espaces de fonctions. En topologie et en analyse fonctionnelle.

\smallskip

3) En chimie: emballement thermique (thermal runaway, jusqu'a l'eclatement, tend vers l'abime) ou non emballement (compacit\'e: objet ne s'abime pas)

\smallskip

4) Problemes variationnels:

Par exemple, regardons le cas regulier, sans singularit\'e:

On a vu dans le print,"Cas d'existence de solutions d'EDP": que la probleme variationnel:

$$ \mu = \inf \{ ||\nabla u||_2, u\in \dot H_1^2(\Omega), \int_{\Omega} Ve^u=1, \} $$
 a une solution positive $ u >0 $ avec la condition $  b|\Omega| < 1 $.

$$ \Delta u = \lambda V e^u, u=0, \,\, {\rm au\,\, bord},\,\, \lambda >0, $$

a) En supposant: $ 0 <a \leq V\leq b <+\infty $ avec la fonction propre, le coefficient, le multiplicateur de lagrange verifie $ \lambda \times a \leq C <+\infty$

b) En supposant $ V $ Lipschitz, on a la compacit\'e. et en fait $ \lambda \not \to 0 $ uniform\'ement. Donc, on a: $ \mu >0 $ uniform\'ement en $ u \in \dot H_1^2(\Omega) $ avec $ \int_{\Omega} Ve^u = 1 $.

Donc avec la condition $ b|\Omega| <1$:

$$ \forall u \in \dot H_1^2(\Omega), \int_{\Omega} Ve^u=1 \Rightarrow ||\nabla u||_2 \geq \mu >0.$$

Donc, la contrapos\'ee donne:

$$ ||\nabla u ||_2 < \mu \Rightarrow \int_{\Omega} Ve^u \not = 1. $$

On va  voir qu'on peut avoir une in\'egalit\'e du type Moser-Trudinger, et les constantes, et $ V $ determinent cette nouvelle in\'egalit\'e.

Soit: $ u\in \dot H_1^2(\Omega), u\not \equiv 0, $ et $ v=\frac{\mu}{2} \frac{u}{||\nabla u||_2} $. Alors, $ ||\nabla v||_2 =\mu/2<\mu $, donc, $ \int_{\Omega} V e^v \not = 1 $.

Supposons qu'il existe $ u_0 $ et $ u_1 $ non nuls tels que $ \int_{\Omega} Ve^{v_0} <1, \int_{\Omega} V e^{v_1} >1 $, avec $ v_0, v_1 $, les fonctions construites  a partir de $ u_0,  u_1$, comme $ u $ et $ v $. On a: $  v_0=\frac{\mu}{2} \frac{u_0}{||\nabla u_0||_2},  v_1=\frac{\mu}{2} \frac{u_1}{||\nabla u_1||_2} $.

Alors, on considere le chemin $ w_t=tv_0+(1-t) v_1$, alors, $ ||\nabla w_t||\leq \mu/2 <\mu $, d'ou, $ \int_{\Omega} Ve^{w_t} \not = 1, \forall \, t \in [0,1] $, or par le theoreme des valeurs intermdiaires en considerant la fonction continue de $ t, g(t)=\int_{\Omega} V e^{w_t} dx $, on aura un $ t_0 $ tel que $ g(t_0)=1=\int_{\Omega} Ve^{w_{t_0}} $, ce qui est constardictoire.

c) Donc, on a soit tout le temps $ \int_{\Omega} Ve^v<1$ ou tout le temps, $ \int_{\Omega} Ve^v > 1 $.

Or, en considerant les fonctions $ h\leq 0, h\not \equiv 0 $, on $ \int_{\Omega} Ve^h \leq b|\Omega|<1$. (Par exemple si $ u \in \dot H_1^2(\Omega),  u\not \equiv 0 $, $ |\nabla (|u|)|=|\nabla u| $, on peut prendre $ h=\frac{\mu}{2}\frac{-|u|}{||\nabla u||_2} $)

d) Donc, on a tout le temps:

$$ \int_{\Omega} Ve^v < 1. $$

e) Soit, 

$$ \nu=\sup \{ \int_{\Omega} V e^v, v=\frac{\mu}{2}\frac{u}{||\nabla u||_2}, u \in \dot H_1^2(\Omega), u\not \equiv 0 \}, $$

Par la compacit\'e de l'injection de Moser Trudinger, ce $ \sup $
 est atteint: $ ||\nabla v||_2\leq \liminf_i ||\nabla v_i||_2 =\mu/2 <\mu $ et la compacit\'e de l'injection de Moser-Tridinger: $\int_{\Omega} Ve^{v_i} \to \int_{\Omega} Ve^v =\nu \leq 1 $. Alors, $ v \in \dot H_1^2(\Omega), ||\nabla v||_2 \leq \mu/2 < \mu $, d'ou, $ \int_{\Omega} Ve^v \not =1$. Donc $ \nu <1 $.

f) Finalement:

$$ \forall \, a, b, A >0, \, b|\Omega|<1, \exists \,\mu >0, \exists \, 0 < \nu <1, \,\, \forall u \in \dot H_1^2(\Omega)-\{0\}, \,\, \int_{\Omega} V e^{\frac{\mu}{2} \frac{u}{||\nabla u||_2}} dx \leq \nu < 1.$$

C'est une in\'egalit\'e du type Moser-Trudinger, la consition, $ b|\Omega| <1 $ et les constantes $ \mu>0 , 0 <\nu <1 $ determinent cette in\'egalit\'e. Donc, pour $ a,b,A>0 $ avec la condition $ b|\Omega| < 1 $, on obtient:

$$ \exists \,\mu >0, \exists \, 0 < \nu <1, \,\, \forall u \in \dot H_1^2(\Omega)-\{0\}, \,\, \int_{\Omega} V e^{\frac{\mu}{2} \frac{u}{||\nabla u||_2}} dx \leq \nu < 1.$$

g) On peut faire la meme chose avec, le resultat de compacit\'e avec singualrit\'e au bord. On a une in\'egalit\'e du type Moser-Trudinger avec singularit\'e au bord, avec des constantes $ \mu >0, 0 < \nu <1$ et la condition sur $ b: b \times \int_{\Omega} \frac{1}{|x-x_0|^{2\alpha}} dx <1, \alpha \in ]0,1/2[ $, $ x_0 \in \partial \Omega $ et $ \Omega $ domaine analytique:

$$ \exists \,\mu >0, \exists \, 0 < \nu <1, \,\, \forall u \in \dot H_1^2(\Omega)-\{0\}, \,\, \int_{\Omega} \frac{V}{|x-x_0|^{2\alpha}} e^{\frac{\mu}{2} \frac{u}{||\nabla u||_2}} dx \leq \nu < 1.$$

h) On peut mettre la valeur absolue pour $ u $ et remplacer $ u $ par $ |u| $ dans ces in\'egalit\'es.

Ces fonctions $ u $ font apparaitre l'Eq. de la courbure scalaire prescrite en dimension 2. Donc: pour un $ u \in \dot H_1^2(\Omega) $, $ \int_{\Omega} e^{ku} dx \approx |\Omega|_{g= e^u \delta}^k $ correspond \`a un volume ou une surface, \`a un r\'eel positif, $ k >0 $ pr\'es (en utilisant l'in\'egalit\'e de Holder par exemple). Le terme $ ||\nabla u||_2^2 \approx |\partial \Omega|_{g=e^u\delta}^2 $, (en utilisant l'inegalit\'e de Cauchy-Schwarz par exemple), correspond au perimetre, en passant par les fonctions BV, voir l'article de. O. Druet dans Numdam, 2001-2002: in\'egalit\'es de Sobolev et in\'egalit\'es isop\'erimetriques.

\smallskip

Dans le cas avec singularit\'e au bord, on a des in\'egalit\'es entre surface et perimetre avec la metrique $ g=\frac{e^u}{|x-x_0|^{2\alpha}} \delta, \alpha \in ]0, 1/2[, x_0\in \partial \Omega $ et $ \Omega $ analytique.

\smallskip

Cela s'applique aussi a l'operateur: $ -div (e^{\epsilon |x|^2/2} \nabla) $. La fonctionnelle est: $ \int_{\Omega} e^{\epsilon |x|^2/2} |\nabla u|^2 dx $.

\smallskip

Donc, ces in\'egalit\'es en dimension 2 de Moser-Trudinger, mettent en relation la surface conforme (volume conforme) et le p\'erimetre (conforme): ce sont des in\'egalit\'es isoperimetriques particulieres.

\smallskip

Il y a aussi l'interpretation en physique, en termes d'energies. (Qui correspond aussi au cas ou l'operateur n'est pas n\'ecessairement le laplacien).

\smallskip

On peut prendre comme contrainte: $ \int_{\Omega} Ve^u = k $, $ k \geq 1 $ au lieu de $ 1 $. Dans ce cas, la condition sur $ b $ est: $ b|\Omega| < k $.

\smallskip

En mettant la valeur absolue dans $ u $, on obtient: in\'egalit\'es du type Moser-Trudinger:

\smallskip

Avec, $ a, b, A>0 $, $ V $ et $ b|\Omega|< 1$:

$$ \exists \,\mu >0, \exists \, 0 < \nu <1, \,\, \forall u \in \dot H_1^2(\Omega)-\{0\}, \,\, \int_{\Omega} V e^{\frac{\mu}{2} \frac{|u|}{||\nabla u||_2}} dx \leq \nu < 1.$$

et, avec $ a, b, A >0 $, $ \alpha \in ]0, 1/2[$, $ V $, $ b\int_{\Omega} \frac{1}{|x-x_0|^{2\alpha}} dx < 1$ et $ \Omega $ analytique:

$$ \exists \,\mu >0, \exists \, 0 < \nu <1, \,\, \forall u \in \dot H_1^2(\Omega)-\{0\}, \,\, \int_{\Omega} \frac{V}{|x-x_0|^{2\alpha}} e^{\frac{\mu}{2} \frac{|u|}{||\nabla u||_2}} dx \leq \nu < 1.$$

De meme pour l'article avec $ V $ holderian et avec condition sur la masse: "About Brezis-Merle Problem with Holderian condition" on a une in\'egalit\'e du type Moser-Trudinger.

\smallskip

//////////////////////////////////////////////////////////////////////////////////////////////////////////////////////

\smallskip

{\bf Sur l'article de Brezis-Merle et leurs problemes dans l'introduction et le probleme 1:}

\smallskip

1) Brezis et Merle ont consid\'erer leur probleme dans les $ L^p, 1 < p \leq +\infty $, or ceci doit inclure le cas $ 1 < p <+\infty $ et la regularit\'e des solutions est $ C^0$ dans ce cas, quand ils utilisent les injections de Sobolev.

\smallskip

Brezis et Merle considerent les fonctions $ u_n, u $ dans les $ L^p$ avec, quand il y a injections de Sobolev, la regularit\'e $ C^0$. (ceci pour inclure tous les cas $ p $ en particulier $ 1< p <+\infty$).(et non les deriv\'ees de  $ u $ et de $ u_n $, car les estimations en normes $ C^1 $ peuvent ne pas exister, au sens usuel et au sens de la regularit\'e).(Aux points du bord, on n'a pas un voisinage dans le quel on puisse considerer le taux d'accroissement, donc, $ \nabla u_n $ n'existe pas forc\'ement, mais $ \partial_{\nu} u_n $ existe, de plus, comme on a considerer les coordonn\'ees cartesiennes, ce terme est un terme limite).

Quand on prend un representant $ C^1, v_n $, il faut voir qu'on considere $ (v_n)_{| \bar \Omega} $: le terme de bord dans l'integration par parties est bien un terme limite qui coincide avec $ v_n$ aux points du bord.

\smallskip

Brezis et Merle utilisent le theoreme de dualit\'e de Stampacchia, en particlier dans le corollaire 4, pour borner uniform. les fonctions harmoniques. Ils ne precisent par la convergence des $ \nabla u_n $ dans $ L^q, 1\leq q<2 $: ils utilisent la borne uniforme du th de dualit\'e.

\smallskip

Ceci nous amene au fait suivant: ils ne considerent pas la deriv\'ee normale, car:

\smallskip

a) La formule d'integration par parties est une formule de trace: le terme de bord est une limite(il coincide avec la deriv\'ee normale quand les fonctions sont regulieres, mais c'est tout de meme un terme limite, car, le gradient n'existe pas au voisinage du bord, mais la deriv\'ee normale existe sans que le gradient existe). En coordonn\'ees cartesiennes.

\smallskip

b) Brezis et Merle considerent $ u_n $, au sens strict, et cette fonction est $C^0 $ presque partout.

\smallskip

c)Brezis et Merle considerent aussi $ u_n$ en coordonn\'ees cartesiennes. Sans notion de cartes. Les coordonn\'ees cartesiennes ne fournissent pas beaucoup d'informations sur $ u_n $ au voisinage des points du bord.

\smallskip

d) quand on multiplie par une fonction test: il y a ambiguit\'e dans la definition des points reguliers et non reguliers d'une mesure.

\smallskip

e) En se placant dans un voisinage d'un point du bord, on a "directement" la notion de carte, et ceci implique que dans cette carte l'expression du laplacien change et l'eq; qu'on considere change. Or leur article est ecrit pour le laplacien, d'autre part, l'article de Y.Y.Li-Chanillo, ne fournit rien de nouveau quand on se place au bord. L'article de Y.Y.Li-Chanillo confirme ses arguments de cartes locales et le changement d'expression de l'operateur dans ces nouvelles cartes.

\smallskip

f) Quand on considere les points du bord: on a les cartes du bord: l'eq. change dans cette nouvelle carte et les th de Brezis et Merle ne sont pas utilisables, car c'est le laplacien usuel qui est utlis\'e pour ces theoremes.

\smallskip

g) Brezis et Merle considerent $ u_n, u $, en coordonn\'ees cartesiennes. Y.Y.Li-Chanillo, considerent un operateur global en forme divergence (comme le laplacien). En se placant au voisinage d'un point $ x_0 \in \partial \Omega $, on a la notion de carte locale $ (x_0\in B_0, \phi)$, la fonction consid\'er\'ee est $ u_no\phi$, l'operateur laplacien change dans cette nouvelle carte et il faut definir les points reguliers et non reguliers, defintion qui est ambigue.

\smallskip

h) Ils n'utilisent pas l'integration par parties car: elle fait apparaitre: soit le gradient de la fonction test, soit le laplacien de la fonction test, or ces termes ne sont pas facilement controlables. Il faut etre precis dans l'estimation de ces termes, en etant precis dans le th. de dualit\'e de Stampacchia. Remarquons que, utiliser la regularit\'e $ C^1$ correspond au cas $ p=+\infty $ et exclu les cas $ 1< p <+\infty $. Et leurs theoremes qui ne prennent en compte que $ u_n $ doivent inclure  le cas $ 1< p <+\infty $, ce qui oblige a considerer la regularit\'e $ C^0$ dans les injections de Sobolev. Donc, ils n'utilisent pas la regularit\'e $ C^1 $ mais $ C^0 $: tout cela pour inclure tous les cas $ 1< p \leq +\infty $, pour le principe du maximum de Stampacchia, il suffit que les solutions soient ici, $ W^{1,2}_0(\Omega) \cap C^0(\bar \Omega) $. 

\smallskip

Ecrivons les formules d'integration par parties: (c'est en coordonn\'ees cartesiennes et c'est une formule de trace):

On a le cas:

$$ \int_{\Omega} V_n e^{u_n} \eta dx =\int_{\Omega} -\Delta u_n \eta = \int_{\partial \Omega} \partial_{\nu} u_n \eta + \int_{\Omega} u_n \Delta \eta, $$

et le cas:

$$ \int_{\Omega} V_n e^{u_n} \eta dx =\int_{\Omega} -\Delta u_n \eta = \int_{\partial \Omega} \partial_{\nu} u_n \eta + \int_{\Omega} \nabla u_n \cdot \nabla\eta, $$

et le cas:

$$ \int_{B_R(x_0)} V_n e^{u_n} dx =\int_{B_R(x_0)} -\Delta u_n = \int_{\partial B_R(x_0)} \partial_{\nu} u_n, $$

Or, pris comme ils sont, ces termes ne sont pas controlables avec seulement la borne uniform du th de dualit\'e de Stampacchia. De plus, il faut definir les points reguliers et non reguliers, par rapport aux cartes, definition qui est ambigue. 

\smallskip

La masse totale est la somme d'une mesure sur le bord et d'une mesure de Lebesgue: quelle mesure considrer ? quelle nouvelle mesure doit on definir ? et comment definir les points reguliers et non reguliers ?

\smallskip

La troisieme formule d'integration par parties, il y a un terme du bord qu'on controle pas non plus: celui de $ \partial B_R(x_0) \cap \Omega $.

\smallskip

La formule d'integration par parties est une formule de trace. En coordonn\'ees caresiennes. C'est formule de trace, le terme de bord est un terme limite, qui coincide avec le terme usuel, quand c'est regulier. C'est une formule de trace.

\smallskip

Les coordonn\'ees cartesiennes "ne fournissent rien de nouveau" pour les points du bord. En general, pour les points du bord, il faut utliser des cartes locales: ici il y a une ambiguit\'e. Il faut aussi definir, une mesure et les points reguliers et non reguliers de cette mesure pour les points du bord: ici aussi il y a une ambiguit\'e.

\smallskip

Brezis et Merle utilisent seulement la borne uniforme du th de dualit\'e de Stampacchia et la regularit\'e $ C^0$ dans les injections de Sobolev: Ils veulent considerer que les cas $ L^{\infty} $ et $ L^{\infty}_{loc} $. Donc: regularit\'e $ C^0 $ dans les  injections de Sobolev. Ils ne considerent que les solutions $ u, u_n $ et non les deriv\'ees de $ u $ et de $ u_n $ car les estimations en normes $ C^1 $  peuvent ne pas exister(en inculuant les cas $ 1< p <+\infty$ et $ \nabla u_n $ n'est pas forc\'ement definit meme en prenant un representant $ C^1$).

\smallskip

2) Brezis et Merle avec les articles suivant et ceux de la dimension $ n\geq 3$ pensent qu'il est possible  d'utiliser le blow-up et la methode moving-plane sur le bord, qui donne la compacit\'e. Ceci est verif\'e par l'article de Z-C.Han. et aussi l'article de De Figueiredo-Lions-Nussbaum, et l'article de Suzuki et l'article de W.Chen- Congming.Li: la methode moving-plane sur le bord donne la compacit\'e.

\smallskip

Leurs apporches en dehors du fait de fixer une mesure:

a) Soit il faut utiliser la technique blow-up: faire un eclatement. Ici, ils font l'eclatement \`a l'interieur. 

\smallskip

b) technique moving-plane: compacit\'e au bord par la methode moving-plane.

\smallskip

Pour Brezis-Merle, Chanillo, YY.Li... d\'es qu'on considere les points du bord, il y a la notion de cartes locales qui apparait, on n'a plus le laplacien usuel, on n'a plus la meme equation. Ceci dans le cas o\`u on utilise la r\'egualrit\'e $ C^0 $ donc $ C^1 $, des solutions.

\smallskip

Dans le cas o\`u on considere les solutions dans $ L^{\infty} $, on n'a plus la r\'egularit\'e au voisinage du bord.

\smallskip

/////////////////////////////////////////////////////////////////////////////////////////////////////////////////////////

\smallskip

{\bf Sur l'article donnant l'unicit\'e au.Bull.Sci.Math.de 2009:}

\smallskip

Remarque sur le th. d'Olivier Druet donnant la compacit\'e et le th. de l'article de 2009 au Bull.Sci.Math. donnant la compacit\'e pour $ \epsilon \to 0 $.

\smallskip

1) La preuve de Druet est vraie quand $ 1\geq \epsilon \geq \tilde \epsilon_0 >0, \forall \tilde \epsilon_0 >0 $. 

\smallskip

2) La preuve de l'article au Bull.Sci.Math. est vraie lorsque $ 0 < \epsilon \leq \epsilon_0 $: $ \exists  m_0 >0, \exists \epsilon_0 >0, $ tels que pour $ 0 < \epsilon \leq \epsilon_0 $, on ait: $ \max_M u_{\epsilon} \leq m_0 $.

\smallskip

3) En suite on a l'unicit\'e: $ \exists \epsilon_1 >0, 0 <\epsilon_1 \leq \epsilon_0 $ tel que pour $ 0 <\epsilon \leq \epsilon_1 $, on ait: $ u_{\epsilon} \equiv \epsilon^{(n-2)/4} $.

\smallskip

La remarque est la suivante: on ne peut pas dire que le resultat de Druet est vrai, par rapport au resultat du Bull.Sci.Math. 2009, pour $ 0< \epsilon_1 \leq \epsilon\leq \epsilon_0 $. Car des qu'on parle de $ \epsilon_0 $, cela nous amene au point 2) precedent qui dit qu'on a la compacit\'e pour $ 0 <\epsilon \leq \epsilon_0 $, or la preuve(la methode) de Druet n'est pas vraie pour $ 0 < \epsilon \leq \epsilon_0 $, donc on ne peut pas parler de la preuve de Druet dans ce cas, car elle n'est pas vraie. On ne peut pas distinguer les cas $ 0 <\epsilon \leq \epsilon_1 $ et $ 0 < \epsilon_1 \leq \epsilon \leq \epsilon_0 $. Car des qu'on parle de $ \epsilon_0$ cela fait reference au point 2) precedent qui dit qu'on a la compacit\'e dans ce cas et que la preuve de Druet n'est pas vrai dans ce cas.

\smallskip

On ne peut pas parler de $\epsilon_0$ et $\epsilon_1$ ou de $ \epsilon_0$, et du resultat de Druet, au meme temps, ce n'est pas compatible.

\smallskip

(Par exemple, considerer $ \epsilon_0 >0$ du point 2), c'est prendre en compte les fonctions $ 
u_{\epsilon} $ avec $ \epsilon >0 $ voisin de $ 0 $, ce qui n'est pas le cas du resultat de Druet).

\smallskip

Tout ceci pour dire qu'on a bien un resultat de compacit\'e dans l'article du Bull.Sci.Math. de 2009, sans faire reference au resultat de Druet et sans faire reference \`a l'unicit\'e. L'unicit\'e vient apres.

\smallskip

Si on considere les operateurs: $ L_t(u)= u-(\Delta+t)^{-1}(u^{N-1}), 0 < \epsilon_1 \leq t\leq \epsilon_0 $, $ \Delta=-\nabla^i \nabla_i $,  $ u\in C=\{u\in C^{2,\alpha}(M), u >0, ||u||_{C^{2,\alpha}(M)} \leq 2m_0 \} $, par les estimations elliptiques et homotopie et en $ t=\epsilon_1 $, $(L_{\epsilon_1})'$ est un isomorphisme et l'eq. $ L_{\epsilon_1}(u)=0 $ a une solution unique d'o\`u en $ \epsilon_1 $ le degr\'e $ L_{\epsilon_1} $ est son indice et il est egal \`a $ \pm 1 $:

$$ deg (L_t, C, 0)=deg (L_{\epsilon_0}, C, 0)= deg (L_{\epsilon_1}, C, 0)=\pm 1. $$

Donc, on a des solutions topologiques \`a $ \Delta u+t u= u^{N-1}, u >0 $ pour $ 0 < \epsilon_1 \leq t\leq \epsilon_0 $ et elles sont born\'ees en norme $ C^{2,\alpha} $, $ ||u||_{C^{2,\alpha}(M)} \leq 2m_0 $.

\smallskip

///////////////////////////////////////////////////////////////////////////////////////////////////////////////////////////

\smallskip

{\bf Sur la 2\`eme note aux Comptes Rendus Math. Acad. Sci. Paris. 2006:}

\smallskip

On a l'in\'egalit\'e de Brezis-Gallouet, qui prend en compte les deriv\'ees d'ordre 2, $ H^2_2 $ et la norme $ \sup $. On peut voir dans le dernier th de la note de 2006 (on utilise la borne uniforme des solutions au voisinage du bord, de l'equation ($ \Delta^2 u= u^{p-\epsilon}, p=\frac{(n+4)}{(n-4)}, n\geq 5 $ dans $ \Omega $ et avec les conditions au bord $ u=\Delta u=0$ sur $ \partial \Omega $, sur un ouvert regulier strictement convexe, par exemple une boule ou un ellipsoide)):

$$ \sup_{\Omega} u \times \inf_K u \geq c(K,\Omega,n) \int_{\Omega} (\Delta u)^2 dx \geq \tilde c(K, \Omega, n) >0, $$

C'est une in\'egalit\'e liant les normes essentielles et la norme $ H^2_2(\Omega) $.

\smallskip

///////////////////////////////////////////////////////////////////////////////////////////////////////

\smallskip

{\bf Sur le r\'esultat de compacit\'e locale en dimension 4: eq.de la courbure scalaire prescrite: application: probleme de Cherrier}

\smallskip

C'est sous-entendu, par la citation du livre d'Aubin. 

\smallskip

Dans le cas plat et non plat.

\begin{displaymath}  \left \{ \begin {split} 
      \Delta u_i+\frac{1}{6} R_g u_i  & = V_i u_i^{3}, \,\,u_i >0     \,\, &&\text{dans} \!\!&& M, \\
                 \partial_{\nu} u_i+h u_i  & =\tilde h_i u_i^2, \,\, u_i >0   \,\,             && \text{sur} \!\!&&\partial M.               
\end {split}\right.
\end{displaymath}

ici, c'est la normale exterieure. $ \Delta =-\nabla^i\nabla_i $, $ 0 < a \leq V_i \leq b <+\infty, ||\nabla V_i||\leq A_i \to 0$. 

\smallskip

Question: peut on avoir, courbures scalaires $ \geq 0 $ et courbures moyennes $ \leq 0 $ ? Est ce compatible ? c 'est peut etre vrai, car quand on prend un objet convexe, sa courbure scalaire est $ >0$ cela veut dire que le bord est concave et donc courbure $ <0 $ (en se placant au bord, le bord devient concave), par exemple une boule (au voisinage du bord).

\smallskip

Exemple de solutions du probleme de Cherrier+ exemple concernant la question precedente: on considere la premiere valeur propre du laplacien conforme(ici, il est cens\'e etre coercif), avec condition de Dirichlet. Si, $ f >0 $ est la fonction propre, on prend $ g=f+ \epsilon >0, \epsilon >0 $. Exemple de J.Escobar, premiere valeur propre.

\smallskip

Par le principe du maximum de Hopf (ici l'operateur laplacien est coercif), si $ x_i $ est un point ou la solution $ u_i >0 $ atteint son minimum: $ \partial_{\nu} u_i(x_i) <0 $. Si, $ h< 0 $ et $ -\infty < c\leq \tilde h_i \leq d <0 \Rightarrow  \min_M u_i = u_i(x_i) \geq k >0 $.

Donc, dans le cas plat et non plat:

Dans le cas d'un ouvert born\'e de $ {\mathbb R}^4 $, une boule par exemple (courbure moyenne $ < 0 $):

\begin{displaymath}  \left \{ \begin {split} 
      \Delta u_i  \,\,\,\, \,\, & = V_i u_i^{3}, \,\,u_i >0     \,\, &&\text{dans} \!\!&& \Omega , \\
                  \partial_{\nu} u_i+h u_i  & =\tilde h_i u_i^2, \,\, u_i >0   \,\,             && \text{sur} \!\!&&\partial \Omega.          
\end {split}\right.
\end{displaymath}

alors,

$$ || u_i||_{L^{\infty}(K)} \leq c(a,b,c,d, (A_i), K, \Omega), \,\, \forall K \subset \subset \Omega \subset {\mathbb R}^4, $$

et dans le cas non plat, qaund le principe du maximum est vrai (le minimum est atteint sur le bord):

$$ || u_i||_{L^{\infty}(K)} \leq c(a,b,c,d, (A_i), K, M), \,\, \forall K \subset \subset M. $$

On a: \underbar {Compacit\'e locale pour un probleme de Neumann non-lineaire}.

\bigskip

1) On a donn\'e un exemple de solution de probleme de Cherrier avec $ h <0, \tilde h<0$ sur une boule: $ u_{\epsilon} =g= f+\epsilon, \epsilon >0$. On veut voir si il y a des solutions qui blow-up:

Pour cela il suffit de prendre $ \epsilon \to + \infty $, on a blow-up total:

Pour chaque $ \epsilon >>1 $, tres grand, il existe une solution au probleme de Cherrier, tel que:

$$ \max_{\bar \Omega} u_{\epsilon} \geq \max_{\partial \Omega} u_{\epsilon}=\epsilon \to + \infty, $$

cela veut dire qu'on depasse n'importe quel rang $ A >>1$, tres grand, en ayant une solution du probleme de Cherrier.

Ici, $ 0< a_{\epsilon} \leq V_{\epsilon} \leq b_{\epsilon} <+\infty $ et $ -\infty < c_{\epsilon} \leq \tilde h_{\epsilon} \leq d_{\epsilon} < 0 $, ces nombres peuvent etre tres petits en valeures absolues. Par exemple ici, en dimension 4 sur une boule, $ h<0 $, $ V_{\epsilon}=\frac{\lambda_1 f}{(f+\epsilon)^3},\lambda_1 >0, \tilde h_{\epsilon}=\frac{\partial_{\nu} f +\epsilon h}{\epsilon^2} $.

On peut prendre aussi $ f >0 $ solution de $  \Delta f=\lambda_1 f,\,\, \Omega, \,\, f=1, \,\, \partial \Omega $ pour eviter d'avoir $ V_{\epsilon} = 0 $ sur $ \partial \Omega $.

\smallskip

2) Cela veut dire que:

\smallskip

Pour chaque rang $ A >>1 $ tres grand, il existe une solution du probleme de Cherrier, dont le maximum depasse ce rang $ A >>1$.

\smallskip

Cela veut dire que le blow-up pour le probleme de Cherrier est possible. En plus, on le determine a partir du bord (et de l'interieur): on a \`a la fois blow-up interieur et au bord. 

$$ \forall \epsilon >>1, \exists u_{\epsilon}: \,\, solution \, \, blow-up \, = \, solution \,\, avec \,\, maximum \,\, grand. $$

Par exemple, on peut considerer l'ensemble des solutions $ v $ du probleme de Cherrier, pour chaque $\epsilon >>1 $ fix\'e. Alors, on a la compacit\'e locale tout en ayant un element dont le maximum depasse le rang $ \epsilon >>1$.

\smallskip

On a considerer le cas d'une boule de $ {\mathbb R}^4 $, mais on peut prendre une boule d'une vari\'et\'e Riemannienne $ M $ de dimension 4. Ici , on s'est interess\'e au cas de la dimension 4, on a la meme chose , pour un operateur (laplacien) coercif, pour l'eq. de Yamabe, en dimension 6, par exemple.

\bigskip

{\bf Remarque:} 1) dans le cas des ouverts de l'espace euclidien, quand on considere le laplacien, le principe du maximum s'applique, le minimum est bien atteint sur le bord et on a le lemme de Hopf. Par contre quand on considere un operateur general avec terme linaire, on a la condition "non positive minimum". Donc, il faut faire attention dans le cas general, ce qu'on a dit n'est pas vrai toujours, il faut que le minimum soit atteint \`a l'interieur(par exemple on considere le cas $ R_g\equiv 0 $ ou le laplacien seul et coercif, petites boules geodesiques), mais pour le laplacien sur des ouverts de l'espace euclidien, c'est vrai. En particulier, on a l'estimation locale en dimension 4 pour un ouvert regulier avec $ h <0 $ sur le bord, par exemple une boule. 

\smallskip

Il faut supposer l'operateur conforme coercif.

\bigskip

2) Dans le cas plat, en dimension 4, on a une estimation uniforme $ L^{\infty}_{oc}$, obtenue, par blow-up, a l'interieur ou eclatement, sans concentration de la mesure. En dimension 2, dans Brezis-Merle, ils ont utilis\'e la concentration de la mesure. Ici, en dimension 4, sans concentration de la mesure. 

3)On a: \underbar { {\bf Compacit\'e locale avec condition de Neumann non-lineaire}}.

\bigskip

On regarde maintenant: la compacit\'e locale, et, si il y a blow-up au bord:

\smallskip

Par exemple, considerons le probleme de Cherrier suivant:

\begin{displaymath}  \left \{ \begin {split} 
      \Delta u_i  \,\,\,\, \,\, & = a u_i^{3}, \,\,u_i >0     \,\, &&\text{dans} \!\!&& \Omega , \\
                 \partial_{\nu} u_i-u_i  & =-u_i^2, \,\, u_i >0   \,\,             && \text{sur} \!\!&&\partial \Omega.               
\end {split}\right.
\end{displaymath}

Avec $ \Delta=-\sum \partial_{jj}$, et $  a >0 $.

\smallskip

Soit $ \lambda_1 >0$, la premiere valeur propre du Laplacien avec condition de Dirichlet.

\smallskip
Comme on l'a fait precedemment, avec la deriv\'ee normale et le principe du maximum, en utilisant la deuxieme equation, on a:

$$ \min_{\Omega} u_i >1 . $$

\smallskip

En utlisant, la premiere fonction propre et une integration par parties, comme dans Brezis-Merle. On obtient:

$$ a \cdot (\min_{\Omega} u_i)^2 \leq \lambda_1, $$

On voit alors que le systeme n'a pas de solution si $ a\geq  \lambda_1$.

\smallskip

Placons nous autour de $ \lambda_1 >0$. Alors, soit le probleme de Cherrier n'a pas de solutions autour de $ \lambda_1 $: $ \exists \epsilon_0 >0 $, tel que pour $ a \in ]\lambda_1-\epsilon_0, \lambda_1]$, le probleme de Cherrier n'a pas de solutions. Ou bien, il y a une infinit\'e de solutions du probleme de Cherrier avec $ a=\lambda- \epsilon, \epsilon \to 0 $, et comme ce parametre converge vers $ \lambda_1 $, ou il n'y a pas de solutions au probleme de Cherrier, les solutions blow-up ( sinon on aurait une solution au probleme de Cherrier en $ \lambda_1 $, ce qui n'est pas possible). 

\smallskip

D'apres le resultat qu'on a donn\'e ci-dessus, de compacit\'e locale pour le probleme de Cherrier, il y a compacit\'e interieure et donc blow-up au bord.

\smallskip

Donc, en placant, $ a >0 $, autour de $ \lambda_1 >0$, on a:

\smallskip

a) Soit il n'y a pas de solutions au probleme de Cherrier, sur un intervalle $]\lambda_1-\epsilon_0,\lambda_1], \epsilon_0>0 $. Un resultat de non-existence.

\smallskip

b) Soit on a un exemple de compacit\'e interieure sans compacit\'e globale: blow-up au bord.

\smallskip

Du point de vue Mathemtiques: on a 2 possibilit\'es: c'est possible.

\smallskip

Du point de vue de la physique: la th\'eorie du tout: cela veut dire, que c'est vrai physiquement. Donc, les 2 possibilit\'es peuvent se realiser.

Il faut voir comment le pheneomene peut se realiser, comment il se produit, a partir de quoi.

\smallskip

Ici, on a montr\'e que,(du point de vue des maths), comment on peut les obtenir: il suffit de se placer autour d'un point $ a >0$ o\`u il n' y a pas de solution. Ici on a pris $ a=\lambda_1 >0$.

\smallskip

On montre comment on obtient: cet exemple (compacit\'e interieure avec blow-up au bord), c'est un exemple de la physique, qui se r\'ealise physiquement. Il suffit de se placer autour de $ a>0 $ o\`u il n'y a pas de solution, ici on a pris $ a=\lambda_1 >0 $.

\smallskip

Notons, que pour utiliser le resultat de compacit\'e locale: on a suppos\'e que les solutions existent.

\smallskip

De maniere similaire, on se place autour de $ a>0 $ o\`u la solution n'existe pas (singularit\'e), (pour $ a=\lambda_1 $, on est sur qu'il existe), et on suppose qu'il y des solutions autour de $ a >0 $: donc : on a la compacit\'e locale sans compacit\'e globale: blow-up au bord. C'est un resultat mathematique, c'est aussi, un exemple de phenomene physique.

\smallskip
L'alternative ci-dessus, a) ou b), est mathematique et le fait de la prouver fait que c'est possible mathematiquement. Le a) est possible et physiquement realisable, le b) est possible et physiquement realisable. Dans la theorie du tout: le b) est un exemple de la physique.

\smallskip

Non existence de solutions: singularit\'e, vide (ce phenomene physique est possible), singularit\'e essentielle (pas de solution en un point, mais autour une infinit\'e de solutions, ce phenomene physique est possible).

\smallskip

//////////////////////////

\smallskip

Une consequence de l'in\'egalit\'e: $ (\sup_K u)^{\alpha} \times \inf_M u \leq c $. On regarde par exemple $ n=4 $, alors $ \alpha =1/3 $: On a:

$$ \inf_{B_{2R}(0)} u \geq c_1 \int_{B_R(0)} u^3 dx = c_1 \int_0^R (\int_{\partial B_r(0)} u^3 d\sigma_r) dr, $$

$$ \sup_{B_R(0)} u \times \inf_{B_{2R}(0)} u \geq c_2 \int_{B_R(0)} u^4 dx, $$

donc,

$$ c\geq (\sup_{B_R(0)} u)^{1/3} \times \inf_{B_{2R}(0)} u \geq c_3 (\int_{B_R(0)} u^4 dx)^{1/3} \times (\inf_{B_{2R}(0)} u )^{2/3} \geq $$

$$ \geq c_4 (\int_{B_R(0)} u^4 dx) ^{1/3} \times (\int_0^R (\int_{\partial B_r(0)} u^3 d\sigma_r) dr )^{2/3}, $$

On obtient:

$$ (\int_{B_R(0)} u^4 dx) ^{1/3} \times \left (\int_0^R (\int_{\partial B_r(0)} u^3 d\sigma_r) dr \right )^{2/3} \leq c_5. $$

On a:

$$ [Vol(u^2\cdot g, B_R(0))]^{1/3} \times \left [\int_0^R (Aire(u^2\cdot g, \partial {B_r(0)})) dr \right ]^{2/3} \leq c_6.$$

On a:

$$ [Vol_{(u^2\cdot g)}(B_R(0))]^{1/3} \times \left [\int_0^R (Aire_{(u^2\cdot g)}( \partial {B_r(0)})) dr \right ]^{2/3} \leq c_6.$$

1) Interpretation du point de vue de la g\'eometrie: in\'egalit\'e isop\'erimetrique: relation entre le volume conforme et l'aire conforme.

\smallskip

C'est une in\'egalit\'e isoperimetrique entre volume conforme et aire conforme pour la metrique conforme $ u^2\cdot g $.

\smallskip

Quand $ \alpha =1 $, on a l'in\'egalit\'e optimale $ \sup \times \inf $ qui implique la majoration uniforme locale du volume conforme. Dans le cas du resultat de YY.Li-L.Zhang, $ \alpha =1 $, les interpretations g\'eometrique, Math\'ematiques et Physique, se confondent: ils bornent localement le volume conforme (interpretation g\'eometrique), ils bornent localement le volume conforme (qui est le $ p-$ volume (point de vue math\'ematique) et la $ p-$ energie (point de vue de la physique), avec $ p=\frac{2n}{n-2}=4 $: exposant maximal).

\smallskip

Ici ($\alpha=\frac{1}{3}$), il y a s\'eparation entre l'interpretation g\'eometrique (in\'egalit\'e isop\'erimtrique, volume conforme et aire conforme) et mathematique ($p-$volume, avec $ p=\frac{10}{3}>3 $) et physique ($p-$energie avec $ p=\frac{10}{3}>3$).

\smallskip

Ici, $ n=4 $ et $ \alpha =1/3 $, on a une relation entre le volume conforme local et l'aire conforme locale, c'est une in\'egalit\'e isop\'erimetrique pour la metrique conforme $ u^2\cdot g $.

\smallskip

Il est possible  d'avoir des resultats de ce type pour $ n=3,5 $. On peut avoir des resultats de ce type en dimension $ n\geq 5 $ d\'es qu'on a l'in\'galit\'e $ (\sup)^{\alpha} \times \inf \leq c $.

\smallskip

Lorsque $ \alpha =1 $ et l'operateur conforme est coercif et $ V=1$, Y.Y.Li-L.Zhang ont prouv\'e l'in\'egalit\'e optimale, il borne le volume local dans ce cas.

\smallskip

Ici, on suppose qu'on considere l'eq. suivante ($ R_g =-6 $):

$$ \Delta u-u=u^3, u >0, $$

avec $ \Delta =-\nabla^i \nabla_i $.

\smallskip

l'operateur conforme n'est pas necessairement coercif. Donc, on n'est pas dans le cas Li-Zhang, on n'a pas l'in\'egalit\'e optimale. Si de plus on suppose que le laplacien est coercif, petites boules geodesiques, alors on peut appliquer ce qu'on a dit precedemment, l'in\'egalit\'e isoperimetrique ci-dessus. c'est un cas, ou n'a pas l'in\'egalit\'e optimale, o\`u on peut appliquer ce qu'on a dit avant.

\smallskip

2) Interpretation du point de vue des Math\'ematiques et de la Physique: le $ p-$ volume ou la $ p-$ \'energie est localement uniformement born\'ee ($ p=\frac{10}{3} > 3 $): il existe une constante positive $ c>0 $, telle que pour toute solution $ u >0 $ de l'equation et tout compact $ K $ de $ M $, on ait:

$$ \int_K u^{\frac{10}{3}} dV_g \leq c. $$

Une norme $ L^p, p=\frac{10}{3}>3 $, est localement uniform\'ement born\'ee.

\smallskip

3) une autre interpretation fondamentale, en g\'eom\'etrie conforme: la m\'etrique conforme est de la forme $ g_u=u^2 \cdot g $, en g\'en\'eral, on a, $ g_u= u^{4/(n-2)} \cdot g $, (la longueur d'un segment est: $ dy = u^{2/(n-2)} dx $). La norme $ L^p $, peut s'ecrire comme $ u^p= (u^{2/(n-2)})^{q+1}, q>0 $, on a alors un volume conforme au sens de la convexit\'e, de la fonction convexe, $ t\to t^{q+1}, t>0, q>0 $. Au sens de la convexit\'e. Pour le volume conforme $ u^{2n/(n-2)} $, cela revient a considerer le volume d'un $ n-$cube, on mulitiplie la longueur d'une arrete $ a>0 $, $ n-$fois. Dans le cas d'une norme $ L^p $, cela revient \`a considerer une fonction convexe $ t\to t^{q+1}, t >0, q>0$ (cela est semblable au fait de fixer une variable et considerer le volume restant par rapport au volume maximal, dans un $ n-$cube). C'est un volume conforme au sens des fonctions convexes (sous-critique, au lieu d'avoir le cas critique, on a un volume au sens de la convexit\'e. Par exemple: volume d'un $ n-$ parall\'el\'epip\`ede: $ \underbrace{u^{2/(n-2)} \times \ldots \times u^{2/(n-2)}}_{q+1:fois} \times \underbrace{1\times 1 \times \ldots \times 1}_{n-q-1:fois} $, quand $ q+1 $ est un entier, c'est le theoreme de Fubini ou le theoreme de Fubini-Tonnelli).

Ici, en dimension $ n=4 $, {\bf on a bien un \underbar{volume conforme} au sens de la convexit\'e}, avec $ u^{10/3}=(\sqrt {u^2})^{10/3} $, la fonction convexe est: $ t\to t^{10/3}, t>0 $.

De meme, en dimension $ n=3 $ et en dimension $ n=5 $, {\bf on a bien des volumes conformes au sens de la convexit\'e}, avec:

$$ u^{16/3}=(\sqrt {u^4})^{8/3}, t \to t^{8/3}, t >0,\,\,  {\rm pour} \,\,  n=3. $$

Et,

$$ u^{52/21}=(\sqrt {u^{4/3}})^{26/7}, t\to t^{26/7}, t>0, \,\,{\rm pour} \,\, n=5. $$

Ceci est en relation avec la g\'eom\'etrie convexe, la g\'eom\'etrie conforme et les volumes des convexes. Volume convexe.

(G\'eom\'etrie convexe: "volume mixte", Alexandrov, Fenchel, Jessen, ... Par exemple: sur le volume mixte: si $ B_1 $ est la boule unit\'e, le volume maximal: $ V(K) $, et les volumes mixtes, $ V(K, B_1)=V(\partial K) $ et apres par recurrence, $ V(K, B_1, B_1)= V(\partial K, B_1)= V(\partial \partial K)$ et ainsi de suite, le bord, puis le bord du bord, puis  le bord du bord du bord....).

\smallskip

///////////////////////////////////////////////////////////////////////////////////////////////////////////////////////////////

\smallskip

{\bf Sur les articles de compacit\'e en dimension 2:}

\smallskip

Sur l'article de 2016. J.M.S.U.Tokyo. Avec singularit\'e au bord. et,

Sur l'article de 2022.Elect.Journal.Diff.Eq.

\smallskip

Pour avoir la compacit\'e interieure, il est n\'ecessaire de supposer $ \int_{\Omega} |x|^{-2\alpha} e^u dx \leq C <+\infty $. Pour appliquer le resultat de Brezis-Merle, corollaire 7. Sans cette condition, on ne peut pas parler de compacit\'e interieure.

\smallskip

Compacit\'e interieure:

\smallskip

a)

\smallskip

Etape 1: defintion du probleme: $ u \in W^{1,1}_0(\Omega), V\in L^{\infty}(\Omega), |x|^{-2\alpha} e^u \in L^1(\Omega) $, d'o\`u la regularit\'e.

\smallskip

Etape 2: $ V\geq 0 $ d'o\`u par le principe du maximum, la solution $ u\geq 0 $.

\smallskip

Etape 3: $ \int_{\Omega} |x|^{-2\alpha} e^u dx \leq C <+\infty \Rightarrow \forall x_0 \in \Omega, \forall r >0, \int_{B_r(x_0)} e^u dx \leq C'=C'(\alpha, r, \Omega) <+\infty, \,\, 0 \leq \tilde V = V\cdot |x|^{-2\alpha } \leq \tilde b $. Ici, ca depend du poids, de la nature du probleme.

Etape 4: application du resultat  de compacit\'e interieure de Brezis-Merle, corollaire 7.

\smallskip

D'o\`u la compacit\'e interieure.

\smallskip

Tout ceci pour dire, que la nature du probleme (singularit\'e au bord, poids, l'etape 2 depend de l'etape 1) et la condition de l'etape 3 sont necessaire \`a l'obtention de la compacit\'e interieure, etape 4.

\smallskip

b)

\smallskip

De meme pour l'article au Electronic.Journ.Diff.Eq, (quand $ \epsilon >0$): il y a des etapes pour obtenir la compacit\'e interieure. 

\smallskip

c)

\smallskip

Pour les 2 autres articles de la dimension 2: la compacit\'e interieure est due au resultat de Brezis-Merle. Ici, il s'agit du probleme 1 de Brezis-Merle.

\smallskip

Sur la compacit\'e globale $ L^{\infty}(\Omega) $ en dimension 2: par exemple, on regarde:

$$ ||u_i||_{L^{\infty}(\Omega)}=u_i(x_i) $$

On suppose par exemple que: $ x_i \to x_0 \in \bar \Omega $: 

\smallskip

a) si $ x_0 \in \Omega $ alors on a la compacit\'e par Brezis-Merle.

\smallskip
 
b) si $ x_0 \in \partial \Omega $, on ne sait pas.

\smallskip

1) le fait de ne pas savoir ou se situe $ x_0 $, fait qu'on a considerer et essayer de prouver la compacit\'e gloable sans le resultat de Brezis-Merle. On sait pas si a) ou b) sont vrais, donc, on regarde $ ||u_i||_{L^{\infty}(\Omega)} $ sans le resultat de compacit\'e locale de Brezis-Merle.

\smallskip

2) 

2-1) Quand on a essay\'e de montrer la compact\'e globale, on a supposer par l'absurde que $ \max_{\Omega} u_i \to +\infty \Rightarrow x_0\in \partial \Omega $, par le resultat de Brezis-Merle. Donc on a utilis\'e Brezis-Merle dans la preuve. Mais ce qu'on a essay\'e de montrer , c'est un resultat de compacit\'e globale sans faire reference au resultat de Brezis-Merle au d\'ebut, car on a soit a) soit b) et le fait de ne pas savoir, fait qu'on a essay\'e de prouver un resultat de compacit\'e globale (sans compacit\'e locale de Brezis-Merle, au d\'ebut).

\smallskip

2-2) Par exemple, on pourrait dire: pourquoi on ne considere pas $ L^{\infty}(\bar \omega) $, avec $ \omega $, un voisinage relativement compact de $ \partial \Omega $. En effet, il suffit de considerer directement $ L^{\infty}(\Omega)$, car considerer, $ \omega $, c'est utiliser le resultat de Brezis-Merle, de compacit\'e interieure, ce fait est possible si on se place d\`es le debut sur tout $ \Omega $. Cela veut dire qu'on s'est plac\'e au debut sur tout $ \Omega $: considerer $ \omega $, c'est utiliser Brezis-Merle et cela veut dire qu'on a considerer des le debut $ \Omega $ tout entier. On a bien deux espaces: $ L^{\infty}(\Omega-\bar \omega), L^{\infty}(\bar \omega)$, soit, $ L^{\infty}(\Omega) $.

\smallskip

Donc, il suffit de prendre $ L^{\infty}(\Omega)$.

\smallskip

3) 

3-1) quand on a dit compacit\'e avec contrainte, c'est pour mettre en evidence qu'il y a une contrainte: energie born\'ee. Mais c'est un resultat de compacit\'e, car dans le print "Cas d'existence de solutions d'edp", on a mis en evidence un exemple avec energies non-born\'ees et $ V_i \equiv b_i=constantes >0, b_i \to 0 $, (l'exemple de Brezis-Merle, avec condition de Dirichlet au bord):

$$ u_i(r)=\log \frac{8c_i^2}{(1+c_i^2 r^2)^2}-\log \frac{8c_i^2}{(1+c_i^2)^2}, c_i \to +\infty.$$

3-2) Puis il y a l'exemple de Brezis-Merle avec energies born\'ees, mais blow-up au bord. Donc, il faut rajouter des conditions.

\smallskip

Dans la compacit\'e des solutions de l'eq. de Yamabe. Il n'y a pas de contrainte, car c'est la dimension $ n\geq 3 $, sauf pour la sphere. Soit il faut avoir un exemple avec energies divergentes, soit il faut le prouver directement. Alors, on aura comme consequence, les energies sont born\'ees.

\smallskip

En dimension 2, on a un exemple avec energies divrgentes. Donc: compacit\'e. Avec energies born\'ees. ici, il y a bien compacit\'e, ce n'est pas seulement un crit\`ere de compacit\'e, il y aussi, compacit\'e. 

\smallskip

//////////////////////////////////////////////////////////////////////////////////////////

\smallskip

On regarde pourquoi, il n'y a pas forcement de compacit\'e pour le probleme de Liouville avec singularit\'e au bord quand $ \alpha \in [1/2,1[ $:

\smallskip

Au sens des distributions, au sens faible, au sens usuel:

On ecrit: soustraire le potentiel Newtionnien: $ u= \log r* \dfrac{1}{|x|^{2\alpha}} + O(1) $, donc:

$ \nabla u= \dfrac{1}{r}*\dfrac{1}{|x|^{2\alpha}}+O(1) $, par Giraud, voir le livre d'Aubin, $ \nabla u $, $ 1+(2-2\alpha) = 3-2\alpha \in (1,2] $, est de l'ordre de: $ 1/r^{2\alpha -1} $ si $ \alpha \in (1/2,1) $, $ \log r $, si $ \alpha =1/2 $,  qui est singulier pour $ \alpha \in[1/2,1[$. Donc, $ u \not \in C^1 $. On sait que par les estimations elliptiques que $ u\in C^0$.

\smallskip

Supposons que pour tout potentiel $ V $, Lipschitzien, on peut avoir la compacit\'e $ C^0 $, on sait qu'elle ne peut pas etre $ C^1$. On prend le cas particulier $ V= |x-x_0|^{2\alpha} $, dans ce probleme d'ordre $ 2\alpha $. Pour $ \alpha \in [1/2,1[ $, $ V $ est bien Lipschitzien. Les solutions deviennent $ C^1$ uniform\'ement, ce qui n'est pas vrai, d'apres ce qu'on a dit avant.

\smallskip

Donc, pour $ \alpha \in [1/2,1[ $, on n'a pas forc\'ement la compacit\'e.

\smallskip

Dans le produit de convolution, le potentiel $ V $ est au sens $ L^{\infty} $: $ V\in [0,b] $, il oscille entre les valeures $ 0 $ et $ b >0 $. Donc, ce qu'il faut considerer, c'est $ k_1=\log r * \dfrac{1}{|x|^{2\alpha}} $ et $ k_2=\dfrac{1}{r}* \dfrac{1}{|x|^{2\alpha}}$ , qui oscille entre $ 0 $ et $ +\infty $, $k_2\in[0,+\infty] $, il peut prendre des valeurs infinies.

\smallskip

Par les estimations elliptiques, on a, $ u \in C^0 $, on n'arrive pas \`a prouver qu'il n'est pas $ C^1$, mais en ecrivant avec le potentiel Newtonien, on voit bien que $ u \not \in C^1 $. Le potentiel $ V $ n'agit qu'au sens $ L^{\infty} $.

\smallskip

Ce probleme \`a l'ordre $ 2\alpha, \alpha \in[1/2,1[ $, a comme propri\'et\'e d'avoir les solutions $ u\in C^0$ et $ u \not \in C^1 $. Il faut voir le terme: $ k_2=\dfrac{1}{r}* \dfrac{1}{|x|^{2\alpha}}$: terme principal.

\smallskip

Apres, pour la compacit\'e, on fait agir un potentiel particulier, qui fait que $ u $ devient $ C^1 $, ce qui n'est pas vrai.
 
\smallskip

Donc, pour $ \alpha \in [1/2,1[ $, on n'a pas forc\'ement la compacit\'e.

\smallskip

//////////////////////////////////////////////////////////////////////////////

\smallskip

{\bf Remarque:} La compacit\'e locale, dans Brezis-Merle est li\'ee \`a la notion de $ (\sup, \inf) $, Brezis-Merle prouvent un r\'esultat de concentration-compacit\'e qui implique la compacit\'e locale si les solutions ont un signe constant (par exemple positif ou nul). La condition de Dirichlet et le principe du maximum impliquent la compacit\'e locale. Ceci est li\'e \`a la notion de $ (\sup,\inf) $. Brezis-Merle donnent un exemple de solutions divergentes au voisinage du bord. Ce qui veut dire que la notion de $ (\sup, \inf) $ n'est plus suffisante pour avoir la compacit\'e globale. Donc, la notion de $ (\sup, \inf) $ n'implique pas directement la compacit\'e globale, elle n'est plus suffisante. 

\smallskip

C'est l'etape suivante: compacit\'e globale: apr\`es la notion de $ (\sup,\inf)$. On a:

\smallskip

a) Principe de concentration compacit\'e.

\smallskip

b) positivit\'e des solutions, ou, condition de Dirichlet: $(\sup, \inf) $ $ \Rightarrow $: compacit\'e locale.

\smallskip

c) condition de Dirichlet, $ (\sup,\inf) $: $ \not \Rightarrow $: compacit\'e globale.

\smallskip

On peut s'arreter la: compacit\'e locale et blow-up au bord.

\smallskip

Si on veut la compacit\'e globale. Par exemple:

\smallskip

d) Potentiel $ C1 $ ou Lipschitzien: $ \Rightarrow $: compacit\'e globale.

\smallskip

e) L'in\'egalit\'e $ (\sup,\inf) $ n'est pas en relation avec la compacit\'e globale. La compacit\'e globale n'est pas en relation avec la notion de $ (\sup,\inf)$. (Dans la preuve ces deux notions sont li\'ees, mais, ceci n'est pas vrai directement).

La compacit\'e globale n'a rien a voir avec l'in\'egalit\'e $ (\sup,\inf) $, dans la preuve elles sont li\'es, mais directement, elles ne sont pas li\'es. On ne deduit pas la compacit\'e globale, directement de l'in\'egalit\'e $ (\sup,\inf) $.

\smallskip

////////) Plus pr\'ecis\'ement,  on r\'epete ce qu'on a dit avant, dans le cas o\`u le potentiel est Lipschitzien: on consid\`ere: Brezis-Li-Shafrir et le probl\`eme de Dirichlet: on explique mieux tout ceci:

\smallskip

 Soit $ E=\{ -\Delta u= Ve^u\} $ et $ (P): \sup +\inf \leq c $: alors considerer $ (E)$ : on a $ (P) $ et $ non (P) $.

\smallskip

 Soit $ \tilde E $: le probleme aux limites de Dirichlet: $ \tilde E=\{-\Delta \tilde u= \tilde V e^{\tilde u}, \,\, \tilde u = 0 \, {\rm \, sur \, le \, bord}\}$:

\smallskip

 alors on ne peut pas definir la propri\'et\'e $ (P'): \tilde u: \sup_K \tilde u \leq c $: la compacit\'e locale, \`a partir de Brezis-Li-Shafrir, car:

 \smallskip
 
si on pose le probleme $ (P') $: rien ne dit qu'on ne tombe pas sur  une solution telle que: $ \tilde u =(u_i) $ verifiant $ non (P) : \sup_K u_i \to +\infty $.

\smallskip

Considerer l'equation, inclut le fait de considerer le probl\`eme de Dirichlet. Si on suppose le domaine regulier, alors considerer l'equation de Liouville, inclut le fait de considerer le probl\`eme de Dirichlet. 

\smallskip

Donc, on ne peut pas definir la propri\'et\'e $ (P') $: compacit\'e locale, du probleme aux limites $ \tilde E $, qu'on veut prouver.

Car considerer $ (P') $ \`a partir de $ E $, c'est avoir $ (P) $ et $ non (P)$.

\smallskip

Il faut definir le probleme aux limites $ \tilde E $ et prouver la compacit\'e locale, $ (P') $, en dehors de Brezis-Li-Shafrir.(avec blow-up).

Par exemple, Brezis-Merle, prouvent la compacit\'e locale, sans Brezis-Li-Shafrir, mais ici, c'est sans blow-up.

\smallskip

Ici, la propri\'et\'e $ (P')$: prend toutes les possibilit\'es. Si, on l'a lie \`a Brezis-Li-Shafrir, c'est prendre en compte $ (P) $ et $ non (P) $. Ceci, pour dire que parfois on utilise Brezis-Li-Shafrir, dans des preuves, il est possible de le faire, car, les suites qu'on considere, ne sont pas dans un enemble plus general, qui prend en compte $ (P) $ et $ non (P) $. On utilise Brezis-Li-Shafrir, sans definir un ensemble o\`u il y a toutes les possibilit\'es.

\smallskip

Par contre, quand on veut considerer, le probleme au limites $ \tilde E $, on doit definir la propri\'et\'e $ (P') $, de maniere generale, c'est a dire, prendre en compte toutes les possibilit\'es. Si de plus, on l'a lie \`a Brezis-Li-Shafrir, cela revient \`a considrer $ (P) $ et $ non (P) $. Ce qui veut dire qu'on peut tomber, d\`es le depart, sur des solutions blow-up, sans avoir pu definir la propri\'et\'e $ (P') $: donc, on ne peut pas definir $(P')$ et $ non (P') $, donc, on ne pas les prouver ou utiliser Brezis-Li-Shafrir pour les prouver.

\bigskip

///////) Ceci est un argument qui confirme ce qu'on a dit auparavant, sur le fait que le $ (\sup,\inf) $ n'a rien a voir avec la compacit\'e globale. 

\smallskip

////////////////////////////////////////////////////////////////////////////////////////////////////////////////////////

\smallskip

Sur le 2eme preprint de 2024: "Blow-up analysis...": on a une solution $ v \in \dot H_1^q(M) $, au sens des distributions de :

$$ \Delta v+h v= \sum_{j=1}^k \alpha_j\delta_{x_j}, \alpha_j >0, $$

On a dit qu'en comparant $ v $ et $ \sum_{j=1}^k \alpha_j G(x_j,\cdot) $, on obtient:

$$ v=\sum_{j=1}^k \alpha_j G(x_j,\cdot), \,\, p.p, $$

Soit $ u= v-\sum_{j=1}^k \alpha_j G(x_j,\cdot) $, alors: $ u\in \dot H_1^q(M), q>1 $, est solution au sens des distributions de:

$$ \Delta u+h u=0,$$

on a, $ \forall \phi \in C^{\infty}_c(M) $, par une integration par parties:

$$ \int_M u(\Delta \phi+h\phi)=0= \int_M (\nabla u\cdot \nabla \phi+h u\phi)=0, $$

Comme $ q>1$, par des theoremes de densite dans $ \dot H_1^{q'}(M)$, et une partition de l'unit\'e, on se ramene a des boules de ${\mathbb R}^n$, on peut prende $ \phi\in C^2_c(M)$.
\smallskip

On veut utiliser le th. de regularit\'e d'Agmon: "The $L^p$ approach to the Dirichlet problem":

\smallskip

Pour le bord, dans des cartes, on se ramene, \`a des demi-boules, $ u $ est localement solutions de, au voisinage de points du bord:

$$ \int_{B^+} g^{ij}(\partial_i u \partial_j\phi+h u\phi) \sqrt {|g|} dx=0, \forall \phi\in C^2_c(B^+), $$

On utilise l'argument de: Brezis-Marcus-Ponce, on prenant $ \phi =\psi \gamma_n $, avec $ \psi \in C^2(B^+)$, nulle au voisinage de $ \partial^+B^+$ et nulle sur $ \partial^- B^+=B^+\cap {\mathbb R}^n_+$. La fonction $ \gamma_n $ est du type Brezis-Marcus-Ponce, s'annule avant d'atteinde la partie de $ {\mathbb R}^n_+$, et fonction elementaire de $ d(x,\partial^-B^+), 0 <d(x, \partial^-B^+) \leq 1/n $.

On a: pour tout $\psi\in C^2_0(B^+),$ comme dans le th. 6.2 d'Agmon:

$$ 0=\int_{B^+} g^{ij}(\partial_i u \partial_j(\psi\gamma_n)+h u(\psi \gamma_n)) \sqrt {|g|} dx \to \int_{B^+} g^{ij}(\partial_i u \partial_j\psi+h u\psi) \sqrt {|g|} dx,$$

La fonction $ \psi $ est a support dans $ B^+$, on la prolonge par $ 0 $ sur tout $ M $, puis on revient \`a l'integration dans $ M $: puis on utlise Stokes, ($ trace(u)=0$), puis on revient dans la boule $ B^+$, pour appliquer le th. 6.2, d'Agmon:

$$ 0=\int_{B^+} g^{ij}(\partial_i u \partial_j\psi+h u\psi) \sqrt {|g|} dx =\int_{M}(\nabla u \cdot \nabla \psi+h u\psi) dV_g= $$

$$= \int_{M} (u \Delta \psi+h u\psi) dV_g =\int_{B^+} (u \partial_i (g^{ij} \sqrt{|g|}\partial_j \psi)+ \tilde h u\psi dx =0. $$

On peut alors appliquer le th.6.2 d'Agmon pour avoir $ u\in H_2^p, p >n $ et par les injections de Sobolev on a $ u \in C^1$, puis on a globalement $ u $ presque partout egale a une fonction $ C^1$, on peut alors utiliser la coercivit\'e de l'operateur dans $ \dot H_1^2 $, pour avoir $ u=0. p.p$.

\smallskip

{\bf Remarque:} On peut dire qu'il y a equivalence entre, Brezis-Merle (dimension $ n=2 $) et le probleme de Dirichlet, pos\'e en dimension $ n\geq 3 $: En effet (raisonnement approximatif), la condition de Brezis Merle $ \int_ {\Omega} e^u dx \leq C $ implique que $ e^u $ "converge" \`a l'infini, par exemple quand $ x\to \partial \Omega $, ce qui veut dire que $ u(x)_{x\to \partial \Omega} \to -\infty $. On peut alors considerer la formulation de Brezis Merle: $ -\Delta u= Ve^u $ dans $\Omega $ et $ u=-\infty $ sur $ \partial \Omega $, (l'analyse blow-up avec singularit\'es, ou le point 3, du th\'eor\`eme de l'alternative de Brezis-Merle, se r\'ealise avec la condition limite $ u=-\infty $ sur le bord). 

On a alors l'analyse blow-up  en dimension $ n=2 $ (Brezis-Merle. ici, on utilise la notion de localement conform\'ement plat, pour ramener le probl\`eme pos\'e sur une surface Riemennienne, au cas d'un probl\`eme pos\'e sur un ouvert de $ {\mathbb R}^2 $) et on a l'analyse blow-up en dimension $ n\geq 3 $ (avec la condition d'interpolation).

\smallskip

Le proc\'ed\'e de concentration-compacit\'e est fait, soit pour l'analyse blow-up et la compacit\'e, soit, pour essayer de trouver des solutions variationnelles, et, topologiques, par le degr\'e topologique de Leray-Schauder. En dimension 2, quand le potentiel est regulier, on n'a pas d'exemples de solutions blow-up, directement. Alors que en dimension $ n\geq 3 $: sous-critique tendant vers le critique (sur des domaines \'etoil\'es), on a des exemples de solutions blow-up. Donc: on peut dire qu'en dimension $ n\geq 3 $, il n'est pas n\'ec\'essaire de chercher des solutions par le degr\'e topologique, car on a des contre-exemples directement. Donc, ce qui importe, c'est de caracteriser le ph\'enomene de blow-up et les criteres de compacit\'e, c'est l'objet du preprint de 2024: "Blow-up analysis...".

\smallskip

-Dans la formulation variationnelle: on utilise la technique du sous-critique tendant vers le critique, or, il y a des exemples (quand le domaine est \'etoil\'e), de solutions qui blow-up. Donc, on a blow-up, ou divergence des solutions dans la formulation variationnelle. Donc, il n'est pas possible d'utiliser la formulation variationnelle pour avoir l'existence de solutions (directement).

\smallskip

-Pour les solutions topologiques, il est n\'ecessaire d'avoir la compacit\'e, or, comme on l'a dit pr\'ecedemment, il y a des exemples de solutions qui blow-up ou divergent. Donc, ce n'est pas possible (directement), d'utiliser la methode topologique pour avoir l'existence de solutions. Dans la formulation topologique aussi, ce n'est pas possible directement, d'avoir des solutions.

\smallskip

Dans l'analyse blow-up: c'est essayer de trouver l'equilibre entre la gravit\'e (conditions sur les solutions $ u >0 $) et les reactions physiques ou chimiques ou biologiques (dus au potentiels $ V $), cet equilibre est li\'e \`a la masse $ m = f(V,u) $ (en dimension $ n\geq 3 $, il faut determiner aussi, ce qu'on appelle la masse). Par exemple, l'effondrement d'une etoile sur elle meme, c'est un trou noir, la gravit\'e prend le dessus. Un critere de compacit\'e, ce sont des conditions pour avoir cet equilibre ou non-effondrement.

\smallskip

La mani\`ere la plus directe, de faire l'analyse blow-up, est de determiner les solutions, les expliciter, puis chercher, les points o\`u il y a blow-up ou non blow-up. En general, il est difficile d'expliciter les solutions, donc, on suppose qu'elles existent, (du point de vus de la physique, de la chimie, de la biologie il y a existence de solutions, puisqu'on essaie de les observer ou les manipuler), puis on fait l'analyse blow-up math\'ematique et on d\'etermine, les critieres de compacit\'e.

\smallskip

C'est pr\'ecisement parce qu'il est difficile d'expliciter, voire, impossible d'expliciter, des solutions, qu'on fait l'analyse blow-up math\'ematique.

\smallskip

Dans le pr\'eprint de 2024: "Blow-up analysis...", on d\'etermine les conditions pour faire l'analyse blow-up math\'ematique:

\smallskip

1) On a la condition d'interpolation: condition gravitationelle: 

$$ ||u||_{\infty}^{\frac{1}{N}} ||u||_{N-1}^{1-\frac{1}{N}} \leq C<+\infty, \,\, N=\frac{2n}{n-2}, n\geq 3.$$

2) On a aussi la condition sur les potentiels: 

$$ 0\leq V \leq b<+\infty.$$

3) On d\'etermine la masse de r\'ef\'erence: $ \epsilon_n $: 

$$ \omega_n/2^n \geq \epsilon_n >0. $$

On d\'etermine une condition pour la compacit\'e, de non-effondrement: crit\`ere de compacit\'e:

$$ b\cdot C^N <\epsilon_n. $$

C'est une condition d'equilibre, condition de masse, comme on l'a dit pr\'ec\'edemment.

\smallskip

////////////////////////////////

\smallskip

Quand on consid\`ere l'equation de la courbure scalaire prescrite:

\smallskip

1) Contrairement \`a la dimension 2, en dimension $ n\geq 3 $, la condition gravitationnelle n'est pas une condition sur le volume conforme. Aussi, le fait de considerer la condition au bord de Dirichlet, fait qu'on ne considere pas seulement l'equation de la courbure scalaire prescrite, de plus, les solutions doivent etre $ >0$, or ici, elle s'annulent au bord. Cette condition au bord, fait que les solutions s'annulent au bord, la condition gravitationnelle n'est plus de nature g\'eom\'etrique conforme.

\smallskip

2) Dans la condition gravitationnelle: $ ||u||_{\infty}^{\frac{1}{N}} ||u||_{N-1}^{1-\frac{1}{N}} \leq C<+\infty, \,\, N=\frac{2n}{n-2}, n\geq 3 $, il y a la contrainte topologique li\'ee au terme: $ ||u||_{\infty} $ (norme $ C^0 $ quand les fonctions sont r\'eguli\`eres), et, il y a la contrainte g\'eom\'etrique $ ||u||_{N-1} $, qui apparait dans l'equation. Donc, c'est bien une condition topologique et g\'eom\'etrique, au lieu du volume conforme.

\smallskip

3) L'exemple qu'on donne dans la note de 2006 et l'article associ\'e de 2007: illustre ce ph\'enom\`ene de blow-up et donne tout son sens \`a l'analyse blow-up en mathematique, en dimension $ n\geq 3 $:

\smallskip

3-1) Dans l'exemple: le potentiel $ V \equiv 1 $ (potentiel plat): entre deux constantes strictement positives, $ C^{\infty}$ et $ C^{\omega} $.

\smallskip

3-2) Dans l'exemple: la condition gravitationnelle est vraie.

\smallskip

3-3) Dans l'exemple: l'energie usuelle ou la norme $ H^1 $ est uniform. entre deux constantes strict.positives.

\smallskip

3-4) Dans l'exemple: les solutions divergent (ou blow-up) et tendent vers 0 en dehors des points de concentrations.

\smallskip

3-5) Dans l'exemple, il y a nombre fini de points de concentrations: ph\'enom\`ene de concentration ou blow-up.

\smallskip

L'exmple donne tout son sens au ph\'enom\`ene blow-up en dimension $ n\geq 3 $, et \`a l'analyse blow-up en math\'ematiques, en dimension $ n\geq 3 $.

\smallskip

////////////////////////////////////////////////////////////////////////////////////////////////

\smallskip

{\bf R\'esum\'e sur les th\'eories en questions, les equations associ\'ees et les notions associ\'ees:}

\smallskip

-"Liouville setting": dimension 2.

\smallskip

-"Einstein-Lichnerowicz setting":dimension $ n\geq 3 $.

\smallskip

1) {\bf Les differentes th\'eories et les differents themes:}

\smallskip

-Th\'eorie de Yang-Mills, Equations de Yang-Mills. 

\smallskip

-Th\'eorie conforme des champs de Liouville. Th\'eorie de la gravitation quantique de Liouville. 

\smallskip

-Th\'eorie de la gravitation quantique; en relativit\'e g\'en\'erale($ n=3 $), de Kaluza-Klein ($ n=4$), dans la th\'eorie des cordes $ (n=5,6$) et des supercordes ($ n=9 $). Particules: Axion, graviton. Supersymetrie, sym\'etrie de la garviations (incorpor\'ees dans l'eq d'Einstein) + symetrie quantique(symetrie conforme): Eq. de Yamabe. Eq. d'Einstein-Lichnerowicz (quantique relativiste).

\smallskip

Eq de Yamabe; Eq. de la courbure scalaire prescrite: equations de la supersymetrie: particules de spin entier et de spin demi-entier. 

\smallskip

-Eq. de Schrodinger(quantique non relativiste).

\smallskip

-Modele des interactions des particules. 

\smallskip

-Modele cosmologique. Trous noirs.

\smallskip

-Modele de Liouville en dimension 2: th\'eorie de Liouville. Th\'eorie quantique des champs conforme de Liouville.

\smallskip

-Modele D'Einstein-Lichnerowicz en dimensions $ n\geq 3 $. Th\'eorie conforme des champs: en relativit\'e g\'enerale. Dimensions supplementaires. Ondes gravitationnelles(tenseur de Weyl non nul). Particules: Axion, Graviton. Quantique relativiste. 

\smallskip

2) {\bf Les Equations en questions:}

\smallskip

-Eq. de Liouville. Eq.du type Liouville (Eq. de la courbure scalaire prescrite en dimension 2). 

\smallskip

-Eq. de Yamabe, Eq. de la courbure scalaire prescrite en dimensions $ n\geq 3 $: (Equations de la supersymetrie). Eq. du type courbure scalaire prescrite. Eq.d'Einstein-Lichnerowicz (dans un certain cas, Equation de la brisure de symetrie).

\smallskip

3) {\bf Les differentes notions:}

\smallskip

D-Branes: Eq. avec condition de Dirichlet au bord (vari\'et\'e \`a bord, bord regulier: $ n=2 $ Eq. de Liouville ou de courbure scalaire et $ n\geq 3 $ Eq. d'Einstein-Lichnerowicz, Eq. de Yamabe). 

\smallskip

Masse positive, expansion de l'univers. Operateurs coercif.

\smallskip

4) {\bf Les differentes notions:}

\smallskip

-Enroulement. Torsion. Distrotion. Noeuds.

\smallskip

-Stabilit\'e. Coh\'erence. 

\smallskip

-Effondrement d'espaces: blow-up. Non-effondrement d'espaces: compacit\'e, rigidit\'e.

\smallskip

Le blow-up a plusieurs sens: eclatement: en mathematique et en physique: en physique une corde qui tend \`a eclater, c'est une limite. Espace qui s'effondre. ou en mathematiques: une transformation, rescaling ou changement d'echelle: eclatement math\'ematique. 

\smallskip

//////////////////////////////////////////////////////////////////////////////////////////////////////

\smallskip

a) relativistic and geometric quantum theory (physics and astronomy): eq.Yamabe. eq; de la courbure scalaire prescrite.

\smallskip

b) relativistic and non geometric quantum theory (physics and astronomy): eq; d'Einstein-Lichnerowicz, generale($ \Psi\not = constante, \nabla \Psi \not =0 $).

\smallskip

c) non-relativistic and non geometric quantum theory (physics): eq. de Schrodinger (par exemple ce qu'on a appel\'e eq; du type Yamabe ou du type courbure scalaire prescrite ou avec termes non lineaires, sous-critiques).

\smallskip

//////////////////////////////////////////////////////////////////////////////////////////////////////////////////////////////

\smallskip

\section{Quelques Problemes:}

\smallskip

{\bf Concernant les in\'egalit\'es de Harnack:}

\smallskip

1) Fixant un minorant $ m >0 $, on a une relation entre $ \sup_K v $ et $ m $ pour toute $ v >0 $ relativement a un $ W $ avec des conditions a priori sur $ W $, $ v $ et $ W $ sont li\'es par une equation ou plusieurs equations. Or, en prenant le minorant $ m=\inf_M u >0 $, on a l'inegalit\'e pour $ v>0 $ or $ u>0 $ est d\'eja une solution du probleme, donc $ \sup_K u $ est fonction de $ \inf_M u $ et des autres param\`etres.

\smallskip

2) Une interpretation fondamentale, en g\'eom\'etrie conforme: {\bf \underbar{volume conforme} au sens de la convexit\'e} et l'\'energie au sens de la convexit\'e:

La m\'etrique conforme est de la forme, en g\'en\'eral, on a, $ g_u= u^{4/(n-2)} \cdot g $, (la longueur d'un segment est: $ dy = u^{2/(n-2)} dx $). La norme $ L^p $, peut s'ecrire comme $ u^p= (u^{2/(n-2)})^{q+1}, q>0 $, on a alors un volume conforme au sens de la convexit\'e, de la fonction convexe, $ t\to t^{q+1}, t>0, q>0 $. Au sens de la convexit\'e. Pour le volume conforme $ u^{2n/(n-2)} $, cela revient a considerer le volume d'un $ n-$cube, on mulitiplie la longueur d'une arrete $ a>0 $, $ n-$fois. Dans le cas d'une norme $ L^p $, cela revient \`a considerer une fonction convexe $ t\to t^{q+1}, t >0, q>0$ (cela est semblable au fait de fixer une variable et considerer le volume restant par rapport au volume maximal, dans un $ n-$cube). C'est un volume conforme au sens des fonctions convexes (sous-critique, au lieu d'avoir le cas critique, on a un volume au sens de la convexit\'e. 

Par exemple: volume d'un $ n-$ parall\'el\'epip\`ede: $ \underbrace{u^{2/(n-2)} \times \ldots \times u^{2/(n-2)}}_{q+1:fois} \times \underbrace{1\times 1 \times \ldots \times 1}_{n-q-1:fois} $, quand $ q+1 $ est un entier, c'est le theoreme de Fubini ou le theoreme de Fubini-Tonnelli).

Ceci est en relation avec la g\'eom\'etrie convexe, la g\'eom\'etrie conforme et les volumes des convexes. Volume convexe.

(G\'eom\'etrie convexe: "volume mixte", Alexandrov, Fenchel, Jessen, ... Par exemple: sur le volume mixte: si $ B_1 $ est la boule unit\'e, le volume maximal: $ V(K) $, et les volumes mixtes, $ V(K, B_1)=V(\partial K) $ et apres par recurrence, $ V(K, B_1, B_1)= V(\partial K, B_1)= V(\partial \partial K)$ et ainsi de suite, le bord, puis le bord du bord, puis  le bord du bord du bord....).

\smallskip

Pour l'equation de la courbure scalaire prescrite en dimensions 3 et 4:

\smallskip

a) En dimension 3, on a la dependance explicite du $ \sup $ en fonction du $ \inf $, ici, l'operateur n'est pas coercif et le bord n'est pas necessairement regulier,

\be (\sup_K u)^{1/3} \times \inf_M u \leq c. \ee

On voit alors l'influence de la courbure: scalaire, Ricci, Weyl, et l'influence de la courbure scalaire prescrite $ V $ et de la non coercivit\'e de l'operateur, dans tous les cas: on n'a pas forc\'ement l'in\'egalit\'e optimale, $ \sup \times \inf $, mais on a l'in\'egalit\'e $ \sup^{1/3} \times \inf $. L'energie maximale ou le volume maximal, ne sont pas n\'ecessairement born\'es. Mais on borne localement, le volume conforme au sens de la convexit\'e et l'energie au sens de la convexit\'e.

(Pour l'in\'egalit\'e optimale $ \sup \times \inf $, et, le volume local=energie locale uniform. born\'es, voir, la partie g\'eom\'etrisation).

En dimension $ n=3 $, pour l'equation de la courbure scalaire prescrite, {\bf on a bien un \underbar{volume} \underbar{conforme} au sens de la convexit\'e}, et, l'\'energie au sens de la convexit\'e, avec la fonction convexe, $ t \to t^{8/3}, t >0 $:

\be \int_K u^{16/3} \leq c, \,\, u^{16/3}=(\sqrt {u^4})^{8/3}=(u^2)^{8/3}, \,\, t \to t^{8/3},\, t >0,\,\,  {\rm pour} \,\,  n=3. \ee

\smallskip

b) En dimension 4, on a une estimation a priori. On plus on a des exemples pour ce cas. Dans le cas ou la courbure scalaire prescrite est constante strict.positive. On a une inegalit\'e de Harnck implicite, cette in\'egalit\'e est differente de celle explicite, ce sont pas les memes constantes, $ c $, $ c: m \to c(m) $.

\smallskip

En general, dans le cas d'une vari\'et\'e Riemannienne $ M $ de dimension 4, on a une estimation a priori. on peut se ramener a une estimation du $ \sup $ en fonction du $ \inf $ en ecrivant:

$$ \sup_K u \leq c(a, b, \inf_M u, K, \Omega), $$

pour toute $ u >0 $ solution de l'equation de la courbure prescrite en dimension 4 sur $ M $, relativement a $ V $ (Lipschitzienne) verifiant:

$$ ||\nabla V||_{\infty} \leq \dfrac{3a {\sqrt a}}{32e^2{\sqrt 2}}\inf_M u, $$

et,

$$ 0 < a \leq V \leq b < +\infty, $$

\smallskip

En effet, $ u >0 $ et $ V $ sont li\'es par l'equation, on peut supposer que le gradient de $ V $ et $ \inf_M u $ sont aussi li\'es, ici on a deux equations liant $ u>0 $ et $ V $, une contrainte de plus. On a des exemples quand $ ||\nabla V_i||_{\infty} \leq A_i \to 0 $ et $ \inf_M u_i \geq m >0 $, on a une estimation locale uniforme.

\smallskip

Ici, la constante $ k(a)=\dfrac{32e^2{\sqrt 2}}{3a {\sqrt a}}$ entre $ \min_M u $ et la norme du gradient est plus grande que dans le cas plat(ouvert de $ {\mathbb R}^4 $). Il faut faire la preuve pour $ \min_M u \geq m >0 $, car on doit eliminer des termes en $ o(1) $. Dans le cas d'un ouvert de $ {\mathbb R}^4 $, la constante $ k(a)=\dfrac{8e^2{\sqrt 2}}{3a {\sqrt a}} $. \footnote{ Pour la preuve, comme c'est locale, on prend $ x_0 \in M $, on fait un changement de metrique conforme $ \tilde g = e^f g=\phi^3 g $ de sorte  que $ \tilde Ricci_{x_0} = 0 $, on a (voir le livre d'Aubin) $ f(x_0)=0 $ et $ \phi(x_0)=1 $. On prouve l'estimation en supposant $ \min_M u \geq m >0 $ pour la nouvelle fonction $ v=u/\phi $ pour $ m>0$ (car, voir l'article de la preuve du cas general en dimension 4, il y a des termes $ \tilde Ricci_{y_i} $ et $ \tilde R_{\tilde g}(y_i) $, des termes en $ o(1) $), comme $ \phi(x_0)=1 $, $ \min_{B_r^{\tilde g}(x_0)} v \geq (\min_M u)/2 \geq m/2 >0 $. Comme on fait un changement de metrique conforme, au voisinage de $ x_0 $, tout est de "l'ordre" de $ u $ et $ m $ car $\phi(x_0)=1 $, $ f(x_0)=1 $, l'exponentielle (en norme de matrices) est $ (1+\epsilon)-$ Lipschitzienne, on a chang\'e de metrique, les majorants restent de l'ordre de $ 1+\epsilon $ au voisinage de $ x_0 $ car $ \phi(x_0)=1 $, ($f(x_0)=0$), et la constante $ k(a)$ entre $ \min_M u $ et la norme du gradient de $ V $ est plus grande. 

On prouve alors, en supposant l'inegalit\'e entre le $ \min_M u $ et le gardient de $ V $: $ \min_M u \geq m >0  \Rightarrow \sup_{B_r^{\tilde g}(x_0)} v \leq c \Rightarrow \sup_K u \leq c(a,b,m,K, M) $. Apres on remplace $ m $ par $ \min_M u $.}

\smallskip

{\bf Remarque:} Dans lecas d'un ouvert $\Omega $ de $ {\mathbb R}^4 $, la deuxieme contrainte est "v\'erifi\'ee" (apr\'es blow-up), dans l'article de C.C.Chen et C.S.Lin, si on suppose $ V \in C^2 $ uniform\'ement. Ils obtiennent une inegalit\'e entre $ \nabla V(y_i) $ et $ [u_i(y_i)]^{-1} $, dans ce cas ils obtiennent l'in\'egalit\'e optimale $ \sup \times  \inf $. (les points $ y_i $ sont les points blow-up).

i) Un exemple, ce resultat s'applique lorsque $ V\equiv constante >0 $. On retrouve le r\'esultat de YY.Li et L. Zhang, eux obtiennent l'in\'egalit\'e optimale $ \sup \times \inf $. 

ii) Un autre cas, on l'a dit ci-dessus, c'est celui o\`u $ ||\nabla V_i ||_{\infty} \leq A_i \to 0 $ et $\min_M u_i \geq m >0 $, on obtient une estimation locale uniforme.
 
iii) On a un r\'esultat analogue a celui de Brezis-Merle en dimension 2, en dimension 4: pour deux suites $ (u_i), (V_i) $ solutions de l'equation de la courbure scalaire prescrite en dimension 4, $ 0 < a \leq V_i \leq b $ et $||\nabla V_i||_{\infty} \leq A_i \to 0 $, alors $ u_i\geq m >0 \Rightarrow \sup_K u_i \leq c(a,b, (A_i), m, K, M) $. En particulier, si on prend $ m=\inf_M u_i >0 $, on a l'in\'egalit\'e de Harnack implicite, $ \sup_K u_i \leq c(a,b, (A_i), \inf_M u_i, K, M) $ (La fonction $ m\to c $ est decroissante de $ m >0 $ et change en une autre fonction $ \tilde c $ si on remplace les suites $ (u_i,V_i) $ par une d'autres suites $ (v_i,W_i) $ avec des hypoth\`eses similaires, et on obtient: $ \sup_K v_i \leq \tilde c(\tilde a, \tilde b, (\tilde A_i), \inf_M v_i, K, M) $).

 \footnote{ Concernant l'exponentielle dans le preuve du r\'esultat ci-dessus. On commence par changer la metrique $ d_g \to d_{\tilde g} $, on a un facteur $ (1+\epsilon) $, puis on considere $ \tilde \exp_y $, l'exponentielle pour la metrique $ \tilde g $, construite a partir d'une carte $ (\Omega, \psi) $ normale en $ x_0 $ pour la metrique $ \tilde g $, on passe de $ d_{\tilde g} $ \`a $ ||\psi o \tilde \exp ||_{{\mathbb R}^4}$, on a un majorant $ (1+\epsilon) $, puis on utlise le fait que $ \dfrac{d}{dz} [z\to \tilde \exp_{x_0} (z)]_{|z=0} =id_{{\mathbb R}^4} $, c'est \`a dire que c'est vrai pour la differentielle de $ \psi o \tilde \exp $ et par continuit\'e des differentielle, (en $ y $ et $ z $ assez petits), on a un majorant $ (1+\epsilon) $, puis on utilise la norme matricielle $ \sup $ des matrices et on ecrit que: $ \psi o \tilde \exp_y(z_1)-\psi o \tilde \exp_y(z_2) =\int_0^1 \partial_s [\psi o \tilde \exp_y(s(z_1-z_2)+z_2)] ds=\int_0^1 [(\partial_z(\psi o \tilde \exp_y)(s(z_1-z_2)+z_2))\cdot (z_1-z_2)] ds $ et on obtient; $ ||\psi o\tilde \exp_y(z_1)-\psi o \tilde \exp_y(z_2)||_{{\mathbb R}^4} \leq (1+\epsilon) ||z_1-z_2||_{{\mathbb R^4}} $, uniform\'ement en $ y $ voisin de $ x_0 $ et $ |z_1|\leq r, |z_2|\leq r, r >0 $ assez petit.}  
\smallskip

La formulation pr\'ec\'edente en dimension 4 s'applique lorsqu'on suppose en $ x_0 $ point critique de $ V>0 $ $ V\in C^1(M) $. on fait la preuve avec $ \min_M u \geq m >0 $ pour $ v=u/\phi$ avec $ \phi>0 $ duc changement de metrique conforme telle que $ \tilde Ricci_{\tilde g}(x_0)=0 $, on obtient:

$$ \exists \,c,R >0,\,\, \sup_{B_R(x_0)} v \leq c, $$

Donc,

$$ \exists \, c,R >0,\,\, u(x_0) \leq \sup_{B_R(x_0)} u \leq c,\, \forall\,  u>0 $$

Donc, on obtient une in\'egalit\'e du type:

$$ u(x_0)\leq c(x_0,V, \inf_M u, M, g), $$

pour toute solution $ u >0 $ de l'equation de la courbure scalaire prescrite en dimension 4 relativement \`a un $ V >0 $, $ V\in C^1(M) $ fix\'e ayant le point $ x_0 $ comme point critique, $ \nabla V(x_0)=0 $. Ou bien $ V \in C^{1,\alpha}(M) $, $ 0 < \alpha \leq 1 $, $ 0 < a \leq V \leq b < +\infty $ et $ ||V||_{C^{1,\alpha}}\leq A $ et $ \nabla V(x_0)=0 $, on obtient:

$$ \forall \,\, u >0,\,\, u(x_0)\leq c(a, b, A, \alpha,x_0,\inf_M u, M, g). $$

Il y a des exemples de solutions de l'equation de la courbure prescrite en dimension 4 avec $ x_0 $ point o\`u $ V $ est maximum et $ M $ compacte sans bord dans le livre d'Aubin.

Quand on prend les notations de YY.Li et L.Zhang: $(M,g)=(B_1(0)\subset {\mathbb R}^4, g) $ et $ 0 < V\in C^1(B_1(0)) $  avec $ \nabla V(0)=0 $, on obtient l'estimation a priori dans un voisinage de $ 0 $ en supposant le minimum unifrom\'ement minor\'e et aussi:

\be u(0) \leq c(V, B_1(0), g, \inf_{B_1(0)} u), \ee

pour toute solution $ u >0 $ de l'equation dela courbure scalaire prescrite en dimension 4 relativement a $ V $ avec les conditions precedentes ou $ V\in C^{1,\alpha}(B_1(0))$ avec $ 0 < \alpha \leq 1 $, $ 0 < a \leq V \leq b < +\infty $ et $ ||V||_{C^{1,\alpha}}\leq A $ et $ \nabla V(0)=0 $.

\be \forall \,\, u >0,\,\, u(0)\leq c(a, b, A, \alpha,\inf_{B_1(0)} u, B_1(0), g).\ee

\smallskip

////////////////////////////////////////////////////////////////////////////////////////////////////////////////////////

\smallskip

Pour l'equation de Yamabe en dimensions 5 et 6:

\smallskip

On voit l'effet des courbures scalaire, de Ricci et de Weyl et de la non-coercivit\'e (l'operateur conforme n'est pas necessairement coercif): on n'a pas l'in\'egalit\'e optimale $ \sup \times \inf $, mais on a des in\'egalit\'es interm\'ediaires. L'energie maximale ou le volume maximal ne sont pas necessairement born\'es. Mais, on borne localement, le volume conforme au sens de la convexit\'e et l'energie au sens de la convexit\'e.

\smallskip

On r\'epete ce qu'on a dit auparavant, sur le fait que borner les solutions dans un $ L^p, p >N-1 $ est le fait d'avoir un {\bf \underbar{volume conforme} au sens de la convexit\'e}, et, l'\'energie au sens de la convexit\'e: 

\smallskip

Une interpretation fondamentale, en g\'eom\'etrie conforme: la m\'etrique conforme est de la forme, en g\'en\'eral, on a, $ g_u= u^{4/(n-2)} \cdot g $, (la longueur d'un segment est: $ dy = u^{2/(n-2)} dx $). La norme $ L^p $, peut s'ecrire comme $ u^p= (u^{2/(n-2)})^{q+1}, q>0 $, on a alors un volume conforme au sens de la convexit\'e, de la fonction convexe, $ t\to t^{q+1}, t>0, q>0 $. Au sens de la convexit\'e. Pour le volume conforme $ u^{2n/(n-2)} $, cela revient a considerer le volume d'un $ n-$cube, on mulitiplie la longueur d'une arrete $ a>0 $, $ n-$fois. Dans le cas d'une norme $ L^p $, cela revient \`a considerer une fonction convexe $ t\to t^{q+1}, t >0, q>0$ (cela est semblable au fait de fixer une variable et considerer le volume restant par rapport au volume maximal, dans un $ n-$cube). C'est un volume conforme au sens des fonctions convexes (sous-critique, au lieu d'avoir le cas critique, on a un volume au sens de la convexit\'e. Par exemple: volume d'un $ n-$ parall\'el\'epip\`ede: $ \underbrace{u^{2/(n-2)} \times \ldots \times u^{2/(n-2)}}_{q+1:fois} \times \underbrace{1\times 1 \times \ldots \times 1}_{n-q-1:fois} $, quand $ q+1 $ est un entier, c'est le theoreme de Fubini ou le theoreme de Fubini-Tonnelli).

Ceci est en relation avec la g\'eom\'etrie convexe, la g\'eom\'etrie conforme et les volumes des convexes. Volume convexe.

(G\'eom\'etrie convexe: "volume mixte", Alexandrov, Fenchel, Jessen, ... Par exemple: sur le volume mixte: si $ B_1 $ est la boule unit\'e, le volume maximal: $ V(K) $, et les volumes mixtes, $ V(K, B_1)=V(\partial K) $ et apres par recurrence, $ V(K, B_1, B_1)= V(\partial K, B_1)= V(\partial \partial K)$ et ainsi de suite, le bord, puis le bord du bord, puis  le bord du bord du bord....).

\smallskip

c) En dimension 5, il y a une relation explicite entre $ \sup_K u $ et $ \inf_M u $,

\be (\sup_K u)^{1/7} \times \inf_M u \leq c. \ee

d) En dimension 6, on a une estimation a priori et une relation non explicite entre $ \sup_K u $ et $ \inf_M u $. Pour toute solution $ u >0 $ de l'equation de Yamabe en dimension 6:

\be \sup_K u \leq c(K, M, g, \inf_M u). \ee

En dimension 6, si on veut inclure le cas $ \inf_M u =0 $ au cas $ \inf_M u>0 $ dans l'in\'egalit\'e: il faut considerer l'in\'egalit\'e suivante:

\be \sup_K u \times \inf_M u \leq \tilde c(K, M, g, \inf_M u) \times \inf_M u. \ee

Avec, $ \tilde c: {\mathbb R}_+ \to {\mathbb R}_+^*: \tilde c(m)=c(m)=c(K,M,g,m) $, si $ m >0 $ et $ \tilde c(0)=1 $. La fonction precedente, $ c $, est decroissante de $ m >0 $.

\smallskip

e) Pour l'equation de Yamabe et du type Yamabe en dimension 4, on a l'in\'egalit\'e:

\be (\sup_K u)^{1/3} \times \inf_M u \leq c. \ee

L'energie maximale ou le volume maximal ne sont pas n\'ecessairement born\'es. Mais on borne localement, le volume conforme au sens de la convexit\'e et l'energie au sens de la convexit\'e.

\smallskip

(Pour l'in\'egalit\'e optimale $ \sup \times \inf $ et le volume local=energie locale, uniform.born\'es, voir, la partie g\'eom\'etrisation).

En dimension $ n=4 $, pour l'equation de Yamabe, {\bf on a bien un \underbar{volume conforme} au sens de la convexit\'e}, et, l'\'energie au sens de la convexit\'e, avec la fonction convexe est: $ t\to t^{10/3}, t>0 $.

\be  \int_K u^{10/3} \leq c, \,\, u^{10/3}=(\sqrt {u^2})^{10/3}=u^{10/3},\,\, t\to t^{10/3}, t>0, \,\, {\rm pour} \,\, n=4.\ee

En dimension $ n=5 $, pour l'equation de Yamabe, {\bf on a bien un \underbar{volume conforme} au sens de la convexit\'e}, et, l'\'energie au sens de la convexit\'e,  avec la fonction convexe, $ t\to t^{26/7}, t>0 $:

\be \int_K u^{52/21} \leq c, \,\, u^{52/21}=(\sqrt {u^{4/3}})^{26/7}=(u^{2/3})^{26/7}, \,\, t\to t^{26/7}, t>0, \,\,{\rm pour} \,\, n=5. \ee

\smallskip

/////////////////////////////////////////////////////////////////////////////////////////////////////////////////////////////////////

\smallskip

Concernant l'article: "Some uniform estimates for scalar curvature type equations". En dimension 4, on a en plus de l'estimation a priori une dependance non explicite du $ \sup_K u $ en fonction du $ \inf_{\Omega} u $:

\be \sup_K u \leq c(a, b, A, \alpha, K, \Omega,\inf_{\Omega} u) ,\ee

ou bien, (pour inclure le cas $ \inf_{\Omega} u=0 $),

\be \sup_K u \times \inf_{\Omega} u \leq \tilde c(a, b, A, \alpha, K, \Omega,\inf_{\Omega} u) \times \inf_{\Omega} u ,\ee

pour le theoreme 3

\smallskip

et,

\be \sup_K u \leq c(a, V, K, \Omega,\inf_{\Omega}u). \ee

ou bien, (pour inclure le cas $ \inf_{\Omega} u=0 $),

\be \sup_K u \times \inf_{\Omega} u \leq \tilde c(a, V, K, \Omega,\inf_{\Omega}u) \times \inf_{\Omega} u. \ee

pour le theoreme 4.

\smallskip

////////////////////////////////////////////////////////////////////////////////////////////////////////////////////////////////

\smallskip

{\bf Remarque:} On a bien une fonction $ c $ de l'inf car, pour chaque $ m $ on a une constante  $ c_0 $ majorant de $ \sup_K u $ et l'ensemble des majorants est lui meme minor\'e par $ \sup_K u $, $ c $ est alors l'inf des majorants de $ \sup_K u $, qui est unique. On voit alors que cette correspondance $ m  \to c $ est bien une fonction et aussi la correspondance $ (a, b, K, \Omega, m) \to c $ est une fonction. Apres on remplace $ m $ par $ \inf_{\Omega} u $ comme on l'a dit ci-dessus. On a bien la correspondance $ (a, b, K, \Omega, \inf_{\Omega}u) \to c $ est une fonction. 

\smallskip

De plus la fonction $ m\to c $ est d\'ecroissante en $ m >0 $. On a alors, si $\inf_{\Omega} u \geq m >0 $ (ou $\inf_M u \geq m >0 $, si on considere une vari\'et\'e Riemannienne $ M $), alors $ c (a, b, K, \Omega, \inf_{\Omega} u) \leq c(a, b, K, \Omega, m) <+\infty $, ou $ ( c(a, b, K, M, \inf_M u) \leq c(a,b, K, M, m) $, si on est sur une vari\'et\'e Riemannienne $ M $), on obtient alors l'estimation a priori lorsque l'inf est uniformement minor\'e.

\smallskip

////////////////////////////////////////////////////////////

\smallskip

Par le proc\'ed\'e de g\'eom\'etrisation: on a: $ (u,V), u >0, V>0 $ solutions de:

$$ -\Delta u+\frac{(n-2)}{4(n-1)} R_g u= Vu^{N-1}, u >0, N=\frac{2n}{n-2}, n\geq 3, \quad (*) $$

avec, $ 0< m \leq \lambda + \frac{(n-2)}{4(n-1)} R_g \leq 1/m $.

\smallskip

$ K $ d\'esigne un compact quelconque de la vari\'et\'e Riemannienne $ (M,g) $ non necessairement compacte sans bord:

$$ \sup_K u\times \inf_M u \leq c,$$

et si $ -\Delta-\lambda $ est coercif, alors:

$$ \int_K u^{2n/(n-2)} dV_g \leq c, $$

Ici on est dans le cas positif. Ceci inclut le cas non loc.conf.plat.

\smallskip

/////////////////////////////////////////////////////////////////////

\smallskip

Sur une vari\'et\'e Riemannienne compacte \`a bord r\'eguliere, $ (M,g) $, $ h $ r\'eguli\`ere avec $ -\Delta +h $ coercif dans $ H^1(M) $:

$$ -\Delta u+hu=Vu^{N-1}, \,\, u >0, \,\, {\rm dans} \,\, M, \,\, N=\frac{2n}{n-2}, n\geq 3, $$

avec $ 0 < a \leq V \leq b <+\infty $, $ V $ h\"olderienne. Alors:

$$ \forall K {\rm \,\, compact\,\, de \,\,}  M, x_0\in K, \,\,\sup_M u\times \inf_K u \geq c(\inf_K u) >0, $$

ou alors,

$$ \forall K {\rm \,\, compact\,\, de \,\,}  M, x_0\in K, \,\,\sup_M u \geq \frac{c(\inf_K u)}{\inf_K u} >0, $$

avec $ m \to c(m) $ croissante de $ m >0 $.

\smallskip

///////////////////////////////////////////////////////////////////////////////////////////////////////////////////////

\smallskip

Dans leur articles de 1998. Journal of Differential Geometry. C.C.Chen-C.S.Lin, montrent qu'on n'a pas forc\'ement $ \sup $ major\'e si $ \inf $ est minor\'e, puis, donnent des conditions pour avoir l'in\'egalit\'e optimale $ \sup \times \inf $. Selon certaines hypotheses, les energies peuvent etre born\'ees ou non.

\smallskip

On a aussi les faits suivants:

I) 1) Conditions de platitudes d'ordre $ C^{n-2} $ de Yanyan Li $ \Rightarrow $:

1-1) blow-up,

1-2) blow-up isol\'e,

1-3) blow-up isol\'e simple $ \Rightarrow $,

1-3-1) $ \sup \times \inf <+\infty $,

1-3-2) Energies born\'ees, $ \int u^{2n/(n-2)} <+\infty $.

Voir les articles de Yanyan Li sur la sphere et de C.C.Chen-C.S.Lin, commun.pure.appl.math. 1997.

\smallskip

2) Conditions de platitudes d'odre $ C^{\alpha}, \frac{(n-2)}{2} \leq \alpha <n-2 $ de C.C.Chen-C.S.Lin $ \Rightarrow $:

2-1) blow-up,

2-2) blow-up isol\'e,

2-3) Pas de blow-up isol\'e simple $ \Rightarrow $,

2-3-1) $ \sup \times \inf  \to +\infty $,

2-3-2) Energies born\'ees, $ \int u^{2n/(n-2)} <+\infty $.

Voir les articles de C.C.Chen-C.S.Lin, J.Differential.Geometry.1998.

\smallskip

3) Conditions de platitudes d'odre $ C^{\alpha}, 1\leq \alpha < \frac{(n-2)}{2} $ de C.C.Chen-C.S.Lin $ \Rightarrow $:

3-1) blow-up,

3-2) Energies non born\'ees, $ \int u^{2n/(n-2)}  \to +\infty $.

Voir les articles de C.C.Chen-C.S.Lin, J.Differential.Geometry.1998.

\smallskip

4) {\it Propri\'et\'e intermediaire:} Propri\'et\'e interm\'ediaire entre les points 1) et 2): In\'egalit\'es du type $ (\sup)^{\alpha} \times \inf <+\infty $. Voir travail de Th\`ese et Journal.Funct.Analysis.

\smallskip

II) Concernant les in\'egalit\'es du type $ (\sup)^{\alpha} \times \inf <+\infty $: Il faut considerer $ M_i^{\beta} $ dans les methodes, avec $ M_i=u_i(y_i) $, $ y_i $ les poins blow-ups et r\'esoudre une equation en $ \beta $ du premier ordre li\'ee a la positivit\'e de la fonction $ (+-) $ la fonction auxiliaire. 

Ceci s'applique aux methodes de:

a) Brezis-Li-Shafrir.

b) YY.Li-L.Zhang,

c) C.C.Chen-C.S.Lin.

Voir travail de Th\`ese et Journ.Funct.Analysis. Acta.Math.Scientia. Analysis in Theory and Applications.
\smallskip

III) Assertion minimale: $ \sup <+\infty $ si $ \inf >0 $: c'est vrai en dimension $ n= 6 $ Eq.de Yamabe. et en dimension $ n=4 $ Eq de Yamabe  et Eq de la courbure scalaire prescrite. Comme consequence $ \sup_K u \leq f(\inf_M u) $ pour l'Eq de Yamabe en dimensions 4 et 6, et, $ \sup_K u_i \leq f(\inf_M u_i) $, pour chaque suite $ (u_i) $ relativement a une suite $ (V_i) $ avec $, 0<a\leq V_i \leq b<+\infty $ et $ |\nabla V_i| \leq A_i \to 0 $, pour l'Eq de la courbure scalaire prescrite en dimension 4. Si on a une autre suite $ (v_i) $ relativement \`a une suite $ (W_i) $ Lipschitziennes avec $ 0 <c \leq W_i \leq d <+\infty $, de constantes de Lipschitz $ B_i \to 0 $ alors $ \sup_K v_i \leq g(\inf_M v_i) $. Voir Travail de Th\`ese en dimension 4 (Eq. courbure scalaire prescrite en dimension 4, cas plat) et Journal. Funct. Analysis. et Journal.Math.Anal.Applications.

\smallskip

//////////////////////////////////////////////////////////////////////////////////////////////////////////////////////////////////////

\smallskip

{\bf Consequence du resultat donnant une solution unique: inegalite de Sobolev et inegalite interpolation}

Soit $ (M,g) $ une vari\'et\'e Riemannienne compacte sans bord de dimension $ n \geq 4 $ et de courbure scalaire $ S_g >0 $ partout sur $ M $ et orientable. Alors, voir le livre d'Aubin pour $\epsilon >0 $ assez petit $ \epsilon <S_g $, l'equation  $ 4\dfrac{(n-1)}{n-2} \Delta u+\epsilon u =u^{N-1}, u>0 $ a une solution et le resultat d'unicit\'e dit que cette solution est unique. (Voir le livre d'Aubin, la fonctionnelle associ\'ee \`a cette equation possede un infimum atteint par cette solution unique). On a alors (voir aussi l'article de L.Veron et J.R. Licois et l'article sur Hal de J.Dolbeault, Esteban, Loss), une in\'egalit\'e d'interpolation et une in\'egalit\'e de Sobolev:

$ \exists \, \epsilon_0=\epsilon_0(n,M,g)>0,$ tel que pour $ 0 < \epsilon <\epsilon_0 $ on ait:

$$ \inf_{u\in H^1(M)-\{0\}} \dfrac{4\dfrac{(n-1)}{n-2} \int_M |\nabla u|^2dV_g+\epsilon\int_M u^2dV_g}{(\int_M |u|^N dV_g)^{2/N}}=\epsilon (vol(M))^{2/n}, $$

Qu'on peut ecrire:

\smallskip

Pour tout $ u \in H^1(M), $ on a:

$$ 4\dfrac{(n-1)}{n-2} \int_M |\nabla u|^2dV_g \geq \epsilon\left ( \left ( \int_M |u|^N dV_g \right )^{2/N} (vol(M))^{2/n}- \int_M u^2dV_g \right ), $$

ou encore,

$$ \forall u\in H^1(M), \, 4\dfrac{(n-1)}{n-2} ||\nabla u||_2^2\geq \epsilon \left ( (vol(M))^{2/n}||u||_{2^*}^2-||u||_2^2 \right )\geq 0, $$

avec $ 2^*=N=\dfrac{2n}{n-2} $ pour $ n\geq 4 $.

Finalement: il existe $ C=C(n,M,g) >0 $ telle que:

$$\forall u\in H^1(M),\, ||\nabla u||_2^2\geq C \left ( (vol(M))^{2/n}||u||_{2^*}^2-||u||_2^2 \right )\geq 0, $$

avec $ 2^*=N=\dfrac{2n}{n-2} $ pour $ n\geq 4 $.

On a une interpretation du point de vue de la physique de cette in\'egalit\'e. Pour les in\'egalit\'es isop\'erimetriques, on a une relation entre le perimetre et l'aire, du point de vue de la g\'eometrie. Ici, on a une relation entre l'\'energie et le $ N-$volume pour une distribution Sobolev $ H^1(M) $. Distributions, \'energie (normes $ H^1$) et $ N-$volume g\'en\'er\'e par cette distribution.

\bigskip

\subsection{Sur les problemes du type $ \sup \times \inf $:\\ }

1) Probleme de Li-Zhang:

\bigskip

Soit $(M,g)$ une vari\'et\'e Riemannienne de dimension $ n \geq 5 $ de courbure scalaire $ S_g $, on considere l'equation de Yamabe:

$$ \dfrac{4(n-1)}{n-2} \Delta u + S_g u=u^{N-1},\,\, u>0,\,\, N=2n/(n-2).$$

$\Delta =-\nabla^i(\nabla_i) $ l'operateur de Laplace-Beltrami.

\smallskip

a) Est il vrai d'avoir pour tout compact $ K \subset M $,
 
 $$ \sup_K u \times \inf_M u \leq c(K, M, g) \,? $$
 
b) En particuler pour $ n=6 $ (in\'egalit\'e explicite), est il vrai qu'on a: 
 
 $$ (\sup_K u)^{\alpha} \times \inf_K u \leq c, \,\, 0 < \alpha \leq 1\, ? $$

Voir les points 37) et 38), ce n'est pas forc\'ement possible d'avoir ce type d'in\'egalit\'es en considerant la classe de metrique $ g $ quelconque, avec la notion de platitude. Il faut considerer une classe de metriques plus petite quand on utilise la notion de platitude.
 
 \bigskip
 
 2) Soit $ (M,g) $ une vari\'et\'e Riemannienne de dimension 4 de courbure scalaire $ S_g $, on considere l'equation de la courbure scalaire prescrite $ V $:
 
 $$ \Delta u +\dfrac{1}{6} S_g u= Vu^3,_,\, u >0. $$
 
 $\Delta = -\nabla^i(\nabla_i) $ l'operateur de Laplace-Beltrami. On suppose:
 
 $$ 0 < a \leq V_i \leq b <+\infty,\,\, ||\nabla V_i ||_{\infty} \leq A_i \to 0 $$
 
 \smallskip
 
a) Est il vrai que:
 
 $$ \exists \, \alpha \in (0,1], \,\, (\sup_K u_i)^{\alpha} \times \inf_M u_i \leq c=c(a, b, (A_i)_i, K, M, g) \,  ? $$

 Voir les points 37) et 38)

\smallskip
 
b) Qu'en est il du cas $ (M, g)=(\Omega \subset {\mathbb R}^4, \delta) $ ?

 Voir les points 37) et 38)
 
\bigskip

c) (R\'esultat du type CC.Chen-CS.Lin du cas plat), qu'en t il si on suppose que $ V $ verifie:

$$  0 < a \leq V_i \leq b < +\infty, \,\, ||\nabla V_i ||_{\infty} \leq A, \,\, ||\nabla^2 V_i ||_{\infty} \leq B \, ? $$
 Plus pr\'ecisement, a t on:
 
 $$ \exists \, \alpha \in )0,1],\,\, (\sup_K u)^{\alpha} \times \inf_M u \leq c=c(a, b, A, B, K, M) \, ? $$

voir les preprints avec potentiel $ C^1 $ et $ C^2 $ sur les espaces localement sym\'etriques.

\subsection{Probleme de Brezis-Li-Shafrir:\\}

\smallskip

Soit $ \Omega $ un ouvert de $ {\mathbb R}^2 $, ici, $ \Delta = \partial_{11}+\partial_{22} $, on considere l'equation suivante:

$$ -\Delta u = Ve^u, $$

avec $ V>0 $ et continue sur $ \Omega $.

\smallskip

A t on, pour tout compact $ K \subset \Omega $ :

$$ \sup_K u+ \inf_{\Omega} u \leq c=c(K, V, \Omega) \, \, ? $$

Le probl\`eme est un probl\`eme ouvert en supposant $ V $ continue strict.positive et fix\'ee, sans etre Holderienne et sans etre Lipschitzienne. C'est dans ce sens la qu'il faut  considerer le probl\`eme, sans etre Holderienne et sans etre Lipschitzienne.(Par exemple, Brezis-Li-Shafrir pose la question de cette maniere, concernant le fait de  savoir si ce resultat est vrai si on prend des fonctions Holderiennes: la reponse a \'et\'e donn\'ee, par C.C.Chen-C.S.Lin: Holderienne sans etre Lipschitzienne).

\smallskip

Avec la notion de: $ (P) $ et $ non(P) $, en utilisant l'article de Brezis-Li-Shafrir avec des potentiels: $ V \in C^1, V>0, 0 < c_1 \leq |\nabla V|\leq c_2<+\infty $, on a des contre exemples. Avec: $ non(P): \sup + \inf \to +\infty $, impossibilit\'e d'utiliser le blow-up. On peut dire  que ce probl\`eme est resolu: en partant de l'article de Brezis-Li-Shafrir et de $ non(P) $.

\smallskip

On dit, maintenant, pourquoi c'est peu probable de le prouver:

\smallskip

La notion de continuit\'e est forte, mais moins precise que la notion de fonction Holderienne ou fonction Lipschitzienne. Comme on part de points $ y_i, z_i $, si on connait $ |y_i-z_i| $ alors on connait $ |V_i(y_i)-V_i(z_i)|$, pour ce qui concerne une fonction Holderienne ou fonction Lipschitzienne.

\smallskip

Par contre si on prend une fonction continue seulement, on n'a pas cette relation plus precise entre les points et les fonctions. C'est peu probable de r\'epondre \`a la question. Il se peut qu'avec les notions de th\'eorie de la mesure et de l'integration, c'est possible de r\'epondre \`a la question.

\smallskip

On a aussi: r\'egularit\'e de $ V \Rightarrow $ r\'egularit\'e de la solution $ u \Rightarrow (P): \sup+\inf \leq c $ et $ non(P) $, donc: moins de r\'egualrit\'e de $ V \Rightarrow $ moins de possibilit\'es d'avoir $ (P) \Rightarrow non(P) $ probable avec $ V\in C^0 $.

\smallskip

//////////////////////////////////////////

\smallskip

Pour le cas Holderien, Brezis-Li-Shafrir, le precisent bien, ils disent que le potentiel doit verifier: $ |V(x)-V(y)|\leq A |x-y|^{\alpha }, 0 <\alpha < 1 $. Donc, si on prend $ V $ Lipschitzienne, $ |V(x)-V(y)|\leq A |x-y| $, alors $ \alpha =1 $, donc, on n'a plus $ 0 <\alpha <1 $, ce qui n'est pas possible. Par contre pour le cas continue, c'est plus ambigu, ce n'est pas dit explicitement et donc, une fonction $ C^1 $ est $ C^0 $, on a alors $ non (P) $ ou un contre-exemple \`a la question 3 de Brezis-Li-Shafrir. Donc, le probleme d'etre $ C^0 $ sans etre Lipschitzienne et sans etre Holderienne, est un nouveau probleme, ce n'est pas la question 3 de Brezis-Li-Shafrir. C'est un nouveau probleme. Le $ non (P)$ est une r\'eponse implicite donn\'ee dans l'article de C-C.Chen-C-S.Lin, 1997. dans le cas de la dimension $ n\geq 3 $, qu'on peut appliquer en dimension 2 \`a la question 3 de Brezis-Li-Shafrir. Au cas o\`u on r\'esout cette question, on aura toujours le $ non (P) $, qui dit qu'il y a un contre-exemple.

\smallskip

1-En consid\'erant l'article de C-C.Chen. C-S.Lin en dimension 2: Communications in Analysis and Geomtery. 1998. On peut dire qu'ils ont repondu \`a la question: en prenant 2 fonctions $ K_1, K_2 \in [a,b], a, b >0 $, $ K_1$ du type logarithmique (un poids, par exemple), comme dans le print "A compactness result for an equation with $ C^0 $ condition", $ K_1$ non lipschitzienne en un point non holderienne en un point mais $ C^0 $, qui verifie la condition $ \frac{K_1(y)}{K_1(x)} \leq 1+\frac{B_1}{|\log|x-y||}$, et, $ K_2 \in C^1 $, $ 0 <c_2\leq |\nabla K_2|\leq c_1<+\infty $, alors $ K_2 $ verifie $ \frac{K_2(y)}{K_2(x)} \leq 1+\frac{B_2}{|\log |x-y||} $. On en deduit que le produit, $ K_1\cdot K_2 $ verifie les conditions de CC.Chen-CS.Lin: $ \frac{K_1(y)K_2(y)}{K_1(x)K_2(x)} \leq 1+\frac{B_3}{|\log |x-y||} $ et donc les solutions $ u $ de $ -\Delta u=K_1\cdot K_2 e^u $, satisfont $ \sup u +\inf u \to +\infty $: impossibilit\'e de faire le blow-up, car la preuve de CC.Chen-CS.Lin utilise les rescaling ou changement d'echelle ou blow-up, finalement le prduit $ K_1 \cdot K_2 $ est $ C^0 $ sans etre lipschitzien et sans etre holderien, on a le $ non(P)$ ou un contre exemple dans ce cas aussi: $ C^0 $ sans etre Lipschitzien et sans etre Holderien.

\smallskip

Un potentiel $ K>0 $ qui v\'erifie: $ \frac{K(y)}{K(x)} \leq 1+\frac{B}{|\log|x-y||}, B \geq 0 $, est $ C^0 $. Dans certains cas, on s'arrange pour qu'il ne soit pas Lipschitzien ni Holderien.

\smallskip

//////////////////////////

\smallskip

2-Par exemple, ceci concerne les probl\`emes 1 et 2 de Brezis Merle: on combine l'analyse blow-up de l'article "About Brezis-Merle problem with Holderian condition" et la compacit\'e pour un potentiel du type $ V=K_1\cdot K_2 $ avec $ K_1>0 $, un poids logarithmique non lipschitzien en un point, non holderian en un point, et, $ K_2 \in C^1, 0 \leq K_2 \leq b, 0< c_2 \leq |\nabla K_2|\leq c_1<+\infty $, dans l'article: "A compactness result for an equation with $ C^0 $ condition".

\smallskip

Alors: impossibilit\'e de faire le blow-up ou rescaling ou changement d'echelle: $ \sup_{\Omega} u_i \to +\infty $. Donc, on obtient des contre exemples aux probl\`emes 1 et 2 de Brezis-Merle avec potentiel $ C^0$ non Lipschitzien et non Holderian.

\smallskip

Les contre-exemple aux probl\`emes 1 et 2 de Brezis-Merle, ont \'et\'e explicit\'es avant. Mais ce contre exemple est donn\'e pour des potentiels $ C^0 $ sans etre Lipschitzien sans etre Holderien, au cas o\`u la question se pose. C'est l'apr\'es Brezis-Merle.

\smallskip

///////////////////////

\smallskip

 Concernant le probl\`eme 3 de Brezis-Merle. La question sur la meilleuere constante, $ C_{optimal}(M)$, elle est donn\'e par le resultat de Brezis-Merle lui meme. Soit: $ \tilde C(M)=\inf\{C, \sup_K u \leq C,\,\, {\rm si} \,\, u\geq -M>-\infty \,\, \forall \,\, u \} $, avec $ u $ solution de l'equation, $ -\Delta u=Ve^u, 0<a\leq V(x)\leq b<+\infty $. Alors: $ C_{optimal}(M)=\tilde C(M)$. Les autres questions sont r\'esolu par I.Shafrir. On peut dire que ce probl\`eme 3 est clos.

\smallskip

\subsection{ Problemes du type Brezis-Merle: \\}

\smallskip

Soit $\Omega $ un ouvert regulier de $ {\mathbb R}^2 $, ici, $\Delta = \partial_{11}+\partial_{22} $. On considere les equations suivantes:

$$ -\Delta u - \epsilon (x_1\partial_1 u+x_2 \partial_2 u) = |x-x_0|^{2\beta} V e^u,\,\,\epsilon \in (0,1], \,\, x_0 \in \partial \Omega, $$

avec,

$$ \beta \geq 0, \,\, 0 \leq V \leq b,\,\, |x-x_0|^{2\beta} e^u \in L^1(\Omega), \,\, u\in W_0^{1,1}(\Omega).$$

a) Donner le comportement "blow-up" au bord et un critere de compacit\'e des solutions en $ u $.

\bigskip

b) Qu'en est il si $ V $ est Lipschitzienne ?

\bigskip

c) De meme, donner le comportement "blow-up" et un critere de compacit\'e pour les solutions de:

$$ -\Delta u-\epsilon(x_1\partial_1u +x_2\partial_2 u) = |x-x_0|^{-2\alpha} V e^u, \,\, \epsilon \in (0,1], \,\, x_0\in \partial \Omega. $$

Ici $ \alpha \in (0,1) $ ou $ \alpha \in (0,1/2) $ et,

$$ 0 \leq V \leq b,\,\, |x-x_0|^{-2\alpha} e^u \in L^1(\Omega),\,\, u\in W_0^{1,1}(\Omega) $$

d) Qu'en est il si on suppose $ V $ Lipschitzienne ?

\smallskip

e) Memes questions avec le cas regulier; equation elliptique avec poids holderien et singualtit\'e au bord:

$$ -\Delta u - \epsilon (x_1\partial_1 u+x_2 \partial_2 u) = (1+|x-x_0|^{2\beta}) V e^u,\,\,\epsilon \in (0,1], \,\, x_0 \in \partial \Omega, $$

avec,

$$ 1/2 > \beta \geq 0, \,\, 0 \leq V \leq b,\,\, e^u \in L^1(\Omega), \,\, u\in W_0^{1,1}(\Omega).$$

\end{document}